\documentclass[nosumlimits,twoside]{amsart}
\usepackage{amssymb}
\usepackage{aliascnt}
\usepackage{graphicx}
\usepackage{float}
\usepackage{srcltx}
\usepackage[all, 2cell, cmtip]{xy}
\usepackage{cancel}
\usepackage{version}
\usepackage[T1]{fontenc}
\usepackage{xcolor}
\usepackage{enumerate}
\usepackage[inline]{enumitem}
\usepackage[normalem]{ulem}
\usepackage[2e]{psfrag}
\usepackage{yhmath}
\usepackage{array}
\usepackage{dsfont}
\usepackage{nicefrac}
\usepackage[colorinlistoftodos,prependcaption,textsize=tiny]{todonotes}

\usepackage{mathtools}
\usepackage{etoolbox}
\usepackage{tikz-cd}
\tikzcdset{scale cd/.style={every label/.append style={scale=#1}, cells={nodes={scale=#1}}}}
\usepackage{quiver}
\usepackage[backref=page]{hyperref}
\usepackage[capitalise,noabbrev]{cleveref}
\usepackage[final]{showlabels}
\usepackage{orcidlink}

\makeatletter
\def\namedlabel#1#2{\begingroup
    #2%
    \def\@currentlabel{#2}%
    \phantomsection\label{#1}\endgroup}
\makeatother

\makeatletter
\patchcmd{\@setaddresses}{\indent}{\noindent}{}{}
\patchcmd{\@setaddresses}{\indent}{\noindent}{}{}
\patchcmd{\@setaddresses}{\indent}{\noindent}{}{}
\patchcmd{\@setaddresses}{\indent}{\noindent}{}{}
\makeatother

\newtheorem{theorem}{\sc Theorem}[section]

\newaliascnt{proposition}{theorem}
\newtheorem{proposition}[proposition]{\sc Proposition}
\aliascntresetthe{proposition}

\newaliascnt{notation}{theorem}

\aliascntresetthe{notation}

\newaliascnt{assumption}{theorem}

\aliascntresetthe{assumption}

\newaliascnt{lemma}{theorem}
\newtheorem{lemma}[lemma]{\sc Lemma}
\aliascntresetthe{lemma}

\newaliascnt{corollary}{theorem}
\newtheorem{corollary}[corollary]{\sc Corollary}
\aliascntresetthe{corollary}

\theoremstyle{definition}
\newaliascnt{definition}{theorem}
\newtheorem{definition}[definition]{\sc Definition}
\aliascntresetthe{definition}

\newaliascnt{definitions}{theorem}

\aliascntresetthe{definitions}

\newaliascnt{example}{theorem}
\newtheorem{example}[example]{\sc Example}
\aliascntresetthe{example}

\newaliascnt{examples}{theorem}

\aliascntresetthe{examples}

\theoremstyle{remark}
\newaliascnt{remark}{theorem}
\newtheorem{remark}[remark]{\sc Remark}
\aliascntresetthe{remark}

\newaliascnt{remarks}{theorem}

\aliascntresetthe{remarks}

\newaliascnt{claim}{theorem}

\aliascntresetthe{claim}

\newaliascnt{noname}{theorem}

\aliascntresetthe{noname}

\newaliascnt{question}{theorem}

\aliascntresetthe{question}

\def\rev{\mathrm{rev}}
\def\revsn{\text{-}}
\def\revsnp{\scriptscriptstyle{+}}

\allowdisplaybreaks
\newenvironment{invisible}{{\noindent\sc \colorbox{yellow}{Invisible:}\;}\color{gray}}{\medskip}
\excludeversion{invisible}

\def\pulb{\ar@{}[dr]|(0.2){\mbox{\Large{$\lrcorner$}}}}

\newcommand{\id}{\mathrm{Id}}
\newcommand{\Cc}{\mathcal{C}}

\newcommand{\mm}{\mathfrak{M}}

\newcommand{\Bb}{\mathcal{B}}

\newcommand{\Rr}{\mathcal{R}}
\newcommand{\Ss}{\mathcal{S}}

\newcommand{\ot}{\otimes}

\def\Bialg{{\sf Bialg}}

\def\cHopf{{\sf Hopf}_{\mathrm{coc}}}
\def\Bimon{{\sf Bimon}}
\def\Comon{{\sf Comon}}
\def\Mon{{\sf Mon}}
\def\Grp{{\sf Grp}}
\def\Hopfmon{{\sf Hopf}}
\def\cHopfmon{{\sf Hopf}_{\mathrm{coc}}}
\def\cBimon{{\sf Bimon}_{\mathrm{coc}}}
\def\cComon{{\sf Comon}_{\mathrm{coc}}}

\def\TrHopf{{\sf TrHopf}}
\def\QTrHopf{{\sf QTrHopf}}
\def\MQTrHopf{{\sf MQTrHopf}}
\def\QTrBialg{{\sf QTrBialg}}

\def\initial{\mathsf{0}}
\def\terminal{\mathsf{1}}

\newcommand{\op}{\mathrm{op}}

\begin{document}
% 16T05--Hopf algebras and their applications
% 18E13--Protomodular categories, semi-abelian categories, Mal'tsev categories
% 18M50--Bimonoidal, skew-monoidal, duoidal categories

\title[]{Protomodularity of cocommutative Hopf monoids in duoidal categories and quasitriangular Hopf algebras}

\author[A.~Ardizzoni]{Alessandro Ardizzoni\, \orcidlink{0000-0001-7384-611X}}
\address{%
\parbox[b]{\linewidth}{University of Turin, Department of Mathematics ``G.Peano'', via
Carlo Alberto 10, I-10123 Torino, Italy}}
\email{alessandro.ardizzoni@unito.it}
\urladdr{\url{www.sites.google.com/site/aleardizzonihome}}
\author[L.~Bottegoni]{Lucrezia Bottegoni\, \orcidlink{0009-0005-8172-5605}}
\email{lucrezia.bottegoni@edu.unito.it}
\author[A.~Cigoli]{Alan Cigoli\, \orcidlink{0000-0002-5181-5096}}
\address{%
\parbox[b]{\linewidth}{University of Turin, Department of Mathematics ``G. Peano'', via
Carlo Alberto 10, I-10123 Torino, Italy}}
\email{alan.cigoli@unito.it}
\author[A.~Sciandra]{Andrea Sciandra\, \orcidlink{0009-0008-0447-287X}}
\address{%
\parbox[b]{\linewidth}{Département de Mathématiques, Université Libre de Bruxelles, Boulevard du Triomphe, B-1050 Bruxelles, Belgium}}
\email{andrea.sciandra@ulb.be}
\urladdr{\url{www.andreasciandra.com}}

\begin{abstract}
In this work, we extend the protomodularity of the category of cocommutative Hopf algebras to the quasitriangular setting. Every quasitriangular Hopf algebra admits a minimal quasitriangular Hopf subalgebra and, as we show, can be regarded as a cocommutative bimonoid in the tensor-braided duoidal category of bimodules over it. This leads us to investigate protomodularity in the broader context of Hopf monoids in duoidal categories. 
To this end, we adopt a slight modification of B\"{o}hm's notion of antipode  associated with a reversion, further refining an earlier one due to B\"{o}hm–Lack. This framework allows us to study Hopf monoids in this setting, Galois and co-Galois maps, and the factorization of Hopf monoids. Using these tools, we prove a factorization of points, the Split Short Five Lemma, and the existence of pullbacks of split epimorphisms along arbitrary morphisms in the category of cocommutative Hopf monoids  with monic unit in any tensor-braided duoidal category with a reversion; hence this category is protomodular.
As applications, we recover the protomodularity of cocommutative Hopf algebras in symmetric monoidal categories under mild assumptions, and we obtain that of the coslice category of quasitriangular (resp.\ triangular) Hopf algebras under a fixed subobject; in the triangular case, this can be traced back to a category of generalized internal groups, introduced in the present work. When the fixed subobject is minimal, we infer the protomodularity of the category of quasitriangular Hopf algebras whose minimal quasitriangular Hopf subalgebra is isomorphic to the fixed subobject, which we interpret as the protomodularity of an essential fibre of a functor. As a byproduct, our results extend the double cross product of cocommutative Hopf algebras to the quasitriangular
setting.
%Finally, for a surjective morphism $\pi:(A,\Rr)\to(H,\Rr)$ of quasitriangular Hopf algebras with $H\subseteq A$, we show that the kernel $K$ can be realized as the product of the coinvariant elements $A^{\mathrm{co}H}$ and the minimal quasitriangular Hopf algebra $B:=H_{\Rr}$, and that $A$ is isomorphic as a Hopf algebra to a generalized version $K\bowtie_{B}H$ of the double cross product of $K$ and $H$.
\end{abstract}

\maketitle

\tableofcontents

%It is known that the category $\Bialg_{cc}$ of cocommutative bialgebras is cartesian monoidal with binary product given by the tensor product. A triangular bialgebra $H$ can be thought of as a general version of cocommutative bialgebra as there is an invertible element  $\Rr\in H\otimes H$ such that $\Delta^\mathrm{op}(\cdot)=\Rr\Delta(\cdot)\Rr^{-1}$. In fact cocommutative Hopf algebras are trivially triangular by taking $\Rr=1\otimes 1$. We then ask whether the category   $\TrBialg$ of triangular bialgebras and morphisms of triangular bialgebras is cartesian monoidal. We will see that this is not true but that changing the morphisms we can solve the problem. This leads to introduce the category $\TwTrBialg$ of triangular bialgebras  and twisted  morphisms of triangular bialgebras and its quotient category $\GaTrBialg$ where morphisms are substituted by their gauge classes.

\section{Introduction}

The starting point of this paper is a result proved by Gran, Sterck and Vercruysse in \cite{GrStVe}, stating that the category of cocommutative Hopf algebras over an arbitrary base field, which is equivalent to the category of internal groups in the category of cocommutative coalgebras, is semi-abelian, hence, in particular, protomodular.
Our main aim is to
undertake the extension of  the protomodularity to quasitriangular Hopf algebras. In fact, any cocommutative Hopf algebra is quasitriangular with universal $\mathcal{R}$-matrix $1\otimes 1$. It is well-known that any quasitriangular Hopf algebra $A$ admits a minimal quasitriangular Hopf subalgebra $H$. We show that $A$ can be regarded as a cocommutative bimonoid in the tensor-braided duoidal category of $H$-bimodules. This leads us to investigate protomodularity in the setting of such duoidal categories. A first difficulty is choosing an appropriate notion of antipode to define Hopf monoids in the duoidal setting. For our purposes, we adopt a slight modification of the notion of antipode associated with a reversion, as sketched in \cite{Bohm-Canada}, which refines an earlier, more restrictive formulation given in \cite{Bohm-Lack}. 

The paper is structured as follows. In \cref{sec:duoidal}, we introduce our modified notion of reversion in the context of a duoidal category $(\Cc,\circ, I,\bullet,J)$, which involves a double lax monoidal functor together with four natural transformations obeying suitable compatibility conditions expressed in terms of commutative diagrams. From these, we recover some further commutative diagrams that will allow us to borrow some of the techniques of \cite{Bohm-Lack} also in our context. In particular, we will be able to define the notions of convolution inverse and of antipode for a bimonoid, and we will retrieve expected properties such as the uniqueness of convolution inverses and antipodes, and the compatibility of bimonoid morphisms between Hopf monoids with antipodes. Next, in presence of an antipode, we recover the invertibility of suitable morphisms  called Galois and co-Galois maps taken from \cite{Bohm-Lack}, and we derive consequences such as compatibility of the antipode with multiplication and comultiplication. We then focus on particular instances of these maps, namely the Galois and co-Galois induced morphisms, which will be crucial for protomodularity.  The machinery developed allows us to address factorization of bimonoids and Hopf monoids. Using results from \cite{BD-CrossII}, we completely describe those bimonoids $A$ endowed with bimonoid morphisms $i:K\to A$ and $j:H\to A$ such that $\varphi:=m_A^\circ(i\circ j):K\circ H\to A$ is invertible. More precisely, in \cref{thm:psi}, we transport to $K\circ H$ a unique bimonoid structure of twisted tensor product $K\circ_\psi H$ such that $\varphi$ is a morphism of bimonoids. Then, \cref{thm:psitriang} explores how an antipode shapes this structure; it will be used later to generalize the so-called matched pairs to the quasitriangular setting.

In \cref{sec:tens-br}, we deal with $\bullet$-braided duoidal categories, i.e.\ duoidal categories $(\Cc,\circ,I,\bullet,J)$ such that $(\Cc,\bullet,J)$ is a braided monoidal category whose braiding satisfies a compatibility condition with the interchange law. When such a duoidal category is equipped with a reversion, we show that the category of Hopf monoids becomes monoidal with respect to $\bullet$ (\cref{prop:Hopfbul}). One of the main results of this section is the \textit{factorization of points}, achieved in \cref{thm:main}. More precisely, assume that $\Cc$ has binary intersections and that $\bullet$ preserves them. Let $\pi:A\to H$ be a split epimorphism in $\Bimon \left(\Cc,\circ,\bullet\right) $ (hence a point) whose unit morphism
$u_{H}^{\circ }:I\to H$ is a monomorphism in $\Cc$ and where $A$ is cocommutative. Then $\pi$ has a kernel. If $\left( \mathcal{C}%
,\circ ,I,\bullet ,J\right) $ is moreover equipped with a reversion and $H$ is a Hopf monoid, then $\varphi:=m_{A}^{\circ }\left( k\circ \sigma \right)
:K\circ H\rightarrow A$ is an isomorphism in $\Cc$ and, if also $A$ is a Hopf monoid, then so is $K$, which turns out to be the kernel of $\pi$ in the category of cocommutative Hopf monoids. 
Thus, by the factorization result for Hopf monoids discussed above, $K\circ H$ becomes a Hopf monoid $K\circ_\psi H$ and $\varphi$ is an isomorphism of Hopf monoids. Moreover, in this setting, one can rewrite $\psi$ in terms of a morphism $\chi$ involving the antipode of $H$ and a right action $\triangleleft$ of $K$ on $H$, by using a co-Galois map, see \cref{coro:main}. The Split Short Five Lemma follows immediately (\cref{lem:split-short}). The second main result of the section is the existence, for each split epimorphism $\pi:A\to H$ as in \cref{thm:main} with $H$ Hopf monoid, of the pullback along any morphism $f:H'\to H$ of bimonoids with $H'$ a Hopf monoid (\cref{thm:ddlstab}).
These ingredients allow us to deduce that, for a $\bullet$-braided duoidal category $(\Cc,\circ,I,\bullet,J)$ with a reversion such that $\Cc$ has binary intersections and $\bullet$ preserves them, the category of cocommutative Hopf monoids whose unit morphism $u^{\circ}_{H}:I\to H$ is a monomorphism in $\Cc$ is protomodular (\cref{thm:protomon}).

In \cref{sec:exapp}, we consider examples and applications of the results achieved. First, considering any braided monoidal category as a duoidal category with reversion given by the identity functor, we prove that Hopf monoids in the duoidal setting coincide with standard Hopf monoids (\cref{pro:Hopfmonbraidedcat}), and we recover the protomodularity of the category of cocommutative Hopf monoids in any symmetric monoidal category $\Cc$ that has binary intersections preserved by the monoidal structure. In particular, we recover the protomodularity of the category of cocommutative Hopf algebras. Next, we extend the notion of internal group to an
arbitrary monoidal category $(\Cc,\circ,I)$, where $\Cc$ is endowed with the binary product $\times$ and a  terminal object $\mathsf{1}$. This is obtained by defining an internal group in $(\Cc,\circ,I)$ as a Hopf monoid in the duoidal category $(\Cc,\circ,I,\times,\mathsf{1})$, provided it has a reversion. As a consequence, we deduce that the category of internal groups whose unit is a monomorphism is protomodular (\cref{coro:Grproto}). 
Given a Hopf algebra $B$, the main example for us is the duoidal category of $B$-bimodules with $\circ=\ot_{B}$, $I=B$, $\bullet=\ot_{\Bbbk}$, $J=\Bbbk$. Equipped this duoidal category with a reversion, we prove that the category of Hopf monoids therein coincides with the coslice category $B/\mathsf{Hopf}_{\Bbbk}$ under $B$ (\cref{pro:Bohm}). If $(B,\Ss)$ is a quasitriangular Hopf algebra, then the category of cocommutative Hopf monoids is isomorphic to the coslice category of quasitriangular Hopf algebras under $(B,\Ss)$ (\cref{cor:quasiascocommutativeHopf}). 
Consequently, its full subcategory whose objects $((H,\Rr),u_{H}^{\circ})$ have injective  unit map $u_{H}^{\circ}:(B,\Ss)\to(H,\Rr)$ is protomodular. The same can be deduced in case $(B,\Ss)$ is even triangular, in which case the coslice category coincides with the category of internal groups in the category of cocommutative comonoids (\cref{coro:trHopfasgroups}) in our generalized sense. 
In case $(B,\Ss)$ is minimal quasitriangular, one obtains the protomodularity of the category of quasitriangular Hopf algebras whose minimal quasitriangular Hopf subalgebra is isomorphic to $(B,\Ss)$, with morphisms preserving minimal quasitriangular Hopf subalgebras. This can be interpreted as protomodularity of any essential fibre of a functor (\cref{thm:protomodularityforfibre}). 

We conclude by pointing out that, given a surjective Hopf algebra map $\pi:A\to H$ with $\pi_{\mid H}=\id$ where $(H,\Rr)$ is a quasitriangular Hopf subalgebra of $(A,\Rr)$, the kernel $K$ turns out to be the product of the space of coinvariant elements $A^{\mathrm{co}H}$ and the minimal quasitriangular Hopf algebra $B:=H_{\Rr}$ and $A$ turns out to be isomorphic to $K\ot_{B}H$ whose Hopf algebra structure generalizes the double cross product attached to a matched pair of Hopf algebras (\cref{thm:crospt}).

\section{Results on duoidal categories}
\label{sec:duoidal}

First, we recall the definition of duoidal category.

\begin{definition}[{cf.\ \cite[Definition 6.1]{Aguiar}}]
    A \emph{duoidal category} is a category $\Cc$ equipped with two monoidal structures $\Cc^\circ= (\Cc, \circ, I, a^\circ, l^\circ, r^\circ )$ and $\Cc^\bullet = (\Cc, \bullet, J, a^\bullet, l^\bullet, r^\bullet )$,  related through a natural transformation
\begin{equation}
\label{def:zetaduoidal}
\zeta=(\zeta_{A,B,C,D}: (A\bullet B)\circ (C\bullet D)\to (A\circ C)\bullet (B\circ D))_{A,B,C,D\in\Cc},
\end{equation}
called the \emph{interchange law}, and morphisms
\begin{equation}\label{def:deltamepsilonduoidal}
\Delta^{\bullet}_I:I\to I\bullet I, \quad m^{\circ}_J: J\circ J\to J, \quad \varepsilon^{\bullet}_I=u^{\circ}_{J}:I\to J
\end{equation}
such that
\begin{itemize}
    \item[$i)$] $(J, m^\circ_J, u^\circ_J)$ is a monoid in $\Cc^\circ$ and $(I, \Delta^\bullet_I, \varepsilon^\bullet_I)$ is a comonoid in $\Cc^\bullet$;
    \item[$ii)$] the following associativity conditions hold, for all $A,B,C,D,E,F\in\Cc$,
    \begin{gather}\label{eq:assoc1}
    \begin{split}
\zeta_{A,B, C\circ E, D\circ F}(\id_{A\bullet B}&\circ\zeta_{C,D,E,F})a^\circ_{A\bullet B, C\bullet D, E\bullet F}\\&=(a^\circ_{A,C,E}\bullet a^\circ_{B,D,F})\zeta_{A\circ C, B\circ D,E,F}(\zeta_{A,B,C,D}\circ\id_{E\bullet F})
 \end{split} \\\label{eq:assoc2}
    \begin{split}
(\id_{A\circ D}\bullet \zeta_{B,C,E,F})&\zeta_{A,B\bullet C, D, E\bullet F}(a^\bullet_{A,B,C}\circ a^\bullet_{D,E,F})\\&=a^\bullet_{A\circ D, B\circ E, C\circ F}(\zeta_{A,B,D, E}\bullet \id_{C\circ F})\zeta_{A\bullet B, C, D\bullet E, F};
    \end{split}
 \end{gather}
 \item[$iii)$] the following unitality conditions hold, for all $A,B\in\Cc$
 \begin{gather}\label{eq:unit1}
         (l^\circ_{A}\bullet l^\circ_B)\zeta_{I,I, A,B}(\Delta^{\bullet}_I\circ\id_{A\bullet B})=l^\circ_{A\bullet B};\quad (r^\circ_A\bullet r^\circ_B)\zeta_{A,B,I,I}(\id_{A\bullet B}\circ\Delta^{\bullet}_I)=r^\circ_{A\bullet B};\\\label{eq:unit2}
    l^\bullet_{A\circ B}(m^{\circ}_J\bullet \id_{A\circ B})\zeta_{J,A,J,B}=l^\bullet_A\circ l^\bullet_B;\quad r^\bullet_{A\circ B}(\id_{A\circ B}\bullet m^{\circ}_J)\zeta_{A,J,B,J}=r^\bullet_A\circ r^\bullet_B.
 \end{gather}
\end{itemize}
\end{definition}
A duoidal category will be denoted by $(\Cc,\circ, I,\bullet, J)$ or $(\Cc,\circ,\bullet)$. We write $\circ ^{\rev }$ for the reverse of the tensor product $\circ$. As observed in \cite[\S 4.1]{Bohm-Lack}, see also \cite[\S 4.3]{Street}, for any duoidal category $(\mathcal{C}, \circ, \bullet)$, we can reverse either or both of the monoidal structures, obtaining duoidal categories $(\mathcal{C}, \circ^\rev, \bullet)$, $(\mathcal{C}, \circ, \bullet^\rev)$,  $(\mathcal{C}, \circ^\rev, \bullet^\rev)$.

It has been proved in \cite{Batanin-Markl} that every duoidal category is duoidal equivalent to a strict duoidal category i.e.\ one where both monoidal categories are strict. For this reason, in computations we will usually omit the associativity constraints.

\begin{remark}
Let $(\mathcal{C}, \circ, I, \bullet, J)$ be a duoidal category. In view of \cite[Proposition 6.36]{Aguiar}, there are canonical equivalences of categories
\begin{equation}\label{equivalnecesBimonC}
\Bimon (\mathcal{C},\circ,\bullet)\cong\Comon (\mathsf{Mon}(\mathcal{C}^\circ), \bullet)\cong \mathsf{Mon}(\Comon (\mathcal{C}^\bullet), \circ).
\end{equation}
%In the following, we will denote $\Bimon (\mathcal{C},\circ,\bullet)$ simply by $\Bimon (\mathcal{C})$.
\end{remark}

\subsection{Reversion}
In \cite{Bohm-Canada}, G. B\"ohm sketches a definition of reversion aimed at introducing a notion of antipode in the duoidal setting, without making explicit the compatibility conditions it is required to satisfy. 
Slightly modifying this definition for our
purposes, we introduce the following one:

\begin{definition}
A \emph{reversion }$\left( \left( -\right) ^{\revsn },\phi^\circ ,\phi^{\bullet},\gamma ,\delta,\widetilde{\gamma},\widetilde{\delta}
\right) $ on a duoidal category $\left( \mathcal{C},\circ ,I,\bullet
,J\right) $ consists of a double lax monoidal functor $\left( \left( -\right) ^{\revsn },\phi^\circ,\phi^{\bullet} \right):(\Cc,\circ^\rev,\bullet^\rev)\to (\Cc,\circ,\bullet)$, see \cite[Definition 6.54]{Aguiar}, constructed from a functor $\left( -\right) ^{\revsn }:\Cc\to \Cc$ that underlies two lax monoidal functors
\begin{eqnarray*}
\left( \left( -\right) ^{\revsn },\phi^\circ \right)
&:&\left( \mathcal{C},\circ^\rev \right) \rightarrow \left( \mathcal{C},\circ \right) ,\quad \phi^\circ _{X,Y}:X^{\revsn }\circ Y^{\revsn }\rightarrow
\left( Y\circ X\right) ^{\revsn },\quad \phi^\circ _{0}:I\rightarrow I^{\revsn }, \\
\left( \left( -\right) ^{\revsn },\phi^\bullet \right)
&:&\left( \mathcal{C},\bullet^\rev \right) \rightarrow \left( \mathcal{C},\bullet \right),\quad
\phi _{X,Y}^{\bullet } :X^{\revsn }\bullet Y^{\revsn }\rightarrow \left(
Y\bullet X\right) ^{\revsn },\quad
\phi _{0}^{\bullet } :J\rightarrow J^{\revsn },
\end{eqnarray*}
together with four natural transformations, defined for every $X,Y,Z,T$ in $\mathcal{C},$ 
\begin{gather*}
\xymatrix@C=1.3cm{
\bigl((X \bullet Y^{\revsn}) \circ Z\bigr) \bullet T
  \ar[r]^-{\gamma_{X,Y,Z,T}}
& X \circ \bigl(Z \bullet (Y \circ T)\bigr)
},\quad
\xymatrix@C=1.3cm{
X \bullet \bigl(Y \circ (Z^{\revsn} \bullet T)\bigr)
  \ar[r]^-{\delta_{X,Y,Z,T}}
& \bigl((X \circ Z) \bullet Y\bigr) \circ T
},\\
\xymatrix@C=1.3cm{
X^{\revsn} \bullet \bigl((Y \bullet Z) \circ T\bigr)
  \ar[r]^-{\widetilde{\gamma}_{X,Y,Z,T}}
& Z \circ \bigl((X \circ Y)^{\revsn} \bullet T\bigr)
},\quad \xymatrix@C=1.3cm{
\bigl(X \circ (Y \bullet Z)\bigr) \bullet T^{\revsn}
  \ar[r]^-{\widetilde{\delta}_{X,Y,Z,T}}
& \bigl(X \bullet (Z \circ T)^{\revsn}\bigr) \circ Y
}.
\end{gather*}
Given the previous data, we can define:
\begin{align*}
\psi _{X,Y} :& \xymatrix@C=1.5cm{
X^{\revsn}\bullet Y
  \ar[r]^-{(r_{X^{\revsn}}^{\circ})^{-1}\,\bullet\, Y}
& (X^{\revsn}\circ I)\bullet Y
  \ar[r]^-{((l_{X^{\revsn}}^{\bullet})^{-1}\circ I)\,\bullet\, Y}
& ((J\bullet X^{\revsn})\circ I)\bullet Y
  \ar[r]^-{\gamma_{J,X,I,Y}}
& J\circ\bigl(I\bullet (X\circ Y)\bigr)
}, \\
\varphi _{X,Y} :&\xymatrix@C=1.5cm{
X \bullet Y^{\revsn}
  \ar[r]^-{X \bullet (l_{Y^{\revsn}}^{\circ})^{-1}}
& X \bullet (I \circ Y^{\revsn})
  \ar[r]^-{X \bullet (I \circ (r_{Y^{\revsn}}^{\bullet})^{-1})}
& X \bullet (I \circ (Y^{\revsn} \bullet J))
  \ar[r]^-{\delta_{X,I,Y,J}}
& ((X \circ Y) \bullet I) \circ J
} , \\
\widetilde{\psi }_{X,Y} :&\xymatrix@C=1.5cm{
X^{\revsn} \bullet \bigl((Y \bullet I)\circ J\bigr)
  \ar[r]^-{\widetilde{\gamma}_{X,Y,I,J}}
& I \circ \bigl((X \circ Y)^{\revsn} \bullet J\bigr)
  \ar[r]^-{I \circ r_{(X\circ Y)^{\revsn}}^{\bullet}}
& I \circ (X \circ Y)^{\revsn}
  \ar[r]^-{l_{(X\circ Y)^{\revsn}}^{\circ}}
& (X \circ Y)^{\revsn}
}, \\
\widetilde{\varphi }_{X,Y} :&\xymatrix@C=1.5cm{
\bigl(J \circ (I \bullet X)\bigr) \bullet Y^{\revsn}
  \ar[r]^-{\widetilde{\delta}_{J,I,X,Y}}
& \bigl(J \bullet (X \circ Y)^{\revsn}\bigr) \circ I
  \ar[r]^-{\,l_{(X\circ Y)^{\revsn}}^{\bullet} \circ I\,}
& (X \circ Y)^{\revsn} \circ I
  \ar[r]^-{r_{(X\circ Y)^{\revsn}}^{\circ}}
& (X \circ Y)^{\revsn}}.
\end{align*}
These data are subject to the following axioms:%
\begin{equation}\label{form:psi1}
\begin{gathered}
\xymatrix@C=52pt@R=28pt{
  (((X\bullet Y^{\revsn}) \circ Z) \bullet T) \circ (U\bullet V)
    \ar[r]^{\zeta_{(X\bullet Y^{\revsn}) \circ Z,\,T,\,U,\,V}}
    \ar[d]_{\gamma_{X,Y,Z,T}\circ (U\bullet V)}
  &
  ((X\bullet Y^{\revsn}) \circ Z \circ U) \bullet (T\circ V)
    \ar[d]^{\gamma_{X,Y,Z\circ U,T\circ V}}
  \\
  X\circ\bigl((Z\bullet (Y\circ T)) \circ (U\bullet V)\bigr)
    \ar[r]_{X\circ \zeta_{Z,\,Y\circ T,\,U,\,V}}
  &
  X\circ\bigl((Z\circ U)\bullet (Y\circ T\circ V)\bigr)
}\end{gathered}
\end{equation}
\vspace{-.2cm}
\begin{equation}\label{form:psi2}
\begin{gathered}
\xymatrix@C=52pt@R=28pt{
  \bigl((X\bullet Y\bullet Z^{\revsn}) \circ (T\bullet U)\bigr)\bullet V
    \ar[d]_-{\gamma_{X\bullet Y,\,Z,\,T\bullet U,\,V}}
    \ar[r]^-{\zeta_{X,\,Y\bullet Z^{\revsn},\,T,\,U}\,\bullet\, V}
  &
  (X\circ T)\bullet\bigl((Y\bullet Z^{\revsn})\circ U\bigr)\bullet V
    \ar[d]^-{(X\circ T)\bullet \gamma_{Y,\,Z,\,U,\,V}}
  \\
  (X\bullet Y)\circ\bigl(T\bullet U\bullet (Z\circ V)\bigr)
    \ar[r]^-{\zeta_{X,\,Y,\,T,\,U\bullet (Z\circ V)}}
  &
  (X\circ T)\bullet\bigl(Y\circ (U\bullet (Z\circ V))\bigr)
} \end{gathered}
\end{equation}
\vspace{-.2cm}
\begin{equation}\label{form:varphi1}
\begin{gathered}
\xymatrix@C=52pt@R=28pt{
  (X\bullet Y)\circ\bigl(Z\bullet (T\circ (U^{\revsn}\bullet V))\bigr)
    \ar[r]^{\zeta_{X,\,Y,\,Z,\,T\circ (U^{\revsn}\bullet V)}}
    \ar[d]_{(X\bullet Y)\circ \delta_{Z,\,T,\,U,\,V}}
  &
  (X\circ Z)\bullet\bigl(Y\circ T\circ (U^{\revsn}\bullet V)\bigr)
    \ar[d]^{\delta_{X\circ Z,\,Y\circ T,\,U,\,V}}
  \\
  (X\bullet Y)\circ\bigl((Z\circ U)\bullet T\bigr)\circ V
    \ar[r]_{\zeta_{X,\,Y,\,Z\circ U,\,T}\circ V}
  &
  \bigl((X\circ Z\circ U)\bullet (Y\circ T)\bigr)\circ V
}
\end{gathered}
\end{equation}
\vspace{-.2cm}
\begin{equation}\label{form:varphi2}
\begin{gathered}
 \xymatrix@C=52pt@R=28pt{
  {X\bullet \left( \left( Y\bullet Z\right) \circ \left( T^{\revsn }\bullet
  U\bullet V\right) \right)}
    \ar[r]^{X\bullet \zeta_{Y,Z,T^{\revsn}\bullet U,V}}
    \ar[d]_{\delta_{X,Y\bullet Z,T,U\bullet V}}
  &
  {X\bullet \left( Y\circ \left( T^{\revsn }\bullet U\right) \right)
  \bullet \left( Z\circ V\right)}
    \ar[d]^{\delta_{X,Y,T,U}\bullet \left( Z\circ V\right)}
  \\
  {\left( \left( X\circ T\right) \bullet Y\bullet Z\right)
  \circ \left( U\bullet V\right)}
    \ar[r]_{\zeta_{\left( X\circ T\right) \bullet Y,Z,U,V}}
  &
  {\left( \left( \left( X\circ T\right) \bullet Y\right) \circ U\right)
  \bullet \left( Z\circ V\right)}
} \end{gathered}
\end{equation}
\vspace{-.2cm}
\begin{equation}\label{form:gammadelta}
    \begin{tikzcd}[column sep=huge, row sep=large]
X\bullet \left( (Y\bullet Z^{\revsn}) \circ (T^{\revsn}\bullet U) \right) \bullet V
  \arrow[r,
    "X\bullet \gamma_{Y,Z,T^{\revsn}\bullet U,V}"]
  \arrow[d,
    "\delta_{X,Y\bullet Z^{\revsn},T,U}\bullet V"']
&
X\bullet \left( Y\circ (T^{\revsn}\bullet U \bullet (Z\circ V)) \right)
  \arrow[d,
    "\delta_{X,Y,T,U\bullet (Z\circ V)}"]
\\
\left(\left( (X\circ T)\bullet Y \bullet Z^{\revsn} \right) \circ U\right) \bullet V
  \arrow[r,
    "\gamma_{(X\circ T)\bullet Y,Z,U,V}"']
&
\left( (X\circ T)\bullet Y \right) \circ (U\bullet (Z\circ V))
\end{tikzcd}
\end{equation}
\vspace{-.2cm}
\begin{equation}\label{form:phi0gamma}
\begin{tikzcd}
	{\left( I\circ I\right) \bullet I} & {I\bullet I} \\
	{\left( \left( I\bullet I^\revsn\right) \circ I\right) \bullet I} & {I\circ \left( I\bullet \left( I\circ I\right) \right)}
	\arrow["{ l_{I}^{\circ }\bullet I}", from=1-1, to=1-2]
	\arrow["{\left( \left( I\bullet \phi^\circ _{0}\right) \Delta _{I}^{\bullet }\circ I\right) \bullet I}"', from=1-1, to=2-1]
	\arrow["{\gamma _{I,I,I,I}}", from=2-1, to=2-2]
	\arrow["{l_{I\bullet I}^{\circ }\left( I\circ \left( I\bullet l_{I}^{\circ }\right) \right)}"', from=2-2, to=1-2]
\end{tikzcd}
\end{equation}%
\vspace{-.2cm}
\begin{equation}\label{form:phi0delta}
\begin{tikzcd}
	{I\bullet(I\circ I)} & {I\bullet I} \\
	{I\bullet(I\circ (I^\revsn\bullet I))} & {\left( \left( I\circ I\right) \bullet I\right) \circ I }
	\arrow["{I\bullet r_{I}^{\circ }}", from=1-1, to=1-2]
	\arrow["{I\bullet \left( I\circ \left( \phi^\circ _{0}\bullet I\right) \Delta _{I}^{\bullet }\right)}"', from=1-1, to=2-1]
	\arrow["{\delta _{I,I,I,I}}", from=2-1, to=2-2]
	\arrow["{r_{I\bullet I}^{\circ }\left(\left( r_{I}^{\circ }\bullet I\right) \circ I\right)}"', from=2-2, to=1-2]
\end{tikzcd}
\end{equation}
\vspace{-.2cm}
\begin{equation}\label{form:phipsi2}
\begin{tikzcd}[scale cd=.9]
	((X\bullet(Y^\revsn\circ Z^\revsn))\circ J)\bullet U & ((X\bullet(Y^{\revsn}\circ ({Z^{\revsn}}\bullet J)))\circ J)\bullet U & (((X\circ Z)\bullet Y^{\revsn})\circ J\circ J)\bullet U \\
	((X\bullet(Z\circ Y)^\revsn)\circ J)\bullet U && X\circ Z\circ ((J\circ J)\bullet (Y\circ  U)) \\
	X\circ (J\bullet (Z\circ Y\circ U)) & X\circ Z\circ Y\circ U & X\circ Z\circ (J\bullet Y)\circ (J\bullet   U)
	\arrow[from=1-1, to=1-2, "((X\bullet(Y^{\revsn}\circ (r^\bullet_{Z^{\revsn}})^{-1}))\circ J)\bullet U"{yshift=1ex}]
	\arrow[from=1-1, to=2-1, "((X\bullet\phi^\circ_{Y, Z})\circ J)\bullet U"']
	\arrow[from=1-2, to=1-3, "(\delta_{X,Y^{\revsn},Z,J}\circ J)\bullet U"{yshift=1ex}]
	\arrow[from=1-3, to=2-3, "\gamma_{X\circ Z,Y,J\circ J,U}"]
	\arrow[from=2-1, to=3-1, "\gamma_{X,Z\circ Y,J,U}"']
	\arrow[from=3-1, to=3-2, "X\circ l^\bullet_{Z\circ Y\circ U}"']
	\arrow[from=3-2, to=3-3, "X\circ Z\circ (l_Y^{\bullet})^{-1}\circ (l_U^{\bullet})^{-1}"']
	\arrow[from=3-3, to=2-3, "X\circ Z\circ\zeta_{J,Y,J,U}"']
\end{tikzcd}
\end{equation}
\vspace{-.2cm}
\begin{equation}\label{form:phivarphi2}
    \begin{tikzcd}[scale cd=.9]
	{X\bullet(J\circ((Y^{\revsn}\circ Z^{\revsn})\bullet U))} & {X\bullet(J\circ(((J\bullet Y^{\revsn})\circ Z^{\revsn})\bullet U))} & {X\bullet(J\circ J\circ (Z^{\revsn}\bullet (Y\circ U)))} \\
	{X\bullet(J\circ((Z\circ Y)^{\revsn}\bullet U))} && {((X\circ Z)\bullet (J\circ J))\circ Y\circ U} \\
	{((X\circ Z\circ Y)\bullet J)\circ U} & {X\circ Z\circ Y\circ U} & {(X\bullet J)\circ (Z\bullet J)\circ Y\circ U}
	\arrow["{X\bullet(J\circ ((l_{Y^{\revsn}}^{\bullet})^{-1}\circ Z^{\revsn})\bullet U))}"{yshift=1ex}, from=1-1, to=1-2]
	\arrow["{X\bullet(J\circ(\phi^\circ_{Y,Z}\bullet U))}"', from=1-1, to=2-1]
	\arrow["{X\bullet (J\circ \gamma_{J,Y,Z^{\revsn},U})}"{yshift=1ex}, from=1-2, to=1-3]
	\arrow["{\delta_{X,J\circ J,Z,Y\circ U}}", from=1-3, to=2-3]
	\arrow["{\delta_{X,J,Z\circ Y,U}}"', from=2-1, to=3-1]
	\arrow["{r^\bullet_{X\circ Z\circ Y}\circ U}"', from=3-1, to=3-2]
	\arrow["{(r_X^\bullet)^{-1}\circ (r_Z^\bullet)^{-1}\circ Y\circ U}"', from=3-2, to=3-3]
	\arrow["{\zeta_{X,J,Z,J}\circ Y\circ U}"', from=3-3, to=2-3]
\end{tikzcd}
\end{equation}
\vspace{-.2cm}
\begin{equation}\label{form:varphipsitilde}
  \begin{tikzcd}
	{X^{\revsn}\bullet Y\bullet Z^\revsn} & {(J\circ (I\bullet(X\circ Y)))\bullet Z^\revsn} \\
	{X^{\revsn}\bullet (((X\circ Y)\bullet I)\circ J)} & {(X\circ Y\circ Z)^{\revsn}}
	\arrow["{\psi_{X,Y}\bullet Z^{\revsn}}", from=1-1, to=1-2]
	\arrow["{X^{\revsn}\bullet \varphi_{Y,Z}}"', from=1-1, to=2-1]
	\arrow["{\widetilde{\varphi}_{X\circ Y,Z}}", from=1-2, to=2-2]
	\arrow["{\widetilde{\psi}_{X,Y\circ Z}}", from=2-1, to=2-2]
\end{tikzcd}
\end{equation}
\vspace{-.2cm}
\begin{equation}
\begin{tikzcd}
	{X^{\revsn }\bullet \left(I\circ J\right)} & {X^\revsn} \\
	{X^{\revsn }\bullet \left( \left( I\bullet I\right) \circ J\right)} & {I\circ \left( \left( X\circ I\right) ^{\revsn }\bullet J\right)}
	\arrow["{r_{X^{\revsn }}^{\bullet }\left( X^{\revsn }\bullet l_{J}^{\circ }\right) }", from=1-1, to=1-2]
	\arrow["{X^{\revsn }\bullet \left( \Delta _{I}^{\bullet }\circ J\right) }"', from=1-1, to=2-1]
	\arrow["{\widetilde{\gamma }_{X,I,I,J}}", from=2-1, to=2-2]
	\arrow["{\left( r_{X}^{\circ }\right) ^{\revsn }l_{\left( X\circ I\right) ^{\revsn }}^{\circ }(I\circ r_{\left( X\circ I\right) ^{\revsn }}^{\bullet })}"', from=2-2, to=1-2]
\end{tikzcd}  \label{form:gammatildeIIJ}
\end{equation}
\vspace{-.2cm}
\begin{equation}
\begin{tikzcd}
	{\left( J\circ I \right) \bullet X^{\revsn }} & {X^\revsn} \\
	{\left( J\circ \left( I\bullet I\right) \right) \bullet X^{\revsn }} & {\left( J\bullet \left( I\circ X\right) ^{\revsn }\right) \circ I}
	\arrow["{l_{X^{\revsn }}^{\bullet }\left( r_{J}^{\circ }\bullet X^{\revsn }\right)   }", from=1-1, to=1-2]
	\arrow["{\left( J\circ \Delta _{I}^{\bullet }\right) \bullet X^{\revsn }}"', from=1-1, to=2-1]
	\arrow["{\widetilde{\delta }_{J,I,I,X}}", from=2-1, to=2-2]
	\arrow["{\left( l_{X}^{\circ }\right) ^{\revsn }r_{\left( I\circ X\right) ^{\revsn }}^{\circ }\left( l_{\left( I\circ X\right) ^{\revsn }}^{\bullet }\circ I\right)}"', from=2-2, to=1-2]
\end{tikzcd}  \label{form:deltatildeJII}
\end{equation}
\vspace{-.2cm}
\begin{equation}
\label{form: psi1gamma}
\begin{tikzcd}
	{((X\bullet Y^\revsn\bullet Z^\revsn)\circ T)\bullet U\bullet V} && {((X\bullet Y^\revsn)\circ (T\bullet (Z\circ U)))\bullet V}\\ && {X\circ (T\bullet (Z\circ U)\bullet (Y\circ V))} \\
	{((X\bullet (Z\bullet Y)^\revsn)\circ T)\bullet U\bullet V} && {(X\circ (T\bullet ((Z\bullet Y)\circ (U\bullet V)))}
	\arrow["{\gamma_{X\bullet Y^\revsn,Z,T,U}\bullet V}", from=1-1, to=1-3]
	\arrow["{((X\bullet \phi_{Y,Z}^{\bullet})\circ T)\bullet U\bullet V}"', from=1-1, to=3-1]
	\arrow["{\gamma_{X,Y,T\bullet (Z\circ U),V}}", from=1-3, to=2-3]
	\arrow["{\gamma_{X,Z\bullet Y,T,U\bullet V}}"', from=3-1, to=3-3]
	\arrow["{X\circ(T\bullet \zeta_{Z,Y,U,V})}"', from=3-3, to=2-3]
\end{tikzcd}
\end{equation}
\vspace{-.2cm}
\begin{equation}\label{form:varphi1delta}
\begin{tikzcd}
	{X\bullet Y\bullet (Z\circ (T^{\revsn}\bullet U^{\revsn}\bullet V))} && {X\bullet Y\bullet (Z\circ ((U\bullet T)^{\revsn}\bullet V))} \\ && {((X\circ U)\bullet (Y\circ T)\bullet Z)\circ V} \\
	{X\bullet Y\bullet (Z\circ ((U\bullet T)^{\revsn}\bullet V))} &&  {(((X\bullet Y)\circ (U\bullet T))\circ Z)\bullet V}
	\arrow["{X\bullet\delta_{Y,Z,T,U^{\revsn}\bullet V}}", from=1-1, to=1-3]
	\arrow["{X\bullet Y\bullet (Z\circ (\phi_{T,U}^{\bullet }\bullet V))}"', from=1-1, to=3-1]
	\arrow["{\delta_{X,(Y\circ T)\bullet Z,U,V}}", from=1-3, to=2-3]
	\arrow["{\delta_{X\bullet Y,Z,U\bullet T,V}}"', from=3-1, to=3-3]
	\arrow["{(\zeta_{X,Y,U,T}\bullet Z)\circ V}"', from=3-3, to=2-3]
\end{tikzcd}
\end{equation}
\vspace{-.2cm}
\begin{equation}
\label{form:phipsi2NEW}
 \begin{tikzcd}
	{((X\bullet Y^\revsn)\circ (Z\bullet T^\revsn)\circ U)\bullet V} && {X\circ(((Z\bullet T^{\revsn})\circ U)\bullet(Y\circ V))}\\ && {X\circ Z\circ (U\bullet (T\circ Y\circ V))} \\
	{(((X\circ Z)\bullet (Y^{\revsn}\circ T^{\revsn}))\circ U)\bullet V} && {(((X\circ Z)\bullet (T\circ Y)^{\revsn})\circ U)\bullet V}
	\arrow["{{\gamma_{X, Y,(Z\bullet T^\revsn)\circ U,V}}}", from=1-1, to=1-3]
	\arrow["{(\zeta_{X,Y^{\revsn},Z,T^{\revsn}}\circ U)\bullet V}"', from=1-1, to=3-1]
	\arrow["{X\circ \gamma_{Z,T,U,Y\circ V}}", from=1-3, to=2-3]
	\arrow["{(((X\circ Z)\bullet \phi_{Y,T}^\circ)\circ U)\bullet V}"'{yshift=-1ex}, from=3-1, to=3-3]
	\arrow["{\gamma_{X\circ Z,T\circ Y,U,V}}"', from=3-3, to=2-3]
\end{tikzcd}
\end{equation}
\vspace{-.2cm}
\begin{equation}
\label{form:phivarphi2NEW}
    \begin{tikzcd}
	{X\bullet((Y\circ(Z^{\revsn}\bullet T)\circ(U^{\revsn}\bullet V))} && {((X\circ U)\bullet(Y\circ(Z^{\revsn}\bullet T))\circ V} \\
	&& {((X\circ U\circ Z)\bullet Y)\circ T\circ V} \\
	{X\bullet(Y\circ((Z^{\revsn}\circ  U^{\revsn})\bullet(T\circ V)))} && {X\bullet(Y\circ((U\circ Z)^{\revsn}\bullet(T\circ V)))}
	\arrow["{\delta_{X,Y\circ(Z^{\revsn}\bullet T),U,V}}", from=1-1, to=1-3]
	\arrow["{X\bullet(Y\circ\zeta_{Z^{\revsn}, T, U^{\revsn}, V})}"', from=1-1, to=3-1]
	\arrow["{\delta_{X\circ U,Y,Z,T}\circ V}", from=1-3, to=2-3]
	\arrow["{X\bullet(Y\circ(\phi^\circ_{Z,U}\bullet(T\circ V)))}"', from=3-1, to=3-3]
	\arrow["{\delta_{X,Y,U\circ Z,T\circ V}}"', from=3-3, to=2-3]
\end{tikzcd}
\end{equation}
\vspace{-.2cm}
\begin{equation}
\label{form:phi0gamma-bul}
  \begin{tikzcd}
	{(J\circ I)\bullet J} & J \\
	{((J\bullet J^\revsn)\circ I)\bullet J} & {J\circ (I\bullet (J\circ J))}
	\arrow["{r_J^\bullet (r^\circ_J\bullet J)}"{pos=0.6}, from=1-1, to=1-2]
	\arrow["{((J\bullet\phi^\bullet_0)\Delta_J^\bullet\circ I)\bullet J}"', from=1-1, to=2-1]
	\arrow["{\gamma_{J,J,I,J}}"', from=2-1, to=2-2]
	\arrow["{r^\circ_J(J\circ r^\bullet_I)(J\circ(I\bullet m^\circ_J))}"', from=2-2, to=1-2]
\end{tikzcd}
\end{equation}
\vspace{-.2cm}
\begin{equation}
\label{form:phi0delta-bul}
    \begin{tikzcd}
	{J\bullet (I\circ J)} & J \\
	{J\bullet(I\circ(J^\revsn\bullet J))} & {((J\circ J)\bullet I)\circ J}
	\arrow["{l_J^\bullet (J\bullet l^\circ_J)}"{pos=0.6}, from=1-1, to=1-2]
	\arrow["{J\bullet(I\circ(\phi_0^\bullet\bullet J)\Delta_J^\bullet)}"', from=1-1, to=2-1]
	\arrow["{\delta_{J,I,J,J}}"', from=2-1, to=2-2]
	\arrow["{l_J^\circ(l^\bullet_I\circ J)((m^\circ_J\bullet I)\circ J)}"', from=2-2, to=1-2]
\end{tikzcd}
\end{equation}
\end{definition}

\begin{remark}
Note that the conditions \eqref{form:varphi1}, \eqref{form:varphi2},
\eqref{form:gammadelta}, \eqref{form:phi0delta}, \eqref{form:phivarphi2},
\eqref{form:varphipsitilde},
\eqref{form:deltatildeJII},
\eqref{form:varphi1delta},
\eqref{form:phivarphi2NEW},
\eqref{form:phi0delta-bul}
are just
\eqref{form:psi1}, \eqref{form:psi2},
\eqref{form:gammadelta}, \eqref{form:phi0gamma}, \eqref{form:phipsi2},
\eqref{form:varphipsitilde},
\eqref{form:gammatildeIIJ},
\eqref{form: psi1gamma},
\eqref{form:phipsi2NEW},
\eqref{form:phi0gamma-bul}
written in $\Cc^\mathrm{rev}:=(\Cc,\circ^\mathrm{rev},\bullet^\mathrm{rev})$ once $\delta$ and $\gamma$ are exchanged. For this reason, we can prove some results from one side only, deducing their counterpart for reasons of symmetry. 
\end{remark}

Using $\gamma$ and $\delta$ of a reversion one can also define the following natural transformations that will be useful afterwords.
\begin{align*}
\psi _{X,Y,Z}^{\prime } :& \xymatrix@C=1.5cm{
X \bullet Y^\revsn \bullet Z
  \ar[r]^-{(r_{X\bullet Y^\revsn}^{\circ})^{-1} \bullet Z}
&
((X \bullet Y^\revsn)\circ I)\bullet Z
  \ar[r]^-{\gamma_{X,Y,I,Z}}
&
X \circ (I \bullet (Y \circ Z))
}, \\
\psi _{X,Y,Z}^{\prime \prime } :& \xymatrix@C=1.5cm{
((X \bullet Y^\revsn)\circ J)\bullet Z
  \ar[r]^-{\gamma_{X,Y,J,Z}}
&
X \circ (J \bullet (Y \circ Z))
  \ar[r]^-{X \circ l^{\bullet}_{Y\circ Z}}
&
X \circ Y \circ Z
}, \\
\varphi _{X,Y,Z}^{\prime } :& \xymatrix@C=1.5cm{
X \bullet Y^\revsn \bullet Z
  \ar[r]^-{X \bullet (l^{\circ}_{Y^\revsn \bullet Z})^{-1}}
&
X \bullet (I \circ (Y^\revsn \bullet Z))
  \ar[r]^-{\delta_{X,I,Y,Z}}
&
((X \circ Y)\bullet I)\circ Z
}, \\
\varphi _{X,Y,Z}^{\prime \prime } :& \xymatrix@C=1.5cm{
X \bullet (J \circ (Y^\revsn \bullet Z))
  \ar[r]^-{\delta_{X,J,Y,Z}}
&
((X \circ Y)\bullet J)\circ Z
  \ar[r]^-{r^{\bullet}_{X\circ Y}\circ Z}
&
X \circ Y \circ Z
}.
\end{align*}

\begin{remark}%[\rd{added:2026/04/21}]
B\"ohm's definition of reversion starts from a strong duoidal equivalence $(-)^\revsn$ together with two natural transformations $\gamma$ and $\widetilde{\gamma}$.
For our purposes, we realized that it is sufficient to require that $(-)^\revsn$ is a double lax monoidal functor, at the cost of adding two additional natural transformations $\delta$ and $\widetilde{\delta}$ which, in the case of an equivalence, can be constructed from the other two.
The compatibility conditions we impose are tailored specifically to secure the results needed for protomodularity. They are therefore not meant to be either optimal or exhaustive. For instance, the diagrams \eqref{form:phipsi2} and $\eqref{form:phipsi2NEW}$ are similar, suggesting that some axioms may follow from others or be replaced by simpler ones.

We also point out that one could replace $\widetilde{\gamma}$ by 
$\widetilde{\gamma}_{X,Y,Z,T}:
X \bullet \bigl((Y \bullet Z) \circ T\bigr)
  \to Z \circ \bigl((Y^{\revsn} \circ X) \bullet T\bigr)$ and recover the former as the composition $(Z \circ \bigl(\phi^\circ_{Y,X} \bullet T\bigr))\widetilde{\gamma}_{X^\revsn,Y,Z,T}$. This would make the new $\widetilde{\gamma}$ more similar to $\gamma$, although it is not clear whether this would actually simplify the axioms.
\end{remark}

The following result is an analogue of \cite[Lemma 4.2]{Bohm-Lack}.

\begin{lemma}%[\rd{added:2026/03/06}]
In a duoidal category $(\Cc,\circ,I,\bullet,J)$ with a reversion, the following diagrams commute.
\begin{equation}\label{form:BohmLem4.2}
\begin{tikzcd}[scale cd=.98]
	\left( X^{\revsn }\bullet Y\right) \circ \left( Z^{\revsn }\bullet T\right) & J\circ \left( I\bullet \left( X\circ Y\right) \right) \circ \left(
Z^{\revsn }\bullet T\right)  & J\circ \left( \left( I\circ Z^{\revsn }\right) \bullet
\left( X\circ Y\circ T\right) \right) \\
	\left( X^{\revsn }\circ Z^{\revsn }\right) \bullet \left( Y\circ T\right) && J\circ
\left( Z^{\revsn }\bullet \left( X\circ Y\circ T\right) \right) \\
	\left( Z\circ X\right) ^{\revsn }\bullet \left( Y\circ T\right) & J\circ \left( I\bullet \left(
Z\circ X\circ Y\circ T\right) \right) & J\circ J\circ \left( I\bullet \left( Z\circ X\circ Y\circ T\right) \right)
	\arrow[from=1-1, to=1-2, "\psi _{X,Y}\circ \left( Z^{\revsn }\bullet T\right)"{yshift=1ex}]
	\arrow[from=1-1, to=2-1,"\zeta _{X^{\revsn },Y,Z^{\revsn },T}"']
	\arrow[from=1-2, to=1-3,"J\circ \zeta _{I,X\circ Y,Z^{\revsn },T}"{yshift=1ex}]
	\arrow[from=1-3, to=2-3,"J\circ
\left( l_{Z^{\revsn }}^{\circ }\bullet \left( X\circ Y\circ T\right) \right)"]
	\arrow[from=2-1, to=3-1,"\phi _{X,Z}^{\circ }\bullet \left( Y\circ T\right)"']
	\arrow[from=2-3, to=3-3,"J\circ \psi _{Z,X\circ Y\circ T}"]
	\arrow[from=3-1, to=3-2,"\psi _{Z\circ X,Y\circ T}"']
	\arrow[from=3-3, to=3-2,"m_{J}^{\circ }\circ \left(
I\bullet \left( Z\circ X\circ Y\circ T\right) \right)"{yshift=-1ex}]
\end{tikzcd}
\end{equation}
% \begin{equation}
% \begin{array}{ccccc}
% \left( X^{\revsn }\bullet Y\right) \circ \left( Z^{\revsn }\bullet T\right)  &
% \overset{\psi _{X,Y}\circ \left( Z^{\revsn }\bullet T\right) }{\longrightarrow
% } & J\circ \left( I\bullet \left( X\circ Y\right) \right) \circ \left(
% Z^{\revsn }\bullet T\right)  & \overset{J\circ \zeta _{I,X\circ Y,Z^{\revsn },T}}%
% {\longrightarrow } & J\circ \left( \left( I\circ Z^{\revsn }\right) \bullet
% \left( X\circ Y\circ T\right) \right)  \\
% \zeta _{X^{\revsn },Y,Z^{\revsn },T}\downarrow  &  &  &  & \downarrow J\circ
% \left( l_{Z^{\revsn }}^{\circ }\bullet \left( X\circ Y\circ T\right) \right)
% \\
% \left( X^{\revsn }\circ Z^{\revsn }\right) \bullet \left( Y\circ T\right)  &  &
% &  & J\circ \left( Z^{\revsn }\bullet \left( X\circ Y\circ T\right) \right)
% \\
% \phi _{X,Z}^{\circ }\bullet \left( Y\circ T\right) \downarrow  &  &  &  &
% \downarrow J\circ \psi _{Z,X\circ Y\circ T} \\
% \left( Z\circ X\right) ^{\revsn }\bullet \left( Y\circ T\right)  & \overset{%
% \psi _{Z\circ X,Y\circ T}}{\longrightarrow } & J\circ \left( I\bullet \left(
% Z\circ X\circ Y\circ T\right) \right)  & \overset{m_{J}^{\circ }\circ \left(
% I\bullet \left( Z\circ X\circ Y\circ T\right) \right) }{\longleftarrow } &
% J\circ J\circ \left( I\bullet \left( Z\circ X\circ Y\circ T\right) \right)
% \end{array}
%\end{equation}
\begin{equation}\label{form:BohmLem4.2B}
\begin{tikzcd}[scale cd=.98]
	\left( X\bullet Y^{\revsn}\right) \circ \left( Z\bullet T^{\revsn }\right) & \left( X\bullet Y^{\revsn }\right) \circ \left( \left(
Z\circ T\right) \bullet I\right) \circ J & \left( \left( X\circ Z\circ
T\right) \bullet \left( Y^{\revsn }\circ I\right) \right) \circ J \\
	\left( X\circ Z\right) \bullet \left( Y^{\revsn }\circ T^{\revsn }\right) && \left( \left( X\circ Z\circ T\right) \bullet Y^{\revsn }\right) \circ J \\
	\left( X\circ Z\right) \bullet \left( T\circ Y\right)^{\revsn } & \left( \left( X\circ
Z\circ T\circ Y\right) \bullet I\right) \circ J & \left( \left( X\circ Z\circ T\circ Y\right) \bullet
I\right) \circ J\circ J
	\arrow[from=1-1, to=1-2, "\left( X\bullet Y^{\revsn }\right) \circ \varphi _{Z,T}"{yshift=1ex}]
	\arrow[from=1-1, to=2-1,"\zeta _{X,Y^{\revsn },Z,T^{\revsn }}"']
	\arrow[from=1-2, to=1-3,"\zeta _{X,Y^{\revsn
},Z\circ T,I}\circ J"{yshift=1ex}]
	\arrow[from=1-3, to=2-3,"\left(
\left( X\circ Z\circ T\right) \bullet r_{Y^{\revsn }}^{\circ }\right) \circ J"]
	\arrow[from=2-1, to=3-1,"\left( X\circ Z\right)\bullet\phi _{Y,T}^{\circ }"']
	\arrow[from=2-3, to=3-3,"\varphi _{X\circ Z\circ T,Y}\circ J"]
	\arrow[from=3-1, to=3-2,"\varphi _{X\circ Z,T\circ Y}"']
	\arrow[from=3-3, to=3-2,"\left(
\left( X\circ Z\circ T\circ Y\right)\bullet I\right)\circ m_{J}^{\circ }"{yshift=-1ex}]
\end{tikzcd}
\end{equation}
\end{lemma}

\begin{proof}
Let us prove that \eqref{form:BohmLem4.2} commutes:
\begin{eqnarray*}
&&\left( m_{J}^{\circ }\circ \left( I\bullet \left( Z\circ X\circ Y\circ
T\right) \right) \right) \left( J\circ \underbracket[0.140ex]{\psi _{Z,X\circ Y\circ T}}\right)
\left( J\circ \left( l_{Z^{\revsn }}^{\circ }\bullet \left( X\circ Y\circ
T\right) \right) \right) \left( J\circ \zeta _{I,X\circ Y,Z^{\revsn
},T}\right) \left( \underbracket[0.140ex]{\psi _{X,Y}}\circ \left( Z^{\revsn }\bullet T\right) \right)
\\
&=&
\left( m_{J}^{\circ }\circ \left( I\bullet \left( Z\circ X\circ Y\circ
T\right) \right) \right) \left( J\circ \gamma _{J,Z,I,X\circ Y\circ
T}\right) \left( J\circ \left( \left( \left( l_{Z^{\revsn }}^{\bullet }\right)
^{-1}\circ I\right) \left( r_{Z^{\revsn }}^{\circ }\right) ^{-1}l_{Z^{\revsn
}}^{\circ }\bullet \left( X\circ Y\circ T\right) \right) \right)  \\
&&\underbracket[0.140ex]{\left( J\circ \zeta _{I,X\circ Y,Z^{\revsn },T}\right)
\left( \gamma _{J,X,I,Y}\circ \left( Z^{\revsn }\bullet T\right) \right) }%
\left( \left( \left( \left( l_{X^{\revsn }}^{\bullet }\right) ^{-1}\circ
I\right) \left( r_{X^{\revsn }}^{\circ }\right) ^{-1}\bullet Y\right) \circ
\left( Z^{\revsn }\bullet T\right) \right)\\
&\overset{\eqref{form:psi1}}{=}&
\left( m_{J}^{\circ }\circ \left( I\bullet \left( Z\circ X\circ Y\circ
T\right) \right) \right) \left( J\circ \gamma _{J,Z,I,X\circ Y\circ
T}\right) \left( J\circ \left( \left( \left( l_{Z^{\revsn }}^{\bullet }\right)
^{-1}\circ I\right) \left( r_{Z^{\revsn }}^{\circ }\right) ^{-1}l_{Z^{\revsn
}}^{\circ }\bullet \left( X\circ Y\circ T\right) \right) \right)  \\
&&\gamma _{J,X,I\circ Z^{\revsn },Y\circ T}\zeta _{\left( J\bullet X^{\revsn
}\right) \circ I,Y,Z^{\revsn },T}\left( \left( \left( \left( l_{X^{\revsn
}}^{\bullet }\right) ^{-1}\circ I\right) \left( r_{X^{\revsn }}^{\circ
}\right) ^{-1}\bullet Y\right) \circ \left( Z^{\revsn }\bullet T\right)
\right)\\
&=&
\left( m_{J}^{\circ }\circ \left( I\bullet \left( Z\circ X\circ Y\circ
T\right) \right) \right) \left( J\circ \gamma _{J,Z,I,X\circ Y\circ
T}\right) \gamma _{J,X,\left( J\bullet Z^{\revsn }\right) \circ I,Y\circ T} \\
&&\left( \left( \left( J\bullet X^{\revsn }\right) \circ \left( \left(
l_{Z^{\revsn }}^{\bullet }\right) ^{-1}\circ I\right) \left( r_{Z^{\revsn
}}^{\circ }\right) ^{-1}l_{Z^{\revsn }}^{\circ }\right) \bullet \left( Y\circ
T\right) \right) \left( \left( \left( \left( l_{X^{\revsn }}^{\bullet }\right)
^{-1}\circ I\right) \left( r_{X^{\revsn }}^{\circ }\right) ^{-1}\circ Z^{\revsn
}\right) \bullet \left( Y\circ T\right) \right) \zeta _{X^{\revsn },Y,Z^{\revsn
},T}\\
&=&
\left( m_{J}^{\circ }\circ \left( I\bullet \left( Z\circ X\circ Y\circ
T\right) \right) \right) \underbracket[0.140ex]{\left( J\circ \gamma
_{J,Z,I,X\circ \left( Y\circ T\right) }\right) \gamma _{J,X,\left( J\bullet
Z^{\revsn }\right) \circ I,Y\circ T}} \\
&&\left( \left( \left( l_{X^{\revsn }}^{\bullet }\right) ^{-1}\circ \left(
l_{Z^{\revsn }}^{\bullet }\right) ^{-1}\circ I\right) \left( r_{X^{\revsn }\circ
Z^{\revsn }}^{\circ }\right) ^{-1}\bullet \left( Y\circ T\right) \right) \zeta
_{X^{\revsn },Y,Z^{\revsn },T}\\
&\overset{\eqref{form:phipsi2NEW}}{=}&
\left( m_{J}^{\circ }\circ \left( I\bullet \left( Z\circ X\circ Y\circ
T\right) \right) \right) \gamma _{J\circ J,Z\circ X,I,Y\circ T}\left( \left(
\left( \left( J\circ J\right) \bullet \phi _{X,Z}^{\circ }\right) \circ
I\right) \bullet \left( Y\circ T\right) \right) \\
&&\left( \left( \zeta
_{J,X^{\revsn },J,Z^{\revsn }}\circ I\right) \bullet \left( Y\circ T\right)
\right)  \left( \left( \left( l_{X^{\revsn }}^{\bullet }\right) ^{-1}\circ \left(
l_{Z^{\revsn }}^{\bullet }\right) ^{-1}\circ I\right) \bullet \left( Y\circ
T\right) \right) \left( \left( r_{X^{\revsn }\circ Z^{\revsn }}^{\circ }\right)
^{-1}\bullet \left( Y\circ T\right) \right) \zeta _{X^{\revsn },Y,Z^{\revsn },T} \\
&=&
\gamma _{J,Z\circ X,I,Y\circ T}\left( \left( \left( m_{J}^{\circ }\bullet
\phi _{X,Z}^{\circ }\right) \circ I\right) \bullet \left( Y\circ T\right)
\right) \left( \left( \zeta _{J,X^{\revsn },J,Z^{\revsn }}\circ I\right) \bullet
\left( Y\circ T\right) \right)  \\
&&\left( \left( \left( l_{X^{\revsn }}^{\bullet }\right) ^{-1}\circ \left(
l_{Z^{\revsn }}^{\bullet }\right) ^{-1}\circ I\right) \bullet \left( Y\circ
T\right) \right) \left( \left( r_{X^{\revsn }\circ Z^{\revsn }}^{\circ }\right)
^{-1}\bullet \left( Y\circ T\right) \right) \zeta _{X^{\revsn },Y,Z^{\revsn },T}  \\
&=&\gamma _{J,Z\circ X,I,Y\circ T}\left( \left( \left( J\bullet \phi
_{X,Z}^{\circ }\right) \circ I\right) \bullet \left( Y\circ T\right) \right)
\left( \left( \left( l_{X^{\revsn }\circ Z^{\revsn }}^{\bullet }\right)
^{-1}\circ I\right) \bullet \left( Y\circ T\right) \right) \left( \left(
r_{X^{\revsn }\circ Z^{\revsn }}^{\circ }\right) ^{-1}\bullet \left( Y\circ
T\right) \right) \zeta _{X^{\revsn },Y,Z^{\revsn },T} \\
&=&\underbracket[0.140ex]{\gamma _{J,Z\circ X,I,Y\circ T}\left( \left( \left( l_{\left( Z\circ
X\right) ^{\revsn }}^{\bullet }\right) ^{-1}\circ I\right) \bullet \left(
Y\circ T\right) \right) \left( \left( r_{\left( Z\circ X\right) ^{\revsn
}}^{\circ }\right) ^{-1}\bullet \left( Y\circ T\right) \right)} \left( \phi
_{X,Z}^{\circ }\bullet \left( Y\circ T\right) \right) \zeta _{X^{\revsn
},Y,Z^{\revsn },T} \\
&=&\psi _{Z\circ X,Y\circ T}\left( \phi _{X,Z}^{\circ }\bullet \left( Y\circ
T\right) \right) \zeta _{X^{\revsn },Y,Z^{\revsn },T},
\end{eqnarray*}
i.e. \eqref{form:BohmLem4.2} commutes. The commutativity of \eqref{form:BohmLem4.2B} follows in a similar way by means of \eqref{form:varphi1} and \eqref{form:phivarphi2NEW}.
\begin{invisible}
We compute
\begin{eqnarray*}
&&\left( \left( \left( X\circ Z\circ T\circ Y\right) \bullet I\right) \circ
m_{J}^{\circ }\right) \left( \varphi _{X\circ Z\circ T,Y}\circ J\right)
\left( \left( \left( X\circ Z\circ T\right) \bullet r_{Y^{\revsn }}^{\circ
}\right) \circ J\right) \left( \zeta _{X,Y^{\revsn },Z\circ T,I}\circ J\right)
\left( \left( X\bullet Y^{\revsn }\right) \circ \varphi _{Z,T}\right)  \\
&=&\left[
\begin{array}{c}
\left( \left( \left( X\circ Z\circ T\circ Y\right) \bullet I\right) \circ
m_{J}^{\circ }\right) \left( \delta _{X\circ Z\circ T,I,Y,J}\circ J\right)
\left( \left( \left( X\circ Z\circ T\right) \bullet \left( I\circ \left(
r_{Y^{\revsn }}^{\bullet }\right) ^{-1}\right) \left( l_{Y^{\revsn }}^{\circ
}\right) ^{-1}r_{Y^{\revsn }}^{\circ }\right) \circ J\right)  \\
\underbracket[0.140ex]{\left( \zeta _{X,Y^{\revsn },Z\circ T,I}\circ J\right) \left(
\left( X\bullet Y^{\revsn }\right) \circ \delta _{Z,I,T,J}\right) }\left(
\left( X\bullet Y^{\revsn }\right) \circ \left( Z\bullet \left( I\circ \left(
r_{T^{\revsn }}^{\bullet }\right) ^{-1}\right) \left( l_{T^{\revsn }}^{\circ
}\right) ^{-1}\right) \right)
\end{array}%
\right]  \\
&\overset{\eqref{form:varphi1}}{=}&\left[
\begin{array}{c}
\left( \left( \left( X\circ Z\circ T\circ Y\right) \bullet I\right) \circ
m_{J}^{\circ }\right) \left( \delta _{X\circ Z\circ T,I,Y,J}\circ J\right)
\left( \left( \left( X\circ Z\circ T\right) \bullet \left( I\circ \left(
r_{Y^{\revsn }}^{\bullet }\right) ^{-1}\right) \left( l_{Y^{\revsn }}^{\circ
}\right) ^{-1}r_{Y^{\revsn }}^{\circ }\right) \circ J\right)  \\
\delta _{X\circ Z,Y^{\revsn }\circ I,T,J}\zeta _{X,Y^{\revsn },Z,I\circ \left(
T^{\revsn }\bullet J\right) }\left( \left( X\bullet Y^{\revsn }\right) \circ
\left( Z\bullet \left( I\circ \left( r_{T^{\revsn }}^{\bullet }\right)
^{-1}\right) \left( l_{T^{\revsn }}^{\circ }\right) ^{-1}\right) \right)
\end{array}%
\right]  \\
&=&\left[
\begin{array}{c}
\left( \left( \left( X\circ Z\circ T\circ Y\right) \bullet I\right) \circ
m_{J}^{\circ }\right) \left( \delta _{X\circ Z\circ T,I,Y,J}\circ J\right)
\delta _{X\circ Z,I\circ \left( Y^{\revsn }\bullet J\right) ,T,J} \\
\left( \left( X\circ Z\right) \bullet \left( \left( I\circ \left( r_{Y^{\revsn
}}^{\bullet }\right) ^{-1}\right) \left( l_{Y^{\revsn }}^{\circ }\right)
^{-1}r_{Y^{\revsn }}^{\circ }\circ \left( T^{\revsn }\bullet J\right) \right)
\left( Y^{\revsn }\circ \left( I\circ \left( r_{T^{\revsn }}^{\bullet }\right)
^{-1}\right) \left( l_{T^{\revsn }}^{\circ }\right) ^{-1}\right) \right) \zeta
_{X,Y^{\revsn },Z,T^{\revsn }}%
\end{array}%
\right]  \\
&=&\left[
\begin{array}{c}
\left( \left( \left( X\circ Z\circ T\circ Y\right) \bullet I\right) \circ
m_{J}^{\circ }\right) \underbracket[0.140ex]{\left( \delta _{X\circ Z\circ
T,I,Y,J}\circ J\right) \delta _{X\circ Z,I\circ \left( Y^{\revsn }\bullet
J\right) ,T,J}} \\
\left( \left( X\circ Z\right) \bullet \left( \left( I\circ \left( r_{Y^{\revsn
}}^{\bullet }\right) ^{-1}\circ \left( r_{T^{\revsn }}^{\bullet }\right)
^{-1}\right) \right) \left( \left( l_{Y^{\revsn }\circ T^{\revsn }}^{\circ
}\right) ^{-1}\right) \right) \zeta _{X,Y^{\revsn },Z,T^{\revsn }}%
\end{array}%
\right]  \\
&\overset{\eqref{form:phivarphi2NEW}}{=}&\left[
\begin{array}{c}
\left( \left( \left( X\circ Z\circ T\circ Y\right) \bullet I\right) \circ
m_{J}^{\circ }\right) \delta _{X\circ Z,I,T\circ Y,J\circ J}\left( \left(
X\circ Z\right) \bullet \left( I\circ \left( \phi _{Y,T}^{\circ }\bullet
\left( J\circ J\right) \right) \right) \right) \left( \left( X\circ Z\right)
\bullet \left( I\circ \zeta _{Y^{\revsn },J,T^{\revsn },J}\right) \right)  \\
\left( \left( X\circ Z\right) \bullet \left( \left( I\circ \left( r_{Y^{\revsn
}}^{\bullet }\right) ^{-1}\circ \left( r_{T^{\revsn }}^{\bullet }\right)
^{-1}\right) \right) \left( \left( l_{Y^{\revsn }\circ T^{\revsn }}^{\circ
}\right) ^{-1}\right) \right) \zeta _{X,Y^{\revsn },Z,T^{\revsn }}%
\end{array}%
\right]  \\
&=&\left[
\begin{array}{c}
\delta _{X\circ Z,I,T\circ Y,J}\left( \left( X\circ Z\right) \bullet \left(
I\circ \left( \phi _{Y,T}^{\circ }\bullet J\right) \right) \right)  \\
\left( \left( X\circ Z\right) \bullet \left( \left( I\circ \underbracket[0.140ex]{%
\left( \left( \left( Y^{\revsn }\circ T^{\revsn }\right) \bullet m_{J}^{\circ
}\right) \right) \zeta _{Y^{\revsn },J,T^{\revsn },J}\left( \left( r_{Y^{\revsn
}}^{\bullet }\right) ^{-1}\circ \left( r_{T^{\revsn }}^{\bullet }\right)
^{-1}\right) }\right) \right) \left( \left( l_{Y^{\revsn }\circ T^{\revsn
}}^{\circ }\right) ^{-1}\right) \right) \zeta _{X,Y^{\revsn },Z,T^{\revsn }}%
\end{array}%
\right]  \\
&=&\delta _{X\circ Z,I,T\circ Y,J}\left( \left( X\circ Z\right) \bullet
\left( I\circ \left( \phi _{Y,T}^{\circ }\bullet J\right) \right) \right)
\left( \left( X\circ Z\right) \bullet \left( \left( I\circ \left( r_{Y^{\revsn
}\circ T^{\revsn }}^{\bullet }\right) ^{-1}\right) \right) \left( \left(
l_{Y^{\revsn }\circ T^{\revsn }}^{\circ }\right) ^{-1}\right) \right) \zeta
_{X,Y^{\revsn },Z,T^{\revsn }} \\
&=&\delta _{X\circ Z,I,T\circ Y,J}\left( \left( X\circ Z\right) \bullet
\left( \left( I\circ \left( r_{\left( T\circ Y\right) ^{\revsn }}^{\bullet
}\right) ^{-1}\right) \right) \left( \left( l_{\left( T\circ Y\right) ^{\revsn
}}^{\circ }\right) ^{-1}\right) \right) \left( \left( X\circ Z\right)
\bullet \phi _{Y,T}^{\circ }\right) \zeta _{X,Y^{\revsn },Z,T^{\revsn }} \\
&=&\varphi _{X\circ Z,T\circ Y}\left( \left( X\circ Z\right) \bullet \phi
_{Y,T}^{\circ }\right) \zeta _{X,Y^{\revsn },Z,T^{\revsn }}.
\end{eqnarray*}
\end{invisible}
\end{proof}

\begin{lemma}
In a duoidal category $(\Cc,\circ,I,\bullet,J)$ with a reversion, the following identities are true.
\begin{eqnarray}
\widetilde{\psi }_{X,I}
\left( X^{\revsn }\bullet \left( \Delta _{I}^{\bullet }\circ J\right) \left(
l_{J}^{\circ }\right) ^{-1}\right) &=&((r_X^\circ)^\revsn)^{-1}r_{X^\revsn}^\bullet ,\label{eq:BohmphiR}
\\
{\widetilde{\varphi }%
_{I,X}}\left( \left( J\circ \Delta _{I}^{\bullet }\right) \left(
r_{J}^{\circ }\right) ^{-1}\bullet X^{\revsn }\right) &=&((l_X^\circ)^\revsn)^{-1}l_{X^\revsn}^\bullet ,\label{eq:BohmphiL}\\
\psi _{I,I}\left( \phi _{0}^{\circ}\bullet I\right)
&=&\left(\varepsilon _{I}^{\bullet }\circ \left( I\bullet
l_{I}^{\circ }\right) ^{-1}\right) \left( l_{I\bullet I}^{\circ }\right)
^{-1}, \label{eq:psiphi0},\\
\varphi _{I,I}\left( I\bullet \phi _{0}^{\circ}\right)&=&\left( \left( \left(
r_{I}^{\circ }\right) ^{-1}\bullet I\right) \circ \varepsilon _{I}^{\bullet }\right) \left(
r_{I\bullet I}^{\circ }\right) ^{-1}, \label{eq:phiphi0}\\
\left( J\circ \left( I\bullet
m_{J}^{\circ }\right) \right) \psi _{J,J}\left( \phi
_{0}^{\bullet }\bullet J\right)&=&\left( J\circ \left( r_{I}^{\bullet }\right) ^{-1}\right) \left( \left(
r_{J}^{\circ }\right) ^{-1}\right)r_{J}^{\bullet },\label{eq:psiphi0bul} \\
 \left( \left( m_{J}^{\circ }\bullet I\right) \circ J\right)\varphi _{J,J}\left( J\bullet
\phi _{0}^{\bullet }\right) &=& \left( \left( l_{I}^{\bullet }\right)
^{-1}\circ J\right) \left( l_{J}^{\circ }\right) ^{-1}l_J^\bullet.\label{eq:phiphi0bul}
\end{eqnarray}
%sono l'unitarietà imposta dalla Bohm alla convolution structure che ha definito nei suoi appunti. Erano sostanzialmente già utilizzate nei nostri conti. Le ho estrapolate per metterle in evidenza. 
\end{lemma}
\begin{proof} We prove \eqref{eq:BohmphiR}:
    \begin{align*}
\underbracket[0.140ex]{\widetilde{\psi }_{X,I}}%
\left( X^{\revsn }\bullet \left( \Delta _{I}^{\bullet }\circ J\right) \left(
l_{J}^{\circ }\right) ^{-1}\right) &=l_{\left( X\circ I\right) ^{\revsn
}}^{\circ }\left( I\circ r_{\left( X\circ I\right) ^{\revsn }}^{\bullet
}\right) \underbracket[0.140ex]{\widetilde{\gamma }_{X,I,I,J}\left( X^{\revsn }\bullet
\left( \Delta _{I}^{\bullet }\circ J\right) \right) }\left( X^{\revsn }\bullet
\left( l_{J}^{\circ }\right) ^{-1}\right)\overset{\eqref{form:gammatildeIIJ}}{=}((r_X^\circ)^\revsn)^{-1}r_{X^\revsn}^\bullet.
\end{align*}%
Similarly, by means of \eqref{form:deltatildeJII}, one proves \eqref{eq:BohmphiL}.
\begin{invisible}
    \begin{eqnarray*}
\underbracket[0.140ex]{\widetilde{\varphi }%
_{I,X}}\left( \left( J\circ \Delta _{I}^{\bullet }\right) \left(
r_{J}^{\circ }\right) ^{-1}\bullet X^{\revsn }\right) 
&=& r_{\left( I\circ X\right) ^{\revsn
}}^{\circ }\left( l_{\left( I\circ X\right) ^{\revsn }}^{\bullet }\circ
I\right) \underbracket[0.140ex]{\widetilde{\delta }_{J,I,I,X}\left( \left( J\circ
\Delta _{I}^{\bullet }\right) \bullet X^{\revsn }\right) }\left( \left(
r_{J}^{\circ }\right) ^{-1}\bullet X^{\revsn }\right) \\&\overset{\eqref{form:deltatildeJII}}{=}&((l_X^\circ)^\revsn)^{-1}l_{X^\revsn}^\bullet.
\end{eqnarray*}
\end{invisible}
We prove \eqref{eq:psiphi0}:
\begin{eqnarray*}
\underbracket[0.140ex]{\psi _{I,I}}\left( \phi _{0}^{\circ}\bullet I\right) 
&=&\gamma
_{J,I,I,I}\left( \left( \left( l_{I^{\revsn }}^{\bullet }\right) ^{-1}\circ
I\right) \bullet I\right) \left( \left( r_{I^{\revsn }}^{\circ }\right)
^{-1}\bullet I\right) \left( \phi _{0}^{\circ}\bullet I\right) 
\\
&=&\gamma
_{J,I,I,I}\left( \left( \left( J\bullet \phi _{0}^{\circ}\right) \circ I\right)
\bullet I\right) \left( \left( \underbracket[0.140ex]{\left( l_{I}^{\bullet }\right)
^{-1}}\circ I\right) \bullet I\right) \left( \left( r_{I}^{\circ }\right)
^{-1}\bullet I\right)  \\
&=& \gamma
_{J,I,I,I}\left( \left( \left( J\bullet \phi _{0}^{\circ}\right) \circ I\right)
\bullet I\right) \left( \left( \left( \varepsilon _{I}^{\bullet }\bullet
I\right) \Delta _{I}^{\bullet }\circ I\right) \bullet I\right) \left( \left(
l_{I}^{\circ }\right) ^{-1}\bullet I\right) \\
&=&\left( \varepsilon _{I}^{\bullet }\circ \left( I\bullet (I\circ I)\right) \right) \underbracket[0.140ex]{\gamma _{I,I,I,I}\left( \left( \left( I\bullet
\phi _{0}^{\circ}\right) \Delta _{I}^{\bullet }\circ I\right) \bullet I\right) }%
\left( \left( l_{I}^{\circ }\right) ^{-1}\bullet I\right)\\
&\overset{\eqref{form:phi0gamma}}{=}&\left( \varepsilon _{I}^{\bullet }\circ
\left( I\bullet (I\circ I)\right) \right) \left( I\circ \left( I\bullet
l_{I}^{\circ }\right) ^{-1}\right) \left( l_{I\bullet I}^{\circ }\right)
^{-1}\left( l_{I}^{\circ }\bullet I\right) \left( \left( l_{I}^{\circ
}\right) ^{-1}\bullet I\right) \\
&=&\left(\varepsilon _{I}^{\bullet }\circ \left( I\bullet
l_{I}^{\circ }\right) ^{-1}\right) \left( l_{I\bullet I}^{\circ }\right)
^{-1}.
\end{eqnarray*}%
Similarly one proves \eqref{eq:phiphi0}, \eqref{eq:psiphi0bul} and \eqref{eq:phiphi0bul}
by means of \eqref{form:phi0delta}, \eqref{form:phi0gamma-bul} and \eqref{form:phi0delta-bul}, respectively. 
\begin{invisible}
Let us prove \eqref{eq:phiphi0}:
   \begin{eqnarray*}
\underbracket[0.140ex]{\varphi _{I,I}}\left( I\bullet \phi _{0}^{\circ}\right)
&=&\delta
_{I,I,I,J}\left( I\bullet \left( I\circ \left( r_{I^{\revsn }}^{\bullet
}\right) ^{-1}\right) \right) \left( I\bullet \left( l_{I^{\revsn }}^{\circ
}\right) ^{-1}\right) \left( I\bullet \phi _{0}^{\circ}\right) 
\\
&=&\delta
_{I,I,I,J}\left( I\bullet \left( I\circ \left( \phi _{0}^{\circ}\bullet J\right)
\right) \right) \left( I\bullet \left( I\circ \underbracket[0.140ex]{\left(
r_{I}^{\bullet }\right) ^{-1}}\right) \right) \left( I\bullet \left(
l_{I}^{\circ }\right) ^{-1}\right) \\
&=&\delta
_{I,I,I,J}\left( I\bullet \left( I\circ \left( \phi _{0}^{\circ}\bullet J\right)
\right) \right) \left( I\bullet \left( I\circ \left( I\bullet \varepsilon
_{I}^{\bullet }\right) \Delta _{I}^{\bullet }\right) \right) \left( I\bullet
\left( l_{I}^{\circ }\right) ^{-1}\right) \\
&=&\left( \left( (I\circ I)\bullet I\right) \circ \varepsilon
_{I}^{\bullet }\right) \underbracket[0.140ex]{\delta _{I,I,I,I}\left( I\bullet \left(
I\circ \left( \phi _{0}^{\circ}\bullet I\right) \Delta _{I}^{\bullet }\right)
\right) }\left( I\bullet \left( l_{I}^{\circ }\right) ^{-1}\right) \\
&\overset{\eqref{form:phi0delta}}{=}&\left( \left( (I\circ I)\bullet
I\right) \circ \varepsilon _{I}^{\bullet }\right) \left( \left( \left(
r_{I}^{\circ }\right) ^{-1}\bullet I\right) \circ I\right) \left(
r_{I\bullet I}^{\circ }\right) ^{-1}\left( I\bullet r_{I}^{\circ }\right)
\left( I\bullet \left( l_{I}^{\circ }\right) ^{-1}\right)  \\
&=& \left( \left( \left(
r_{I}^{\circ }\right) ^{-1}\bullet I\right) \circ \varepsilon _{I}^{\bullet }\right) \left(
r_{I\bullet I}^{\circ }\right) ^{-1} .
\end{eqnarray*} 
We now prove \eqref{eq:phiphi0bul}
\begin{eqnarray*}
 \left( \left( m_{J}^{\circ }\bullet I\right) \circ J\right)\underbracket[0.140ex]{\varphi _{J,J}}\left( J\bullet
\phi _{0}^{\bullet }\right) 
&=& \left( \left( m_{J}^{\circ }\bullet I\right) \circ J\right)\delta
_{J,I,J,J}\left( J\bullet \left( I\circ \left( r_{J^{\revsn }}^{\bullet
}\right) ^{-1}\right) \right) \left( J\bullet \left( l_{J^{\revsn }}^{\circ
}\right) ^{-1}\right) \left( J\bullet \phi _{0}^{\bullet }\right)  \\
&=&\left( \left( m_{J}^{\circ }\bullet I\right) \circ J\right)\delta
_{J,I,J,J}\left( J\bullet \left( I\circ \left( \phi _{0}^{\bullet }\bullet
J\right) \right) \right) \left( J\bullet \left( I\circ \left( r_{J}^{\bullet
}\right) ^{-1}\right) \right) \left( J\bullet \left( l_{J}^{\circ }\right)
^{-1}\right)  \\
&=&\underbracket[0.140ex]{\left( \left( m_{J}^{\circ }\bullet I\right) \circ J\right)
\delta _{J,I,J,J}\left( J\bullet \left( I\circ \left( \phi _{0}^{\bullet
}\bullet J\right) \Delta _{J}^{\bullet }\right) \right) \left( J\bullet
\left( l_{J}^{\circ }\right) ^{-1}\right)} \\
&\overset{\eqref{form:phi0delta-bul}}{=}&\left( \left( l_{I}^{\bullet }\right)
^{-1}\circ J\right) \left( l_{J}^{\circ }\right) ^{-1}l_J^\bullet .
\end{eqnarray*}%
We prove \eqref{eq:psiphi0bul} 
\begin{eqnarray*}
\left( J\circ \left( I\bullet
m_{J}^{\circ }\right) \right) \underbracket[0.140ex]{\psi _{J,J}}\left( \phi
_{0}^{\bullet }\bullet J\right) 
&=&\left( J\circ \left( I\bullet m_{J}^{\circ }\right) \right) \gamma
_{J,J,I,J}\left( \left( \left( l_{J^{\revsn }}^{\bullet }\right) ^{-1}\circ
I\right) \bullet J\right) \left( \left( r_{J^{\revsn }}^{\circ }\right)
^{-1}\bullet J\right) \left( \phi _{0}^{\bullet }\bullet J\right) \\
&=&\left( J\circ \left( I\bullet m_{J}^{\circ }\right) \right) \gamma
_{J,J,I,J}\left( \left( \left( J\bullet \phi _{0}^{\bullet }\right) \circ
I\right) \bullet J\right) \left( \left( \left( l_{J}^{\bullet }\right)
^{-1}\circ I\right) \bullet J\right) \left( \left( r_{J}^{\circ }\right)
^{-1}\bullet J\right) \\
&=&\underbracket[0.140ex]{\left( J\circ \left(
I\bullet m_{J}^{\circ }\right) \right) \gamma _{J,J,I,J}\left( \left( \left(
J\bullet \phi _{0}^{\bullet }\right) \Delta _{J}^{\bullet }\circ I\right)
\bullet J\right) \left( \left( r_{J}^{\circ }\right) ^{-1}\bullet J\right)} \\
&\overset{\eqref{form:phi0gamma-bul}}{=}&\left( J\circ \left( r_{I}^{\bullet }\right) ^{-1}\right) \left( \left(
r_{J}^{\circ }\right) ^{-1}\right) r_{J}^{\bullet }.
\end{eqnarray*}
\end{invisible}
\end{proof}

\subsection{The antipode}

To establish a coherent notion of antipode, we need an analogue of the convolution product.  We obtain this by drawing it from \cite[Proposition 6.2]{Bohm-Lack}. Note that we are not going to recover in our setting the entire structure of the convolution category considered in \cite{Bohm-Lack}, but only what is needed to obtain the uniqueness of the convolution inverse, which will ensure the uniqueness of the antipode, its compatibility with multiplication and comultiplication, and the fact that a bimonoid morphism between Hopf monoids is compatible with antipodes.\medskip

For an object $A$ in $\mathsf{Mon}(\Cc^{\circ})$ and an object $C$ in $\mathsf{Comon}(\Cc^{\bullet})$, we can define the following morphisms in $\Cc$:
\begin{eqnarray*}
1_{C,A}^{l} &:=&\left( \left( u_{A}^{\circ }\bullet I\right) \circ J\right)
\left( \Delta _{I}^{\bullet }\circ J\right) \left( l_{J}^{\circ }\right)
^{-1}\varepsilon _{C}^{\bullet }:C\rightarrow \left( A\bullet I\right) \circ
J, \\
1_{C,A}^{r} &:=&\left( J\circ \left( I\bullet u_{A}^{\circ }\right) \right)
\left( J\circ \Delta _{I}^{\bullet }\right) \left( r_{J}^{\circ }\right)
^{-1}\varepsilon _{C}^{\bullet }:C\rightarrow J\circ \left( I\bullet
A\right) .
\end{eqnarray*}
Given $f:C\rightarrow A$ and $g:C\rightarrow A^{\revsn
}$ arbitrary morphisms in $\Cc$, we set
\begin{eqnarray*}
f\star g &:&\xymatrix@C=1cm{
C \ar[r]^-{(f \bullet g)\,\Delta_C^{\bullet}} &
A \bullet A^{\revsn} \ar[r]^-{\varphi_{A,A}} &
\bigl((A \circ A) \bullet I\bigr) \circ J
\ar[r]^-{(m_A^{\circ} \bullet I)\,\circ\, J} &
(A \bullet I) \circ J
}, \\
g\star f &:&\xymatrix@C=1cm{
C \ar[r]^-{(g \bullet f)\,\Delta_C^{\bullet}} &
A^{\revsn} \bullet A \ar[r]^-{\psi_{A,A}} &
J \circ \bigl(I \bullet (A \circ A)\bigr)
\ar[r]^-{J \circ (I \bullet m_A^{\circ})} &
J \circ (I \bullet A)
}.
\end{eqnarray*}
Moreover, given arbitrary morphisms $h:C\to((A\bullet I)\circ J)$ and $h':C\to(J\circ(I\bullet A))$ in $\Cc$, we set
\begin{eqnarray*}
g\star h &:&\xymatrix@R=0cm@C=1cm{
C \ar[r]^-{(g \bullet h)\,\Delta_C^{\bullet}} &
A^{\revsn} \bullet \bigl((A \bullet I)\circ J\bigr)
   \ar[r]^-{\widetilde{\psi}_{A,A}} &
(A \circ A)^{\revsn}
   \ar[r]^-{(m_A^{\circ})^{\revsn}} &
A^{\revsn}
},
\\
h' \star g^{\prime } &:&\xymatrix@R=0cm@C=1cm{
C \ar[r]^-{(h' \bullet g')\,\Delta_C^{\bullet}} &
\bigl(J \circ (I \bullet A)\bigr) \bullet A^{\revsn}
   \ar[r]^-{\widetilde{\varphi}_{A,A}} &
(A \circ A)^{\revsn}
   \ar[r]^-{(m_A^{\circ})^{\revsn}} &
A^{\revsn}
}.
\end{eqnarray*}
% \begin{eqnarray*}
% g\star \left( f\star g^{\prime }\right)  &:&C\overset{\left( g\bullet \left(
% f\star g^{\prime }\right) \right) \Delta _{C}^{\bullet }}{\longrightarrow }%
% A^{\revsn }\bullet \left( \left( A\bullet I\right) \circ J\right) \overset{%
% \widetilde{\psi }_{A,A}}{\longrightarrow }\left( A\circ A\right) ^{\revsn }%
% \overset{\left( m_{A}^{\circ }\right) ^{\revsn }}{\longrightarrow }A^{\revsn },
% \\
% \left( g\star f\right) \star g^{\prime } &:&C\overset{\left( \left( g\star
% f\right) \bullet g^{\prime }\right) \Delta _{C}^{\bullet }}{\longrightarrow }%
% \left( J\circ \left( I\bullet A\right) \right) \bullet A^{\revsn }\overset{%
% \widetilde{\varphi }_{A,A}}{\longrightarrow }\left( A\circ A\right) ^{\revsn }%
% \overset{\left( m_{A}^{\circ }\right) ^{\revsn }}{\longrightarrow }A^{\revsn }.
% \end{eqnarray*}%

Given %an object $C$ in $\mathsf{Comon}(\Cc^{\bullet})$, an object $A$ in $\mathsf{Mon}(\Cc^{\circ})$ and 
morphisms $f:C\rightarrow A$ and $g:C\rightarrow A^{\revsn
}$ as above, we say that $g$ is a convolution inverse of $f$ if $f \star g=1^{l}_{C,A}$ and $g \star f=1^{r}_{C,A}$.
We are now able to recover the uniqueness of convolution inverses, cf. \cite[Proposition 6.2]{Bohm-Lack}.

\begin{lemma}
\label{lem:uniqinv}
\bigskip Let $C$ be an object in $\mathsf{Comon}(\Cc^{\bullet})$ and $A$ be an object in $\mathsf{Mon}(\Cc^{\circ})$. Then, a morphism $
f:C\rightarrow A$ in $\Cc$ may have a unique convolution inverse.
\end{lemma}

\begin{proof}
Consider morphisms $g,g^{\prime }:C\rightarrow A^{\revsn }$ in $\Cc$. 
We compute
\begin{eqnarray*}
g\star \left( f\star g^{\prime }\right)  &=&\left( m_{A}^{\circ }\right) ^{%
\revsn}\widetilde{\psi }_{A,A}\left( g\bullet \left( f\star g^{\prime
}\right) \right) \Delta _{C}^{\bullet } \\
&=&\left( m_{A}^{\circ }\right) ^{\revsn}\widetilde{\psi }_{A,A}\left( A^{%
\revsn}\bullet \left( \left( m_{A}^{\circ }\bullet I\right) \circ J\right)
\right) \left( A^{\revsn}\bullet \varphi _{A,A}\right) \left( g\bullet
\left( f\bullet g^{\prime }\right) \Delta _{C}^{\bullet }\right) \Delta
_{C}^{\bullet } \\
&=&\left( m_{A}^{\circ }\right) ^{\revsn}\left( A\circ m_{A}^{\circ }\right)
^{\revsn}\underbracket[0.140ex]{\widetilde{\psi }_{A,A\circ A}\left(
A^{\revsn }\bullet \varphi _{A,A}\right) }\left( g\bullet \left( f\bullet
g^{\prime }\right) \Delta _{C}^{\bullet }\right) \Delta _{C}^{\bullet } \\
&\overset{\eqref{form:varphipsitilde}}{=}&\left( m_{A}^{\circ }\right) ^{%
\revsn}\left( m_{A}^{\circ }\circ A\right) ^{\revsn}\widetilde{\varphi }%
_{A\circ A,A}\left( \psi _{A,A}\bullet A^{\revsn}\right) \left( \left(
g\bullet f\right) \Delta _{C}^{\bullet }\bullet g^{\prime }\right) \Delta
_{C}^{\bullet } \\
&=&\left( m_{A}^{\circ }\right) ^{\revsn}\widetilde{\varphi }_{A,A}\left(
\left( J\circ \left( I\bullet m_{A}^{\circ }\right) \right) \bullet A^{\revsn%
}\right) \left( \psi _{A,A}\bullet A^{\revsn}\right) \left( \left( g\bullet
f\right) \Delta _{C}^{\bullet }\bullet g^{\prime }\right) \Delta
_{C}^{\bullet } \\
&=&\left( m_{A}^{\circ }\right) ^{\revsn}\widetilde{\varphi }_{A,A}\left(
\left( g\star f\right) \bullet g^{\prime }\right) \Delta _{C}^{\bullet
}=\left( g\star f\right) \star g^{\prime }
\end{eqnarray*}
so that $g\star \left( f\star g^{\prime }\right) =\left( g\star f\right) \star
g^{\prime }.$ Moreover,
\begin{eqnarray*}
g\star 1_{C,A}^{l} &=&\left( m_{A}^{\circ }\right) ^{\revsn }\widetilde{\psi }%
_{A,A}\left( g\bullet 1_{C,A}^{l}\right) \Delta _{C}^{\bullet } \\
&=&\left( m_{A}^{\circ }\right) ^{\revsn }\widetilde{\psi }_{A,A}\left(
A^{\revsn }\bullet \left( \left( u_{A}^{\circ }\bullet I\right) \circ J\right)
\right) \left( g\bullet \left( \Delta _{I}^{\bullet }\circ J\right) \left(
l_{J}^{\circ }\right) ^{-1}\right) \left( C\bullet \varepsilon _{C}^{\bullet
}\right) \Delta _{C}^{\bullet } \\
&=&\left( m_{A}^{\circ }\right) ^{\revsn }\left( A\circ u_{A}^{\circ }\right)
^{\revsn }\underbracket[0.140ex]{\widetilde{\psi }_{A,I}\left( A^\revsn\bullet \left( \Delta _{I}^{\bullet
}\circ J\right) \left( l_{J}^{\circ }\right) ^{-1}\right)}(g\bullet J) \left(
r_{C}^{\bullet }\right) ^{-1} \\
&\overset{\eqref{eq:BohmphiR}}=&\left( r_{A}^{\circ }\right) ^{\revsn }((r_A^\circ)^\revsn)^{-1}r_{A^\revsn}^\bullet \left( r_{A^{\revsn }}^{\bullet }\right)
^{-1}g =g.
\end{eqnarray*}%
Similarly, by means of \eqref{eq:BohmphiL} one proves that $1_{C,A}^{r}\star g=g.$
\begin{invisible}
\begin{eqnarray*}
1_{C,A}^{r}\star g &=&\left( m_{A}^{\circ }\right) ^{\revsn }\widetilde{\varphi
}_{A,A}\left( 1_{C,A}^{r}\bullet g\right) \Delta _{C}^{\bullet } \\
&=&\left( m_{A}^{\circ }\right) ^{\revsn }\widetilde{\varphi }_{A,A}\left(
\left( J\circ \left( I\bullet u_{A}^{\circ }\right) \right) \bullet A^{\revsn
}\right) \left( \left( J\circ \Delta _{I}^{\bullet }\right) \left(
r_{J}^{\circ }\right) ^{-1}\bullet g\right) \left( \varepsilon _{C}^{\bullet
}\bullet C\right) \Delta _{C}^{\bullet } \\
&=&\left( m_{A}^{\circ }\right) ^{\revsn }\left( u_{A}^{\circ }\circ A\right)
^{\revsn }\underbracket[0.140ex]{\widetilde{\varphi }_{I,A}\left( \left( J\circ \Delta _{I}^{\bullet
}\right) \left( r_{J}^{\circ }\right) ^{-1}\bullet A^\revsn\right)} (J\bullet g)\left(
l_{C}^{\bullet }\right) ^{-1} \\
&\overset{\eqref{eq:BohmphiL}}=&\left( l_{A}^{\circ }\right) ^{\revsn }((l_A^\circ)^\revsn)^{-1}l_{A^\revsn}^\bullet \left( l_{A^{\revsn
}}^{\bullet }\right) ^{-1}g =g.
\end{eqnarray*}%
so that $g\star 1_{C,A}^{l}=g=1_{C,A}^{r}\star g.$ \end{invisible}
As a consequence, if $f\star
g^{\prime }=1_{C,A}^{l}$ and $g\star f=1_{C,A}^{r},$ we obtain
$g=g\star 1_{C,A}^{l}=g\star \left( f\star g^{\prime }\right) =\left( g\star
f\right) \star g^{\prime }=1_{C,A}^{r}\star g^{\prime }=g^{\prime }$.
\end{proof}

Following \cite[Theorem 7.2]{Bohm-Lack}, we give the following definition of antipode for a bimonoid in a duoidal category $(\Cc,\circ,\bullet)$.

\begin{definition}
\label{def:Hopfmon}
Let $\left( \mathcal{C},\circ ,I,\bullet
,J\right) $ be a duoidal category  with a reversion and $\left( H,m_{H}^{\circ },u_{H}^{\circ },\Delta
_{H}^{\bullet },\varepsilon _{H}^{\bullet }\right) $ be an object in $\Bimon(\Cc,\circ,\bullet)$.
A morphism $%
\Sigma _{H}:H\rightarrow H^{\revsn }$
%we set
%\begin{eqnarray*}
%\mathrm{Id}_{H}\star \Sigma _{H} &:=&\left( H\overset{\left( H\bullet \Sigma
%_{H}\right) \Delta _{H}^{\bullet }}{\longrightarrow }H\bullet H^{\revsn }%
%\overset{\varphi _{H,H}}{\longrightarrow }\left( \left( H\circ H\right)
%\bullet I\right) \circ J\overset{\left( m_{H}^{\circ }\bullet I\right) \circ
%J}{\longrightarrow }\left( H\bullet I\right) \circ J\right) , \\
%\Sigma _{H}\star \mathrm{Id}_{H} &:=&\left( H\overset{\left( \Sigma
%_{H}\bullet H\right) \Delta _{H}^{\bullet }}{\longrightarrow }H^{\revsn
%}\bullet H\overset{\psi _{H,H}}{\longrightarrow }J\circ \left( I\bullet
%\left( H\circ H\right) \right) \overset{J\circ \left( I\bullet m_{H}^{\circ
%}\right) }{\longrightarrow }J\circ \left( I\bullet H\right) \right) .
%\end{eqnarray*}%
is an
\emph{antipode} for $H$ if it is a convolution inverse of $\id_H$ i.e.\ if it
obeys $\mathrm{Id}_{H}\star \Sigma _{H}=1_{H,H}^{l}$ and $\Sigma _{H}\star \mathrm{Id}_{H}=1_{H,H}^{r}$ that read as:
\begin{equation}\label{eq:antipode}
  (\left( m_{H}^{\circ }\bullet I\right) \circ
J)\varphi_{H,H}(H\bullet\Sigma_{H})\Delta^{\bullet}_{H}=1_{H,H}^{l},\qquad (J\circ \left( I\bullet m_{H}^{\circ
}\right))\psi_{H,H}(\Sigma_{H}\bullet H)\Delta_{H}^{\bullet}=1_{H,H}^{r}.
\end{equation}
%\begin{eqnarray*}
%\mathrm{Id}_{H}\star \Sigma _{H} &=&\left( H\overset{\varepsilon
%_{H}^{\bullet }}{\longrightarrow }J\overset{\left( l_{J}^{\circ }\right)
%^{-1}}{\longrightarrow }I\circ J\overset{\Delta _{I}^{\bullet }\circ J}{%
%\longrightarrow }\left( I\bullet I\right) \circ J\overset{\left(
%u_{H}^{\circ }\bullet I\right) \circ J}{\longrightarrow }\left( H\bullet
%I\right) \circ J\right) , \\
%\Sigma _{H}\star \mathrm{Id}_{H} &=&\left( H\overset{\varepsilon
%_{H}^{\bullet }}{\longrightarrow }J\overset{\left( r_{J}^{\circ }\right)
%^{-1}}{\longrightarrow }J\circ I\overset{J\circ \Delta _{I}^{\bullet }}{%
%\longrightarrow }J\circ \left( I\bullet I\right) \overset{J\circ \left(
%I\bullet u_{H}^{\circ }\right) }{\longrightarrow }J\circ \left( I\bullet
%H\right) \right) .
%\end{eqnarray*}%
In this case, we say that $\left( H,m_{H}^{\circ },u_{H}^{\circ
},\Delta _{H}^{\bullet },\varepsilon _{H}^{\bullet },\Sigma _{H}\right) $ is
a \emph{Hopf monoid} in $\left( \mathcal{C},\circ ,I,\bullet ,J\right) $. 

\cref{lem:uniqinv} ensures the uniqueness of the antipode.

Hopf monoids constitute a full subcategory of $\Bimon(\Cc,\circ,\bullet)$ that we denote by $\Hopfmon(\Cc,\circ,\bullet)$.
\end{definition}

As expected, unit objects are examples of Hopf monoids.

\begin{lemma}
\label{lem:IHopf} $I$ is a Hopf monoid with antipode $\phi _{0}^{\circ}$ and $J$ is a Hopf monoid with antipode $\phi _{0}^{\bullet }$.
\end{lemma}

\begin{proof}
We compute%
\begin{eqnarray*}
&&\phi _{0}^{\circ}\star \mathrm{Id}_{I} =\left( J\circ \left( I\bullet m_{I}^{\circ
}\right) \right) \underbracket[0.140ex]{\psi _{I,I}\left( \phi _{0}^{\circ}\bullet I\right)} \Delta
_{I}^{\bullet }\overset{\eqref{eq:psiphi0}}=
\left( J\circ \left( I\bullet m_{I}^{\circ
}\right) \right) \left(\varepsilon _{I}^{\bullet }\circ \left( I\bullet
l_{I}^{\circ }\right) ^{-1}\right) \left( l_{I\bullet I}^{\circ }\right)
^{-1} \Delta
_{I}^{\bullet }\\
&=&%\left( \varepsilon _{I}^{\bullet }\circ \left( I\bullet I\right) \right)
%\left( l_{I\bullet I}^{\circ }\right) ^{-1}\Delta _{I}^{\bullet } =
\left( \varepsilon _{I}^{\bullet }\circ \left( I\bullet I\right) \right)
\left( I\circ \Delta _{I}^{\bullet }\right) \left( l_{I}^{\circ }\right)^{-1} =
\left( J\circ \Delta _{I}^{\bullet }\right) \left( \varepsilon
_{I}^{\bullet }\circ I\right) \left( r_{I}^{\circ }\right) ^{-1} 
=%\left( J\circ \Delta _{I}^{\bullet }\right) \left( r_{J}^{\circ }\right)
%^{-1}\varepsilon _{I}^{\bullet }=
\left( J\circ \left( I\bullet u_{I}^{\circ
}\right) \right) \left( J\circ \Delta _{I}^{\bullet }\right) \left(
r_{J}^{\circ }\right) ^{-1}\varepsilon _{I}^{\bullet }=1_{I,I}^{r}.
\end{eqnarray*}%
Similarly, by means of \eqref{eq:phiphi0}, one proves that $\mathrm{Id}_{I}\star \phi _{0}^{\circ}=1_{I,I}^{l}$, and, by means of \eqref{eq:phiphi0bul} and \eqref{eq:psiphi0bul} one checks that
$J$ is a Hopf monoid with antipode $\phi _{0}^{\bullet }$.
\begin{invisible}
We have
\begin{eqnarray*}
\mathrm{Id}_{I}\star \phi _{0}^{\circ} &=&\left( \left( m_{I}^{\circ }\bullet
I\right) \circ J\right) \underbracket[0.140ex]{\varphi _{I,I}\left( I\bullet \phi _{0}^{\circ}\right)}
\Delta _{I}^{\bullet } \\
&\overset{\eqref{eq:phiphi0}}=&\left( \left( m_{I}^{\circ }\bullet
I\right) \circ J\right) \left( \left( \left(
r_{I}^{\circ }\right) ^{-1}\bullet I\right) \circ \varepsilon _{I}^{\bullet }\right) \left(
r_{I\bullet I}^{\circ }\right) ^{-1}
\Delta _{I}^{\bullet }\\
&=&\left( \left( I\bullet I\right) \circ \varepsilon _{I}^{\bullet }\right)
\left( r_{I\bullet I}^{\circ }\right) ^{-1}\Delta _{I}^{\bullet } \\
&=&\left( \left( I\bullet I\right) \circ \varepsilon _{I}^{\bullet }\right)
\left( \Delta _{I}^{\bullet }\circ I\right) \left( r_{I}^{\circ }\right)
^{-1} \\
&=&\left( \Delta _{I}^{\bullet }\circ J\right) \left( I\circ \varepsilon
_{I}^{\bullet }\right) \left( l_{I}^{\circ }\right) ^{-1} \\
&=&\left( \Delta _{I}^{\bullet }\circ J\right) \left( l_{J}^{\circ }\right)
^{-1}\varepsilon _{I}^{\bullet } \\
&=&\left( \left( u_{I}^{\circ }\bullet I\right) \circ J\right) \left( \Delta
_{I}^{\bullet }\circ J\right) \left( l_{J}^{\circ }\right) ^{-1}\varepsilon
_{I}^{\bullet }=1_{I,I}^{l}.
\end{eqnarray*}
Let us prove that $J$ is a Hopf monoid with antipode $\phi _{0}^{\bullet }$.
We compute
\begin{eqnarray*}
\mathrm{Id}_{J}\star \phi _{0}^{\bullet } &=&\underbracket[0.140ex]{\left( \left( m_{J}^{\circ
}\bullet I\right) \circ J\right) \varphi _{J,J}\left( J\bullet
\phi _{0}^{\bullet }\right)} \Delta _{J}^{\bullet } \\
&\overset{\eqref{eq:phiphi0bul}}=&\left( \left( l_{I}^{\bullet }\right)
^{-1}\circ J\right) \left( l_{J}^{\circ }\right) ^{-1}l_J^\bullet \Delta _{J}^{\bullet }\\
&=&\left( \left( l_{I}^{\bullet }\right)
^{-1}\circ J\right) \left( l_{J}^{\circ }\right) ^{-1}
=\left( \left(
\varepsilon _{I}^{\bullet }\bullet I\right) \circ J\right) \left( \Delta
_{I}^{\bullet }\circ J\right) \left( l_{J}^{\circ }\right) ^{-1}\\&=&\left(
\left( u_{J}^{\circ }\bullet I\right) \circ J\right) \left( \Delta
_{I}^{\bullet }\circ J\right) \left( l_{J}^{\circ }\right) ^{-1}\varepsilon
_{J}^{\bullet }=1_{J,J}^l.
\end{eqnarray*}%
Moreover, we have
\begin{eqnarray*}
\phi _{0}^{\bullet }\star \mathrm{Id}_{J} &=&\underbracket[0.140ex]{\left( J\circ \left( I\bullet
m_{J}^{\circ }\right) \right) \psi _{J,J}\left( \phi
_{0}^{\bullet }\bullet J\right)} \Delta _{J}^{\bullet } \\
&\overset{\eqref{eq:psiphi0bul}}=&\left( J\circ \left( r_{I}^{\bullet }\right) ^{-1}\right) \left( \left(
r_{J}^{\circ }\right) ^{-1}\right)r_{J}^{\bullet } \Delta _{J}^{\bullet } \\
&=&\left( J\circ \left( r_{I}^{\bullet }\right) ^{-1}\right) \left( \left(
r_{J}^{\circ }\right) ^{-1}\right) =\left( J\circ \left( I\bullet
\varepsilon _{I}^{\bullet }\right) \right) \left( J\circ \Delta
_{I}^{\bullet }\right) \left( \left( r_{J}^{\circ }\right) ^{-1}\right)\\
&=&\left( J\circ \left( I\bullet u_{J}^{\circ }\right) \right) \left( J\circ
\Delta _{I}^{\bullet }\right) \left( \left( r_{J}^{\circ }\right)
^{-1}\right) \varepsilon _{J}^{\bullet }=1_{J,J}^r
\end{eqnarray*}\end{invisible}
\end{proof}

In the following result, we obtain some useful properties of the convolution product.

\begin{lemma}\label{lem:convcomp}
Let $B$ be an object in $\mathsf{Bimon}(\Cc,\circ,\bullet)$ equipped with a morphism $\Sigma
_{B}:B\rightarrow B^{\revsn }$. The following statements hold:
\begin{itemize}
    \item[1)] If $f:B\rightarrow A$ is a morphism in $\mathsf{Mon}(\Cc^{\circ})$, then
\[
f^{\revsn
}\Sigma _{B}\star f=\left( J\circ \left( I\bullet f\right) \right) \left(
\Sigma _{B}\star \mathrm{Id}_{B}\right) ,\quad  f\star f^{\revsn }\Sigma _{B}=\left(
\left( f\bullet I\right) \circ J\right) \left( \mathrm{Id}_{B}\star \Sigma
_{B}\right) ,
\]
and
$\left( J\circ \left( I\bullet f\right) \right)
1_{B,B}^{r}=1_{B,A}^{r}$, $\left( \left( f\bullet I\right) \circ J\right)
1_{B,B}^{l}=1_{B,A}^{l}.$

As a consequence, if $\Sigma _{B}$ is an antipode, $f^{\revsn }\Sigma _{B}\ $is the convolution inverse of $f.$

\item[2)] If $f:C\rightarrow B$ is a morphism in $\mathsf{Comon}(\Cc^{\bullet})$, then
\[
\Sigma _{B}f\star f=\left( \Sigma _{B}\star \mathrm{Id}_{B}\right) f,\quad f\star
\Sigma _{B}f=\left( \mathrm{Id}_{B}\star \Sigma _{B}\right) f,\quad %
1_{B,B}^{r}f=1_{C,B}^{r},\quad 1_{B,B}^{l}f=1_{C,B}^{l}.
\]
As a consequence,
if $\Sigma _{B}$ is an antipode, we get that $\Sigma _{B}f\ $is the
convolution inverse of $f.$
\end{itemize}
\end{lemma}

\begin{proof}
1) By definition of convolution product, naturality of $\psi_{X,Y}$ and $\varphi _{X,Y}$ and multiplicativity of $f$ one easily checks the first two identities, while the remaining two follow by unitality of $f$.
\begin{invisible}
 We compute%
\begin{eqnarray*}
f^{\revsn }\Sigma _{B}\star f &=&\left( J\circ \left( I\bullet m_{A}^{\circ
}\right) \right) \psi _{A,A}\left( f^{\revsn }\Sigma _{B}\bullet f\right)
\Delta _{B}^{\bullet } \\
&=&\left( J\circ \left( I\bullet m_{A}^{\circ }\right) \right) \psi
_{A,A}\left( f^{\revsn }\bullet f\right) \left( \Sigma _{B}\bullet B\right)
\Delta _{B}^{\bullet } \\
&=&\left( J\circ \left( I\bullet m_{A}^{\circ }\right) \right) \left( J\circ
\left( I\bullet \left( f\circ f\right) \right) \right) \psi _{B,B}\left(
\Sigma _{B}\bullet B\right) \Delta _{B}^{\bullet } \\
&=&\left( J\circ \left( I\bullet f\right) \right) \left( J\circ \left(
I\bullet m_{B}^{\circ }\right) \right) \psi _{B,B}\left( \Sigma _{B}\bullet
B\right) \Delta _{B}^{\bullet } \\
&=&\left( J\circ \left( I\bullet f\right) \right) \left( \Sigma _{B}\star
\mathrm{Id}_{B}\right) ,
\end{eqnarray*}  
and also 
\begin{eqnarray*}
\left( J\circ \left( I\bullet f\right) \right) 1_{B,B}^{r} &=&\left( J\circ
\left( I\bullet f\right) \right) \left( J\circ \left( I\bullet u_{B}^{\circ
}\right) \right) \left( J\circ \Delta _{I}^{\bullet }\right) \left(
r_{J}^{\circ }\right) ^{-1}\varepsilon _{B}^{\bullet } \\
&=&\left( J\circ \left( I\bullet u_{A}^{\circ }\right) \right) \left( J\circ
\Delta _{I}^{\bullet }\right) \left( r_{J}^{\circ }\right) ^{-1}\varepsilon
_{B}^{\bullet }=1_{B,A}^{r}.
\end{eqnarray*} 
We compute
 \begin{eqnarray*}
f\star f^{\revsn }\Sigma _{B} &=&\left( \left( m_{A}^{\circ }\bullet I\right)
\circ J\right) \varphi _{A,A}\left( f\bullet f^{\revsn }\Sigma _{B}\right)
\Delta _{B}^{\bullet } \\
&=&\left( \left( m_{A}^{\circ }\bullet I\right) \circ J\right) \varphi
_{A,A}\left( f\bullet f^{\revsn }\right) \left( B\bullet \Sigma _{B}\right)
\Delta _{B}^{\bullet } \\
&=&\left( \left( m_{A}^{\circ }\bullet I\right) \circ J\right) \left( \left(
\left( f\circ f\right) \bullet I\right) \circ J\right) \varphi _{B,B}\left(
B\bullet \Sigma _{B}\right) \Delta _{B}^{\bullet } \\
&=&\left( \left( f\bullet I\right) \circ J\right) \left( \left( m_{B}^{\circ
}\bullet I\right) \circ J\right) \varphi _{B,B}\left( B\bullet \Sigma
_{B}\right) \Delta _{B}^{\bullet } \\
&=&\left( \left( f\bullet I\right) \circ J\right) \left( \mathrm{Id}_{B}\star
\Sigma _{B}\right)
\end{eqnarray*}%
and
\begin{eqnarray*}
\left( \left( f\bullet I\right) \circ J\right) 1_{B,B}^{l} &=&\left( \left(
f\bullet I\right) \circ J\right) \left( \left( u_{B}^{\circ }\bullet
I\right) \circ J\right) \left( \Delta _{I}^{\bullet }\circ J\right) \left(
l_{J}^{\circ }\right) ^{-1}\varepsilon _{B}^{\bullet } \\
&=&\left( \left( u_{A}^{\circ }\bullet I\right) \circ J\right) \left( \Delta
_{I}^{\bullet }\circ J\right) \left( l_{J}^{\circ }\right) ^{-1}\varepsilon
_{B}^{\bullet }=1_{B,A}^{l}.
\end{eqnarray*}
\end{invisible}
If $\Sigma _{B}$ is an antipode, we get $f^{\revsn }\Sigma _{B}\star f=\left(
J\circ \left( I\bullet f\right) \right) \left( \Sigma _{B}\star \mathrm{Id}%
_{B}\right) =\left( J\circ \left( I\bullet f\right) \right)
1_{B,B}^{r}=1_{B,A}^{r}$ and similarly $f\star f^{\revsn }\Sigma _{B}=\left(
\left( f\bullet I\right) \circ J\right) \left( \mathrm{Id}_{B}\star \Sigma
_{B}\right) =\left( \left( f\bullet I\right) \circ J\right)
1_{B,B}^{l}=1_{B,A}^{l}.$ This means that $f^{\revsn }\Sigma _{B}\ $is the
convolution inverse of $f.$

2)  By definition of convolution product and comultiplicativity of $f$ one easily checks the first two identities, while the remaining two follow by counitality of $f$.
\begin{invisible}
We compute
\begin{eqnarray*}
\Sigma _{B}f\star f &=&\left( J\circ \left( I\bullet m_{B}^{\circ }\right)
\right) \psi _{B,B}\left( \Sigma _{B}f\bullet f\right) \Delta _{C}^{\bullet
}=\left( J\circ \left( I\bullet m_{B}^{\circ }\right) \right) \psi
_{B,B}\left( \Sigma _{B}\bullet B\right) \Delta _{B}^{\bullet }f=\left(
\Sigma _{B}\star \mathrm{Id}_{B}\right) f, \\
1_{B,B}^{r}f &=&\left( J\circ \left( I\bullet u_{B}^{\circ }\right) \right)
\left( J\circ \Delta _{I}^{\bullet }\right) \left( r_{J}^{\circ }\right)
^{-1}\varepsilon _{B}^{\bullet }f=\left( J\circ \left( I\bullet u_{B}^{\circ
}\right) \right) \left( J\circ \Delta _{I}^{\bullet }\right) \left(
r_{J}^{\circ }\right) ^{-1}\varepsilon _{C}^{\bullet }=1_{C,B}^{r},\\
f\star \Sigma _{B}f &=&\left( \left( m_{B}^{\circ }\bullet I\right) \circ
J\right) \varphi _{B,B}\left( f\bullet \Sigma _{B}f\right) \Delta
_{C}^{\bullet }=\left( \left( m_{B}^{\circ }\bullet I\right) \circ J\right)
\varphi _{B,B}\left( B\bullet \Sigma _{B}\right) \Delta _{B}^{\bullet
}f=\left( \mathrm{Id}_{B}\star \Sigma _{B}\right) f, \\
1_{B,B}^{l}f &=&\left( \left( u_{B}^{\circ }\bullet I\right) \circ J\right)
\left( \Delta _{I}^{\bullet }\circ J\right) \left( l_{J}^{\circ }\right)
^{-1}\varepsilon _{B}^{\bullet }f=\left( \left( u_{B}^{\circ }\bullet
I\right) \circ J\right) \left( \Delta _{I}^{\bullet }\circ J\right) \left(
l_{J}^{\circ }\right) ^{-1}\varepsilon _{C}^{\bullet }=1_{C,B}^{l}.
\end{eqnarray*}%
\end{invisible}
If $\Sigma _{B}$ is an antipode, we get $\Sigma _{B}f\star f=\left( \Sigma
_{B}\star \mathrm{Id}_{B}\right) f=1_{B,B}^{r}f=1_{C,B}^{r}$ and similarly $f\star
\Sigma _{B}f=\left( \mathrm{Id}_{B}\star \Sigma _{B}\right)
f=1_{B,B}^{l}f=1_{C,B}^{l}.$ This means that $\Sigma _{B}f\ $is the
convolution inverse of $f.$
\end{proof}

We now show that any bimonoid morphism between Hopf monoids preserves antipodes.

\begin{corollary}
\label{coro:convcomp}
 Let $B$ and $H$ be in $\mathsf{Hopf}(\Cc,\circ,\bullet)$. If $f:B\rightarrow H$ is a morphism in $\mathsf{Bimon}(\Cc,\circ,\bullet)$, we have $\Sigma _{H}f=f^{\revsn }\Sigma _{B}.$    
\end{corollary}

\begin{proof}
 By \cref{lem:convcomp}, $f^{\revsn }\Sigma _{B}\ $and $\Sigma _{H}f$ are both convolution inverses of $f.$ By the uniqueness of the convolution inverse established in \cref{lem:uniqinv}, we conclude.
\end{proof}

Another outcome of \cref{lem:convcomp} is that, sometimes, the antipode transfers between bimonoids.

\begin{lemma}%[\rd{upgrade:27/2/2026}]
\label{lem:retract}
Let $H$ and $H'$ be bimonoids in  a duoidal category $(\Cc,\circ,I,\bullet,J)$. If $H$ is a Hopf monoid, then
so is $H'$ in each of the following cases:
\begin{itemize}
    \item[$i)$] there is a morphism $s:H'\to H$  in $\Bimon(\Cc,\circ,\bullet)$ and a morphism $\Sigma_{H'}$ in $\Cc$
    such that $s^\revsn\Sigma_{H'}=\Sigma_{H}s$ and both $\left( s\bullet I\right) \circ J$ and $J\circ \left( I\bullet s\right)$ are monomorphisms in $\Cc$;
\item[$ii)$] there is a morphism $s:H\to H'$  in $\Bimon(\Cc,\circ,\bullet)$ and a morphism $\Sigma_{H'}$ in $\Cc$
    such that $s^\revsn\Sigma_{H}=\Sigma_{H'}s$ and $s$ is an epimorphism in $\Cc$;  
    \item[$iii)$] $H'$ is a retract of $H$ in $\Bimon(\Cc,\circ,\bullet)$.
\end{itemize}
\end{lemma}

\begin{proof}
$i)$ Since $s:H'\to H$ is a morphism of bimonoids, by \cref{lem:convcomp}, we have%
\begin{align*}
\left( \left( s\bullet I\right) \circ J\right) \left( \mathrm{Id}_{H'}\star
\Sigma _{H'}\right)  &=s\star s^{\revsn }\Sigma _{H'}=s\star \Sigma _{H}s=\left(
\mathrm{Id}_{H}\star \Sigma _{H}\right) s=1_{H,H}^{l}s=1_{H',H}^{l}=\left(
\left( s\bullet I\right) \circ J\right) 1_{H',H'}^{l}, \\
\left( J\circ \left( I\bullet s\right) \right) \left( \Sigma _{H'}\star
\mathrm{Id}_{H'}\right)  &=s^{\revsn }\Sigma _{H'}\star s=\Sigma _{H}s\star
s=\left( \Sigma _{H}\star \mathrm{Id}_{H}\right)
s=1_{H,H}^{r}s=1_{H',H}^{r}=\left( J\circ \left( I\bullet s\right) \right)
1_{H',H'}^{r}.
\end{align*}
Since $\left( s\bullet I\right) \circ J$ and $J\circ \left( I\bullet s\right)$ are monomorphisms in $\Cc$, we get $\mathrm{Id}_{H'}\star \Sigma _{H'}=1_{H',H'}^{l}$ and $\Sigma _{H'}\star
\mathrm{Id}_{H'}=1_{H',H'}^{r}.$
This means that $\Sigma _{H'}$ is an antipode.

$ii)$ It follows similarly, again using \cref{lem:convcomp}.
\begin{invisible}
Since $s:H\to H'$ is a morphism of bimonoids, by \cref{lem:convcomp}, we have
	\begin{align*}
	 \left( \mathrm{Id}_{H'}\star
		\Sigma _{H'}\right)s  &=s\star \Sigma _{H'}s=s\star s^{\revsn}\Sigma _{H}=	\left( \left( s\bullet I\right) \circ J\right)\left(
		\mathrm{Id}_{H}\star \Sigma _{H}\right) = 	\left( \left( s\bullet I\right) \circ J\right) 1_{H,H}^{l}=1_{H,H'}^{l}= 1_{H',H'}^{l}s, \\
		 \left( \Sigma _{H'}\star
		\mathrm{Id}_{H'}\right) s &=\Sigma _{H'}s\star
		s=s^{\revsn }\Sigma _{H}\star s=\left( J\circ \left( I\bullet s\right) \right)\left( \Sigma _{H}\star \mathrm{Id}_{H}\right)
		=\left( J\circ \left( I\bullet s\right) \right)1_{H,H}^{r}=1_{H,H'}^{r}=1_{H',H'}^{r}s.
	\end{align*}
By assumption $s$ is an epimorphism in $\Cc$, so we get $\mathrm{Id}_{H'}\star \Sigma _{H'}=1_{H',H'}^{l}$ and $\Sigma _{H'}\star \mathrm{Id}_{H'}=1_{H',H'}^{r}.$ This means that $\Sigma _{H'}$ is an antipode for $H'$.    
\end{invisible}

$iii)$ By assumption there are morphisms $s:H'\to H$ and $p:H\to H'$ in $\mathsf{Bimon}(\Cc,\circ,\bullet)$ with $ps=\id_{H'}$. Then $sp:H\to H$ commutes with the antipode by \cref{coro:convcomp}, i.e.\ $(sp)^\revsn\Sigma_H=\Sigma_Hsp$. Hence, if we set $\Sigma_{H'}:=p^\revsn\Sigma_Hs$, we get $s^\revsn\Sigma_{H'}=s^\revsn p^\revsn\Sigma_Hs=(sp)^\revsn\Sigma_Hs=\Sigma_Hsps=\Sigma_{H}s$. Moreover, both $\left( s\bullet I\right) \circ J$ and $J\circ \left( I\bullet s\right)$ are split-monomorphisms in $\Cc$. Thus, we can apply the item $i)$.
\end{proof}

\begin{corollary}
    The category $\mathsf{Hopf}(\Cc,\circ,\bullet)$ is a replete subcategory of $\Bimon(\Cc,\circ,\bullet)$.
\end{corollary}

\begin{proof}
    We have to prove that, given an isomorphism $f:H\to H'$ in $\mathsf{Bimon}(\Cc,\circ,\bullet)$ with $H$ in $\mathsf{Hopf}(\Cc,\circ,\bullet)$, then $H'$ is in $\mathsf{Hopf}(\Cc,\circ,\bullet)$. This follows by applying $iii)$ of \cref{lem:retract}.
\end{proof}

\subsection{Galois and co-Galois maps}

In the present subsection, we introduce the notions of Galois map and co-Galois map, which will be essential in what follows. These are inspired by \cite{Bohm-Lack}.

\begin{definition}
\label{def:Galmap} Let $\left( \mathcal{C},\circ
,I,\bullet ,J\right) $ be  a duoidal category with a reversion and let $\left( H,m_{H}^{\circ },u_{H}^{\circ
},\Delta _{H}^{\bullet },\varepsilon _{H}^{\bullet }\right) $ be a bimonoid.
Given a morphism $\Sigma _{H}:H\rightarrow H^{\revsn }$, we
define:

\begin{itemize}
\item the \emph{Galois maps} (cf. \cite[\S 3.5 and Eq. (7.6)]{Bohm-Lack}) as the natural transformations%
\begin{align*}
\varkappa _{P,X,Q}^{H} &: \xymatrix@C=3.5em{
 (P \bullet X) \circ Q 
   \ar[r]^-{(P \bullet X)\circ \rho_Q^{\bullet}} &
 (P \bullet X)\circ (H \bullet Q)
   \ar[r]^-{\zeta_{P,X,H,Q}} &
 (P \circ H)\bullet (X \circ Q)
   \ar[r]^-{\mu_P^{\circ}\bullet (X\circ Q)} &
 P \bullet (X \circ Q)
} \\
\overline{\varkappa }_{P,X,Q}^{H} &:\xymatrix@C=3.5em{
 P \bullet (X \circ Q)
   \ar[r]^-{P \bullet \bigl(X \circ (\Sigma_H \bullet Q)\,\rho_Q^{\bullet}\bigr)} &
 P \bullet \bigl(X \circ (H^{\revsn} \bullet Q)\bigr)
   \ar[r]^-{\delta_{P,X,H,Q}} &
 \bigl((P \circ H)\bullet X\bigr) \circ Q
   \ar[r]^-{(\mu_P^{\circ}\bullet X)\circ Q} &
 (P \bullet X)\circ Q
}
\end{align*}%
for any object $X$, for any right $H$-module $\left( P,\mu _{P}^{\circ
}\right) $ and left $H$-comodule $\left( Q,\rho _{Q}^{\bullet }\right) ;$

\item the \emph{co-Galois maps} (cf. \cite[\S 3.6 and Eq. (7.5)]{Bohm-Lack}) as the natural transformations
\begin{align*}
\varsigma _{P,X,Q}^{H} &:\xymatrix@C=3.5em{
 P \circ (X \bullet Q)
   \ar[r]^-{\rho_P^{\bullet} \circ (X \bullet Q)} &
 (P \bullet H) \circ (X \bullet Q)
   \ar[r]^-{\zeta_{P,H,X,Q}} &
 (P \circ X) \bullet (H \circ Q)
   \ar[r]^-{(P \circ X)\bullet \mu_Q^{\circ}} &
 (P \circ X)\bullet Q
}\\
\overline{\varsigma}_{P,X,Q}^{H} &: \xymatrix@C=3.5em{
 (P \circ X)\bullet Q
   \ar[r]^-{\bigl((P \bullet \Sigma_H)\,\rho_P^{\bullet} \circ X\bigr)\bullet Q} &
 \bigl((P \bullet H^{\revsn}) \circ X\bigr)\bullet Q
   \ar[r]^-{\gamma_{P,H,X,Q}} &
 P \circ \bigl(X \bullet (H \circ Q)\bigr)
   \ar[r]^-{P \circ (X \bullet \mu_Q^{\circ})} &
 P \circ (X \bullet Q)
}
\end{align*}%
for any object $X$, for any right $H$-comodule $\left( P,\rho _{P}^{\bullet
}\right) $ and left $H$-module $\left( Q,\mu _{Q}^{\circ }\right) $.
\end{itemize}
\end{definition}
\begin{invisible}
Let us check that $\varkappa _{P,X,Q}^{H}$ and $\overline{\varkappa }%
_{P,X,Q}^{H}$ are natural for any object $X$, for any right $H$-module $%
\left( P,\mu _{P}^{\circ }\right) $ and left $H$-comodule $\left( Q,\rho
_{Q}^{\bullet }\right) $. Let $\chi :X\rightarrow X^{\prime },$ $p:\left(
P,\mu _{P}^{\circ }\right) \rightarrow \left( P^{\prime },\mu _{P^{\prime
}}^{\circ }\right) $ and $q:\left( Q,\rho _{Q}^{\bullet }\right) \rightarrow
\left( Q^{\prime },\rho _{Q^{\prime }}^{\bullet }\right) $ be morphisms. Then%
\begin{eqnarray*}
\left( p\bullet \left( \chi \circ q\right) \right) \varkappa _{P,X,Q}^{H}
&=&\left( p\bullet \left( \chi \circ q\right) \right) \left( \mu _{P}^{\circ
}\bullet \left( X\circ Q\right) \right) \zeta _{P,X,H,Q}\left( \left(
P\bullet X\right) \circ \rho _{Q}^{\bullet }\right)  \\
&=&\left( \mu _{P^{\prime }}^{\circ }\bullet \left( X^{\prime }\circ
Q^{\prime }\right) \right) \left( \left( p\circ H\right) \bullet \left( \chi
\circ q\right) \right) \zeta _{P,X,H,Q}\left( \left( P\bullet X\right) \circ
\rho _{Q}^{\bullet }\right)  \\
&=&\left( \mu _{P^{\prime }}^{\circ }\bullet \left( X^{\prime }\circ
Q^{\prime }\right) \right) \zeta _{P^{\prime },X^{\prime },H,Q^{\prime
}}\left( \left( p\bullet \chi \right) \circ \left( H\bullet q\right) \right)
\left( \left( P\bullet X\right) \circ \rho _{Q}^{\bullet }\right)  \\
&=&\left( \mu _{P^{\prime }}^{\circ }\bullet \left( X^{\prime }\circ
Q^{\prime }\right) \right) \zeta _{P^{\prime },X^{\prime },H,Q^{\prime
}}\left( \left( P^{\prime }\bullet X^{\prime }\right) \circ \rho _{Q^{\prime
}}^{\bullet }\right) \left( \left( p\bullet \chi \right) \circ q\right)  \\
&=&\varkappa _{P^{\prime },X^{\prime },Q^{\prime }}^{H}\left( \left(
p\bullet \chi \right) \circ q\right)
\end{eqnarray*}%
and
\begin{eqnarray*}
\left( \left( p\bullet \chi \right) \circ q\right) \overline{\varkappa }%
_{P,X,Q}^{H} &=&\left( \left( p\bullet \chi \right) \circ q\right) \left(
\left( \mu _{P}^{\circ }\bullet X\right) \circ Q\right) \delta
_{P,X,H,Q}\left( P\bullet \left( X\circ \left( \Sigma _{H}\bullet Q\right)
\rho _{Q}^{\bullet }\right) \right)  \\
&=&\left( \left( \mu _{P^{\prime }}^{\circ }\bullet X^{\prime }\right) \circ
Q^{\prime }\right) \left( \left( \left( p\circ H\right) \bullet \chi \right)
\circ q\right) \delta _{P,X,H,Q}\left( P\bullet \left( X\circ \left( \Sigma
_{H}\bullet Q\right) \rho _{Q}^{\bullet }\right) \right)  \\
&=&\left( \left( \mu _{P^{\prime }}^{\circ }\bullet X^{\prime }\right) \circ
Q^{\prime }\right) \delta _{P^{\prime },X^{\prime },H,Q^{\prime }}\left(
p\bullet \left( \chi \circ \left( H^{\revsn }\bullet q\right) \right) \right)
\left( P\bullet \left( X\circ \left( \Sigma _{H}\bullet Q\right) \rho
_{Q}^{\bullet }\right) \right)  \\
&=&\left( \left( \mu _{P^{\prime }}^{\circ }\bullet X^{\prime }\right) \circ
Q^{\prime }\right) \delta _{P^{\prime },X^{\prime },H,Q^{\prime }}\left(
P^{\prime }\bullet \left( X^{\prime }\circ \left( \Sigma _{H}\bullet
Q^{\prime }\right) \rho _{Q^{\prime }}^{\bullet }\right) \right) \left(
p\bullet \left( \chi \circ q\right) \right)  \\
&=&\overline{\varkappa }_{P^{\prime },X^{\prime },Q^{\prime }}^{H}\left(
p\bullet \left( \chi \circ q\right) \right)
\end{eqnarray*}%
By symmetric argument one prove the naturality of co-Galois maps.
\end{invisible}
Note that the co-Galois maps are just the Galois maps written in $\Cc^\mathrm{rev}:=(\Cc,\circ^\mathrm{rev},\bullet^\mathrm{rev})$. Therefore, the aforementioned symmetry principle, which allows us to prove results from only one side, still holds.

\begin{remark}
\label{rmk:Galunit} We may set $\Sigma _{I}:=\phi^\circ
_{0}:I\rightarrow I^{\revsn }.$ With this choice, we can always consider $\varkappa
_{P,X,Q}^{I}$ and $\overline{\varkappa }_{P,X,Q}^{I}.$ Moreover, in this
case, we necessarily have $\mu _{P}^{\circ }=r_{P}^{\circ }$ so that the
corresponding $P$ is just an object in $\Cc$.

Similarly, we can consider $\varsigma _{P,X,Q}^{I}$ and $\overline{\varsigma}%
_{P,X,Q}^{I}$ where $\mu _{Q}^{\circ }=l_{Q}^{\circ }$ and hence $Q$ is just
an object in $\Cc$.
\end{remark}

We aim to prove that the existence of an antipode ensures the invertibility of both the Galois and co-Galois maps. This will be achieved by means of some technical lemmas.

\begin{lemma}
\label{lem:inv} We have the identities $\varkappa _{X,Y,I}^{I}=\left(
X\bullet \left( r_{Y}^{\circ }\right) ^{-1}\right) r_{X\bullet Y}^{\circ }$
and $\varsigma _{I,X,Y}^{I}=\left( \left( l_{X}^{\circ }\right) ^{-1}\bullet
Y\right) l_{X\bullet Y}^{\circ }$ so that $\varkappa _{X,Y,I}^{I}$ and $%
\varsigma _{I,X,Y}^{I}$ are both invertible.
\end{lemma}

\begin{proof}
Using the right $I$-module $\left( X,r_{X}^{\circ }\right) ,$ we compute%
\begin{eqnarray*}
\varkappa _{X,Y,I}^{I} &=&\left( r_{X}^{\circ }\bullet \left( Y\circ
I\right) \right) \underbracket[0.140ex]{\zeta _{X,Y,I,I}\left( \left( X\bullet Y\right) \circ
\Delta _{I}^{\bullet }\right)}  %\\&=&%\left( X\bullet \left( r_{Y}^{\circ }\right) ^{-1}\right) \underbracket[0.140ex]{%
%\left( r_{X}^{\circ }\bullet r_{Y}^{\circ }\right) \zeta _{X,Y,I,I}\left(
%\left( X\bullet Y\right) \circ \Delta _{I}^{\bullet }\right) }
\overset{\eqref{eq:unit1}}=\left(
X\bullet \left( r_{Y}^{\circ }\right) ^{-1}\right) r_{X\bullet Y}^{\circ }.
\end{eqnarray*}%
The proof for $\varsigma _{I,X,Y}^{I}$ follows from a symmetric argument.
\begin{invisible}
\begin{eqnarray*}
\varsigma _{I,X,Y}^{I} &=&\left( \left( I\circ X\right) \bullet l_{Y}^{\circ
}\right) \zeta _{I,I,X,Y}\left( \Delta _{I}^{\bullet }\circ \left( X\bullet
Y\right) \right)  \\
&=&\left( \left( l_{X}^{\circ }\right) ^{-1}\bullet Y\right) \underbracket[0.140ex]{%
\left( l_{X}^{\circ }\bullet l_{Y}^{\circ }\right) \zeta _{I,I,X,Y}\left(
\Delta _{I}^{\bullet }\circ \left( X\bullet Y\right) \right) }\overset{\eqref{eq:unit1}}{=}\left( \left(
l_{X}^{\circ }\right) ^{-1}\bullet Y\right) l_{X\bullet Y}^{\circ }
\end{eqnarray*}\end{invisible}
\end{proof}

\begin{lemma}%[\rd{added:2026/02/21}]
The following squares are serially commutative.
{\small \begin{equation}
\xymatrix{(X\bullet Y)\circ (H\bullet Z)\circ Q\ar[d]_-{ \zeta _{X,Y,H,Z}\circ Q}\ar@<.5ex>[r]^-{\left( X\bullet Y\right) \circ \varkappa
_{H,Z,Q}^{H}}&\ar@<.6ex>[l]^-{\left( X\bullet Y\right) \circ    \overline{\varkappa}
_{H,Z,Q}^{H}} (X\bullet Y)\circ (H\bullet (Z\circ Q))\ar[d]^-{\zeta _{X,Y,H,Z\circ Q}}\\
((X\circ H)\bullet(Y\circ Z))\circ Q\ar@<.5ex>[r]^-{\varkappa _{X\circ H,Y\circ Z,Q}^{H}} & (X\circ H)\bullet(Y\circ Z\circ Q)\ar@<.6ex>[l]^-{\overline{\varkappa} _{X\circ H,Y\circ Z,Q}^{H}}
} \label{from:varkappazeta}
\quad \xymatrix{
  {\left( P\bullet X\bullet Y\right) \circ \left( H\bullet Z\right)}
    \ar@<0.6ex>[r]^{\varkappa_{P,X\bullet Y,H\bullet Z}^{H}}
    \ar[d]_{\zeta _{P\bullet X,Y,H,Z}}
  &
  {P\bullet \left( \left( X\bullet Y\right) \circ \left( H\bullet Z\right) \right)}
    \ar[d]^{P\bullet \zeta_{X,Y,H,Z}}\ar@<0.6ex>[l]^{\overline{\varkappa }_{P,X\bullet Y,H\bullet Z}^{H}}
  \\
  {\left( \left( P\bullet X\right) \circ H\right) \bullet \left( Y\circ Z\right)}
    \ar@<0.6ex>[r]^{\varkappa _{P,X,H}^{H}\bullet \left( Y\circ Z\right)}
  &
  {P\bullet \left( X\circ H\right) \bullet \left( Y\circ Z\right)}\ar@<0.6ex>[l]^{\overline{\varkappa }_{P,X,H}^{H}\bullet \left( Y\circ Z\right)}
}   \end{equation}
\begin{equation}\label{from:varsigmazeta}
\xymatrix{
P\circ (X\bullet H)\circ (Y\bullet Z)
  \ar@<0.6ex>[r]^{\varsigma^{H}_{P,X,H}\circ (Y\bullet Z)}
    \ar[d]_{P\circ \zeta_{X,H,Y,Z}}
&
((P\circ X)\bullet H)\circ (Y\bullet Z)\ar@<0.6ex>[l]^{\overline{\varsigma}^{H}_{P,X,H}\circ (Y\bullet Z)}
  \ar[d]^{\zeta_{P\circ X,H,Y,Z}}
\\
P\circ\big((X\circ Y)\bullet (H\circ Z)\big)
  \ar@<0.6ex>[r]^{\varsigma^{H}_{P,X\circ Y,H\circ Z}}
  &
(P\circ X\circ Y)\bullet (H\circ Z)\ar@<0.6ex>[l]^{\overline{\varsigma}^{H}_{P,X\circ Y,H\circ Z}}
}\quad
\xymatrix{
 (X \bullet H)\circ (Y \bullet Z \bullet Q)
   \ar@<0.6ex>[r]^{\varsigma^{H}_{X\bullet H,\,Y\bullet Z,\,Q}}
      \ar[d]_{\zeta_{X,H,Y,Z\bullet Q}}
 &
 ((X\bullet H)\circ (Y\bullet Z)) \bullet Q\ar@<0.6ex>[l]^{\overline{\varsigma}^{H}_{X\bullet H,\,Y\bullet Z,\,Q}}
   \ar[d]^{\zeta_{X,H,Y,Z}\bullet Q}
 \\
 (X\circ Y)\bullet (H\circ (Z\bullet Q))
   \ar@<0.5ex>[r]^{(X\circ Y)\bullet \varsigma^{H}_{H,Z,Q}}
   &
 (X\circ Y)\bullet (H\circ Z)\bullet Q
 \ar@<0.5ex>[l]^{(X\circ Y)\bullet \overline{\varsigma}^{H}_{H,Z,Q}}}
\end{equation}
}
\end{lemma}

\begin{proof}
Concerning the first square, by definition of $\varkappa$, naturality of $\zeta$ and \eqref{eq:assoc1} one proves the commutativity of the part involving $\varkappa$, while by definition of $\overline{\varkappa }$, naturality of $\zeta$ and \eqref{form:varphi1} one proves the commutativity of the part involving $\overline{\varkappa }$.
\begin{invisible}
We compute%
\begin{eqnarray*}
&&\varkappa _{X\circ H,Y\circ Z,Q}^{H}\left( \zeta _{X,Y,H,Z}\circ Q\right)\\
&=&\left( \left( X\circ m_{H}^{\circ }\right) \bullet \left( Y\circ Z\circ
Q\right) \right) \zeta _{X\circ H,Y\circ Z,H,Q}\left( \left( \left( X\circ
H\right) \bullet \left( Y\circ Z\right) \right) \circ \rho _{Q}^{\bullet
}\right) \left( \zeta _{X,Y,H,Z}\circ Q\right)  \\
&=&\left( \left( X\circ m_{H}^{\circ }\right) \bullet \left( Y\circ Z\circ
Q\right) \right) \underbracket[0.140ex]{\zeta _{X\circ H,Y\circ Z,H,Q}\left( \zeta
_{X,Y,H,Z}\circ \left( H\bullet Q\right) \right) }\left( \left( X\bullet
Y\right) \circ \left( H\bullet Z\right) \circ \rho _{Q}^{\bullet }\right)  \\
&\overset{\eqref{eq:assoc1}}=&\left( \left( X\circ m_{H}^{\circ }\right) \bullet \left( Y\circ Z\circ
Q\right) \right) \zeta _{X,Y,H\circ H,Z\circ Q}\left( \left( X\bullet
Y\right) \circ \zeta _{H,Z,H,Q}\right) \left( \left( X\bullet Y\right) \circ
\left( H\bullet Z\right) \circ \rho _{Q}^{\bullet }\right)  \\
&=&\zeta _{X,Y,H,Z\circ Q}\left( \left( X\bullet Y\right) \circ \left(
m_{H}^{\circ }\bullet \left( Z\circ Q\right) \right) \right) \left( \left(
X\bullet Y\right) \circ \zeta _{H,Z,H,Q}\right) \left( \left( X\bullet
Y\right) \circ \left( H\bullet Z\right) \circ \rho _{Q}^{\bullet }\right)  \\
&=&\zeta _{X,Y,H,Z\circ Q}\left( \left( X\bullet Y\right) \circ \varkappa
_{H,Z,Q}^{H}\right)
\end{eqnarray*}%
and% 
\begin{eqnarray*}
&&\overline{\varkappa }_{X\circ H,Y\circ Z,Q}^{H}\zeta _{X,Y,H,Z\circ Q}\\
&=&\left( \left( \left( X\circ m_{H}^{\circ }\right) \bullet \left( Y\circ
Z\right) \right) \circ Q\right) \delta _{X\circ H,Y\circ Z,H,Q}\left( \left(
X\circ H\right) \bullet \left( Y\circ Z\circ \left( \Sigma _{H}\bullet
Q\right) \rho _{Q}^{\bullet }\right) \right) \zeta _{X,Y,H,Z\circ Q} \\
&=&\left( \left( \left( X\circ m_{H}^{\circ }\right) \bullet \left( Y\circ
Z\right) \right) \circ Q\right) \underbracket[0.140ex]{\delta _{X\circ H,Y\circ
Z,H,Q}\zeta _{X,Y,H,Z\circ \left( H^{\revsn }\bullet Q\right) }}\left( \left(
X\bullet Y\right) \circ \left( H\bullet \left( Z\circ \left( \Sigma
_{H}\bullet Q\right) \rho _{Q}^{\bullet }\right) \right) \right)  \\
&\overset{\eqref{form:varphi1} }{=}&\left( \left( \left(
X\circ m_{H}^{\circ }\right) \bullet \left( Y\circ Z\right) \right) \circ
Q\right) \left( \zeta _{X,Y,H\circ H,Z}\circ Q\right) \left( \left( X\bullet
Y\right) \circ \delta _{H,Z,H,Q}\right) \left( \left( X\bullet Y\right)
\circ \left( H\bullet \left( Z\circ \left( \Sigma _{H}\bullet Q\right) \rho
_{Q}^{\bullet }\right) \right) \right)  \\
&=&\left( \zeta _{X,Y,H,Z}\circ Q\right) \left( \left( X\bullet Y\right)
\circ \left( m_{H}^{\circ }\bullet Z\right) \circ Q\right) \left( \left(
X\bullet Y\right) \circ \delta _{H,Z,H,Q}\right) \left( \left( X\bullet
Y\right) \circ \left( H\bullet \left( Z\circ \left( \Sigma _{H}\bullet
Q\right) \rho _{Q}^{\bullet }\right) \right) \right)  \\
&=&\left( \zeta _{X,Y,H,Z}\circ Q\right) \left( \left( X\bullet Y\right)
\circ \overline{\varkappa }_{H,Z,Q}^{H}\right)
\end{eqnarray*}%
\end{invisible}
The remaining squares are treated in a similar way by using
\eqref{eq:assoc1},\eqref{eq:assoc2}, \eqref{form:psi1}, \eqref{form:psi2} and \eqref{form:varphi2}.
\begin{invisible}
 Similarly%
\begin{eqnarray*}
\left( P\bullet \zeta _{X,Y,H,Z}\right) \varkappa _{P,X\bullet Y,H\bullet
Z}^{H} &=&\left( P\bullet \zeta _{X,Y,H,Z}\right) \left( \mu _{P}^{\circ
}\bullet \left( \left( X\bullet Y\right) \circ \left( H\bullet Z\right)
\right) \right) \zeta _{P,X\bullet Y,H,H\bullet Z}\left( \left( P\bullet
X\bullet Y\right) \circ \left( \Delta _{H}^{\bullet }\bullet Z\right)
\right)  \\
&=&\left( \mu _{P}^{\circ }\bullet \left( X\circ H\right) \bullet \left(
Y\circ Z\right) \right) \underbracket[0.140ex]{\left( P\bullet \zeta _{X,Y,H,Z}\right)
\zeta _{P,X\bullet Y,H,H\bullet Z}}\left( \left( P\bullet X\bullet Y\right)
\circ \left( \Delta _{H}^{\bullet }\bullet Z\right) \right)  \\
&\overset{\eqref{eq:assoc2}}=&\left( \mu _{P}^{\circ }\bullet \left( X\circ H\right) \bullet \left(
Y\circ Z\right) \right) \left( \zeta _{P,X,H,H}\bullet \left( Y\circ
Z\right) \right) \zeta _{P\bullet X,Y,H\bullet H,Z}\left( \left( P\bullet
X\bullet Y\right) \circ \left( \Delta _{H}^{\bullet }\bullet Z\right)
\right)  \\
&=&\left( \mu _{P}^{\circ }\bullet \left( X\circ H\right) \bullet \left(
Y\circ Z\right) \right) \left( \zeta _{P,X,H,H}\bullet \left( Y\circ
Z\right) \right) \left( \left( \left( P\bullet X\right) \circ \Delta
_{H}^{\bullet }\right) \bullet \left( Y\circ Z\right) \right) \zeta
_{P\bullet X,Y,H,Z} \\
&=&\left( \varkappa _{P,X,H}^{H}\bullet \left( Y\circ Z\right) \right) \zeta
_{P\bullet X,Y,H,Z}
\end{eqnarray*}%
and%
\begin{eqnarray*}
\zeta _{P\bullet X,Y,H,Z}\overline{\varkappa }_{P,X\bullet Y,H\bullet Z}^{H}
&=&\zeta _{P\bullet X,Y,H,Z}\left( \left( \mu _{P}^{\circ }\bullet X\bullet
Y\right) \circ \left( H\bullet Z\right) \right) \delta _{P,X\bullet
Y,H,H\bullet Z}\left( P\bullet \left( \left( X\bullet Y\right) \circ \left(
\Sigma _{H}\bullet H\bullet Z\right) \left( \Delta _{H}^{\bullet }\bullet
Z\right) \right) \right)  \\
&=&\left( \left( \left( \mu _{P}^{\circ }\bullet X\right) \circ H\right)
\bullet \left( Y\circ Z\right) \right) \underbracket[0.140ex]{\zeta _{\left( P\circ
H\right) \bullet X,Y,H,Z}\delta _{P,X\bullet Y,H,H\bullet Z}}\\ & &\left( P\bullet
\left( \left( X\bullet Y\right) \circ \left( \Sigma _{H}\bullet H\bullet
Z\right) \left( \Delta _{H}^{\bullet }\bullet Z\right) \right) \right)  \\
&\overset{\eqref{form:varphi2}}=&\left( \left( \left( \mu _{P}^{\circ }\bullet X\right) \circ H\right)
\bullet \left( Y\circ Z\right) \right) \left( \delta _{P,X,H,H}\bullet
\left( Y\circ Z\right) \right) \left( P\bullet \zeta _{X,Y,H^{\revsn }\bullet
H,Z}\right) \\ & &\left( P\bullet \left( \left( X\bullet Y\right) \circ \left(
\left( \Sigma _{H}\bullet H\right) \Delta _{H}^{\bullet }\bullet Z\right)
\right) \right)  \\
&=&\left( \left( \left( \mu _{P}^{\circ }\bullet X\right) \circ H\right)
\bullet \left( Y\circ Z\right) \right) \left( \delta _{P,X,H,H}\bullet
\left( Y\circ Z\right) \right) \\& &\left( P\bullet \left( \left( X\circ \left(
\Sigma _{H}\bullet H\right) \Delta _{H}^{\bullet }\right) \bullet \left(
Y\circ Z\right) \right) \right) \left( P\bullet \zeta _{X,Y,H,Z}\right)  \\
&=&\overline{\varkappa }_{P,X,H}^{H}\bullet \left( Y\circ Z\right) \left(
P\bullet \zeta _{X,Y,H,Z}\right).
\end{eqnarray*}%

We also compute%
\begin{eqnarray*}
&&\varsigma _{P,X\circ Y,H\circ Z}^{H}\left( P\circ \zeta _{X,H,Y,Z}\right)
\\
&=&\left( \left( P\circ X\circ Y\right) \bullet \left( m_{H}^{\circ }\circ
Z\right) \right) \zeta _{P,H,X\circ Y,H\circ Z}\left( \rho _{P}^{\bullet
}\circ \left( \left( X\circ Y\right) \bullet \left( H\circ Z\right) \right)
\right) \left( P\circ \zeta _{X,H,Y,Z}\right)  \\
&=&\left( \left( P\circ X\circ Y\right) \bullet \left( m_{H}^{\circ }\circ
Z\right) \right) \underbracket[0.140ex]{\zeta _{P,H,X\circ Y,H\circ Z}\left( \left(
P\bullet H\right) \circ \zeta _{X,H,Y,Z}\right) }\left( \rho _{P}^{\bullet
}\circ \left( X\bullet H\right) \circ \left( Y\bullet Z\right) \right)  \\
&\overset{\eqref{eq:assoc1}}=&\left( \left( P\circ X\circ Y\right) \bullet \left( m_{H}^{\circ }\circ
Z\right) \right) \zeta _{P\circ X,H\circ H,Y,Z}\left( \zeta _{P,H,X,H}\circ
\left( Y\bullet Z\right) \right) \left( \left( \rho _{P}^{\bullet }\circ
\left( X\bullet H\right) \right) \circ \left( Y\bullet Z\right) \right)  \\
&=&\zeta _{P\circ X,H,Y,Z}\left( \left( \left( P\circ X\right) \bullet
m_{H}^{\circ }\right) \circ \left( Y\bullet Z\right) \right) \left( \zeta
_{P,H,X,H}\circ \left( Y\bullet Z\right) \right) \left( \left( \rho
_{P}^{\bullet }\circ \left( X\bullet H\right) \right) \circ \left( Y\bullet
Z\right) \right)  \\
&=&\zeta _{P\circ X,H,Y,Z}\left( \varsigma _{P,X,H}^{H}\circ \left( Y\bullet
Z\right) \right)
\end{eqnarray*}%
and%
\begin{eqnarray*}
\overline{\varsigma }_{P,X\circ Y,H\circ Z}^{H}\zeta _{P\circ X,H,Y,Z}
&=&\left( P\circ \left( \left( X\circ Y\right) \bullet \left( m_{H}^{\circ
}\circ Z\right) \right) \right) \gamma _{P,H,X\circ Y,H\circ Z}\left( \left(
\left( P\bullet \Sigma _{H}\right) \rho _{P}^{\bullet }\circ X\circ Y\right)
\bullet \left( H\circ Z\right) \right) \zeta _{P\circ X,H,Y,Z} \\
&=&\left( P\circ \left( \left( X\circ Y\right) \bullet \left( m_{H}^{\circ
}\circ Z\right) \right) \right) \underbracket[0.140ex]{\gamma _{P,H,X\circ Y,H\circ
Z}\zeta _{\left( P\bullet H^{\revsn }\right) \circ X,H,Y,Z}}\\ & &\left( \left(
\left( \left( P\bullet \Sigma _{H}\right) \rho _{P}^{\bullet }\circ X\right)
\bullet H\right) \circ \left( Y\bullet Z\right) \right)  \\
&\overset{\eqref{form:psi1}}=&\left( P\circ \left( \left( X\circ Y\right) \bullet \left( m_{H}^{\circ
}\circ Z\right) \right) \right) \left( P\circ \zeta _{X,H\circ H,Y,Z}\right)
\left( \gamma _{P,H,X,H}\circ \left( Y\bullet Z\right) \right)\\ & & \left( \left(
\left( \left( P\bullet \Sigma _{H}\right) \rho _{P}^{\bullet }\circ X\right)
\bullet H\right) \circ \left( Y\bullet Z\right) \right)  \\
&=&\left( P\circ \zeta _{X,H,Y,Z}\right) \left( P\circ \left( X\bullet
m_{H}^{\circ }\right) \circ \left( Y\bullet Z\right) \right) \left( \gamma
_{P,H,X,H}\circ \left( Y\bullet Z\right) \right)\\ & & \left( \left( \left( \left(
P\bullet \Sigma _{H}\right) \rho _{P}^{\bullet }\circ X\right) \bullet
H\right) \circ \left( Y\bullet Z\right) \right)  \\
&=&\left( P\circ \zeta _{X,H,Y,Z}\right) \left( \overline{\varsigma }%
_{P,X,H}^{H}\circ \left( Y\bullet Z\right) \right).
\end{eqnarray*}%

We also compute%
\begin{eqnarray*}
&&\left( \zeta _{X,H,Y,Z}\bullet Q\right) \varsigma _{X\bullet H,Y\bullet
Z,Q}^{H} \\
&=&\left( \zeta _{X,H,Y,Z}\bullet Q\right) \left( \left( \left( X\bullet
H\right) \circ \left( Y\bullet Z\right) \right) \bullet \mu _{Q}^{\circ
}\right) \zeta _{X\bullet H,H,Y\bullet Z,Q}\left( \left( X\bullet \Delta
_{H}^{\bullet }\right) \circ \left( Y\bullet Z\bullet Q\right) \right)  \\
&=&\left( \left( X\circ Y\right) \bullet \left( H\circ Z\right) \bullet \mu
_{Q}^{\circ }\right) \underbracket[0.140ex]{\left( \zeta _{X,H,Y,Z}\bullet \left(
H\circ Q\right) \right) \zeta _{X\bullet H,H,Y\bullet Z,Q}}\left( \left(
X\bullet \Delta _{H}^{\bullet }\right) \circ \left( Y\bullet Z\bullet
Q\right) \right)  \\
&\overset{\eqref{eq:assoc2}}=&\left( \left( X\circ Y\right) \bullet \left( H\circ Z\right) \bullet \mu
_{Q}^{\circ }\right) \left( \left( X\circ Y\right) \bullet \zeta
_{H,H,Z,Q}\right) \zeta _{X,H\bullet H,Y,Z\bullet Q}\left( \left( X\bullet
\Delta _{H}^{\bullet }\right) \circ \left( Y\bullet Z\bullet Q\right)
\right)  \\
&=&\left( \left( X\circ Y\right) \bullet \left( H\circ Z\right) \bullet \mu
_{Q}^{\circ }\right) \left( \left( X\circ Y\right) \bullet \zeta
_{H,H,Z,Q}\right) \left( \left( X\circ Y\right) \bullet \left( \Delta
_{H}^{\bullet }\circ \left( Z\bullet Q\right) \right) \right) \zeta
_{X,H,Y,Z\bullet Q} \\
&=&\left( \left( X\circ Y\right) \bullet \varsigma _{H,Z,Q}^{H}\right) \zeta
_{X,H,Y,Z\bullet Q}
\end{eqnarray*}

and%
\begin{eqnarray*}
\zeta _{X,H,Y,Z\bullet Q}\overline{\varsigma }_{X\bullet H,Y\bullet Z,Q}^{H}
&=&\zeta _{X,H,Y,Z\bullet Q}\left( \left( X\bullet H\right) \circ \left(
Y\bullet Z\bullet \mu _{Q}^{\circ }\right) \right) \gamma _{X\bullet
H,H,Y\bullet Z,Q}\\ & &\left( \left( \left( X\bullet H\bullet \Sigma _{H}\right)
\left( X\bullet \Delta _{H}^{\bullet }\right) \circ \left( Y\bullet Z\right)
\right) \bullet Q\right)  \\
&=&\left( \left( X\circ Y\right) \bullet \left( H\circ \left( Z\bullet \mu
_{Q}^{\circ }\right) \right) \right) \underbracket[0.140ex]{\zeta _{X,H,Y,Z\bullet
\left( H\circ Q\right) }\gamma _{X\bullet H,H,Y\bullet Z,Q}}\\ & &\left( \left(
\left( X\bullet \left( H\bullet \Sigma _{H}\right) \Delta _{H}^{\bullet
}\right) \circ \left( Y\bullet Z\right) \right) \bullet Q\right)  \\
&\overset{\eqref{form:psi2}}=&\left( \left( X\circ Y\right) \bullet \left( H\circ \left( Z\bullet \mu
_{Q}^{\circ }\right) \right) \right) \left( \left( X\circ Y\right) \bullet
\gamma _{H,H,Z,Q}\right)\\ & & \left( \zeta _{X,H\bullet H^{\revsn },Y,Z}\bullet
Q\right) \left( \left( \left( X\bullet \left( H\bullet \Sigma _{H}\right)
\Delta _{H}^{\bullet }\right) \circ \left( Y\bullet Z\right) \right) \bullet
Q\right)  \\
&=&\left( \left( X\circ Y\right) \bullet \left( H\circ \left( Z\bullet \mu
_{Q}^{\circ }\right) \right) \right) \left( \left( X\circ Y\right) \bullet
\gamma _{H,H,Z,Q}\right)\\ & & \left( \left( X\circ Y\right) \bullet \left( \left(
H\bullet \Sigma _{H}\right) \Delta _{H}^{\bullet }\circ Z\right) \bullet
Q\right) \left( \zeta _{X,H,Y,Z}\bullet Q\right)  \\
&=&\left( \left( X\circ Y\right) \bullet \overline{\varsigma }%
_{H,Z,Q}^{H}\right) \left( \zeta _{X,H,Y,Z}\bullet Q\right).
\end{eqnarray*}%
\end{invisible}
\end{proof}

\begin{lemma}%[\rd{added:2026/03/30}]
Let $\left( \mathcal{C},\circ ,I,\bullet ,J\right) $ be a duoidal category
with a reversion and let $\left( H,m_{H}^{\circ },u_{H}^{\circ },\Delta
_{H}^{\bullet },\varepsilon _{H}^{\bullet }\right) $ be a bimonoid. Let $%
\Sigma _{H}:H\rightarrow H^{\text{-}}$ be a morphism in $\Cc$ and consider the Galois and co-Galois
maps of \cref{def:Galmap}.\\
For any right $H$-module $\left( P,\mu
_{P}^{\circ }\right) $ and left $H$-comodule $\left( Q,\rho _{Q}^{\bullet
}\right) $, we have%
\begin{eqnarray}
\varkappa _{H,I,Q}^{H}\left( \left( u_{H}^{\circ }\bullet I\right) \Delta
_{I}^{\bullet }\circ Q\right) \left( l_{Q}^{\circ }\right) ^{-1} &=&\left(
H\bullet \left( l_{Q}^{\circ }\right) ^{-1}\right) \rho _{Q}^{\bullet },
\label{form:Galunit} \\
r_{P}^{\bullet }\left( P\bullet m_{J}^{\circ }\left( J\circ \varepsilon
_{H}^{\bullet }\right) \right) \varkappa _{P,J,H}^{H} &=&\mu _{P}^{\circ
}\left( r_{P}^{\bullet }\circ H\right) ,  \label{form:Galcounit}\end{eqnarray}
For any right $H$-comodule $(P,\rho_P^\bullet)$ and  left $H$-module $(Q,\mu_Q^\circ)$, we have%
\begin{eqnarray}
\varsigma _{P,I,H}^{H}\left( P\circ \left( I\bullet u_{H}^{\circ }\right)
\Delta _{I}^{\bullet }\right) \left( r_{P}^{\circ }\right) ^{-1} &=&\left(
\left( r_{P}^{\circ }\right) ^{-1}\bullet H\right) \rho _{P}^{\bullet },
\label{form:coGalunit} \\
l_{Q}^{\bullet }\left( m_{J}^{\circ }\left( \varepsilon _{H}^{\bullet }\circ
J\right) \bullet Q\right) \varsigma _{H,J,Q}^{H} &=&\mu _{Q}^{\circ }\left(
H\circ l_{Q}^{\bullet }\right) .  \label{form:coGalcounit}
\end{eqnarray}
\end{lemma}

\begin{proof}
By definition of $\varkappa$, naturality of $\zeta$, unitality of $m_{H}^{\circ }$ and \eqref{eq:unit1}, one easily checks \eqref{form:Galunit}, while by definition of $\varkappa$, naturality of $\zeta$, counitality of $\Delta_{H}^{\bullet}$ and \eqref{eq:unit2}, one gets \eqref{form:Galcounit}.
\begin{invisible}
We compute%
\begin{eqnarray*}
&&\varkappa _{H,I,Q}^{H}\left( \left( u_{H}^{\circ }\bullet I\right) \Delta
_{I}^{\bullet }\circ Q\right) \left( l_{Q}^{\circ }\right) ^{-1} \\
&=&\left( m_{H}^{\circ }\bullet \left( I\circ Q\right) \right) \zeta
_{H,I,H,Q}\left( \left( H\bullet I\right) \circ \rho _{Q}^{\bullet }\right)
\left( \left( u_{H}^{\circ }\bullet I\right) \Delta _{I}^{\bullet }\circ
Q\right) \left( l_{Q}^{\circ }\right) ^{-1} \\
&=&\left( m_{H}^{\circ }\left( u_{H}^{\circ }\circ H\right) \bullet \left(
I\circ Q\right) \right) \zeta _{I,I,H,Q}\left( \left( I\bullet I\right)
\circ \rho _{Q}^{\bullet }\right) \left( \Delta _{I}^{\bullet }\circ
Q\right) \left( l_{Q}^{\circ }\right) ^{-1} \\
&=&\left( H\bullet \left( l_{Q}^{\circ }\right) ^{-1}\right) \left(
l_{H}^{\circ }\bullet l_{Q}^{\circ }\right) \zeta _{I,I,H,Q}\left( \Delta
_{I}^{\bullet }\circ \left( H\bullet Q\right) \right) \left( I\circ \rho
_{Q}^{\bullet }\right) \left( l_{Q}^{\circ }\right) ^{-1} \\
&=&\left( H\bullet \left( l_{Q}^{\circ }\right) ^{-1}\right) %
\underbracket[0.140ex]{\left( l_{H}^{\circ }\bullet l_{Q}^{\circ }\right)
\zeta _{I,I,H,Q}\left( \Delta _{I}^{\bullet }\circ \left( H\bullet Q\right)
\right) \left( l_{H\bullet Q}^{\circ }\right) ^{-1}}\rho _{Q}^{\bullet } \overset{\eqref{eq:unit1}}=\left( H\bullet \left( l_{Q}^{\circ }\right) ^{-1}\right) \rho
_{Q}^{\bullet }
\end{eqnarray*}%
and%
\begin{eqnarray*}
&&r_{P}^{\bullet }\left( P\bullet m_{J}^{\circ }\left( J\circ \varepsilon
_{H}^{\bullet }\right) \right) \varkappa _{P,J,H}^{H} \\
&=&r_{P}^{\bullet }\left( P\bullet m_{J}^{\circ }\left( J\circ \varepsilon
_{H}^{\bullet }\right) \right) \left( \mu _{P}^{\circ }\bullet \left( J\circ
H\right) \right) \zeta _{P,J,H,H}\left( \left( P\bullet J\right) \circ
\Delta _{H}^{\bullet }\right)  \\
&=&r_{P}^{\bullet }\left( P\bullet m_{J}^{\circ }\right) \left( \mu
_{P}^{\circ }\bullet \left( J\circ J\right) \right) \zeta _{P,J,H,J}\left(
\left( P\bullet J\right) \circ \left( H\bullet \varepsilon _{H}^{\bullet
}\right) \Delta _{H}^{\bullet }\right)  \\
&=&r_{P}^{\bullet }\left( \mu _{P}^{\circ }\bullet J\right) \left( \left(
P\circ H\right) \bullet m_{J}^{\circ }\right) \zeta _{P,J,H,J}\left( \left(
r_{P}^{\bullet }\right) ^{-1}\circ \left( r_{H}^{\bullet }\right)
^{-1}\right) \left( r_{P}^{\bullet }\circ H\right)  \\
&=&\mu _{P}^{\circ }\underbracket[0.140ex]{r_{P\circ H}^{\bullet }\left(
\left( P\circ H\right) \bullet m_{J}^{\circ }\right) \zeta _{P,J,H,J}\left(
\left( r_{P}^{\bullet }\right) ^{-1}\circ \left( r_{H}^{\bullet }\right)
^{-1}\right) }\left( r_{P}^{\bullet }\circ H\right)  \overset{\eqref{eq:unit2}}=\mu _{P}^{\circ }\left( r_{P}^{\bullet }\circ H\right) .
\end{eqnarray*}
\end{invisible}
Similarly, one proves the remaining equalities.
\begin{invisible}
We compute%
\begin{eqnarray*}
&&\varsigma _{P,I,H}^{H}\left( P\circ \left( I\bullet u_{H}^{\circ }\right)
\Delta _{I}^{\bullet }\right) \left( r_{P}^{\circ }\right) ^{-1}\\&=&\left(
\left( P\circ I\right) \bullet m_{H}^{\circ }\right) \zeta _{P,H,I,H}\left(
\rho _{P}^{\bullet }\circ \left( I\bullet H\right) \right) \left( P\circ
\left( I\bullet u_{H}^{\circ }\right) \Delta _{I}^{\bullet }\right) \left(
r_{P}^{\circ }\right) ^{-1} \\
&=&\left( \left( P\circ I\right) \bullet m_{H}^{\circ }\left( H\circ
u_{H}^{\circ }\right) \right) \zeta _{P,H,I,I}\left( \left( P\bullet
H\right) \circ \Delta _{I}^{\bullet }\right) \left( \rho _{P}^{\bullet
}\circ I\right) \left( r_{P}^{\circ }\right) ^{-1} \\
&=&\left( \left( r_{P}^{\circ }\right) ^{-1}\bullet H\right) %
\underbracket[0.140ex]{\left( r_{P}^{\circ }\bullet r_{H}^{\circ }\right)
\zeta _{P,H,I,I}\left( \left( P\bullet H\right) \circ \Delta _{I}^{\bullet
}\right) \left( r_{P\bullet H}^{\circ }\right) ^{-1}}\rho _{P}^{\bullet
}\overset{\eqref{eq:unit1}}=\left( \left( r_{P}^{\circ }\right) ^{-1}\bullet H\right) \rho
_{P}^{\bullet }
\end{eqnarray*}%
\begin{eqnarray*}
l_{Q}^{\bullet }\left( m_{J}^{\circ }\left( \varepsilon _{H}^{\bullet }\circ
J\right) \bullet Q\right) \varsigma _{H,J,Q}^{H} &=&l_{Q}^{\bullet }\left(
m_{J}^{\circ }\left( \varepsilon _{H}^{\bullet }\circ J\right) \bullet
Q\right) \left( \left( H\circ J\right) \bullet \mu _{Q}^{\circ }\right)
\zeta _{H,H,J,Q}\left( \Delta _{H}^{\bullet }\circ \left( J\bullet Q\right)
\right)  \\
&=&l_{Q}^{\bullet }\left( J\bullet \mu _{Q}^{\circ }\right) \left(
m_{J}^{\circ }\bullet \left( H\circ Q\right) \right) \left( \left(
\varepsilon _{H}^{\bullet }\circ J\right) \bullet \left( H\circ Q\right)
\right) \zeta _{H,H,J,Q}\left( \Delta _{H}^{\bullet }\circ \left( J\bullet
Q\right) \right)  \\
&=&\mu _{Q}^{\circ }l_{H\circ Q}^{\bullet }\left( m_{J}^{\circ }\bullet
\left( H\circ Q\right) \right) \zeta _{J,H,J,Q}\left( \left( \varepsilon
_{H}^{\bullet }\bullet H\right) \circ \left( J\bullet Q\right) \right)
\left( \Delta _{H}^{\bullet }\circ \left( J\bullet Q\right) \right)  \\
&=&\mu _{Q}^{\circ }l_{H\circ Q}^{\bullet }\left( m_{J}^{\circ }\bullet
\left( H\circ Q\right) \right) \zeta _{J,H,J,Q}\left( \left( l_{H}^{\bullet
}\right) ^{-1}\circ \left( l_{Q}^{\bullet }\right) ^{-1}\right) \left(
H\circ l_{Q}^{\bullet }\right) \overset{\eqref{eq:unit2}}=\mu _{Q}^{\circ }\left( H\circ
l_{Q}^{\bullet }\right)
\end{eqnarray*}
\end{invisible}
\end{proof}

The following result establishes preliminary conditions for the invertibility of the
Galois maps.

\begin{proposition}%[\rd{revised:2026/03/30}]
\label{prop:cogal2} Let $\left( \mathcal{C},\circ ,I,\bullet ,J\right) $ be
a duoidal category with a reversion and let $\left( H,m_{H}^{\circ
},u_{H}^{\circ },\Delta _{H}^{\bullet },\varepsilon _{H}^{\bullet }\right) $
be a bimonoid. Let $\Sigma _{H}:H\rightarrow H^{\text{-} }$ be a morphism in $\Cc$
and consider the Galois maps\textbf{\ }$\varkappa _{P,X,Q}^{H}$ and $%
\overline{\varkappa }_{P,X,Q}^{H}$ of \cref{def:Galmap} for any object $X$,
right $H$-module $\left( P,\mu _{P}^{\circ }\right) $ and left $H$-comodule $%
\left( Q,\rho _{Q}^{\bullet }\right)$.
%\as{[E' forse meglio inserire il "per ogni $P$, per ogni $Q$" in 1) e 2)?]}
The following statements hold:

1) $\overline{\varkappa }_{P,X,Q}^{H}\varkappa _{P,X,Q}^{H}=\mathrm{Id}$ if
and only if $\overline{\varkappa }_{H,I,Q}^{H}\varkappa _{H,I,Q}^{H}=\mathrm{%
Id}$ if and only if
\begin{equation}
\overline{\varkappa }_{H,I,Q}^{H}\varkappa _{H,I,Q}^{H}\left( \left(
u_{H}^{\circ }\bullet I\right) \Delta _{I}^{\bullet }\circ Q\right) = \left( u_{H}^{\circ }\bullet I\right)
\Delta _{I}^{\bullet }\circ Q.
\label{eq:GalinvGal}
\end{equation}%

2) $\varkappa _{P,X,Q}^{H}\overline{\varkappa }_{P,X,Q}^{H}=\mathrm{Id}$ if
and only if $\varkappa _{P,J,H}^{H}\overline{\varkappa }_{P,J,H}^{H}=\mathrm{%
Id}$ if and only if%
\begin{equation}
\left( P\bullet m_{J}^{\circ }\left( J\circ \varepsilon
_{H}^{\bullet }\right) \right) \varkappa _{P,J,H}^{H}\overline{\varkappa }%
_{P,J,H}^{H}= P\bullet m_{J}^{\circ }\left( J\circ
\varepsilon _{H}^{\bullet }\right) .  \label{eq:GalGalinv}
\end{equation}
\end{proposition}

\begin{proof}
1) Clearly $\overline{\varkappa }_{P,X,Q}^{H}\varkappa _{P,X,Q}^{H}=\mathrm{%
Id}$ for every $(P,\mu^{\circ}_{P}),X$ implies $\overline{\varkappa }_{H,I,Q}^{H}\varkappa
_{H,I,Q}^{H}=\mathrm{Id.}$ If the latter equality is true$\mathrm{,}$ it is
clear that \eqref{eq:GalinvGal} is true. If \eqref{eq:GalinvGal} is true, we
get $\overline{\varkappa }_{P,X,Q}^{H}\varkappa _{P,X,Q}^{H}=\mathrm{Id}$ as
\begin{align*}
\overline{\varkappa }_{P,X,Q}^{H}\varkappa _{P,X,Q}^{H} &=\overline{%
\varkappa }_{P,X,Q}^{H}\left( \mu _{P}^{\circ }\bullet \left( X\circ
Q\right) \right) \zeta _{P,X,H,Q}\left( \left( P\bullet X\right) \circ \rho
_{Q}^{\bullet }\right)  \\
&=\left( \left( \mu _{P}^{\circ }\bullet X\right) \circ Q\right) \overline{%
\varkappa }_{P\circ H,X,Q}^{H}\zeta _{P,X,H,Q}\left( \left( P\bullet
X\right) \circ \rho _{Q}^{\bullet }\right)  \\
&=\left( \left( \mu _{P}^{\circ }\bullet r_{X}^{\circ }\right) \circ
Q\right) \left( \left( \left( P\circ H\right) \bullet \left( r_{X}^{\circ
}\right) ^{-1}\right) \circ Q\right) \overline{\varkappa }_{P\circ
H,X,Q}^{H}\zeta _{P,X,H,Q}\left( \left( P\bullet X\right) \circ \rho
_{Q}^{\bullet }\right)  \\
&=\left( \left( \mu _{P}^{\circ }\bullet r_{X}^{\circ }\right) \circ
Q\right) \overline{\varkappa }_{P\circ H,X\circ I,Q}^{H}\left( \left( P\circ
H\right) \bullet \left( \left( r_{X}^{\circ }\right) ^{-1}\circ Q\right)
\right) \zeta _{P,X,H,Q}\left( \left( P\bullet X\right) \circ \rho
_{Q}^{\bullet }\right)  \\
&=\left( \left( \mu _{P}^{\circ }\bullet r_{X}^{\circ }\right) \circ
Q\right) \overline{\varkappa }_{P\circ H,X\circ I,Q}^{H}\left( \left( P\circ
H\right) \bullet \left( X\circ \left( l_{Q}^{\circ }\right) ^{-1}\right)
\right) \zeta _{P,X,H,Q}\left( \left( P\bullet X\right) \circ \rho
_{Q}^{\bullet }\right)  \\
&=\left( \left( \mu _{P}^{\circ }\bullet r_{X}^{\circ }\right) \circ
Q\right) \underbracket[0.140ex]{\overline{\varkappa }_{P\circ H,X\circ
I,Q}^{H}\zeta _{P,X,H,I\circ Q}}\left( \left( P\bullet X\right) \circ (H\bullet\left(
l_{Q}^{\circ }\right) ^{-1})\rho _{Q}^{\bullet }\right)  \\
&\overset{\mathclap{\eqref{from:varkappazeta}}}{=}\left( \left( \mu _{P}^{\circ
}\bullet r_{X}^{\circ }\right) \circ Q\right) \left( \zeta _{P,X,H,I}\circ
Q\right) \left( \left( P\bullet X\right) \circ \overline{\varkappa }_{H,I,Q}^{H}\underbracket[0.140ex]{%
\left( H\bullet \left( l_{Q}^{\circ
}\right) ^{-1}\right) \rho _{Q}^{\bullet }}\right)  \\
&\overset{\mathclap{\eqref{form:Galunit}}}{=}\left( \left( \mu _{P}^{\circ }\bullet
r_{X}^{\circ }\right) \circ Q\right) \left( \zeta _{P,X,H,I}\circ Q\right)
\left( \left( P\bullet X\right) \circ \underbracket[0.140ex]{\overline{\varkappa }%
_{H,I,Q}^{H}\varkappa _{H,I,Q}^{H}\left( \left( u_{H}^{\circ }\bullet
I\right) \Delta _{I}^{\bullet }\circ Q\right) }\left( l_{Q}^{\circ }\right)
^{-1}\right)  \\
&\overset{\mathclap{\eqref{eq:GalinvGal}}}{=}\left( \left( \mu _{P}^{\circ }\bullet
r_{X}^{\circ }\right) \circ Q\right) \left( \zeta _{P,X,H,I}\circ Q\right)
\left( \left( P\bullet X\right) \circ \left( \left( u_{H}^{\circ }\bullet
I\right) \Delta _{I}^{\bullet }\circ Q\right) \left( l_{Q}^{\circ }\right)
^{-1}\right)  \\
&=\left( \left( \mu _{P}^{\circ }\left( P\circ u_{H}^{\circ }\right)
\bullet r_{X}^{\circ }\right) \circ Q\right) \left( \zeta _{P,X,I,I}\circ
Q\right) \left( \left( P\bullet X\right) \circ \Delta _{I}^{\bullet }\circ
Q\right) \left( \left( P\bullet X\right) \circ \left( l_{Q}^{\circ }\right)
^{-1}\right)  \\
&=\underbracket[0.140ex]{\left( \left( r_{P}^{\circ }\bullet r_{X}^{\circ }\right) \circ Q\right)
\left( \zeta _{P,X,I,I}\circ Q\right) \left( \left( P\bullet X\right) \circ
\Delta _{I}^{\bullet }\circ Q\right)} \left( \left( P\bullet X\right) \circ
\left( l_{Q}^{\circ }\right) ^{-1}\right)  \\
&\overset{\mathclap{\eqref{eq:unit1}}}{=}\left( r_{P\bullet X}^{\circ }\circ Q\right) \left( \left( P\bullet
X\right) \circ \left( l_{Q}^{\circ }\right) ^{-1}\right) =\mathrm{Id}.
\end{align*}

2) If $\varkappa _{P,X,Q}^{H}\overline{\varkappa }_{P,X,Q}^{H}=\mathrm{Id}$
for every $(Q,\rho^{\bullet}_{Q}),X$ then in particular $\varkappa _{P,J,H}^{H}\overline{\varkappa }%
_{P,J,H}^{H}=\mathrm{Id.}$ The latter implies \eqref{eq:GalGalinv}. If we
assume \eqref{eq:GalGalinv}, we get $\varkappa _{P,X,Q}^{H}\overline{\varkappa
}_{P,X,Q}^{H}=\mathrm{Id}$ by the following computation:
\begin{align*}
\varkappa _{P,X,Q}^{H}&\overline{\varkappa }_{P,X,Q}^{H} =\left( \mu
_{P}^{\circ }\bullet \left( X\circ Q\right) \right) \zeta _{P,X,H,Q}\left(
\left( P\bullet X\right) \circ \rho _{Q}^{\bullet }\right) \overline{%
\varkappa }_{P,X,Q}^{H} \\
&=\left( \mu _{P}^{\circ }\bullet \left( X\circ Q\right) \right) \zeta
_{P,X,H,Q}\overline{\varkappa }_{P,X,H\bullet Q}^{H}\left( P\bullet \left(
X\circ \rho _{Q}^{\bullet }\right) \right)  \\
&=\left( \mu _{P}^{\circ }\bullet \left( X\circ Q\right) \right) \zeta
_{P,X,H,Q}\overline{\varkappa }_{P,X,H\bullet Q}^{H}\left( P\bullet \left(
l_{X}^{\bullet }\circ \left( H\bullet Q\right) \right) \right) \left(
P\bullet \left( \left( l_{X}^{\bullet }\right) ^{-1}\circ \rho _{Q}^{\bullet
}\right) \right)  \\
&=\left( \mu _{P}^{\circ }\bullet \left( X\circ Q\right) \right) \zeta
_{P,X,H,Q}\left( \left( P\bullet l_{X}^{\bullet }\right) \circ \left(
H\bullet Q\right) \right) \overline{\varkappa }_{P,J\bullet X,H\bullet
Q}^{H}\left( P\bullet \left( \left( l_{X}^{\bullet }\right) ^{-1}\circ \rho
_{Q}^{\bullet }\right) \right)  \\
&=\left( \mu _{P}^{\circ }\bullet \left( X\circ Q\right) \right) \zeta
_{P,X,H,Q}\left( \left( r_{P}^{\bullet }\bullet X\right) \circ \left(
H\bullet Q\right) \right) \overline{\varkappa }_{P,J\bullet X,H\bullet
Q}^{H}\left( P\bullet \left( \left( l_{X}^{\bullet }\right) ^{-1}\circ \rho
_{Q}^{\bullet }\right) \right)  \\
&=\left( \mu _{P}^{\circ }\left( r_{P}^{\bullet }\circ H\right) \bullet
\left( X\circ Q\right) \right) \underbracket[0.140ex]{\zeta _{P\bullet
J,X,H,Q}\overline{\varkappa }_{P,J\bullet X,H\bullet Q}^{H}}\left( P\bullet
\left( \left( l_{X}^{\bullet }\right) ^{-1}\circ \rho _{Q}^{\bullet }\right)
\right)  \\
&\overset{\mathclap{\eqref{from:varkappazeta}}}{=}\left( \underbracket[0.140ex]{\mu
_{P}^{\circ }\left( r_{P}^{\bullet }\circ H\right) }\bullet \left( X\circ
Q\right) \right) \left( {\overline{\varkappa
}_{P,J,H}^{H}}\bullet \left( X\circ Q\right) \right) \left( P\bullet \zeta
_{J,X,H,Q}\right) \left( P\bullet \left( \left( l_{X}^{\bullet }\right)
^{-1}\circ \rho _{Q}^{\bullet }\right) \right)  \\
&\overset{\mathclap{\eqref{form:Galcounit}}}{=}\left( r_{P}^{\bullet }\underbracket[0.140ex]{\left(
P\bullet m_{J}^{\circ }\left( J\circ \varepsilon _{H}^{\bullet }\right)
\right) \varkappa _{P,J,H}^{H}\overline{\varkappa }_{P,J,H}^{H}}\bullet
\left( X\circ Q\right) \right) \left( P\bullet \zeta _{J,X,H,Q}\right)
\left( P\bullet \left( \left( l_{X}^{\bullet }\right) ^{-1}\circ \rho
_{Q}^{\bullet }\right) \right)  \\
&\overset{\mathclap{\eqref{eq:GalGalinv}}}{=}\left( r_{P}^{\bullet }\left( P\bullet
m_{J}^{\circ }\left( J\circ \varepsilon _{H}^{\bullet }\right) \right)
\bullet \left( X\circ Q\right) \right) \left( P\bullet \zeta
_{J,X,H,Q}\right) \left( P\bullet \left( \left( l_{X}^{\bullet }\right)
^{-1}\circ \rho _{Q}^{\bullet }\right) \right)  \\
&=\left( \left( r_{P}^{\bullet }\left( P\bullet m_{J}^{\circ }\right)
\right) \bullet \left( X\circ Q\right) \right) \left( P\bullet \zeta
_{J,X,J,Q}\right) \left( P\bullet \left( \left( J\bullet X\right) \circ
\left( \varepsilon _{H}^{\bullet }\bullet Q\right) \right) \right) \left(
P\bullet \left( \left( l_{X}^{\bullet }\right) ^{-1}\circ \rho _{Q}^{\bullet
}\right) \right)  \\
&=\left( r_{P}^{\bullet }\bullet \left( X\circ Q\right) \right) %
\underbracket[0.140ex]{\left( P\bullet m_{J}^{\circ }\bullet \left( X\circ
Q\right) \right) \left( P\bullet \zeta _{J,X,J,Q}\right) \left( P\bullet
\left( \left( l_{X}^{\bullet }\right) ^{-1}\circ \left( l_{Q}^{\bullet
}\right) ^{-1}\right) \right) } \\
&\overset{\mathclap{\eqref{eq:unit2}}}=\left( r_{P}^{\bullet }\bullet \left( X\circ Q\right) \right) \left(
P\bullet \left( l_{X\circ Q}^{\bullet }\right) ^{-1}\right) =\mathrm{Id}.\qedhere
\end{align*}
\end{proof}

A symmetric argument yields the following result, providing preliminary conditions for the invertibility of the co-Galois maps.

\begin{proposition}%[\rd{revised:2026/03/30}]
\label{prop:cogal1} Let $\left( \mathcal{C},\circ ,I,\bullet ,J\right) $ be
a duoidal category with a reversion and let $\left( H,m_{H}^{\circ
},u_{H}^{\circ },\Delta _{H}^{\bullet },\varepsilon _{H}^{\bullet }\right) $
be a bimonoid. Let $\Sigma _{H}:H\rightarrow H^{\text{-} }$ be a morphism in $\Cc$
and consider the co-Galois maps\textbf{\ }$\varsigma _{P,X,Q}^{H}$ and $%
\overline{\varsigma}_{P,X,Q}^{H}$ of \cref{def:Galmap} for any object $X$,
right $H$-comodule $\left( P,\rho _{P}^{\bullet }\right) $ and left $H$%
-module $\left( Q,\mu _{Q}^{\circ }\right) .$

1) $\overline{\varsigma }_{P,X,Q}^{H}\varsigma _{P,X,Q}^{H}=\mathrm{Id}$ if
and only if $\overline{\varsigma }_{P,I,H}^{H}\varsigma _{P,I,H}^{H}=\mathrm{%
Id}$ if and only if
\begin{equation}
\overline{\varsigma }_{P,I,H}^{H}\varsigma _{P,I,H}^{H}\left( P\circ \left(
I\bullet u_{H}^{\circ }\right) \Delta _{I}^{\bullet }\right) = P\circ \left( I\bullet u_{H}^{\circ
}\right) \Delta _{I}^{\bullet }.
\label{eq:coGalinvGal}
\end{equation}

2) $\varsigma _{P,X,Q}^{H}\overline{\varsigma }_{P,X,Q}^{H}=\mathrm{Id}$ if
and only if $\varsigma _{H,J,Q}^{H}\overline{\varsigma }_{H,J,Q}^{H}=\mathrm{%
Id}$ if and only if%
\begin{equation}
\left( m_{J}^{\circ }\left( \varepsilon _{H}^{\bullet }\circ
J\right) \bullet Q\right) \varsigma _{H,J,Q}^{H}\overline{\varsigma }%
_{H,J,Q}^{H}= m_{J}^{\circ }\left( \varepsilon
_{H}^{\bullet }\circ J\right) \bullet Q.  \label{eq:coGalGalinv}
\end{equation}
\end{proposition}

\begin{invisible}
\begin{proof}
1) Clearly $\overline{\varsigma }_{P,X,Q}^{H}\varsigma _{P,X,Q}^{H}=\mathrm{%
Id}$ for every $\left( Q,\mu _{Q}^{\circ }\right),X$ implies $\overline{\varsigma }_{P,I,H}^{H}\varsigma
_{P,I,H}^{H}=\mathrm{Id.}$ If the latter equality is true$\mathrm{,}$ it is
clear that \eqref{eq:coGalinvGal} is true. If \eqref{eq:coGalinvGal} is true, we
get $\overline{\varsigma }_{P,X,Q}^{H}\varsigma _{P,X,Q}^{H}=\mathrm{Id}$ by
the following computation:%
\begin{eqnarray*}
\overline{\varsigma }_{P,X,Q}^{H}\varsigma _{P,X,Q}^{H} &=&\overline{%
\varsigma }_{P,X,Q}^{H}\left( \left( P\circ X\right) \bullet \mu _{Q}^{\circ
}\right) \zeta _{P,H,X,Q}\left( \rho _{P}^{\bullet }\circ \left( X\bullet
Q\right) \right)  \\
&=&\left( P\circ \left( X\bullet \mu _{Q}^{\circ }\right) \right) \overline{%
\varsigma }_{P,X,H\circ Q}^{H}\zeta _{P,H,X,Q}\left( \rho _{P}^{\bullet
}\circ \left( X\bullet Q\right) \right)  \\
&=&\left( P\circ \left( l_{X}^{\circ }\bullet \mu _{Q}^{\circ }\right)
\right) \left( P\circ \left( \left( l_{X}^{\circ }\right) ^{-1}\bullet
\left( H\circ Q\right) \right) \right) \overline{\zeta }_{P,X,H\circ
Q}^{H}\zeta _{P,H,X,Q}\left( \rho _{P}^{\bullet }\circ \left( X\bullet
Q\right) \right)  \\
&=&\left( P\circ \left( l_{X}^{\circ }\bullet \mu _{Q}^{\circ }\right)
\right) \overline{\varsigma }_{P,I\circ X,H\circ Q}^{H}\left( \left( P\circ
\left( l_{X}^{\circ }\right) ^{-1}\right) \bullet \left( H\circ Q\right)
\right) \zeta _{P,H,X,Q}\left( \rho _{P}^{\bullet }\circ \left( X\bullet
Q\right) \right)  \\
&=&\left( P\circ \left( l_{X}^{\circ }\bullet \mu _{Q}^{\circ }\right)
\right) \overline{\varsigma }_{P,I\circ X,H\circ Q}^{H}\left( \left( \left(
r_{P}^{\circ }\right) ^{-1}\circ X\right) \bullet \left( H\circ Q\right)
\right) \zeta _{P,H,X,Q}\left( \rho _{P}^{\bullet }\circ \left( X\bullet
Q\right) \right)  \\
&=&\left( P\circ \left( l_{X}^{\circ }\bullet \mu _{Q}^{\circ }\right)
\right) \underbracket[0.140ex]{\overline{\varsigma }_{P,I\circ X,H\circ
Q}^{H}\zeta _{P\circ I,H,X,Q}}\left( \left( \left( r_{P}^{\circ }\right)
^{-1}\bullet H\right) \rho _{P}^{\bullet }\circ \left( X\bullet Q\right)
\right)  \\
&\overset{\eqref{from:varsigmazeta}}{=}&\left( P\circ \left( l_{X}^{\circ }\bullet
\mu _{Q}^{\circ }\right) \right) \left( P\circ \zeta _{I,H,X,Q}\right)
\left( \overline{\varsigma }_{P,I,H}^{H}\circ \left( X\bullet Q\right)
\right) \left( \underbracket[0.140ex]{\left( \left( r_{P}^{\circ }\right) ^{-1}\bullet
H\right) \rho _{P}^{\bullet }}\circ \left( X\bullet Q\right) \right)  \\
&\overset{\eqref{form:coGalunit}}{=}&\left( P\circ \left( l_{X}^{\circ
}\bullet \mu _{Q}^{\circ }\right) \right) \left( P\circ \zeta
_{I,H,X,Q}\right) \left( \underbracket[0.140ex]{\overline{\varsigma }%
_{P,I,H}^{H}\varsigma _{P,I,H}^{H}\left( P\circ \left( I\bullet u_{H}^{\circ
}\right) \Delta _{I}^{\bullet }\right)} \left( r_{P}^{\circ }\right) ^{-1}%
\circ \left( X\bullet Q\right) \right)  \\
&\overset{\eqref{eq:coGalinvGal}}{=}&\left( P\circ \left( l_{X}^{\circ
}\bullet \mu _{Q}^{\circ }\right) \right) \left( P\circ \zeta
_{I,H,X,Q}\right) \left( \left( \left( P\circ \left( I\bullet u_{H}^{\circ
}\right) \Delta _{I}^{\bullet }\right) \left( r_{P}^{\circ }\right)
^{-1}\right) \circ \left( X\bullet Q\right) \right)  \\
&=&\left( P\circ \left( l_{X}^{\circ }\bullet \mu _{Q}^{\circ }\left(
u_{H}^{\circ }\circ Q\right) \right) \right) \left( P\circ \zeta
_{I,I,X,Q}\right) \left( \left( \left( P\circ \Delta _{I}^{\bullet }\right)
\left( r_{P}^{\circ }\right) ^{-1}\right) \circ \left( X\bullet Q\right)
\right)  \\
&=&\underbracket[0.140ex]{\left( P\circ \left( l_{X}^{\circ }\bullet
l_{Q}^{\circ }\right) \right) \left( P\circ \zeta _{I,I,X,Q}\right) \left(
P\circ \Delta _{I}^{\bullet }\circ \left( X\bullet Q\right) \right) }\left(
\left( r_{P}^{\circ }\right) ^{-1}\circ \left( X\bullet Q\right) \right)  \\
&\overset{\eqref{eq:unit1}}{=}&\left( P\circ l_{X\bullet Q}^{\circ }\right) \left( \left( r_{P}^{\circ
}\right) ^{-1}\circ \left( X\bullet Q\right) \right) =\mathrm{Id}.
\end{eqnarray*}

2) Clearly $\varsigma _{P,X,Q}^{H}\overline{\varsigma }_{P,X,Q}^{H}=\mathrm{%
Id}$ implies $\varsigma _{H,J,Q}^{H}\overline{\varsigma }_{H,J,Q}^{H}=%
\mathrm{Id}$ which in turn implies \eqref{eq:coGalGalinv}. On the other
hand, if we assume the latter, we obtain
\begin{eqnarray*}
\varsigma _{P,X,Q}^{H}\overline{\varsigma }_{P,X,Q}^{H} &=&\left( \left(
P\circ X\right) \bullet \mu _{Q}^{\circ }\right) \zeta _{P,H,X,Q}\left( \rho
_{P}^{\bullet }\circ \left( X\bullet Q\right) \right) \overline{\varsigma }%
_{P,X,Q}^{H} \\
&=&\left( \left( P\circ X\right) \bullet \mu _{Q}^{\circ }\right) \zeta
_{P,H,X,Q}\overline{\varsigma }_{P\bullet H,X,Q}^{H}\left( \left( \rho
_{P}^{\bullet }\circ X\right) \bullet Q\right)  \\
&=&\left( \left( P\circ X\right) \bullet \mu _{Q}^{\circ }\right) \zeta
_{P,H,X,Q}\overline{\varsigma }_{P\bullet H,X,Q}^{H}\left( \left( \left(
P\bullet H\right) \circ r_{X}^{\bullet }\right) \bullet Q\right) \left(
\left( \rho _{P}^{\bullet }\circ \left( r_{X}^{\bullet }\right) ^{-1}\right)
\bullet Q\right)  \\
&=&\left( \left( P\circ X\right) \bullet \mu _{Q}^{\circ }\right) \zeta
_{P,H,X,Q}\left( \left( P\bullet H\right) \circ \left( r_{X}^{\bullet
}\bullet Q\right) \right) \overline{\varsigma }_{P\bullet H,X\bullet
J,Q}^{H}\left( \left( \rho _{P}^{\bullet }\circ \left( r_{X}^{\bullet
}\right) ^{-1}\right) \bullet Q\right)  \\
&=&\left( \left( P\circ X\right) \bullet \mu _{Q}^{\circ }\right) \zeta
_{P,H,X,Q}\left( \left( P\bullet H\right) \circ \left( X\bullet
l_{Q}^{\bullet }\right) \right) \overline{\varsigma }_{P\bullet H,X\bullet
J,Q}^{H}\left( \left( \rho _{P}^{\bullet }\circ \left( r_{X}^{\bullet
}\right) ^{-1}\right) \bullet Q\right)  \\
&=&\left( \left( P\circ X\right) \bullet \mu _{Q}^{\circ }\right) \left(
\left( P\circ X\right) \bullet \left( H\circ l_{Q}^{\bullet }\right) \right) %
\underbracket[0.140ex]{\zeta _{P,H,X,J\bullet Q}\overline{\varsigma
}_{P\bullet H,X\bullet J,Q}^{H}}\left( \left( \rho _{P}^{\bullet }\circ
\left( r_{X}^{\bullet }\right) ^{-1}\right) \bullet Q\right)  \\
&\overset{\eqref{from:varsigmazeta}}{=}&\left( \left( P\circ X\right) \bullet
\underbracket[0.140ex]{\mu _{Q}^{\circ }\left( H\circ l_{Q}^{\bullet }\right) }\right)
\left( \left( P\circ X\right) \bullet \overline{\varsigma }%
_{H,J,Q}^{H}\right) \left( \zeta _{P,H,X,J}\bullet Q\right) \left( \left(
\rho _{P}^{\bullet }\circ \left( r_{X}^{\bullet }\right) ^{-1}\right)
\bullet Q\right)  \\
&\overset{\eqref{form:coGalcounit}}{=}&\left( \left( P\circ X\right) \bullet
l_{Q}^{\bullet }\underbracket[0.140ex]{\left( m_{J}^{\circ }\left( \varepsilon
_{H}^{\bullet }\circ J\right) \bullet Q\right) \varsigma _{H,J,Q}^{H}%
\overline{\varsigma }_{H,J,Q}^{H}}\right) \left( \zeta _{P,H,X,J}\bullet
Q\right) \left( \left( \rho _{P}^{\bullet }\circ \left( r_{X}^{\bullet
}\right) ^{-1}\right) \bullet Q\right)  \\
&\overset{\eqref{eq:coGalGalinv}}{=}&\left( \left( P\circ X\right) \bullet
l_{Q}^{\bullet }\left( m_{J}^{\circ }\left( \varepsilon _{H}^{\bullet }\circ
J\right) \bullet Q\right) \right) \left( \zeta _{P,H,X,J}\bullet Q\right)
\left( \left( \rho _{P}^{\bullet }\circ \left( r_{X}^{\bullet }\right)
^{-1}\right) \bullet Q\right)  \\
&=&\left( \left( P\circ X\right) \bullet l_{Q}^{\bullet }\right) \left(
\left( P\circ X\right) \bullet m_{J}^{\circ }\bullet Q\right) \left( \left(
P\circ X\right) \bullet \left( \varepsilon _{H}^{\bullet }\circ J\right)
\bullet Q\right) \left( \zeta _{P,H,X,J}\bullet Q\right) \left( \left( \rho
_{P}^{\bullet }\circ \left( r_{X}^{\bullet }\right) ^{-1}\right) \bullet
Q\right)  \\
&=&\left( \left( P\circ X\right) \bullet l_{Q}^{\bullet }\right) \left(
\left( P\circ X\right) \bullet m_{J}^{\circ }\bullet Q\right) \left( \zeta
_{P,J,X,J}\bullet Q\right) \left( \left( \left( P\bullet \varepsilon
_{H}^{\bullet }\right) \circ \left( X\bullet J\right) \right) \bullet
Q\right) \left( \left( \rho _{P}^{\bullet }\circ \left( r_{X}^{\bullet
}\right) ^{-1}\right) \bullet Q\right)  \\
&=&\left( \left( P\circ X\right) \bullet l_{Q}^{\bullet }\right) \left(
\left( P\circ X\right) \bullet m_{J}^{\circ }\bullet Q\right) \left( \zeta
_{P,J,X,J}\bullet Q\right) \left( \left( \left( r_{P}^{\bullet }\right)
^{-1}\circ \left( r_{X}^{\bullet }\right) ^{-1}\right) \bullet Q\right)  \\
&\overset{\eqref{eq:unit2}}{=}&\left( \left( P\circ X\right) \bullet l_{Q}^{\bullet }\right) \left(
r_{P\circ Q}^{\bullet }\right) ^{-1}=\mathrm{Id}.
\end{eqnarray*}
\end{proof}
\end{invisible}

A first easy consequence of the previous two propositions is the following corollary.

\begin{corollary}%[\rd{extrapolated:2026/04/27}]
\label{coro:cogal1}
In a duoidal category $\left( \mathcal{C},\circ
,I,\bullet ,J\right) $  with a reversion, we have $\overline{\varkappa }%
_{X,Y,I}^{I}=\left( \varkappa _{X,Y,I}^{I}\right) ^{-1}$ and $\overline{%
\varsigma }_{I,X,Y}^{I}=\left( \varsigma _{I,X,Y}^{I}\right) ^{-1}$ for any $%
X,Y$ in $\Cc$.
\end{corollary}

\begin{proof}
By using \cref{lem:inv} in the last
equality, we get
\[\overline{\varkappa }_{I,I,I}^{I}=\left( \left( r_{I}^{\circ }\bullet
I\right) \circ I\right) \underbracket[0.140ex]{\delta _{I,I,I,I}\left( I\bullet \left( I\circ
\left( \phi^\circ _{0}\bullet I\right) \Delta _{I}^{\bullet }\right) \right)}
\overset{\eqref{form:phi0delta}}=\left( r_{I\bullet I}^{\circ }\right) ^{-1}\left( I\bullet r_{I}^{\circ
}\right) =\left( \varkappa _{I,I,I}^{I}\right) ^{-1}.\]
In particular $\overline{\varkappa }%
_{I,I,I}^{I}\varkappa _{I,I,I}^{I}=\mathrm{Id.}$ By \cref{prop:cogal2} this
implies $\overline{\varkappa }_{X,Y,I}^{I}\varkappa _{X,Y,I}^{I}=\mathrm{Id}$
for any right $I$-module $X$ (i.e.\ just an object) and any object $Y.$ 
Since, again by \cref{lem:inv}, we know that $\varkappa _{X,Y,I}^{I}$ is always
invertible, we get that $\overline{\varkappa }_{X,Y,I}^{I}=\left( \varkappa
_{X,Y,I}^{I}\right) ^{-1}$. Similarly, by means of \eqref{form:phi0gamma}, \cref{lem:inv} and \cref{prop:cogal1}, one gets the other equality.
\begin{invisible}
   We have
\[\overline{\varsigma }_{I,I,I}^{I}=\left( I\circ \left( I\bullet
l_{I}^{\circ }\right) \right) \gamma _{I,I,I,I}\left( \left( \left( I\bullet
\phi^\circ _{0}\right) \Delta _{I}^{\bullet }\circ I\right) \bullet I\right)
\overset{\eqref{form:phi0gamma}}=\left( l_{I\bullet I}^{\circ }\right) ^{-1}\left( l_{I}^{\circ }\bullet
Y\right) =\left( \varsigma _{I,I,I}^{I}\right) ^{-1}\] 
where the last
equality follows by \cref{lem:inv}. In particular $\overline{\varsigma }%
_{I,I,I}^{I}\varsigma _{I,I,I}^{I}=\mathrm{Id.}$ Then, by \cref{prop:cogal1}%
, we have that $\overline{\varsigma }_{I,X,Y}^{I}\varsigma _{I,X,Y}^{I}=%
\mathrm{Id}$ for every object $X$ and left $I$-module (whence just an
object) $Y.$ Since, again by \cref{lem:inv}, we know that $\varsigma _{I,X,Y}^{I}$
is always invertible, we get that $\overline{\varsigma }_{I,X,Y}^{I}=\left(
\varsigma _{I,X,Y}^{I}\right) ^{-1}$ for any $X,Y.$ 
We already noted that
$\varsigma _{I,X,Y}^{I}=\left( \left( l_{X}^{\circ }\right) ^{-1}\bullet
Y\right) l_{X\bullet Y}^{\circ }$
is invertible and hence $\overline{\varsigma }_{I,X,Y}^{I}\varsigma
_{I,X,Y}^{I}=\mathrm{Id}$ implies $\overline{\varsigma }_{I,X,Y}^{I}=\left(
\varsigma _{I,X,Y}^{I}\right) ^{-1}=\left( l_{X\bullet Y}^{\circ }\right)
^{-1}\left( l_{X}^{\circ }\bullet Y\right) .$ 
\end{invisible}
\end{proof}

\cref{lem:inv} and \cref{coro:cogal1} allow us to rewrite $\psi _{X,Y,Z}^{\prime }$ and $%
\varphi _{X,Y,Z}^{\prime }$ in a different way that will be useful. To this aim, we first define the following natural transformations, for any objects $X,Y$ in $\mathcal{C}$:
\begin{align*}
\xi _{X,Y} :& \xymatrix{
X\bullet \left( J\circ \left( I\bullet Y\right) \right)
  \ar[rr]^-{X\bullet \left( J\circ \left( \left( \phi^\circ_0\bullet I\right)
      \Delta _{I}^{\bullet }\bullet Y\right) \right)}&&
X\bullet \left( J\circ \left( I^{\revsn }\bullet I\bullet Y\right) \right)
  \ar[r]^-{\varphi _{X,I,I\bullet Y}^{\prime \prime }}&
X\circ I\circ \left( I\bullet Y\right)
  \ar[r]^-{r_{X}^{\circ }\circ \left( I\bullet Y\right)}&
X\circ \left( I\bullet Y\right)} , \\
\xi _{X,Y}^{\prime } :& \xymatrix{
\left( \left( X\bullet I\right) \circ J\right) \bullet Y
  \ar[rr]^-{\left( \left( X\bullet \left( I\bullet \phi^\circ_0\right)
      \Delta _{I}^{\bullet }\right) \circ J\right) \bullet Y}&&
\left( \left( X\bullet I\bullet I^{\revsn }\right) \circ J\right) \bullet Y
  \ar[r]^-{\psi _{X\bullet I,I,Y}^{\prime \prime }}&\left( X\bullet I\right) \circ I\circ Y
  \ar[r]^-{\left( X\bullet I\right) \circ l_{Y}^{\circ }}&
\left( X\bullet I\right) \circ Y},\\
\nu_{X,Y}:&\xymatrix{((X\bullet I)\circ J\circ J)\bullet Y\ar[rr]^-{((X\bullet \phi_0^\circ)\circ m_J^\circ)\bullet Y}&&
((X\bullet I^\revsn)\circ J)\bullet Y\ar[r]^-{\psi_{X,I,Y}''}&X\circ I\circ Y\ar[r]^-{X\circ l_Y^\circ}&X\circ Y},\\
\nu_{X,Y}':&\xymatrix{X\bullet(J\circ J\circ (I\bullet Y))\ar[rr]^-{X\bullet (m_J^\circ\circ(\phi_0^\circ\bullet Y))}&&
X\bullet(J\circ (I^\revsn\bullet Y))\ar[r]^-{\varphi_{X,I,Y}''}&X\circ I\circ Y\ar[r]^-{r_X^\circ\circ Y}&X\circ Y}.
\end{align*}

\begin{lemma}%[\rd{revised:2027/04/27}]
\label{lem:rev}
In a duoidal category $\left( \mathcal{C},\circ
,I,\bullet ,J\right) $  with a reversion, for any objects $X,Y, Z$ in $\Cc$, 
we have:%
\begin{align}
\psi _{X,Y,Z}^{\prime } &=\xymatrix{
  X \bullet Y^{\revsn} \bullet Z
    \ar[r]^-{X \bullet \psi_{Y,Z}} &
  X \bullet (J \circ (I \bullet (Y \circ Z)))
    \ar[r]^-{\xi_{X,Y\circ Z}} &
  X \circ (I \bullet (Y \circ Z))
},  \label{form:psi1-rev} \\
\varphi _{X,Y,Z}^{\prime } &= \xymatrix{
  X \bullet Y^{\revsn} \bullet Z
    \ar[r]^-{{\varphi_{X,Y} \bullet Z}} &
  (((X \circ Y) \bullet I) \circ J) \bullet Z
    \ar[r]^-{{\xi'_{X\circ Y, Z}}} &
  ((X \circ Y) \bullet I) \circ Z
},  \label{form:varphi1-rev}\\
\psi _{X,Y,Z}^{\prime \prime }&=\xymatrix{\left( \left( X\bullet Y^{\revsn }\right) \circ J\right) \bullet Z\ar[rr]^-{\left( \varphi _{X,Y}\circ J\right) \bullet Z}&&\left(
\left( \left( X\circ Y\right) \bullet I\right) \circ J\circ J\right) \bullet
Z\ar[r]^-{\nu_{X\circ Y,Z}}&X\circ Y\circ Z},  \label{form:psi2-rev}\\
\varphi _{X,Y,Z}^{\prime \prime }&=\xymatrix{ X\bullet \left( J\circ \left( Y^{\revsn }\bullet Z\right) \right)
 \ar[rr]^-{X\bullet \left( J\circ \psi _{Y,Z}\right) }&&X\bullet \left( J\circ J\circ \left( I\bullet \left( Y\circ Z\right)
\right) \right)\ar[r]^-{\nu_{X,Y\circ Z}'}& X\circ Y\circ Z }.  \label{form:varphi2-rev}
\end{align}%
\end{lemma}

\begin{proof} For the first two equalities, we need the identities in  \cref{coro:cogal1}.
We check \eqref{form:psi1-rev}:
\begin{eqnarray*}
&&\xi _{X,Y\circ Z}\left( X\bullet \psi _{Y,Z}\right)  \\
&=&\left( r_{X}^{\circ }\circ \left( I\bullet \left( Y\circ Z\right) \right)
\right) \varphi _{X,I,I\bullet \left( Y\circ Z\right) }^{\prime \prime
}\left( X\bullet \left( J\circ \left( \left( \phi^\circ_0\bullet I\right)
\Delta _{I}^{\bullet }\bullet \left( Y\circ Z\right) \right) \right) \right)
\left( X\bullet \psi _{Y,Z}\right)  \\
&=&
\left( r_{X}^{\circ }\circ \left( I\bullet \left( Y\circ Z\right) \right)
\right) \left( r_{X\circ I}^{\bullet }\circ \left( I\bullet \left( Y\circ
Z\right) \right) \right) \delta _{X,J,I,I\bullet \left( Y\circ Z\right)
}\left( X\bullet \left( J\circ \left( \left( \phi^\circ_0\bullet I\right)
\Delta _{I}^{\bullet }\bullet \left( Y\circ Z\right) \right) \right) \right)
\\
&&\left( X\bullet \gamma _{J,Y,I,Z}\right) \left( X\bullet \left( \left(
l_{Y^{\revsn }}^{\bullet }\right) ^{-1}\circ I\right) \bullet Z\right) \left(
X\bullet \left( r_{Y^{\revsn }}^{\circ }\right) ^{-1}\bullet Z\right)\\
&=&
\left( r_{X}^{\bullet }\circ \left( I\bullet \left( Y\circ Z\right) \right)
\right) \left( \left( r_{X}^{\circ }\bullet J\right) \circ \left( I\bullet
\left( Y\circ Z\right) \right) \right) \underbracket[0.140ex]{\delta _{X,J,I,I\bullet
\left( Y\circ Z\right) }\left( X\bullet \gamma _{J,Y,I^{\revsn }\bullet
I,Z}\right) } \\
&&\left( X\bullet \left( \left( l_{Y^{\revsn }}^{\bullet }\right) ^{-1}\circ
\left( \phi^\circ_0\bullet I\right) \Delta _{I}^{\bullet }\right) \bullet
Z\right) \left( X\bullet \left( r_{Y^{\revsn }}^{\circ }\right) ^{-1}\bullet
Z\right)\\
&\overset{\eqref{form:gammadelta}}{=}&
\left( r_{X}^{\bullet }\circ \left( I\bullet \left( Y\circ Z\right) \right)
\right) \left( \left( r_{X}^{\circ }\bullet J\right) \circ \left( I\bullet
\left( Y\circ Z\right) \right) \right) \gamma _{\left( X\circ I\right)
\bullet J,Y,I,Z}\left( \delta _{X,J\bullet Y^{\revsn },I,I}\bullet
Z\right)  \\
&&\left( X\bullet \left( \left( l_{Y^{\revsn }}^{\bullet }\right) ^{-1}\circ
\left( \phi^\circ_0\bullet I\right) \Delta _{I}^{\bullet }\right) \bullet
Z\right) \left( X\bullet \left( r_{Y^{\revsn }}^{\circ }\right) ^{-1}\bullet
Z\right)\\
&=&
\gamma _{X,Y,I,Z}\left( \left( \left( r_{X}^{\bullet }\bullet Y^{\revsn
}\right) \circ I\right) \bullet Z\right) \left( \left( \left( r_{X}^{\circ
}\bullet \left( l_{Y^{\revsn }}^{\bullet }\right) ^{-1 }\right) \circ
I\right) \bullet Z\right) \left( \delta _{X,Y^{\revsn },I,I}\bullet
Z\right)  \\
&&\left( X\bullet \left( Y^{\revsn }\circ \left( \phi^\circ_0\bullet I\right)
\Delta _{I}^{\bullet }\right) \bullet Z\right) \left( X\bullet \left(
r_{Y^{\revsn }}^{\circ }\right) ^{-1}\bullet Z\right)\\
&=&\gamma _{X,Y,I,Z}\left( \left( \left( r_{X}^{\circ }\bullet Y^{\revsn
}\right) \circ I\right) \bullet Z\right) \left( \delta _{X,Y^{\revsn },I,I}\bullet Z\right) \left( X\bullet \left( Y^{\revsn }\circ \left(\phi^\circ_0\bullet I\right) \Delta _{I}^{\bullet }\right) \bullet Z\right) \left(
X\bullet \left( r_{Y^{\revsn }}^{\circ }\right) ^{-1}\bullet Z\right)  \\
&=&\gamma _{X,Y,I,Z}\left( \left( \underbracket[0.140ex]{\left( \left(
r_{X}^{\circ }\bullet Y^{\revsn }\right) \circ I\right) \delta _{X,Y^{\revsn
},I,I}\left( X\bullet \left( Y^{\revsn }\circ \left( \phi^\circ_0\bullet
I\right) \Delta _{I}^{\bullet }\right) \right) }\left( X\bullet \left(
r_{Y^{\revsn }}^{\circ }\right) ^{-1}\right) \right) \bullet Z\right)  \\
&=&\gamma _{X,Y,I,Z}\left( \left( \overline{\varkappa }_{X,Y^{\revsn
},I}^{I}\left( X\bullet \left( r_{Y^{\revsn }}^{\circ }\right) ^{-1}\right)
\right) \bullet Z\right)  \\
&=&\gamma _{X,Y,I,Z}\left( \left( \left( X\bullet r_{Y^{\revsn }}^{\circ
}\right) \underbracket[0.140ex]{\varkappa _{X,Y^{\revsn },I}^{I}}\right) ^{-1}\bullet
Z\right)  
=\gamma _{X,Y,I,Z}\left( \left( r_{X\bullet Y^{\revsn }}^{\circ }\right)
^{-1}\bullet Z\right)  
=\psi _{X,Y,Z}^{\prime }.
\end{eqnarray*}

The proof of \eqref{form:varphi1-rev} follows by symmetric argument.
\begin{invisible}
Let us check that $\varphi
_{X,Y,Z}^{\prime }=\xi _{X\circ Y,Z}^{\prime }\left( \varphi _{X,Y}\bullet
Z\right) :$%
\begin{eqnarray*}
&&\xi _{X\circ Y,Z}^{\prime }\left( \varphi _{X,Y}\bullet Z\right)  \\
&=&\left( \left( \left( X\circ Y\right) \bullet I\right) \circ l_{Z}^{\circ
}\right) \psi _{\left( X\circ Y\right) \bullet I,I,Z}^{\prime \prime }\left(
\left( \left( \left( X\circ Y\right) \bullet \left( I\bullet \phi^\circ_0\right) \Delta _{I}^{\bullet }\right) \circ J\right) \bullet Z\right)
\left( \varphi _{X,Y}\bullet Z\right)  \\
&=&\left[
\begin{array}{c}
\left( \left( \left( X\circ Y\right) \bullet I\right) \circ l_{Z}^{\circ
}\right) \left( \left( \left( X\circ Y\right) \bullet I\right) \circ
l_{I\circ Z}^{\bullet }\right) \gamma _{\left( X\circ Y\right) \bullet
I,I,J,Z}\left( \left( \left( \left( X\circ Y\right) \bullet \left( I\bullet
\phi^\circ_0\right) \Delta _{I}^{\bullet }\right) \circ J\right) \bullet
Z\right)  \\
\left( \delta _{X,I,Y,J}\bullet Z\right) \left( X\bullet \left( I\circ
\left( r_{Y^{\revsn }}^{\bullet }\right) ^{-1}\right) \bullet Z\right) \left(
X\bullet \left( l_{Y^{\revsn }}^{\circ }\right) ^{-1}\bullet Z\right)
\end{array}%
\right]  \\
&=&\left[
\begin{array}{c}
\left( \left( \left( X\circ Y\right) \bullet I\right) \circ l_{Z}^{\circ
}\right) \left( \left( \left( X\circ Y\right) \bullet I\right) \circ
l_{I\circ Z}^{\bullet }\right) \underbracket[0.140ex]{\gamma _{\left( X\circ Y\right)
\bullet I,I,J,Z}\left( \delta _{X,I\bullet I^{\revsn },Y,J}\bullet Z\right) }
\\
\left( \left( X\bullet \left( \left( I\bullet \phi^\circ_0\right) \Delta
_{I}^{\bullet }\circ \left( r_{Y^{\revsn }}^{\bullet }\right) ^{-1}\right)
\right) \bullet Z\right) \left( X\bullet \left( l_{Y^{\revsn }}^{\circ
}\right) ^{-1}\bullet Z\right)
\end{array}%
\right]  \\
&\overset{\eqref{form:gammadelta}}{=}&\left[
\begin{array}{c}
\left( \left( \left( X\circ Y\right) \bullet I\right) \circ l_{Z}^{\circ
}\right) \left( \left( \left( X\circ Y\right) \bullet I\right) \circ
l_{I\circ Z}^{\bullet }\right) \delta _{X,I,Y,J\bullet \left( I\circ
Z\right) }\left( X\bullet \gamma _{I,I,Y^{\revsn }\bullet J,Z}\right)  \\
\left( \left( X\bullet \left( \left( I\bullet \phi^\circ_0\right) \Delta
_{I}^{\bullet }\circ \left( r_{Y^{\revsn }}^{\bullet }\right) ^{-1}\right)
\right) \bullet Z\right) \left( X\bullet \left( l_{Y^{\revsn }}^{\circ
}\right) ^{-1}\bullet Z\right)
\end{array}%
\right]  \\
&=&\left[
\begin{array}{c}
\delta _{X,I,Y,Z}\left( X\bullet \left( I\circ \left( Y^{\revsn }\bullet
l_{Z}^{\circ }\right) \right) \right) \left( X\bullet \left( I\circ \left(
Y^{\revsn }\bullet l_{I\circ Z}^{\bullet }\right) \right) \right) \left(
X\bullet \gamma _{I,I,Y^{\revsn }\bullet J,Z}\right)  \\
\left( \left( X\bullet \left( \left( I\bullet \phi^\circ_0\right) \Delta
_{I}^{\bullet }\circ \left( r_{Y^{\revsn }}^{\bullet }\right) ^{-1}\right)
\right) \bullet Z\right) \left( X\bullet \left( l_{Y^{\revsn }}^{\circ
}\right) ^{-1}\bullet Z\right)
\end{array}%
\right]  \\
&=&\left[
\begin{array}{c}
\delta _{X,I,Y,Z}\left( X\bullet \left( I\circ \left( Y^{\revsn }\bullet
l_{Z}^{\circ }\right) \right) \right) \left( X\bullet \left( I\circ \left(
r_{Y^{\revsn }}^{\bullet }\bullet \left( I\circ Z\right) \right) \right)
\right) \left( X\bullet \gamma _{I,I,Y^{\revsn }\bullet J,Z}\right)  \\
\left( \left( X\bullet \left( \left( I\bullet \phi^\circ_0\right) \Delta
_{I}^{\bullet }\circ \left( r_{Y^{\revsn }}^{\bullet }\right) ^{-1}\right)
\right) \bullet Z\right) \left( X\bullet \left( l_{Y^{\revsn }}^{\circ
}\right) ^{-1}\bullet Z\right)
\end{array}%
\right]  \\
&=&\delta _{X,I,Y,Z}\left( X\bullet \left( I\circ \left( Y^{\revsn }\bullet
l_{Z}^{\circ }\right) \right) \right) \left( X\bullet \gamma _{I,I,Y^{\revsn
},Z}\right) \left( \left( X\bullet \left( \left( I\bullet \phi^\circ_0\right)
\Delta _{I}^{\bullet }\circ Y^{\revsn }\right) \right) \bullet Z\right) \left(
X\bullet \left( l_{Y^{\revsn }}^{\circ }\right) ^{-1}\bullet Z\right)  \\
&=&\delta _{X,I,Y,Z}\left( X\bullet \underbracket[0.140ex]{\overline{\varsigma }%
_{I,Y^{\revsn },Z}^{I}}\right) \left( X\bullet \underbracket[0.140ex]{\left( l_{Y^{\revsn
}}^{\circ }\right) ^{-1}\bullet Z}\right)  \\
&=&\delta _{X,I,Y,Z}\left( X\bullet \left( l_{Y^{\revsn }\bullet Z}^{\circ
}\right) ^{-1}\right) =\varphi _{X,Y,Z}^{\prime }.
\end{eqnarray*}\end{invisible}

Let us prove \eqref{form:psi2-rev}: 
\begin{eqnarray*}
&&\nu _{X\circ Y,Z}\left( \left( \varphi _{X,Y}\circ J\right) \bullet
Z\right) =\left( \left( X\circ Y\right) \circ l_{Z}^{\circ }\right) \psi
_{X\circ Y,I,Z}^{\prime \prime }\left( \left( \left( \left( X\circ Y\right)
\bullet \phi _{0}^{\circ }\right) \circ m_{J}^{\circ }\right) \bullet
Z\right) \left( \left( \varphi _{X,Y}\circ J\right) \bullet Z\right)  \\
&=&
\left( \left( X\circ Y\right) \circ l_{Z}^{\circ }\right) \psi _{X\circ
Y,I,Z}^{\prime \prime }\left( \left( \left( \left( X\circ Y\right) \bullet
\phi _{0}^{\circ }\right) \circ m_{J}^{\circ }\right) \bullet Z\right)  \\
&&\left( \left( \delta _{X,I,Y,J}\circ J\right) \bullet Z\right) \left( \left(
\left( X\bullet \left( I\circ \left( r_{Y^{\revsn}}^{\bullet }\right)
^{-1}\right) \right) \circ J\right) \bullet Z\right) \left( \left( \left(
X\bullet \left( l_{Y^{\revsn}}^{\circ }\right) ^{-1}\right) \circ J\right)
\bullet Z\right)\\
&=&
\left( \left( X\circ Y\right) \circ l_{Z}^{\circ }l_{I\circ Z}^{\bullet
}\right) \gamma _{X\circ Y,I,J,Z}\left( \left( \left( \left( X\circ Y\right)
\bullet I^{\revsn}\right) \circ m_{J}^{\circ }\right) \bullet Z\right)
\left( \left( \left( \left( X\circ Y\right) \bullet \phi _{0}^{\circ
}\right) \circ J\circ J\right) \bullet Z\right)  \\
&&\left( \left( \delta _{X,I,Y,J}\circ J\right) \bullet Z\right) \left( \left(
\left( X\bullet \left( I\circ \left( r_{Y^{\revsn}}^{\bullet }\right)
^{-1}\right) \right) \circ J\right) \bullet Z\right) \left( \left( \left(
X\bullet \left( l_{Y^{\revsn}}^{\circ }\right) ^{-1}\right) \circ J\right)
\bullet Z\right)\\
&=&
\left( X\circ Y\circ l_{Z}^{\circ }l_{I\circ Z}^{\bullet }\right) \left(
X\circ Y\circ \left( m_{J}^{\circ }\bullet \left( I\circ Z\right) \right)
\right) \gamma _{X\circ Y,I,J\circ J,Z}\left( \left( \delta _{X,I^{\revsn%
},Y,J}\circ J\right) \bullet Z\right)  \\
&&\left( \left( \left( X\bullet \left( \phi _{0}^{\circ }\circ \left( Y^{\revsn%
}\bullet J\right) \right) \right) \circ J\right) \bullet Z\right) \left(
\left( \left( X\bullet \left( I\circ \left( r_{Y^{\revsn}}^{\bullet }\right)
^{-1}\right) \right) \circ J\right) \bullet Z\right) \left( \left( \left(
X\bullet \left( l_{Y^{\revsn}}^{\circ }\right) ^{-1}\right) \circ J\right)
\bullet Z\right)\\
&=&
\left( X\circ Y\circ l_{Z}^{\circ }l_{I\circ Z}^{\bullet }\right) \left(
X\circ Y\circ \left( m_{J}^{\circ }\bullet \left( I\circ Z\right) \right)
\right) \gamma _{X\circ Y,I,J\circ J,Z}\left( \left( \delta _{X,I^{\revsn%
},Y,J}\circ J\right) \bullet Z\right)  \\
&&\left( \left( \left( X\bullet \left( I^{\revsn}\circ \left( r_{Y^{\revsn%
}}^{\bullet }\right) ^{-1}\right) \right) \circ J\right) \bullet Z\right)
\left( \left( \left( X\bullet \left( \phi _{0}^{\circ }\circ Y^{\revsn%
}\right) \right) \circ J\right) \bullet Z\right) \left( \left( \left(
X\bullet \left( l_{Y^{\revsn}}^{\circ }\right) ^{-1}\right) \circ J\right)
\bullet Z\right)\\
&\overset{\eqref{form:phipsi2}}{=}&
\left( X\circ Y\circ l_{Z}^{\circ }\underbracket[0.140ex]{l_{I\circ
Z}^{\bullet }}\right) \left( X\circ Y\circ \left( %
\underbracket[0.140ex]{m_{J}^{\circ }\bullet \left( I\circ Z\right) }\right)
\right) \left( X\circ Y\circ \underbracket[0.140ex]{\zeta _{J,I,J,Z}\left(
\left( l_{I}^{\bullet }\right) ^{-1}\circ \left( l_{Z}^{\bullet }\right)
^{-1}\right) }\right)  \\
&&\left( X\circ l_{Y\circ I\circ Z}^{\bullet }\right) \gamma _{X,Y\circ
I,J, Z}\left( \left( \left( X\bullet \phi _{I,Y}^{\circ }\right) \circ J\right)
\bullet Z\right) \left( \left( \left( X\bullet \left( \phi _{0}^{\circ
}\circ Y^{\revsn}\right) \left( l_{Y^{\revsn}}^{\circ }\right) ^{-1}\right)
\circ J\right) \bullet Z\right)\\
&\overset{\eqref{eq:unit2}}{=}&\left( X\circ l_{Y\circ Z}^{\bullet }\right) \left( X\circ \left(
J\bullet \left( I\circ l_{Z}^{\circ }\right) \right) \right) \gamma
_{X,Y\circ I,J,Z}\left( \left( \left( X\bullet \underbracket[0.140ex]{\phi
_{I,Y}^{\circ }\left(\phi _{0}^{\circ }\circ Y^{\revsn }\right)\left( l_{Y^{%
\revsn}}^{\circ }\right) ^{-1}}\right) \circ J\right) \bullet Z\right)  \\
&\overset{\mathtt{lax}}{=}&\left( X\circ l_{Y\circ Z}^{\bullet }\right)
\left( X\circ \left( J\bullet \left( r_{Y}^{\circ }\circ Z\right) \right)
\right) \gamma _{X,Y\circ I,J,Z}\left( \left( \left( X\bullet \left( \left(
r_{Y}^{\circ }\right) ^{\revsn}\right) ^{-1}\right) \circ J\right) \bullet Z\right)
\\
&=&\left( X\circ l_{Y\circ Z}^{\bullet }\right) \gamma _{X,Y,J,Z}\left(
\left( \left( X\bullet \left( r_{Y}^{\circ }\right) ^{\revsn}\right) \circ
J\right) \bullet Z\right) \left( \left( \left( X\bullet \left( \left(
r^{\circ }_{Y}\right) ^{\revsn}\right) ^{-1}\right) \circ J\right) \bullet Z\right)
=\psi _{X,Y,Z}^{\prime \prime }.
\end{eqnarray*}
The proof of \eqref{form:varphi2-rev} follows by symmetric argument.
\begin{invisible}
Let us prove \eqref{form:varphi2-rev}:%
\begin{eqnarray*}
&&\nu_{X,Y\circ Z}'\left( J\circ \psi _{Y,Z}\right)=\\&=&\left( r_{X}^{\circ }r_{X\circ I}^{\bullet }\circ Y\circ Z\right) \delta
_{X,J,I,Y\circ Z}\left( X\bullet \left( J\circ \left( \phi^\circ _{0}\bullet
\left( Y\circ Z\right) \right) \right) \right) \left( X\bullet \left(
m_{J}^{\circ }\circ \left( I\bullet \left( Y\circ Z\right) \right) \right)
\right) \left( X\bullet \left( J\circ \psi _{Y,Z}\right) \right)  \\
&=&\left(
\begin{array}{c}
\left( r_{X}^{\circ }r_{X\circ I}^{\bullet }\circ Y\circ Z\right) \delta
_{X,J,I,Y\circ Z}\left( X\bullet \left( J\circ \left( \phi^\circ _{0}\bullet
\left( Y\circ Z\right) \right) \right) \right) \left( X\bullet \left(
m_{J}^{\circ }\circ \left( I\bullet \left( Y\circ Z\right) \right) \right)
\right)  \\
\left( X\bullet \left( J\circ \gamma _{J,Y,I,X}\right) \right) \left(
X\bullet \left( J\circ \left( \left( \left( l_{Y^{\revsn }}^{\bullet }\right)
^{-1}\circ I\right) \bullet Z\right) \right) \right) \left( X\bullet \left(
J\circ \left( \left( r_{Y^{\revsn }}^{\circ }\right) ^{-1}\bullet Z\right)
\right) \right)
\end{array}%
\right)  \\
&=&\left(
\begin{array}{c}
\left( r_{X}^{\circ }r_{X\circ I}^{\bullet }\circ Y\circ Z\right) \delta
_{X,J,I,Y\circ Z}\left( X\bullet \left( m_{J}^{\circ }\circ \left( I^{\revsn
}\bullet \left( Y\circ Z\right) \right) \right) \right) \left( X\bullet
\left( J\circ J\circ \left( \phi^\circ _{0}\bullet \left( Y\circ Z\right) \right)
\right) \right)  \\
\left( X\bullet \left( J\circ \gamma _{J,Y,I,X}\right) \right) \left(
X\bullet \left( J\circ \left( \left( \left( l_{Y^{\revsn }}^{\bullet }\right)
^{-1}\circ I\right) \bullet Z\right) \right) \right) \left( X\bullet \left(
J\circ \left( \left( r_{Y^{\revsn }}^{\circ }\right) ^{-1}\bullet Z\right)
\right) \right)
\end{array}%
\right)  \\
&=&\left(
\begin{array}{c}
\left( r_{X}^{\circ }r_{X\circ I}^{\bullet }\circ Y\circ Z\right) \left(
\left( \left( X\circ I\right) \bullet m_{J}^{\circ }\right) \circ Y\circ
Z\right) \delta _{X,J\circ J,I,Y\circ Z}\left( X\bullet \left( J\circ \gamma
_{J,Y,I^{\revsn },X}\right) \right)  \\
\left( X\bullet \left( J\circ \left( \left( \left( J\bullet Y^{\revsn }\right)
\circ \phi^\circ _{0}\right) \bullet Z\right) \right) \right) \left( X\bullet
\left( J\circ \left( \left( \left( l_{Y^{\revsn }}^{\bullet }\right)
^{-1}\circ I\right) \bullet Z\right) \right) \right) \left( X\bullet \left(
J\circ \left( \left( r_{Y^{\revsn }}^{\circ }\right) ^{-1}\bullet Z\right)
\right) \right)
\end{array}%
\right)  \\
&=&\left(
\begin{array}{c}
\left( r_{X}^{\circ }r_{X\circ I}^{\bullet }\circ Y\circ Z\right) \left(
\left( \left( X\circ I\right) \bullet m_{J}^{\circ }\right) \circ Y\circ
Z\right) \delta _{X,J\circ J,I,Y\circ Z}\left( X\bullet \left( J\circ \gamma
_{J,Y,I^{\revsn },X}\right) \right)  \\
\left( X\bullet \left( J\circ \left( \left( \left( l_{Y^{\revsn }}^{\bullet
}\right) ^{-1}\circ I^{\revsn }\right) \bullet Z\right) \right) \right) \left(
X\bullet \left( J\circ \left( \left( Y^{\revsn }\circ \phi^\circ _{0}\right) \bullet
Z\right) \right) \right) \left( X\bullet \left( J\circ \left( \left(
r_{Y^{\revsn }}^{\circ }\right) ^{-1}\bullet Z\right) \right) \right)
\end{array}%
\right)  \\
&\overset{\eqref{form:phivarphi2}}{=}&\left(
\begin{array}{c}
\left( r_{X}^{\circ }\underbracket[0.140ex]{r_{X\circ I}^{\bullet }}\circ Y\circ
Z\right) \left( \left( \underbracket[0.140ex]{\left( X\circ I\right) \bullet
m_{J}^{\circ }}\right) \circ Y\circ Z\right) \left( \underbracket[0.140ex]{\zeta
_{X,J,I,J}}\circ Y\circ Z\right) \left( \underbracket[0.140ex]{\left( r_{X}^{\bullet
}\right) ^{-1}\circ \left( r_{I}^{\bullet }\right) ^{-1}}\circ Y\circ
Z\right)  \\
\left( r_{X\circ I\circ Y}^{\bullet }\circ Z\right) \delta _{X,J,I\circ
Y,Z}\left( X\bullet \left( J\circ \left( \phi _{Y,I}^{\circ}\bullet Z\right)
\right) \right) \left( X\bullet \left( J\circ \left( \left( Y^{\revsn }\circ
\phi^\circ _{0}\right) \bullet Z\right) \right) \right) \left( X\bullet \left(
J\circ \left( \left( r_{Y^{\revsn }}^{\circ }\right) ^{-1}\bullet Z\right)
\right) \right)
\end{array}%
\right)  \\
&\overset{\eqref{eq:unit2}}{=}&\left( r_{X}^{\circ }\circ Y\circ Z\right) \left( r_{X\circ I\circ
Y}^{\bullet }\circ Z\right) \delta _{X,J,I\circ Y,Z}\left( X\bullet \left(
J\circ \left( \phi _{Y,I}^{\circ}\bullet Z\right) \right) \right) \left(
X\bullet \left( J\circ \left( \left( Y^{\revsn }\circ \phi^\circ _{0}\right) \bullet
Z\right) \right) \right) \left( X\bullet \left( J\circ \left( \left(
r_{Y^{\revsn }}^{\circ }\right) ^{-1}\bullet Z\right) \right) \right)  \\
&=&\left( r_{X\circ Y}^{\bullet }\circ Z\right) \left( \left( \left(
r_{X}^{\circ }\circ Y\right) \bullet J\right) \circ Z\right) \delta
_{X,J,I\circ Y,Z}\left( X\bullet \left( J\circ \left( \phi
_{Y,I}^{\circ}\bullet Z\right) \right) \right) \left( X\bullet \left( J\circ
\left( \left( Y^{\revsn }\circ \phi^\circ _{0}\right) \bullet Z\right) \right)
\right)\\&& \left( X\bullet \left( J\circ \left( \left( r_{Y^{\revsn }}^{\circ
}\right) ^{-1}\bullet Z\right) \right) \right)  \\
&=&\left( r_{X\circ Y}^{\bullet }\circ Z\right) \left( \left( \left( X\circ
l_{Y}^{\circ }\right) \bullet J\right) \circ Z\right) \delta _{X,J,I\circ
Y,Z}\left( X\bullet \left( J\circ \left( \phi _{Y,I}^{\circ}\bullet Z\right)
\right) \right) \left( X\bullet \left( J\circ \left( \left( Y^{\revsn }\circ
\phi^\circ _{0}\right) \bullet Z\right) \right) \right) \\&&\left( X\bullet \left(
J\circ \left( \left( r_{Y^{\revsn }}^{\circ }\right) ^{-1}\bullet Z\right)
\right) \right)  \\
&=&\left( r_{X\circ Y}^{\bullet }\circ Z\right) \delta _{X,J,Y,Z}\left(
X\bullet \left( J\circ \left( \underbracket[0.140ex]{\left( l_{Y}^{\circ }\right)
^{\revsn }}\bullet Z\right) \right) \right) \left( X\bullet \left( J\circ
\left( \underbracket[0.140ex]{\phi _{Y,I}^{\circ}}\bullet Z\right) \right) \right)
\left( X\bullet \left( J\circ \left( \left( \underbracket[0.140ex]{Y^{\revsn }\circ \phi^\circ
_{0}}\right) \bullet Z\right) \right) \right)\\&& \left( X\bullet \left( J\circ
\left( \left( r_{Y^{\revsn }}^{\circ }\right) ^{-1}\bullet Z\right) \right)
\right)  \\
&=&\varphi _{X,Y,Z}^{\prime \prime }\left( X\bullet \left( J\circ \left(
\underbracket[0.140ex]{\left( l_{Y}^{\circ }\right) ^{\revsn }}\bullet Z\right) \right)
\right) \left( X\bullet \left( J\circ \left( \underbracket[0.140ex]{\phi _{Y,I}^{\circ}}%
\bullet Z\right) \right) \right) \left( X\bullet \left( J\circ \left( \left(
\underbracket[0.140ex]{Y^{\revsn }\circ \phi^\circ _{0}}\right) \bullet Z\right) \right)
\right) \left( X\bullet \left( J\circ \left( \left( r_{Y^{\revsn }}^{\circ
}\right) ^{-1}\bullet Z\right) \right) \right)  \\
&\overset{\text{lax}}{=}&\varphi _{X,Y,Z}^{\prime \prime }\left( X\bullet
\left( J\circ \left( r_{Y^{\revsn }}^{\circ }\bullet Z\right) \right) \right)
\left( X\bullet \left( J\circ \left( \left( r_{Y^{\revsn }}^{\circ }\right)
^{-1}\bullet Z\right) \right) \right)  \\
&=&\varphi _{X,Y,Z}^{\prime \prime }.
\end{eqnarray*}
\end{invisible}
\end{proof}

We also need the following identities.

\begin{lemma}%[\rd{added:2027/04/28}]
\label{lem:rev2}
In a duoidal category $\left( \mathcal{C},\circ
,I,\bullet ,J\right) $  with a reversion, for any object $X$  in $\Cc$, the following identities hold:%
\begin{align}
\label{form:xidelta}
 \xi _{X,I}\left( X\bullet \left( J\circ \Delta _{I}^{\bullet }\right)
\right)&=\left( X\circ \Delta
_{I}^{\bullet }\right) \left( r_{X}^{\bullet }\circ I\right)
\overline{\varkappa} _{X,J,I}^{I},\\
\label{form:xi1delta}
 \xi _{I,X}^{\prime }\left( \left( \Delta _{I}^{\bullet }\circ J\right)
\bullet X\right)  &= \left( \Delta _{I}^{\bullet }\circ X\right) \left( I\circ l_{X}^{\bullet
}\right) \overline{\varsigma }_{I,J,X}^{I},\\
\nu _{I,X}\left( \left( \Delta _{I}^{\bullet }\circ J\circ J\right)
\bullet X\right)&=\left( I\circ l_{X}^{\bullet }\right) \overline{\varsigma }%
_{I,J,X}^{I}\left( \left( I\circ m_{J}^{\circ }\right) \bullet X\right),\label{form:nudelta}\\
\nu _{X,I}^{\prime }\left(X\bullet \left( J\circ J\circ \Delta
_{I}^{\bullet }\right) \right)&=\left( r_{X}^{\bullet }\circ I\right) \overline{\varkappa }%
_{X,J,I}^{I}\left( X\bullet \left( m_{J}^{\circ }\circ I\right) \right).\label{form:nu1delta}
\end{align}
\end{lemma}

\begin{proof}
By definition of $\xi _{X,I}$, definition and naturality of $\varphi^{\prime \prime }$, and by definition of $\overline{\varkappa }$, one gets \eqref{form:xidelta}.
\begin{invisible}
We now check \eqref{form:xidelta}:
\begin{eqnarray*}
\xi _{X,I}\left( X\bullet \left( J\circ \Delta _{I}^{\bullet }\right)
\right)  &=&\left( r_{X}^{\circ }\circ \left( I\bullet I\right) \right) \varphi
_{X,I,I\bullet I}^{\prime \prime }\left( X\bullet \left( J\circ \left(
\left( \Sigma _{I}\bullet I\right) \Delta _{I}^{\bullet }\bullet I\right)
\right) \right) \left( X\bullet \left( J\circ \Delta _{I}^{\bullet }\right)
\right)  \\
&=&\left( r_{X}^{\circ }\circ \left( I\bullet I\right) \right) \varphi
_{X,I,I\bullet I}^{\prime \prime }\left( X\bullet \left( J\circ \left(
I^\revsn\bullet \Delta _{I}^{\bullet }\right) \right) \right) \left( X\bullet
\left( J\circ \left( \Sigma _{I}\bullet I\right) \Delta _{I}^{\bullet
}\right) \right)  \\
&=&\left( X\circ \Delta _{I}^{\bullet }\right) \left( r_{X}^{\circ }\circ
I\right) \varphi _{X,I,I}^{\prime \prime }\left( X\bullet \left( J\circ
\left( \Sigma _{I}\bullet I\right) \Delta _{I}^{\bullet }\right) \right)  \\
&=&\left( X\circ \Delta _{I}^{\bullet }\right) \left( r_{X}^{\bullet }\circ
I\right) \left( \left( r_{X}^{\circ }\bullet J\right) \circ I\right) \delta
_{X,J,I,I}\left( X\bullet \left( J\circ \left( \Sigma _{I}\bullet I\right)
\Delta _{I}^{\bullet }\right) \right)  \\
&=&\left( X\circ \Delta _{I}^{\bullet }\right) \left( r_{X}^{\bullet }\circ
I\right) \overline{\varkappa }_{X,J,I}^{I}.
\end{eqnarray*}%
\end{invisible}
The proof of \eqref{form:xi1delta}, \eqref{form:nudelta} and \eqref{form:nu1delta} follow by similar arguments.
\begin{invisible}
We now check \eqref{form:xi1delta}:
\begin{eqnarray*}
\xi _{I,X}^{\prime }\left( \left( \Delta _{I}^{\bullet }\circ J\right)
\bullet X\right)  &=&\left( \left( I\bullet I\right) \circ l_{X}^{\circ
}\right) \psi _{I\bullet I,I,X}^{\prime \prime }\left( \left( \left(
I\bullet \left( I\bullet \Sigma _{I}\right) \Delta _{I}^{\bullet }\right)
\circ J\right) \bullet X\right) \left( \left( \Delta _{I}^{\bullet }\circ
J\right) \bullet X\right)  \\
&=&\left( \left( I\bullet I\right) \circ l_{X}^{\circ }\right) \psi
_{I\bullet I,I,X}^{\prime \prime }\left( \left( \left( \Delta _{I}^{\bullet
}\bullet I^\revsn\right) \circ J\right) \bullet X\right) \left( \left( \left(
I\bullet \Sigma _{I}\right) \Delta _{I}^{\bullet }\circ J\right) \bullet
X\right)  \\
&=&\left( \Delta _{I}^{\bullet }\circ X\right) \left( I\circ l_{X}^{\circ
}\right) \psi _{I,I,X}^{\prime \prime }\left( \left( \left( I\bullet \Sigma
_{I}\right) \Delta _{I}^{\bullet }\circ J\right) \bullet X\right)  \\
&=&\left( \Delta _{I}^{\bullet }\circ X\right) \left( I\circ l_{X}^{\circ
}\right) \left( I\circ l_{I\circ X}^{\bullet }\right) \gamma
_{I,I,J,X}\left( \left( \left( I\bullet \Sigma _{I}\right) \Delta
_{I}^{\bullet }\circ J\right) \bullet X\right)  \\
&=&\left( \Delta _{I}^{\bullet }\circ X\right) \left( I\circ l_{X}^{\bullet
}\right) \left( I\circ \left( J\bullet l_{X}^{\circ }\right) \right) \gamma
_{I,I,J,X}\left( \left( \left( I\bullet \Sigma _{I}\right) \Delta
_{I}^{\bullet }\circ J\right) \bullet X\right)  \\
&=&\left( \Delta _{I}^{\bullet }\circ X\right) \left( I\circ l_{X}^{\bullet
}\right) \overline{\varsigma }_{I,J,X}^{I}.
\end{eqnarray*}
We check \eqref{form:nudelta}:
\begin{eqnarray*}
&&\nu _{I,Q}\left( \left( \Delta _{I}^{\bullet }\circ J\circ J\right)
\bullet Q\right)  \\
&=&\left( I\circ l_{Q}^{\circ }l_{I\circ Q}^{\bullet }\right) \gamma
_{I,I,J,Q}\left( \left( \left( I\bullet \phi _{0}^{\circ }\right) \circ
J\right) \bullet Q\right) \left( \left( \left( I\bullet I\right) \circ
m_{J}^{\circ }\right) \bullet Q\right) \left( \left( \Delta _{I}^{\bullet
}\circ J\circ J\right) \bullet Q\right)  \\
&=&\left( I\circ l_{Q}^{\bullet }\right) \underbracket[0.140ex]{\left(
I\circ \left( J\bullet l_{Q}^{\circ }\right) \right) \gamma _{I,I,J,Q}\left(
\left( \left( I\bullet \phi^\circ_0\right) \Delta _{I}^{\bullet }\circ
J\right) \bullet Q\right) }\left( \left( I\circ m_{J}^{\circ }\right)
\bullet Q\right)  \\
&=&\left( I\circ l_{Q}^{\bullet }\right) \overline{\varsigma }%
_{I,J,Q}^{I}\left( \left( I\circ m_{J}^{\circ }\right) \bullet Q\right).
\end{eqnarray*}
We check \eqref{form:nu1delta}:
\begin{eqnarray*}
&&\nu _{P,I}^{\prime }\left( P\bullet \left( J\circ J\circ \Delta
_{I}^{\bullet }\right) \right)  \\
&=&\left( r_{P}^{\circ }r_{P\circ I}^{\bullet }\circ I\right) \delta
_{P,J,I,I}\left( P\bullet \left( m_{J}^{\circ }\circ \left( \phi _{0}^{\circ
}\bullet I\right) \Delta _{I}^{\bullet }\right) \right)  \\
&=&\left( r_{P}^{\bullet }\circ I\right) \underbracket[0.140ex]{\left(
\left( r_{P}^{\circ }\bullet J\right) \circ I\right) \delta _{P,J,I,I}\left(
P\bullet \left( J\circ \left( \phi^\circ_0\bullet I\right) \Delta
_{I}^{\bullet }\right) \right) }\left( P\bullet \left( m_{J}^{\circ }\circ
I\right) \right)  \\
&=&\left( r_{P}^{\bullet }\circ I\right) \overline{\varkappa }%
_{P,J,I}^{I}\left( P\bullet \left( m_{J}^{\circ }\circ I\right) \right)
\end{eqnarray*}.
\end{invisible}    
\end{proof}

We are now ready to show that, when $H$ is a Hopf monoid, the Galois and co-Galois maps are invertible and the inverses are given explicitly in terms of the antipode.

\begin{theorem}%[\rd{upgrade:2026/03/30}]
\label{thm:Galinv} Let $\left( \mathcal{C},\circ ,I,\bullet ,J\right) $ be a
duoidal category with a reversion and let $\left( H,m_{H}^{\circ
},u_{H}^{\circ },\Delta _{H}^{\bullet },\varepsilon _{H}^{\bullet }\right) $
be a bimonoid with a morphism $\Sigma _{H}:H\rightarrow H^\revsn.$ Then

\begin{itemize}
\item[$1)$] $\overline{\varkappa }_{P,X,Q}^{H}\varkappa _{P,X,Q}^{H}=\mathrm{%
Id}$ if and only if $\overline{\varkappa }_{H,I,H}^{H}\varkappa _{H,I,H}^{H}=%
\mathrm{Id}$ if and only if $\mathrm{Id}_{H}\star \Sigma _{H}=1_{H,H}^{l};$

\item[$2)$] $\varkappa _{P,X,Q}^{H}\overline{\varkappa }_{P,X,Q}^{H}=\mathrm{%
Id}$ if $\Sigma _{H}\star \mathrm{Id}_{H}=1_{H,H}^{r};$

\item[$3)$] $\overline{\varsigma }_{P,X,Q}^{H}\varsigma _{P,X,Q}^{H}=\mathrm{%
Id}$  if and only if $\overline{\varsigma }_{H,I,H}^{H}\varsigma
_{H,I,H}^{H}=\mathrm{Id}$ if and only if $\Sigma _{H}\star \mathrm{Id}%
_{H}=1_{H,H}^{r};$

\item[$4)$] $\varsigma _{P,X,Q}^{H}\overline{\varsigma }_{P,X,Q}^{H}=\mathrm{%
Id}$ if $\mathrm{Id}_{H}\star \Sigma _{H}=1_{H,H}^{l}.$
\end{itemize}

As a consequence, both the Galois and the co-Galois maps are invertible with
$\overline{\varkappa }_{P,X,Q}^{H}=\left( \varkappa _{P,X,Q}^{H}\right) ^{-1}
$ and $\overline{\varsigma }_{P,X,Q}^{H}=\left( \varsigma
_{P,X,Q}^{H}\right) ^{-1}$ if and only if $\Sigma _{H}$ is an antipode.
\end{theorem}

\begin{proof}
We prove points 1) and 2) concerning the Galois maps.

Recall that by \cref{coro:cogal1}, we have $\overline{\varsigma }%
_{I,X,Y}^{I}=\left( \varsigma _{I,X,Y}^{I}\right) ^{-1}$ for any $X,Y$ in $\mathcal{C}$.

$1)$\ If $\mathrm{Id}_{H}\star \Sigma _{H}=1_{H,H}^{l},$ we prove \eqref{eq:GalinvGal} as follows%
\begin{eqnarray*}
&&\overline{\varkappa }_{H,I,Q}^{H}\underbracket[0.140ex]{\varkappa _{H,I,Q}^{H}\left(
\left( u_{H}^{\circ }\bullet I\right) \Delta _{I}^{\bullet }\circ Q\right)
\left( l_{Q}^{\circ }\right) ^{-1}}\overset{\eqref{form:Galunit}}{=}\overline{%
\varkappa }_{H,I,Q}^{H}\left( H\bullet \left( l_{Q}^{\circ }\right)
^{-1}\right) \rho _{Q}^{\bullet } \\
&=&\left( \left( m_{H}^{\circ }\bullet I\right) \circ Q\right) \delta
_{H,I,H,Q}\left( H\bullet \left( I\circ \left( \Sigma _{H}\bullet Q\right)
\rho _{Q}^{\bullet }\right) \right) \left( H\bullet \left( l_{Q}^{\circ
}\right) ^{-1}\right) \rho _{Q}^{\bullet } \\
&=&\left( \left( m_{H}^{\circ }\bullet I\right) \circ Q\right) \delta
_{H,I,H,Q}\left( H\bullet \left( l_{H^{\revsn}\bullet Q}^{\circ }\right)
^{-1}\right) \left( H\bullet \left( \Sigma _{H}\bullet Q\right) \rho
_{Q}^{\bullet }\right) \rho _{Q}^{\bullet } \\
&=&\left( \left( m_{H}^{\circ }\bullet I\right) \circ Q\right) \varphi
_{H,H,Q}^{\prime }\left( H\bullet \left( \Sigma _{H}\bullet Q\right) \rho
_{Q}^{\bullet }\right) \rho _{Q}^{\bullet } \\
&=&\left( \left( m_{H}^{\circ }\bullet I\right) \circ Q\right) %
\underbracket[0.140ex]{\varphi _{H,H,Q}^{\prime }}\left( \left( H\bullet
\Sigma _{H}\right) \Delta _{H}^{\bullet }\bullet Q\right) \rho _{Q}^{\bullet
} \\
&\overset{\eqref{form:varphi1-rev}}{=}&\left( \left( m_{H}^{\circ }\bullet
I\right) \circ Q\right) \xi _{H\circ H,Q}^{\prime }\left( \varphi
_{H,H}\bullet Q\right) \left( \left( H\bullet \Sigma _{H}\right) \Delta
_{H}^{\bullet }\bullet Q\right) \rho _{Q}^{\bullet } \\
&=&\xi _{H,Q}^{\prime }\left( \left( \underbracket[0.140ex]{\left(
m_{H}^{\circ }\bullet I\right) \circ J}\right) \bullet Q\right) \left( %
\underbracket[0.140ex]{\varphi _{H,H}}\bullet Q\right) \left( %
\underbracket[0.140ex]{\left( H\bullet \Sigma _{H}\right) \Delta
_{H}^{\bullet }}\bullet Q\right) \rho _{Q}^{\bullet } \\
&=&\xi _{H,Q}^{\prime }\left( \left( \mathrm{Id}_{H}\star \Sigma _{H}\right)
\bullet Q\right) \rho _{Q}^{\bullet } \\
&\overset{\eqref{eq:antipode}}=&\xi _{H,Q}^{\prime }\left( \left( \left( u_{H}^{\circ }\bullet I\right)
\circ J\right) \bullet Q\right) \left( \left( \Delta _{I}^{\bullet }\circ
J\right) \bullet Q\right) \left( \left( l_{J}^{\circ }\right)
^{-1}\varepsilon _{H}^{\bullet }\bullet Q\right) \rho _{Q}^{\bullet } \\
&=&\left( \left( u_{H}^{\circ }\bullet I\right) \circ Q\right) %
\underbracket[0.140ex]{\xi _{I,Q}^{\prime }\left( \left( \Delta
_{I}^{\bullet }\circ J\right) \bullet Q\right) }\left( \left( l_{J}^{\circ
}\right) ^{-1}\bullet Q\right) \left( l_{Q}^{\bullet }\right) ^{-1} \\
&\overset{\eqref{form:xi1delta}}{=}&\left( \left( u_{H}^{\circ }\bullet
I\right) \circ Q\right) \left( \Delta _{I}^{\bullet }\circ Q\right) \left(
I\circ l_{Q}^{\bullet }\right) \underbracket[0.140ex]{\overline{\varsigma
}_{I,J,Q}^{I}\left( \left( l_{J}^{\circ }\right) ^{-1}\bullet Q\right) }%
\left( l_{Q}^{\bullet }\right) ^{-1} \\
&=&\left( \left( u_{H}^{\circ }\bullet I\right) \Delta _{I}^{\bullet }\circ
Q\right) \left( I\circ l_{Q}^{\bullet }\right) \left( l_{J\bullet Q}^{\circ
}\right) ^{-1}\left( l_{Q}^{\bullet }\right) ^{-1} 
=\left( \left( u_{H}^{\circ }\bullet I\right) \Delta _{I}^{\bullet }\circ
Q\right) \left( l_{Q}^{\circ }\right) ^{-1}.
\end{eqnarray*}%
By \cref{prop:cogal2}, this implies $\overline{\varkappa }_{P,X,Q}^{H}%
\varkappa _{P,X,Q}^{H}=\mathrm{Id.}$ Conversely, if $\overline{\varkappa }%
_{P,X,Q}^{H}\varkappa _{P,X,Q}^{H}=\mathrm{Id}$ is true for every $P,X,Q,$
in particular $\overline{\varkappa }_{H,I,H}^{H}\varkappa _{H,I,H}^{H}=%
\mathrm{Id.}$ If the latter holds, then%
\begin{eqnarray*}
1_{H,H}^{l} &=&\left( \left( u_{H}^{\circ }\bullet I\right) \Delta
_{I}^{\bullet }\circ J\right) \left( l_{J}^{\circ }\right) ^{-1}\varepsilon
_{H}^{\bullet } \\
&=&\left( \left( H\bullet I\right) \circ \varepsilon _{H}^{\bullet }\right)
\left( \left( u_{H}^{\circ }\bullet I\right) \Delta _{I}^{\bullet }\circ
H\right) \left( l_{H}^{\circ }\right) ^{-1} \\
&=&\left( \left( H\bullet I\right) \circ \varepsilon _{H}^{\bullet }\right)
\overline{\varkappa }_{H,I,H}^{H}\underbracket[0.140ex]{\varkappa _{H,I,H}^{H}\left(
\left( u_{H}^{\circ }\bullet I\right) \Delta _{I}^{\bullet }\circ H\right)
\left( l_{H}^{\circ }\right) ^{-1}} \\
&\overset{\eqref{form:Galunit}}{=}&\left( \left( H\bullet I\right) \circ
\varepsilon _{H}^{\bullet }\right) \overline{\varkappa }_{H,I,H}^{H}\left(
H\bullet \left( l_{H}^{\circ }\right) ^{-1}\right) \Delta _{H}^{\bullet } \\
&=&\left( \left( H\bullet I\right) \circ \varepsilon _{H}^{\bullet }\right)
\left( \left( m_{H}^{\circ }\bullet I\right) \circ H\right) \delta
_{H,I,H,H}\left( H\bullet \left( I\circ \left( \Sigma _{H}\bullet H\right)
\Delta _{H}^{\bullet }\right) \right) \left( H\bullet \left( l_{H}^{\circ
}\right) ^{-1}\right) \Delta _{H}^{\bullet } \\
&=&\left( \left( m_{H}^{\circ }\bullet I\right) \circ J\right) \delta
_{H,I,H,J}\left( H\bullet \left( I\circ \left( \Sigma _{H}\bullet
\varepsilon _{H}^{\bullet }\right) \Delta _{H}^{\bullet }\right) \right)
\left( H\bullet \left( l_{H}^{\circ }\right) ^{-1}\right) \Delta
_{H}^{\bullet } \\
&=&\left( \left( m_{H}^{\circ }\bullet I\right) \circ J\right) \delta
_{H,I,H,J}\left( H\bullet \left( I\circ \left( \Sigma _{H}\bullet J\right)
\left( r_{H}^{\bullet }\right) ^{-1}\right) \right) \left( H\bullet \left(
l_{H}^{\circ }\right) ^{-1}\right) \Delta _{H}^{\bullet } \\
&=&\left( \left( m_{H}^{\circ }\bullet I\right) \circ J\right) \delta
_{H,I,H,J}\left( H\bullet \left( I\circ \left( r_{H^{\revsn}}^{\bullet
}\right) ^{-1}\right) \left( l_{H^{\revsn}}^{\circ }\right) ^{-1}\right)
\left( H\bullet \Sigma _{H}\right) \Delta _{H}^{\bullet } \\
&=&\left( \left( m_{H}^{\circ }\bullet I\right) \circ J\right) \varphi
_{H,H}\left( H\bullet \Sigma _{H}\right) \Delta _{H}^{\bullet }=\mathrm{Id}%
_{H}\star \Sigma _{H}
\end{eqnarray*}

$2)$\ Now, note that, by \cref{coro:cogal1}, we know that $\overline{\varkappa }%
_{P,J,I}^{I}=\left( \varkappa _{P,J,I}^{I}\right) ^{-1}$. Thus, if $\Sigma
_{H}\star \mathrm{Id}_{H}=1_{H,H}^{r},$ we get \eqref{eq:GalGalinv} as follows%
\begin{eqnarray*}
&&\underbracket[0.140ex]{r_{P}^{\bullet }\left( P\bullet m_{J}^{\circ
}\left( J\circ \varepsilon _{H}^{\bullet }\right) \right) \varkappa
_{P,J,H}^{H}}\overline{\varkappa }_{P,J,H}^{H} \\
&\overset{\eqref{form:Galcounit}}{=}&\mu _{P}^{\circ }\left( r_{P}^{\bullet
}\circ H\right) \overline{\varkappa }_{P,J,H}^{H} \\
&=&\mu _{P}^{\circ }\left( r_{P}^{\bullet }\circ H\right) \left( \left( \mu
_{P}^{\circ }\bullet J\right) \circ H\right) \delta _{P,J,H,H}\left(
P\bullet \left( J\circ \left( \Sigma _{H}\bullet H\right) \Delta
_{H}^{\bullet }\right) \right)  \\
&=&\mu _{P}^{\circ }\left( \mu _{P}^{\circ }\circ H\right) \left( r_{P\circ
H}^{\bullet }\circ H\right) \delta _{P,J,H,H}\left( P\bullet \left( J\circ
\left( \Sigma _{H}\bullet H\right) \Delta _{H}^{\bullet }\right) \right)  \\
&=&\mu _{P}^{\circ }\left( P\circ m_{H}^{\circ }\right) %
\underbracket[0.140ex]{\varphi _{P,H,H}^{\prime \prime }}\left( P\bullet
\left( J\circ \left( \Sigma _{H}\bullet H\right) \Delta _{H}^{\bullet
}\right) \right)  \\
&\overset{\eqref{form:varphi2-rev}}{=}&\mu _{P}^{\circ }\left( P\circ
m_{H}^{\circ }\right) \nu _{P,H\circ H}^{\prime }\left( P\bullet \left(
J\circ \psi _{H,H}\right) \right) \left( P\bullet \left( J\circ \left(
\Sigma _{H}\bullet H\right) \Delta _{H}^{\bullet }\right) \right)  \\
&=&\mu _{P}^{\circ }\nu _{P,H}^{\prime }\left( P\bullet \left( J\circ J\circ
\left( I\bullet m_{H}^{\circ }\right) \right) \right) \left( P\bullet \left(
J\circ \psi _{H,H}\right) \right) \left( P\bullet \left( J\circ \left(
\Sigma _{H}\bullet H\right) \Delta _{H}^{\bullet }\right) \right)  \\
&=&\mu _{P}^{\circ }\nu _{P,H}^{\prime }\left( P\bullet \left( J\circ \left(
\Sigma _{H}\star \mathrm{Id}_{H}\right) \right) \right)  \\
&\overset{\eqref{eq:antipode}}{=}&\mu _{P}^{\circ }\nu _{P,H}^{\prime
}\left( P\bullet \left( J\circ J\circ \left( I\bullet u_{H}^{\circ }\right)
\right) \right) \left( P\bullet \left( J\circ J\circ \Delta _{I}^{\bullet
}\right) \right) \left( P\bullet \left( J\circ \left( r_{J}^{\circ }\right)
^{-1}\right) \right) \left( P\bullet \left( J\circ \varepsilon _{H}^{\bullet
}\right) \right)  \\
&=&\mu _{P}^{\circ }\left( P\circ u_{H}^{\circ }\right) \underbracket[0.140ex]{\nu
_{P,I}^{\prime }\left( P\bullet \left( J\circ J\circ \Delta _{I}^{\bullet
}\right) \right) }\left( P\bullet \left( r_{J\circ J}^{\circ }\right)
^{-1}\right) \left( P\bullet \left( J\circ \varepsilon _{H}^{\bullet
}\right) \right)  \\
&\overset{\eqref{form:nu1delta}}=&r_{P}^{\circ }\left( r_{P}^{\bullet }\circ I\right) \overline{\varkappa }%
_{P,J,I}^{I}\left( P\bullet \left( m_{J}^{\circ }\circ I\right) \right)
\left( P\bullet \left( r_{J\circ J}^{\circ }\right) ^{-1}\right) \left(
P\bullet \left( J\circ \varepsilon _{H}^{\bullet }\right) \right)  %\\
%&=&r_{P}^{\bullet }r_{P\bullet J}^{\circ }\underbracket[0.140ex]{\left(
%\left( r_{P}^{\circ }\bullet %J\right) \circ I\right) \delta _{P,J,I,I}\left(
%P\bullet \left( J\circ \left( \phi^\circ_0\bullet I\right) \Delta
%_{I}^{\bullet }\right) \right) %}\left( P\bullet \left( r_{J}^{\circ }\right)
%^{-1}\right) \left( P\bullet m_{J}^{\circ }\left( J\circ \varepsilon
%_{H}^{\bullet }\right) \right)  
\\
&=&r_{P}^{\bullet }r_{P\bullet J}^{\circ }\overline{\varkappa }%
_{P,J,I}^{I}\left( P\bullet \left( r_{J}^{\circ }\right) ^{-1}\right) \left(
P\bullet m_{J}^{\circ }\left( J\circ \varepsilon _{H}^{\bullet }\right)
\right)  \\
&=&r_{P}^{\bullet }\underbracket[0.140ex]{r_{P\bullet J}^{\circ }\left(
\varkappa _{P,J,I}^{I}\right) ^{-1}\left( P\bullet \left( r_{J}^{\circ
}\right) ^{-1}\right) }\left( P\bullet m_{J}^{\circ }\left( J\circ
\varepsilon _{H}^{\bullet }\right) \right) =r_{P}^{\bullet }\left( P\bullet
m_{J}^{\circ }\left( J\circ \varepsilon _{H}^{\bullet }\right) \right)
\end{eqnarray*}%
where in the last step we used the formula of $\varkappa _{P,J,I}^{I}$
obtained in \cref{lem:inv}. By \cref{prop:cogal2}, this implies that $\varkappa
_{P,X,Q}^{H}\overline{\varkappa }_{P,X,Q}^{H}=\mathrm{Id.}$

In a similar way, by means of \eqref{form:coGalunit},
\eqref{form:psi1-rev}, \eqref{eq:antipode} and \eqref{form:xidelta}
one proves $3)$, while by means of \eqref{form:coGalcounit}, \eqref{form:psi2-rev}, \eqref{eq:antipode} and \eqref{form:nudelta} 
one gets $4).$ The last assertion is clear.
\begin{invisible}
3) Recall that by \cref{coro:cogal1}, we have that $\overline{\varkappa }%
_{X,Y,I}^{I}=\left( \varkappa _{X,Y,I}^{I}\right) ^{-1}$.

If $\Sigma _{H}\star \mathrm{Id}_{H}=1_{H,H}^{r},$ we get%
\begin{eqnarray*}
&&\overline{\varsigma }_{P,I,H}^{H}\underbracket[0.140ex]{\varsigma _{P,I,H}^{H}\left(
P\circ \left( I\bullet u_{H}^{\circ }\right) \Delta _{I}^{\bullet }\right)
\left( r_{P}^{\circ }\right) ^{-1}} \\
&\overset{\eqref{form:coGalunit}}{=}&\overline{\varsigma }_{P,I,H}^{H}\left(
\left( r_{P}^{\circ }\right) ^{-1}\bullet Z\right) \rho _{P}^{\bullet } \\
&=&\left( P\circ \left( I\bullet m_{H}^{\circ }\right) \right) \gamma
_{P,H,I,H}\left( \left( \left( P\bullet \Sigma _{H}\right) \rho
_{P}^{\bullet }\circ I\right) \bullet H\right) \left( \left( r_{P}^{\circ
}\right) ^{-1}\bullet Z\right) \rho _{P}^{\bullet } \\
&=&\left( P\circ \left( I\bullet m_{H}^{\circ }\right) \right) \gamma
_{P,H,I,H}\left( \left( r_{P\bullet H^{\revsn}}^{\circ }\right)
^{-1}\bullet Z\right) \left( \left( P\bullet \Sigma _{H}\right) \rho
_{P}^{\bullet }\bullet H\right) \rho _{P}^{\bullet } \\
&=&\left( P\circ \left( I\bullet m_{H}^{\circ }\right) \right) \psi
_{P,H,H}^{\prime }\left( P\bullet \Sigma _{H}\bullet H\right) \left( \rho
_{P}^{\bullet }\bullet H\right) \rho _{P}^{\bullet } \\
&=&\left( P\circ \left( I\bullet m_{H}^{\circ }\right) \right) %
\underbracket[0.140ex]{\psi _{P,H,H}^{\prime }}\left( P\bullet \left( \Sigma
_{H}\bullet H\right) \Delta _{H}^{\bullet }\right) \rho _{P}^{\bullet } \\
&\overset{\eqref{form:psi1-rev}}{=}&\underbracket[0.140ex]{\left( P\circ
\left( I\bullet m_{H}^{\circ }\right) \right) \xi _{P,H\circ H}}\left(
P\bullet \psi _{H,H}\right) \left( P\bullet \left( \Sigma _{H}\bullet
H\right) \Delta _{H}^{\bullet }\right) \rho _{P}^{\bullet } \\
&=&\xi _{P,H}\left( P\bullet \left( \underbracket[0.140ex]{J\circ \left(
I\bullet m_{H}^{\circ }\right) }\right) \right) \left( P\bullet %
\underbracket[0.140ex]{\psi _{H,H}}\right) \left( P\bullet %
\underbracket[0.140ex]{\left( \Sigma _{H}\bullet H\right) \Delta
_{H}^{\bullet }}\right) \rho _{P}^{\bullet } \\
&\overset{\text{antip.}}{=}&\xi _{P,H}\left( P\bullet \left( \Sigma _{H}\star
\mathrm{Id}_{H}\right) \right) \rho _{P}^{\bullet } \\
&\overset{\eqref{eq:antipode}}=&\xi _{P,H}\left( P\bullet \left( J\circ \left( I\bullet u_{H}^{\circ
}\right) \right) \right) \left( P\bullet \left( \left( J\circ \Delta
_{I}^{\bullet }\right) \left( r_{J}^{\circ }\right) ^{-1}\right) \right)
\left( P\bullet \varepsilon _{H}^{\bullet }\right) \rho _{P}^{\bullet } \\
&=&\left( P\circ \left( I\bullet u_{H}^{\circ }\right) \right) %
\underbracket[0.140ex]{\xi _{P,I}\left( P\bullet \left( J\circ \Delta
_{I}^{\bullet }\right) \right) }\left( P\bullet \left( r_{J}^{\circ }\right)
^{-1}\right) \left( P\bullet \varepsilon _{H}^{\bullet }\right) \rho
_{P}^{\bullet } \\
&\overset{\eqref{form:xidelta}}{=}&\left( P\circ \left( I\bullet
u_{H}^{\circ }\right) \Delta _{I}^{\bullet }\right) \left( r_{P}^{\bullet
}\circ I\right) \underbracket[0.140ex]{\left( \varkappa _{P,J,I}^{I}\right)
^{-1}}\left( P\bullet \left( r_{J}^{\circ }\right) ^{-1}\right) \left(
P\bullet \varepsilon _{H}^{\bullet }\right) \rho _{P}^{\bullet } \\
&=&\left( P\circ \left( I\bullet u_{H}^{\circ }\right) \Delta _{I}^{\bullet
}\right) \left( r_{P}^{\bullet }\circ I\right) \left( r_{P\bullet J}^{\circ
}\right) ^{-1}\left( P\bullet r_{J}^{\circ }\right) \left( P\bullet \left(
r_{J}^{\circ }\right) ^{-1}\right) \left( P\bullet \varepsilon _{H}^{\bullet
}\right) \rho _{P}^{\bullet } \\
&=&\left( P\circ \left( I\bullet u_{H}^{\circ }\right) \Delta _{I}^{\bullet
}\right) \left( r_{P}^{\bullet }\circ I\right) \left( r_{P\bullet J}^{\circ
}\right) ^{-1}\left( P\bullet \varepsilon _{H}^{\bullet }\right) \rho
_{P}^{\bullet } \\
&=&\left( P\circ \left( I\bullet u_{H}^{\circ }\right) \Delta _{I}^{\bullet
}\right) \left( r_{P}^{\circ }\right) ^{-1}r_{P}^{\bullet }\left( P\bullet
\varepsilon _{H}^{\bullet }\right) \rho _{P}^{\bullet }=\left( P\circ \left(
I\bullet u_{H}^{\circ }\right) \Delta _{I}^{\bullet }\right) \left(
r_{P}^{\circ }\right) ^{-1}.
\end{eqnarray*}%
By \cref{prop:cogal1}, we get that $\overline{\varsigma }_{P,X,Q}^{H}%
\varsigma _{P,X,Q}^{H}=\mathrm{Id.}$

We now see the converse. Thus, if $\overline{\varsigma }_{P,X,Q}^{H}%
\varsigma _{P,X,Q}^{H}=\mathrm{Id,}$ then $\overline{\varsigma }%
_{H,I,Q}^{H}\varsigma _{H,I,H}^{H}=\mathrm{Id}$ and so
\begin{eqnarray*}
1_{H,H}^{r} &=&\left( J\circ \left( I\bullet u_{H}^{\circ }\right) \Delta
_{I}^{\bullet }\right) \left( r_{J}^{\circ }\right) ^{-1}\varepsilon
_{H}^{\bullet } \\
&=&\left( \left( \varepsilon _{H}^{\bullet }\circ I\right) \bullet H\right)
\left( H\circ \left( I\bullet u_{H}^{\circ }\right) \Delta _{I}^{\bullet
}\right) \left( r_{H}^{\circ }\right) ^{-1} \\
&=&\left( \left( \varepsilon _{H}^{\bullet }\circ I\right) \bullet H\right)
\overline{\varsigma }_{H,I,H}^{H}\underbracket[0.140ex]{\varsigma _{H,I,H}^{H}\left(
H\circ \left( I\bullet u_{H}^{\circ }\right) \Delta _{I}^{\bullet }\right)
\left( r_{H}^{\circ }\right) ^{-1}} \\
&\overset{\eqref{form:coGalunit}}{=}&\left( \left( \varepsilon _{H}^{\bullet
}\circ I\right) \bullet H\right) \overline{\varsigma }_{H,I,H}^{H}\left(
\left( r_{H}^{\circ }\right) ^{-1}\bullet H\right) \Delta _{H}^{\bullet } \\
&=&\left( \left( \varepsilon _{H}^{\bullet }\circ I\right) \bullet H\right)
\left( H\circ \left( I\bullet m_{H}^{\circ }\right) \right) \gamma
_{H,H,I,H}\left( \left( \left( H\bullet \Sigma _{H}\right) \Delta
_{H}^{\bullet }\circ I\right) \bullet H\right) \left( \left( r_{H}^{\circ
}\right) ^{-1}\bullet H\right) \Delta _{H}^{\bullet } \\
&=&\left( J\circ \left( I\bullet m_{H}^{\circ }\right) \right) \gamma
_{J,H,I,H}\left( \left( \left( \varepsilon _{H}^{\bullet }\bullet \Sigma
_{H}\right) \Delta _{H}^{\bullet }\circ I\right) \bullet H\right) \left(
\left( r_{H}^{\circ }\right) ^{-1}\bullet H\right) \Delta _{H}^{\bullet } \\
&=&\left( J\circ \left( I\bullet m_{H}^{\circ }\right) \right) \gamma
_{J,H,I,H}\left( \left( \left( J\bullet \Sigma _{H}\right) \left(
l_{H}^{\bullet }\right) ^{-1}\circ I\right) \bullet H\right) \left( \left(
r_{H}^{\circ }\right) ^{-1}\bullet H\right) \Delta _{H}^{\bullet } \\
&=&\left( J\circ \left( I\bullet m_{H}^{\circ }\right) \right) \gamma
_{J,H,I,H}\left( \left( \left( l_{H^{\revsn}}^{\bullet }\right) ^{-1}\circ
I\right) \left( r_{H^{\revsn}}^{\circ }\right) ^{-1}\bullet H\right)
\left( \Sigma _{H}\bullet H\right) \Delta _{H}^{\bullet } \\
&=&\left( J\circ \left( I\bullet m_{H}^{\circ }\right) \right) \psi
_{H,H}\left( \Sigma _{H}\bullet H\right) \Delta _{H}^{\bullet }=\Sigma
_{H}\star \mathrm{Id}_{H}.
\end{eqnarray*}

4) Note that, by \cref{coro:cogal1}, we know that $\overline{\varsigma }%
_{I,J,Q}^{I}=\left( \varsigma _{I,J,Q}^{I}\right) ^{-1}$. As a consequence,
if we assume $\mathrm{Id}_{H}\star \Sigma _{H}=1_{H,H}^{l},$ we get%
\begin{eqnarray*}
&&\underbracket[0.140ex]{l_{Q}^{\bullet }\left( m_{J}^{\circ }\left(
\varepsilon _{H}^{\bullet }\circ J\right) \bullet Q\right) \varsigma
_{H,J,Q}^{H}}\overline{\varsigma }_{H,J,Q}^{H} \\
&\overset{\eqref{form:coGalcounit}}{=}&\mu _{Q}^{\circ }\left( H\circ
l_{Q}^{\bullet }\right) \overline{\varsigma }_{H,J,Q}^{H} \\
&=&\mu _{Q}^{\circ }\left( H\circ l_{Q}^{\bullet }\right) \left( H\circ
\left( J\bullet \mu _{Q}^{\circ }\right) \right) \gamma _{H,H,J,Q}\left(
\left( \left( H\bullet \Sigma _{H}\right) \Delta _{H}^{\bullet }\circ
J\right) \bullet Q\right)  \\
&=&\mu _{Q}^{\circ }\left( H\circ \mu _{Q}^{\circ }\right) \left( H\circ
l_{H\circ Q}^{\bullet }\right) \gamma _{H,H,J,Q}\left( \left( \left(
H\bullet \Sigma _{H}\right) \Delta _{H}^{\bullet }\circ J\right) \bullet
Q\right)  \\
&=&\mu _{Q}^{\circ }\left( H\circ \mu _{Q}^{\circ }\right) %
\underbracket[0.140ex]{\psi _{H,H,Q}^{\prime \prime }}\left( \left( \left(
H\bullet \Sigma _{H}\right) \Delta _{H}^{\bullet }\circ J\right) \bullet
Q\right)  \\
&\overset{\eqref{form:psi2-rev}}{=}&\mu _{Q}^{\circ }\left( m_{H}^{\circ
}\circ Q\right) \nu _{H\circ H,Q}\left( \left( \varphi _{H,H}\circ J\right)
\bullet Q\right) \left( \left( \left( H\bullet \Sigma _{H}\right) \Delta
_{H}^{\bullet }\circ J\right) \bullet Q\right)  \\
&=&\mu _{Q}^{\circ }\nu _{H,Q}\left( \left( \underbracket[0.140ex]{\left(
m_{H}^{\circ }\bullet I\right) \circ J}\circ J\right) \bullet Q\right)
\left( \left( \underbracket[0.140ex]{\varphi _{H,H}}\circ J\right) \bullet
Q\right) \left( \left( \underbracket[0.140ex]{\left( H\bullet \Sigma
_{H}\right) \Delta _{H}^{\bullet }}\circ J\right) \bullet Q\right)  \\
&=&\mu _{Q}^{\circ }\nu _{H,Q}\left( \left( \underbracket[0.140ex]{%
\mathrm{Id}_{H}\star \Sigma _{H}}\right) \bullet Q\right)  \\
&\overset{\eqref{eq:antipode}}{=}&\mu _{Q}^{\circ }\nu _{H,Q}\left( \left(
\left( u_{H}^{\circ }\bullet I\right) \circ J\circ J\right) \bullet Q\right)
\left( \left( \Delta _{I}^{\bullet }\circ J\circ J\right) \bullet Q\right)
\left( \left( \left( l_{J}^{\circ }\right) ^{-1}\varepsilon _{H}^{\bullet
}\circ J\right) \bullet Q\right)  \\
&=&\mu _{Q}^{\circ }\left( u_{H}^{\circ }\circ Q\right) \underbracket[0.140ex]{\nu
_{I,Q}\left( \left( \Delta _{I}^{\bullet }\circ J\circ J\right) \bullet
Q\right) }\left( \left( \left( l_{J}^{\circ }\right) ^{-1}\circ J\right)
\bullet Q\right) \left( \left( \varepsilon _{H}^{\bullet }\circ J\right)
\bullet Q\right)  \\
&\overset{\eqref{form:nudelta}}=&l_{Q}^{\circ }\left( I\circ l_{Q}^{\bullet }\right) \overline{\varsigma }%
_{I,J,Q}^{I}\left( \left( I\circ m_{J}^{\circ }\right) \bullet Q\right)
\left( \left( l_{J\circ J}^{\circ }\right) ^{-1}\bullet Q\right) \left(
\left( \varepsilon _{H}^{\bullet }\circ J\right) \bullet Q\right)  \\
&=&l_{Q}^{\bullet }l_{J\bullet Q}^{\circ }\overline{\varsigma }%
_{I,J,Q}^{I}\left( \left( l_{J}^{\circ }\right) ^{-1}\bullet Q\right) \left(
m_{J}^{\circ }\left( \varepsilon _{H}^{\bullet }\circ J\right) \bullet
Q\right)  \\
&=&l_{Q}^{\bullet }\underbracket[0.140ex]{l_{J\bullet Q}^{\circ }\left(
\varsigma _{I,J,Q}^{I}\right) ^{-1}\left( \left( l_{J}^{\circ }\right)
^{-1}\bullet Q\right) }\left( m_{J}^{\circ }\left( \varepsilon _{H}^{\bullet
}\circ J\right) \bullet Q\right) =l_{Q}^{\bullet }\left( m_{J}^{\circ
}\left( \varepsilon _{H}^{\bullet }\circ J\right) \bullet Q\right)
\end{eqnarray*}
where in the last step we used the formula of $\left( \varsigma
_{I,J,Q}^{I}\right) ^{-1}$ obtained in \cref{lem:inv}.

By \cref{prop:cogal1}, we get that $\varsigma _{P,X,Q}^{H}\overline{%
\varsigma }_{P,X,Q}^{H}=\mathrm{Id.}$
\end{invisible}
\end{proof}

By employing \cref{lem:rev} and \cref{lem:rev2}, we now prove the following analogue of \cite[Lemma 4.3]{Bohm-Lack}, which will play a key role later on.

\begin{lemma}
\label{lem:Bohm4.3}
Define the morphisms $\vartheta_{X,Y,Z}:X^{\revsn }\bullet \left( J\circ \left( I\bullet Y\right)
\right) \bullet Z\to J\circ \left( I\bullet Y\bullet \left(
X\circ Z\right) \right)$ and $\vartheta_{X,Y,Z}':X\bullet((Y\bullet I)\circ J)\bullet Z^{\revsn}\to((X\circ Z)\bullet Y\bullet I)\circ J$ as the following compositions
{\small\[\begin{tikzcd}
	X^{\revsn }\bullet \left( J\circ \left( I\bullet Y\right)
\right) \bullet Z & \left( X^{\revsn }\circ \left( I\bullet Y\right) \right) \bullet Z \\
	J\circ \left( I\bullet Y\bullet \left(
X\circ Z\right) \right) & \left( \left( J\bullet X^{\revsn
}\right) \circ \left( I\bullet Y\right) \right) \bullet Z
	\arrow[from=1-1, to=1-2, "\xi _{X^{\revsn },Y}\bullet Z"]
	\arrow[from=1-1, to=2-1, dotted, "\vartheta_{X,Y,Z}"']
	\arrow[from=1-2, to=2-2,"\left( \left( l_{X^{\revsn }}^{\bullet }\right) ^{-1}\circ \left( I\bullet
Y\right) \right) \bullet Z"]
	\arrow[from=2-2, to=2-1,"\gamma
_{J,X,I\bullet Y,Z}"]
\end{tikzcd}\quad
\begin{tikzcd}
	X\bullet \left( \left( Y\bullet I\right)
\circ J\right) \bullet Z^{\revsn } & X\bullet \left( \left( Y\bullet I\right) \circ Z^{\revsn
}\right)  \\
	\left( \left(
X\circ Z\right) \bullet Y\bullet I\right) \circ J & X\bullet
\left( \left( Y\bullet I\right) \circ \left( Z^{\revsn }\bullet J\right)
\right)
	\arrow[from=1-1, to=1-2, "X\bullet \xi _{Y,Z^{\revsn }}^{\prime}"]
	\arrow[from=1-1, to=2-1, dotted, "\vartheta'_{X,Y,Z}"']
	\arrow[from=1-2, to=2-2,"X\bullet \left( \left( Y\bullet I\right) \circ \left(
r_{Z^{\revsn }}^{\bullet }\right) ^{-1}\right)"]
	\arrow[from=2-2, to=2-1,"\delta _{X,Y\bullet I,Z,J}"]
\end{tikzcd}\]}
Then, the following diagrams commute:
\begin{equation}\label{form:vartheta}
\begin{tikzcd}
	X^{\revsn}\bullet Y^{\revsn}\bullet Z\bullet T && X^{\revsn }\bullet \left( J\circ \left(
I\bullet \left( Y\circ Z\right) \right) \right) \bullet T \\
	&& J\circ \left( I\bullet \left( Y\circ
Z\right) \bullet \left( X\circ T\right) \right) \\
	\left( Y\bullet X\right) ^{\revsn
}\bullet Z\bullet T && J\circ \left( I\bullet \left( \left( Y\bullet X\right) \circ \left( Z\bullet
T\right) \right) \right)
	\arrow[from=1-1, to=1-3,"X^{\revsn }\bullet
\psi _{Y,Z}\bullet T"]
	\arrow[from=1-1, to=3-1, "\phi _{X,Y}^{\bullet
}\bullet Z\bullet T"']
	\arrow[from=1-3, to=2-3,"\vartheta
_{X,Y\circ Z,T}"]
	\arrow[from=3-1, to=3-3,"\psi _{Y\bullet X,Z\bullet T}"']
	\arrow[from=3-3, to=2-3,"J\circ \left( I\bullet \zeta
_{Y,X,Z,T}\right) "']
\end{tikzcd}
\end{equation}
% \begin{equation}
% \begin{tikzcd}[column sep=huge, row sep=large]
% X^{\revsn}\bullet Y^{\revsn}\bullet Z \bullet T
%   \arrow[r,
%     "X^{\revsn }\bullet
% \psi _{Y,Z}\bullet T"]
%   \arrow[d,
%     "\phi _{X,Y}^{\bullet
% }\bullet Z\bullet T"']
% &
% X^{\revsn }\bullet \left( J\circ \left(
% I\bullet \left( Y\circ Z\right) \right) \right) \bullet T
%   \arrow[r,
%     "\vartheta
% _{X,Y\circ Z,T}"]
% &
% J\circ \left( I\bullet \left( Y\circ
% Z\right) \bullet \left( X\circ T\right) \right)
% \\
% \left( Y\bullet X\right) ^{\revsn
% }\bullet Z\bullet T
%   \arrow[r,
%     "\psi _{Y\bullet X,Z\bullet T}"']
% &
% J\circ \left( I\bullet \left( \left( Y\bullet X\right) \circ \left( Z\bullet
% T\right) \right) \right)
%   \arrow[ur,
%     "J\circ \left( I\bullet \zeta
% _{Y,X,Z,T}\right) "']
% \end{tikzcd}
% \end{equation}
\begin{equation}\label{form:vartheta1}
\begin{tikzcd}
	X\bullet Y\bullet Z^{\revsn} \bullet T^{\revsn} && X\bullet( \left( \left( Y\circ
Z\right) \bullet I\right) \circ J)\bullet T^{\revsn } \\
	&& \left( \left( X\circ T\right)
\bullet \left( Y\circ Z\right)  \bullet I\right) \circ J \\
	X\bullet Y\bullet \left( T\bullet
Z\right) ^{\revsn } && \left( \left( \left( X\bullet Y\right) \circ \left( T\bullet Z\right)
\right) \bullet I\right) \circ J
	\arrow[from=1-1, to=1-3,"X\bullet
\varphi _{Y,Z}\bullet T^{\revsn }"]
	\arrow[from=1-1, to=3-1,"X\bullet Y\bullet\phi _{Z,T}^{\bullet
}"']
	\arrow[from=1-3, to=2-3,"\vartheta'
_{X,Y\circ Z,T}"]
	\arrow[from=3-1, to=3-3,"\varphi _{X\bullet Y,T\bullet Z}"']
	\arrow[from=3-3, to=2-3,"(\zeta
_{X,Y,T,Z}\bullet I)\circ J"']
\end{tikzcd}
\end{equation}
% \begin{equation}
% \begin{tikzcd}[column sep=huge, row sep=large]
% X\bullet Y\bullet Z^{\revsn} \bullet T^{\revsn}
%   \arrow[r,
%     "X\bullet
% \varphi _{Y,Z}\bullet T^{\revsn }"]
%   \arrow[d,
%     "X\bullet Y\bullet\phi _{Z,T}^{\bullet
% }"']
% &
% X\bullet \left( \left( Y\circ
% Z\right) \bullet I\right) \circ J\bullet T^{\revsn }
%   \arrow[r,
%     "\vartheta'
% _{X,Y\circ Z,T}"]
% &
% \left( \left( X\circ T\right)
% \bullet \left( Y\circ Z\right)  \bullet I\right) \circ J
% \\
% X\bullet Y\bullet \left( T\bullet
% Z\right) ^{\revsn }
%   \arrow[r,
%     "\varphi _{X\bullet Y,T\bullet Z}"']
% &
% \left( \left( \left( X\bullet Y\right) \circ \left( T\bullet Z\right)
% \right) \bullet I\right) \circ J
%   \arrow[ur,
%     "(\zeta
% _{X,Y,T,Z}\bullet I)\circ J"']
% \end{tikzcd}
% \end{equation}
Moreover, the following equalities hold:
\begin{equation}\label{form:varthetaI}
\vartheta _{X,I,X}\left( X^{\revsn }\bullet \left( J\circ \Delta _{I}^{\bullet
}\right) \bullet X\right) =\left( J\circ \left( \Delta _{I}^{\bullet
}\bullet \left( X\circ X\right) \right) \right) \psi _{X,X}\left( r_{X^{\revsn
}}^{\bullet }\left( X^{\revsn }\bullet r_{J}^{\circ }\right)
\bullet X\right)
\end{equation}%
\begin{equation}
\label{form:vartheta1I}
\vartheta _{X,I,X}^{\prime }\left( X\bullet \left( \Delta _{I}^{\bullet
}\circ J\right) \bullet X^{\revsn }\right) =\left( \left( \left( X\circ
X\right) \bullet \Delta _{I}^{\bullet }\right) \circ J\right) \varphi
_{X,X}\left( X\bullet l_{X^{\revsn }}^{\bullet }\left( l_{J}^{\circ
}\bullet X^{\revsn }\right) \right)
\end{equation}%
\end{lemma}

\begin{proof}
Let us prove \eqref{form:vartheta}:
\begin{eqnarray*}
&&\vartheta _{X,Y\circ Z,T}\left( X^{\revsn }\bullet \psi _{Y,Z}\bullet
T\right)  \\
&=&\gamma _{J,X,I\bullet \left( Y\circ Z\right) ,T}\left( \left( \left(
l_{X^{\revsn }}^{\bullet }\right) ^{-1}\circ \left( I\bullet \left( Y\circ
Z\right) \right) \right) \bullet T\right) \left( \underbracket[0.140ex]{\xi _{X^{\revsn
},Y\circ Z}}\bullet T\right) \left( \underbracket[0.140ex]{X^{\revsn }\bullet \psi _{Y,Z}}%
\bullet T\right)  \\
&\overset{\eqref{form:psi1-rev} }{=}&\gamma _{J,X,I\bullet
\left( Y\circ Z\right) ,T}\left( \left( \left( l_{X^{\revsn }}^{\bullet
}\right) ^{-1}\circ \left( I\bullet \left( Y\circ Z\right) \right) \right)
\bullet T\right) \left( \psi _{X^{\revsn },Y,Z}^{\prime }\bullet T\right)  \\
&=&\gamma _{J,X,I\bullet \left( Y\circ Z\right) ,T}\left( \psi _{J\bullet
X^{\revsn },Y,Z}^{\prime }\bullet T\right) \left( \left( l_{X^{\revsn
}}^{\bullet }\right) ^{-1}\bullet Y^{\revsn }\bullet Z\bullet T\right) \\
&=&\underbracket[0.140ex]{\gamma _{J,X,I\bullet \left( Y\circ Z\right) ,T}\left( \gamma
_{J\bullet X^{\revsn },Y,I,Z}\bullet T\right) }\left( \left( r_{J\bullet
X^{\revsn }\bullet Y^{\revsn }}^{\circ }\right) ^{-1}\bullet Z\bullet T\right)
\left( \left( l_{X^{\revsn }\bullet Y^{\revsn }}^{\bullet }\right) ^{-1}\bullet
Z\bullet T\right)  \\
&\overset{\eqref{form: psi1gamma}}{=}&\left( J\circ \left( I\bullet \zeta
_{Y,X,Z,T}\right) \right) \gamma _{J,Y\bullet X,I,Z\bullet T}\left( \left(
\left( J\bullet \phi _{X,Y}^{\bullet }\right) \circ I\right) \bullet
Z\bullet T\right)\\
&&\left( \left( \left( l_{X^{\revsn }\bullet Y^{\revsn
}}^{\bullet }\right) ^{-1}\circ I\right) \bullet Z\bullet T\right) \left(
\left( r_{X^{\revsn }\bullet Y^{\revsn }}^{\circ }\right) ^{-1}\bullet Z\bullet
T\right)  \\
&=&\left( J\circ \left( I\bullet \zeta _{Y,X,Z,T}\right) \right) \gamma
_{J,Y\bullet X,I,Z\bullet T}\left( \left( \left( l_{\left( Y\bullet X\right)
^{\revsn }}^{\bullet }\right) ^{-1}\circ I\right) \bullet Z\bullet T\right)\\
&&
\left( \left( r_{\left( Y\bullet X\right) ^{\revsn }}^{\circ }\right)
^{-1}\bullet Z\bullet T\right) \left( \phi _{X,Y}^{\bullet }\bullet Z\bullet
T\right)  \\
&=&\left( J\circ \left( I\bullet \zeta _{Y,X,Z,T}\right) \right) \psi
_{Y\bullet X,Z\bullet T}\left( \phi _{X,Y}^{\bullet }\bullet Z\bullet
T\right).
\end{eqnarray*}
Similarly, one proves \eqref{form:vartheta1} by using \eqref{form:varphi1-rev}. 
\begin{invisible}
\begin{eqnarray*}
&&\vartheta _{X,Y\circ Z,T}^{\prime }\left( X\bullet \varphi _{Y,Z}\bullet
T^{\revsn }\right)  \\
&=&\delta _{X,\left( Y\circ Z\right) \bullet I,T,J}\left( X\bullet \left(
\left( \left( Y\circ Z\right) \bullet I\right) \circ \left( r_{T^{\revsn
}}^{\bullet }\right) ^{-1}\right) \right) \left( X\bullet \underbracket[0.140ex]{\xi
_{Y\circ Z,T^{\revsn }}^{\prime }}\right) \left( X\bullet \underbracket[0.140ex]{\varphi
_{Y,Z}\bullet T^{\revsn }}\right)  \\
&\overset{\eqref{form:varphi1-rev} }{=}&\delta _{X,\left( Y\circ
Z\right) \bullet I,T,J}\left( X\bullet \left( \left( \left( Y\circ Z\right)
\bullet I\right) \circ \left( r_{T^{\revsn }}^{\bullet }\right) ^{-1}\right)
\right) \left( X\bullet \varphi _{Y,Z,T^{\revsn }}^{\prime }\right)  \\
&=&\underbracket[0.140ex]{\delta _{X,\left( Y\circ Z\right) \bullet I,T,J}\left(
X\bullet \varphi _{Y,Z,T^{\revsn }\bullet J}^{\prime }\right) \left( X\bullet
Y\bullet Z^{\revsn }\bullet \left( r_{T^{\revsn }}^{\bullet }\right)
^{-1}\right) } \\
&=&\delta _{X,\left( Y\circ Z\right) \bullet I,T,J}\left( X\bullet \delta
_{Y,I,Z,T^{\revsn }\bullet J}\right) \left( Y\bullet X\bullet \left(
l_{Z^{\revsn }\bullet T^{\revsn }\bullet J}^{\circ }\right) ^{-1}\right) \left(
X\bullet Y\bullet Z^{\revsn }\bullet \left( r_{T^{\revsn }}^{\bullet }\right)
^{-1}\right)  \\
&=&\underbracket[0.140ex]{\delta _{X,\left( Y\circ Z\right) \bullet I,T,J}\left(
X\bullet \delta _{Y,I,Z,T^{\revsn }\bullet J}\right) }\left( Y\bullet X\bullet
\left( l_{Z^{\revsn }\bullet T^{\revsn }\bullet J}^{\circ }\right) ^{-1}\right)
\left( X\bullet Y\bullet \left( r_{Z^{\revsn }\bullet T^{\revsn }}^{\bullet
}\right) ^{-1}\right)  \\
&\overset{\eqref{form:varphi1delta}}{=}&\left( \left( \zeta
_{X,Y,T,Z}\bullet I\right) \circ J\right) \delta _{X\bullet Y,I,T\bullet
Z,J}\left( X\bullet Y\bullet \left( I\circ \left( \phi _{Z,T}^{\bullet
}\bullet J\right) \right) \right) \left( Y\bullet X\bullet \left( l_{Z^{\revsn
}\bullet T^{\revsn }\bullet J}^{\circ }\right) ^{-1}\right) \left( X\bullet
Y\bullet \left( r_{Z^{\revsn }\bullet T^{\revsn }}^{\bullet }\right)
^{-1}\right)  \\
&=&\left( \left( \zeta _{X,Y,T,Z}\bullet I\right) \circ J\right) \delta
_{X\bullet Y,I,T\bullet Z,J}\left( X\bullet Y\bullet \left( I\circ \left(
\phi _{Z,T}^{\bullet }\bullet J\right) \right) \right) \left( X\bullet
Y\bullet \left( I\circ \left( r_{Z^{\revsn }\bullet T^{\revsn }}^{\bullet
}\right) ^{-1}\right) \right) \left( X\bullet Y\bullet \left( l_{Z^{\revsn
}\bullet T^{\revsn }}^{\circ }\right) ^{-1}\right)  \\
&=&\left( \left( \zeta _{X,Y,T,Z}\bullet I\right) \circ J\right) \delta
_{X\bullet Y,I,T\bullet Z,J}\left( X\bullet Y\bullet \left( I\circ \left(
r_{\left( T\bullet Z\right) ^{\revsn }}^{\bullet }\right) ^{-1}\right) \right)
\left( X\bullet Y\bullet \left( l_{\left( T\bullet Z\right) ^{\revsn }}^{\circ
}\right) ^{-1}\right) \left( X\bullet Y\bullet \phi _{Z,T}^{\bullet }\right)
\\
&=&\left( \left( \zeta _{X,Y,T,Z}\bullet I\right) \circ J\right) \varphi
_{X\bullet Y,T\bullet Z}\left( X\bullet Y\bullet \phi _{Z,T}^{\bullet
}\right)
\end{eqnarray*}
\end{invisible}
We now compute
\begin{align*}
\vartheta _{X,I,X}\left( X^{\revsn }\bullet \left( J\circ \Delta _{I}^{\bullet
}\right) \bullet X\right)  &=\gamma _{J,X,I\bullet I,X}\left( \left( \left(
l_{X^{\revsn }}^{\bullet }\right) ^{-1}\circ \left( I\bullet I\right) \right)
\bullet X\right) \left( \underbracket[0.140ex]{\xi _{X^{\revsn },I}}\bullet X\right)
\left( \underbracket[0.140ex]{X^{\revsn }\bullet \left( J\circ \Delta _{I}^{\bullet
}\right) }\bullet X\right)  \\
&\overset{\mathclap{\eqref{form:xidelta}}}{=}\gamma _{J,X,I\bullet I,X}\left(
\left( \left( l_{X^{\revsn }}^{\bullet }\right) ^{-1}\circ \left( I\bullet
I\right) \right) \bullet X\right) \left( \left( X^{\revsn }\circ \Delta
_{I}^{\bullet }\right) \left( r_{X^{\revsn }}^{\bullet }\circ I\right)
\overline{\varkappa }_{X^{\revsn },J,I}^{I}\bullet X\right)  \\
&=\left( J\circ \left( \Delta _{I}^{\bullet }\bullet \left( X\circ X\right)
\right) \right) \gamma _{J,X,I,X}\left( \left( \left( l_{X^{\revsn }}^{\bullet
}\right) ^{-1}\circ I\right) \bullet X\right) \left( \left( r_{X^{\revsn
}}^{\bullet }\circ I\right) \overline{\varkappa }_{X^{\revsn },J,I}^{I}\bullet
X\right)  \\
&=\left( J\circ \left( \Delta _{I}^{\bullet }\bullet \left( X\circ X\right)
\right) \right) \psi _{X,X}\left( r_{X^{\revsn }}^{\circ }\bullet X\right)
\left( \left( r_{X^{\revsn }}^{\bullet }\circ I\right) \left( \underbracket[0.140ex]{%
\varkappa _{X^{\revsn },J,I}^{I}}\right) ^{-1}\bullet X\right)  \\
&=\left( J\circ \left( \Delta _{I}^{\bullet }\bullet \left( X\circ X\right)
\right) \right) \psi _{X,X}\left( r_{X^{\revsn }}^{\bullet }r_{X^{\revsn
}\bullet J}^{\circ }\left( \left( X^{\revsn }\bullet \left( r_{J}^{\circ
}\right) ^{-1}\right) r_{X^{\revsn }\bullet J}^{\circ }\right) ^{-1}\bullet
X\right)  \\
&=\left( J\circ \left( \Delta _{I}^{\bullet }\bullet \left( X\circ X\right)
\right) \right) \psi _{X,X}\left( r_{X^{\revsn }}^{\bullet } \left(
X^{\revsn }\bullet r_{J}^{\circ }\right) \bullet X\right).
\end{align*}
Similarly, one proves \eqref{form:vartheta1I} by using \eqref{form:xi1delta}.
\begin{invisible}
\begin{eqnarray*}
\vartheta _{X,I,X}^{\prime }\left( X\bullet \left( \Delta _{I}^{\bullet
}\circ J\right) \bullet X^{\revsn }\right)  &=&\delta _{X,I\bullet
I,X,J}\left( X\bullet \left( \left( I\bullet I\right) \circ \left(
r_{X^{\revsn }}^{\bullet }\right) ^{-1}\right) \right) \left( X\bullet
\underbracket[0.140ex]{\xi _{I,X^{\revsn }}^{\prime }}\right) \left( X\bullet \underbracket[0.140ex]{%
\left( \Delta _{I}^{\bullet }\circ J\right) \bullet X^{\revsn }}\right)  \\
&\overset{\eqref{form:xi1delta} }{=}&\delta _{X,I\bullet I,X,J}\left(
X\bullet \left( \left( I\bullet I\right) \circ \left( r_{X^{\revsn }}^{\bullet
}\right) ^{-1}\right) \right) \left( X\bullet \left( \Delta _{I}^{\bullet
}\circ X^{\revsn }\right) \left( I\circ l_{X^{\revsn }}^{\bullet }\right)
\overline{\varsigma }_{I,J,X^{\revsn }}^{I}\right)  \\
&=&\left( \left( \left( X\circ X\right) \bullet \Delta _{I}^{\bullet
}\right) \circ J\right) \underbracket[0.140ex]{\delta _{X,I,X,J}\left( X\bullet \left(
I\circ \left( r_{X^{\revsn }}^{\bullet }\right) ^{-1}\right) \right) }\left(
X\bullet \left( I\circ l_{X^{\revsn }}^{\bullet }\right) \overline{\varsigma }%
_{I,J,X^{\revsn }}^{I}\right)  \\
&=&\left( \left( \left( X\circ X\right) \bullet \Delta _{I}^{\bullet
}\right) \circ J\right) \varphi _{X,X}\left( X\bullet l_{X^{\revsn }}^{\circ
}\left( I\circ l_{X^{\revsn }}^{\bullet }\right) \overline{\varsigma }%
_{I,J,X^{\revsn }}^{I}\right)  \\
&=&\left( \left( \left( X\circ X\right) \bullet \Delta _{I}^{\bullet
}\right) \circ J\right) \varphi _{X,X}\left( X\bullet l_{X^{\revsn }}^{\bullet
}l_{J\bullet X^{\revsn }}^{\circ }\left( \underbracket[0.140ex]{\varsigma _{I,J,X^{\revsn
}}^{I}}\right) ^{-1}\right)  \\
&=&\left( \left( \left( X\circ X\right) \bullet \Delta _{I}^{\bullet
}\right) \circ J\right) \varphi _{X,X}\left( X\bullet l_{X^{\revsn }}^{\bullet
}l_{J\bullet X^{\revsn }}^{\circ }\left( \left( \left( l_{J}^{\circ }\right)
^{-1}\bullet X^{\revsn }\right) l_{J\bullet X^{\revsn }}^{\circ }\right)
^{-1}\right)  \\
&=&\left( \left( \left( X\circ X\right) \bullet \Delta _{I}^{\bullet
}\right) \circ J\right) \varphi _{X,X}\left( X\bullet l_{X^{\revsn }}^{\bullet
} \left( l_{J}^{\circ }\bullet X^{\revsn }\right) \right)
\end{eqnarray*}%
\end{invisible}
\end{proof}

We note in passing that \cref{lem:Bohm4.3} allows us to obtain the following result, which is an analogue of %\cite[Equation (7.3)]{Bohm-Lack}.
\cite[Theorem 7.5]{Bohm-Lack}. We state it for sake of completeness, to show that our assumptions are enough to obtain these expected properties of the antipode, even though we do not make use of them.

\begin{lemma}%[\rd{added:2026/03/06}]
Let $H$ be a Hopf monoid in a duoidal category with a reversion.
Then, the following equalities hold:%
\begin{align}
\left( m_{H}^{\circ }\right) ^{\revsn }\phi _{H,H}^{\circ }\left( \Sigma
_{H}\circ \Sigma _{H}\right) &=\Sigma _{H}m_{H}^{\circ }.
\label{form:antipode-mult}\\
\phi
_{H,H}^{\bullet }\left( \Sigma _{H}\bullet \Sigma _{H}\right) \Delta
_{H}^{\bullet }&=\left( \Delta _{H}^{\bullet }\right) ^{\revsn }\Sigma _{H}.  \label{form:antipcomult}\\
\Sigma_{H}u^\circ_H &=(u^\circ_H)^\revsn\phi^\circ_0\label{form:antip-unit}
\\
(\varepsilon^\bullet_H)^\revsn \Sigma_H &=\phi^\bullet_0 \varepsilon^\bullet_H.
\label{form:antip-counit}
\end{align}
\end{lemma}

\begin{proof}
Since \cref{lem:convcomp} entails that $\Sigma _{H}m_{H}^{\circ }$ is the
convolution inverse of the morphism of $\bullet$-comonoids $m_{H}^{\circ },$ by
uniqueness of convolution inverse, it suffices to prove that $\left(
m_{H}^{\circ }\right) ^{\revsn }\phi _{H,H}^{\circ }\left( \Sigma _{H}\circ
\Sigma _{H}\right) \star m_{H}^{\circ }=1_{H\circ H,H}^{r}$. Since we have \eqref{form:BohmLem4.2}, the equation \eqref{form:antipode-mult} is obtained following the
lines of \cite[Proposition 6.5]{Bohm-Lack}.

\begin{invisible}
By means of this equality one gets%
\begin{eqnarray*}
&&\left( m_{A}^{\circ }\right) ^{\revsn }\phi _{A,A}^{\circ }\left( \Sigma
_{A}\circ \Sigma _{A}\right) \star m_{A}^{\circ } \\
&=&\left( J\circ \left( I\bullet m_{A}^{\circ }\right) \right) \psi
_{A,A}\left( \left( m_{A}^{\circ }\right) ^{\revsn }\phi _{A,A}^{\circ }\left(
\Sigma _{A}\circ \Sigma _{A}\right) \bullet m_{A}^{\circ }\right) \Delta
_{A\circ A}^{\bullet } \\
&=&\left( J\circ \left( I\bullet m_{A}^{\circ }\right) \right) \psi
_{A,A}\left( \left( m_{A}^{\circ }\right) ^{\revsn }\bullet m_{A}^{\circ
}\right) \left( \phi _{A,A}^{\circ }\left( \Sigma _{A}\circ \Sigma
_{A}\right) \bullet \left( A\circ A\right) \right) \Delta _{A\circ
A}^{\bullet } \\
&=&\left( J\circ \left( I\bullet m_{A}^{\circ }\left( m_{A}^{\circ }\circ
m_{A}^{\circ }\right) \right) \right) \psi _{A\circ A,A\circ A}\left( \phi
_{A,A}^{\circ }\bullet \left( A\circ A\right) \right) \left( \left( \Sigma
_{A}\circ \Sigma _{A}\right) \bullet \left( A\circ A\right) \right) \zeta
_{A,A,A,A}\left( \Delta _{A}^{\bullet }\circ \Delta _{A}^{\bullet }\right)
\\
&=&\left( J\circ \left( I\bullet m_{A}^{\circ }\left( A\circ m_{A}^{\circ
}\right) \left( A\circ m_{A}^{\circ }\circ A\right) \right) \right)
\underbracket[0.140ex]{\psi _{A\circ A,A\circ A}\left( \phi _{A,A}^{\circ }\bullet
\left( A\circ A\right) \right) \zeta _{A^{\revsn },A,A^{\revsn },A}}\left(
\left( \Sigma _{A}\bullet A\right) \circ \left( \Sigma _{A}\bullet A\right)
\right) \left( \Delta _{A}^{\bullet }\circ \Delta _{A}^{\bullet }\right)  \\
&\overset{\eqref{form:BohmLem4.2}}{=}&\left[
\begin{array}{c}
\left( J\circ \left( I\bullet m_{A}^{\circ }\left( A\circ m_{A}^{\circ
}\right) \left( A\circ m_{A}^{\circ }\circ A\right) \right) \right) \left(
m_{J}^{\circ }\circ \left( I\bullet \left( A\circ A\circ A\circ A\right)
\right) \right)  \\
\left( J\circ \psi _{A,A\circ A\circ A}\right) \left( J\circ \left(
l_{A^{\revsn }}^{\circ }\bullet \left( A\circ A\circ A\right) \right) \right)
\left( J\circ \zeta _{I,A\circ A,A^{\revsn },A}\right) \left( \psi _{A,A}\circ
\left( A^{\revsn }\bullet A\right) \right)  \\
\left( \left( \Sigma _{A}\bullet A\right) \circ \left( \Sigma _{A}\bullet
A\right) \right) \left( \Delta _{A}^{\bullet }\circ \Delta _{A}^{\bullet
}\right)
\end{array}%
\right]  \\
&=&\left[
\begin{array}{c}
\left( J\circ \left( I\bullet m_{A}^{\circ }\left( A\circ m_{A}^{\circ
}\right) \right) \right) \left( m_{J}^{\circ }\circ \left( I\bullet \left(
A\circ A\circ A\right) \right) \right)  \\
\left( J\circ \psi _{A,A\circ A}\right) \left( J\circ \left( l_{A^{\revsn
}}^{\circ }\bullet \left( A\circ A\right) \right) \right) \left( J\circ
\zeta _{I,A,A^{\revsn },A}\right)  \\
\left( \left( J\circ \left( I\bullet m_{A}^{\circ }\right) \right) \psi
_{A,A}\left( \Sigma _{A}\bullet A\right) \Delta _{A}^{\bullet }\circ \left(
\Sigma _{A}\bullet A\right) \Delta _{A}^{\bullet }\right)
\end{array}%
\right]  \\
&=&\left[
\begin{array}{c}
\left( J\circ \left( I\bullet m_{A}^{\circ }\left( A\circ m_{A}^{\circ
}\right) \right) \right) \left( m_{J}^{\circ }\circ \left( I\bullet \left(
A\circ A\circ A\right) \right) \right)  \\
\left( J\circ \psi _{A,A\circ A}\right) \left( J\circ \left( l_{A^{\revsn
}}^{\circ }\bullet \left( A\circ A\right) \right) \right) \left( J\circ
\zeta _{I,A,A^{\revsn },A}\right)  \\
\left( \left( \Sigma _{A}\star \mathrm{Id}_{A}\right) \circ \left( \Sigma
_{A}\bullet A\right) \Delta _{A}^{\bullet }\right)
\end{array}%
\right]  \\
&=&\left[
\begin{array}{c}
\left( J\circ \left( I\bullet m_{A}^{\circ }\left( A\circ m_{A}^{\circ
}\right) \right) \right) \left( m_{J}^{\circ }\circ \left( I\bullet \left(
A\circ A\circ A\right) \right) \right)  \\
\left( J\circ \psi _{A,A\circ A}\right) \left( J\circ \left( l_{A^{\revsn
}}^{\circ }\bullet \left( A\circ A\right) \right) \right) \left( J\circ
\zeta _{I,A,A^{\revsn },A}\right)  \\
\left( \left( J\circ \left( I\bullet u_{A}^{\circ }\right) \right) \left(
J\circ \Delta _{I}^{\bullet }\right) \left( r_{J}^{\circ }\right)
^{-1}\varepsilon _{A}^{\bullet }\circ \left( \Sigma _{A}\bullet A\right)
\Delta _{A}^{\bullet }\right)
\end{array}%
\right]  \\
&=&\left[
\begin{array}{c}
\left( J\circ \left( I\bullet m_{A}^{\circ }\left( A\circ m_{A}^{\circ
}\right) \right) \right) \left( m_{J}^{\circ }\circ \left( I\bullet \left(
A\circ u_{A}^{\circ }\circ A\right) \right) \right)  \\
\left( J\circ \psi _{A,I\circ A}\right) \left( J\circ \left( l_{A^{\revsn
}}^{\circ }\bullet \left( I\circ A\right) \right) \right) \left( J\circ
\zeta _{I,I,A^{\revsn },A}\left( \Delta _{I}^{\bullet }\circ \left( A^{\revsn
}\bullet A\right) \right) \right)  \\
\left( \left( r_{J}^{\circ }\right) ^{-1}\varepsilon _{A}^{\bullet }\circ
\left( \Sigma _{A}\bullet A\right) \Delta _{A}^{\bullet }\right)
\end{array}%
\right]  \\
&=&\left[
\begin{array}{c}
\left( J\circ \left( I\bullet m_{A}^{\circ }\right) \right) \left(
m_{J}^{\circ }\circ \left( I\bullet \left( A\circ l_{A}^{\circ }\right)
\right) \right) \left( J\circ \psi _{A,I\circ A}\right)  \\
\left( J\circ \left( l_{A^{\revsn }}^{\circ }\bullet \left( I\circ A\right)
\right) \right) \left( J\circ \zeta _{I,I,A^{\revsn },A}\left( \Delta
_{I}^{\bullet }\circ \left( A^{\revsn }\bullet A\right) \right) \right) \left(
\left( r_{J}^{\circ }\right) ^{-1}\circ \left( A^{\revsn }\bullet A\right)
\right) \left( \varepsilon _{A}^{\bullet }\circ \left( \Sigma _{A}\bullet
A\right) \Delta _{A}^{\bullet }\right)
\end{array}%
\right]  \\
&=&\left[
\begin{array}{c}
\left( m_{J}^{\circ }\circ \left( I\bullet A\right) \right) \left( J\circ
J\circ \left( I\bullet m_{A}^{\circ }\right) \right) \left( J\circ \psi
_{A,A}\right)  \\
\left( J\circ \underbracket[0.140ex]{\left( l_{A^{\revsn }}^{\circ }\bullet l_{A}^{\circ
}\right) \zeta _{I,I,A^{\revsn },A}\left( \Delta _{I}^{\bullet }\circ \left(
A^{\revsn }\bullet A\right) \right) \left( l_{A^{\revsn }\bullet A}^{\circ
}\right) ^{-1}}\right) \left( \varepsilon _{A}^{\bullet }\circ \left( \Sigma
_{A}\bullet A\right) \Delta _{A}^{\bullet }\right)
\end{array}%
\right]  \\
&\overset{\eqref{eq:unit1}}{=}&\left( m_{J}^{\circ }\circ \left( I\bullet A\right) \right) \left( J\circ
J\circ \left( I\bullet m_{A}^{\circ }\right) \right) \left( J\circ \psi
_{A,A}\right) \left( \varepsilon _{A}^{\bullet }\circ \left( \Sigma
_{A}\bullet A\right) \Delta _{A}^{\bullet }\right)  \\
&=&\left( m_{J}^{\circ }\circ \left( I\bullet A\right) \right) \left( J\circ
\left( \Sigma _{A}\star \mathrm{Id}_{A}\right) \right)  \\
&=&\left( m_{J}^{\circ }\circ \left( I\bullet A\right) \right) \left(
\varepsilon _{A}^{\bullet }\circ \left( J\circ \left( I\bullet u_{A}^{\circ
}\right) \right) \left( J\circ \Delta _{I}^{\bullet }\right) \left(
r_{J}^{\circ }\right) ^{-1}\varepsilon _{A}^{\bullet }\right)  \\
&=&\left( m_{J}^{\circ }\circ \left( I\bullet A\right) \right) \left( J\circ
J\circ \left( I\bullet u_{A}^{\circ }\right) \right) \left( J\circ J\circ
\Delta _{I}^{\bullet }\right) \left( J\circ \left( r_{J}^{\circ }\right)
^{-1}\right) \left( \varepsilon _{A}^{\bullet }\circ \varepsilon
_{A}^{\bullet }\right)  \\
&=&\left( J\circ \left( I\bullet u_{A}^{\circ }\right) \right) \left( J\circ
\Delta _{I}^{\bullet }\right) \left( m_{J}^{\circ }\circ I\right) \left(
r_{J\circ J}^{\circ }\right) ^{-1}\left( \varepsilon _{A}^{\bullet }\circ
\varepsilon _{A}^{\bullet }\right)  \\
&=&\left( J\circ \left( I\bullet u_{A}^{\circ }\right) \right) \left( J\circ
\Delta _{I}^{\bullet }\right) \left( r_{J}^{\circ }\right) ^{-1}m_{J}^{\circ
}\left( \varepsilon _{A}^{\bullet }\circ \varepsilon _{A}^{\bullet }\right)
\\
&=&\left( J\circ \left( I\bullet u_{A}^{\circ }\right) \right) \left( J\circ
\Delta _{I}^{\bullet }\right) \left( r_{J}^{\circ }\right) ^{-1}\varepsilon
_{A\circ A}^{\bullet }=1_{A\circ A,A}^{r}.
\end{eqnarray*}%
It remains to prove the commutativity of the diagram above.
\end{invisible}

Since \cref{lem:convcomp} entails that $(\Delta_{H}^{\bullet})^{\revsn}\Sigma_{H}$ is the convolution inverse of the morphism of $\circ$-monoids $\Delta_{H}^{\bullet},$ by
uniqueness of convolution inverse, it suffices to prove $\left( \phi
_{H,H}^{\bullet }\left( \Sigma _{H}\bullet \Sigma _{H}\right) \Delta
_{H}^{\bullet }\right) \star \Delta _{H}^{\bullet }=1_{H,H\bullet H}^{r}$. Since we have \eqref{form:vartheta} and \eqref{form:varthetaI}, the equation \eqref{form:antipcomult} is obtained following the lines of \cite[Proposition 6.6]{Bohm-Lack}.

\begin{invisible}
By \cref{lem:convcomp}, we know that $\left( \Delta _{H}^{\bullet }\right)
^{\revsn }\Sigma _{H}$ is the convolution inverse of the morphism of $\circ $%
-monoids $\Delta _{H}^{\bullet }:H\rightarrow H\bullet H.$ By uniqueness of
the convolution inverse it suffices to check that $\left( \phi
_{H,H}^{\bullet }\left( \Sigma _{H}\bullet \Sigma _{H}\right) \Delta
_{H}^{\bullet }\right) \star \Delta _{H}^{\bullet }=1_{H,H\bullet H}^{r}.$
Indeed, we have
\begin{eqnarray*}
&&\left( \phi _{H,H}^{\bullet }\left( \Sigma _{H}\bullet \Sigma _{H}\right)
\Delta _{H}^{\bullet }\right) \star \Delta _{H}^{\bullet } \\
&=&\left( J\circ \left( I\bullet m_{H\bullet H}^{\circ }\right) \right) \psi
_{H\bullet H,H\bullet H}\left( \phi _{H,H}^{\bullet }\left( \Sigma
_{H}\bullet \Sigma _{H}\right) \Delta _{H}^{\bullet }\bullet \Delta
_{H}^{\bullet }\right) \Delta _{H}^{\bullet } \\
&=&\left( J\circ \left( I\bullet m_{H}^{\circ }\bullet m_{H}^{\circ }\right)
\right) \underbracket[0.140ex]{\left( J\circ \left( I\bullet \zeta _{H,H,H,H}\right)
\right) \psi _{H\bullet H,H\bullet H}\left( \phi _{H,H}^{\bullet }\bullet
H\bullet H\right) }\left( H^{\revsn }\bullet \left( \Sigma _{H}\bullet
H\right) \Delta _{H}^{\bullet }\bullet H\right) \left( \Sigma _{H}\bullet
\Delta _{H}^{\bullet }\right) \Delta _{H}^{\bullet } \\
&\overset{\eqref{form:vartheta}}{=}&\left( J\circ \left(
I\bullet m_{H}^{\circ }\bullet m_{H}^{\circ }\right) \right) \vartheta
_{H,H\circ H,H}\left( H^{\revsn }\bullet \psi _{H,H}\bullet H\right) \left(
H^{\revsn }\bullet \left( \Sigma _{H}\bullet H\right) \Delta _{H}^{\bullet
}\bullet H\right) \left( \Sigma _{H}\bullet \Delta _{H}^{\bullet }\right)
\Delta _{H}^{\bullet } \\
&=&\left( J\circ \left( I\bullet H\bullet m_{H}^{\circ }\right) \right)
\vartheta _{H,H,H}\left( H^{\revsn }\bullet \left( J\circ \left( I\bullet
m_{H}^{\circ }\right) \right) \bullet H\right) \left( H^{\revsn }\bullet \psi
_{H,H}\bullet H\right) \left( H^{\revsn }\bullet \left( \Sigma _{H}\bullet
H\right) \Delta _{H}^{\bullet }\bullet H\right) \left( \Sigma _{H}\bullet
\Delta _{H}^{\bullet }\right) \Delta _{H}^{\bullet } \\
&=&\left( J\circ \left( I\bullet H\bullet m_{H}^{\circ }\right) \right)
\vartheta _{H,H,H}\left( H^{\revsn }\bullet \left( \Sigma _{H}\star \mathrm{Id}%
_{H}\right) \bullet H\right) \left( \Sigma _{H}\bullet \Delta _{H}^{\bullet
}\right) \Delta _{H}^{\bullet } \\
&=&\left( J\circ \left( I\bullet H\bullet m_{H}^{\circ }\right) \right)
\vartheta _{H,H,H}\left( H^{\revsn }\bullet \left( J\circ \left( I\bullet
u_{H}^{\circ }\right) \right) \bullet H\right) \left( H^{\revsn }\bullet
\left( J\circ \Delta _{I}^{\bullet }\right) \left( r_{J}^{\circ }\right)
^{-1}\bullet H\right) \left( H^{\revsn }\bullet \varepsilon _{H}^{\bullet
}\bullet H\right) \left( \Sigma _{H}\bullet \Delta _{H}^{\bullet }\right)
\Delta _{H}^{\bullet } \\
&=&\left( J\circ \left( I\bullet u_{H}^{\circ }\bullet m_{H}^{\circ }\right)
\right) \underbracket[0.140ex]{\vartheta _{H,I,H}\left( H^{\revsn }\bullet \left( J\circ
\Delta _{I}^{\bullet }\right) \bullet H\right) }\left( H^{\revsn }\bullet
\left( r_{J}^{\circ }\right) ^{-1}\bullet H\right) \left( H^{\revsn }\bullet
\left( l_{H}^{\bullet }\right) ^{-1}\right) \left( \Sigma _{H}\bullet
H\right) \Delta _{H}^{\bullet } \\
&\overset{\eqref{form:varthetaI}}{=}&\left( J\circ \left( I\bullet
u_{H}^{\circ }\bullet m_{H}^{\circ }\right) \right) \left( J\circ \left(
\Delta _{I}^{\bullet }\bullet \left( H\circ H\right) \right) \right) \psi
_{H,H}\left( r_{H^{\revsn }}^{\bullet }\left( \left( H^{\revsn }\bullet
r_{J}^{\circ }\right) \right) \bullet H\right) \left( H^{\revsn }\bullet
\left( r_{J}^{\circ }\right) ^{-1}\bullet H\right) \left( (r_{H^{\revsn
}}^{\bullet })^{-1}\bullet H\right) \left( \Sigma _{H}\bullet H\right) \Delta
_{H}^{\bullet } \\
&=&\left( J\circ \left( \left( I\bullet u_{H}^{\circ }\right) \Delta
_{I}^{\bullet }\bullet H\right) \right) \left( J\circ \left( I\bullet
m_{H}^{\circ }\right) \right) \psi _{H,H}\left( \Sigma _{H}\bullet H\right)
\Delta _{H}^{\bullet } \\
&=&\left( J\circ \left( \left( I\bullet u_{H}^{\circ }\right) \Delta
_{I}^{\bullet }\bullet H\right) \right) \left( \Sigma _{H}\star \mathrm{Id}%
_{H}\right)  \\
&=&\left( J\circ \left( \left( I\bullet u_{H}^{\circ }\right) \Delta
_{I}^{\bullet }\bullet H\right) \right) \left( J\circ \left( I\bullet
u_{H}^{\circ }\right) \Delta _{I}^{\bullet }\right) \left( r_{J}^{\circ
}\right) ^{-1}\varepsilon _{H}^{\bullet } \\
&=&\left( J\circ \left( I\bullet \left( u_{H}^{\circ }\bullet u_{H}^{\circ
}\right) \Delta _{I}^{\bullet }\right) \right) \left( J\circ \Delta
_{I}^{\bullet }\right) \left( r_{J}^{\circ }\right) ^{-1}\varepsilon
_{H}^{\bullet } \\
&=&\left( J\circ \left( I\bullet u_{H\bullet H}^{\circ }\right) \right)
\left( J\circ \Delta _{I}^{\bullet }\right) \left( r_{J}^{\circ }\right)
^{-1}\varepsilon _{H}^{\bullet }=1_{H,H\bullet H}^{r}.
\end{eqnarray*}
\end{invisible}

The equations \eqref{form:antip-unit} and \eqref{form:antip-counit} follow from \cref{coro:convcomp}, as $I$, $J$ are in $\mathsf{Hopf}(\Cc,\circ,\bullet)$ by \cref{lem:IHopf} and $u^\circ_H$, $\varepsilon^\bullet_H$ are morphisms in $\mathsf{Bimon}(\Cc,\circ,\bullet)$.
\end{proof}

\subsection{Galois and co‑Galois induced morphisms}

In our discussion, the following morphisms induced by the Galois and co-Galois maps will play a central role.

\begin{definition}
Let $(\Cc,\circ, I,\bullet,J)$ be a duoidal category. Consider an object $H$ in $\Bimon (\Cc,\circ,\bullet) $.\\
Given a
morphism $g:C\rightarrow H$  in $\Comon(\Cc^\bullet)$, we define the following morphism in $\Cc$:
\begin{eqnarray*}
g^{r} &:=&\left( \left( C\circ \varepsilon _{H}^{\bullet }\right) \bullet
m_{H}^{\circ }\left( g\circ H\right) \right) \Delta _{C\circ H}^{\bullet
}:C\circ H\rightarrow \left( C\circ J\right) \bullet H.
%\\g^{X}&:=&\rd{\varsigma^H_{C,X,H}=} \left( \left( C\circ X\right) \bullet m_{H}^{\circ }\right) \zeta
%_{C,H,X,H}\left( \left( C\bullet g\right) \Delta _{C}^{\bullet }\circ \left(
%X\bullet H\right) \right) :C\circ \left( X\bullet H\right) \rightarrow
%\left( C\circ X\right) \bullet H.
\end{eqnarray*}
Given a morphism  $g:H\rightarrow A$ in $\Mon(\Cc^\circ)$, we define the following morphism in $\Cc$:
\begin{eqnarray*}
g_{r}&:=&m_{A\bullet H}^{\circ }\left( \left( A\bullet u_{H}^{\circ }\right)
\circ \left( g\bullet H\right) \Delta _{H}^{\bullet }\right) :\left(
A\bullet I\right) \circ H\rightarrow A\bullet H.
%\\g_{X}&:=&\rd{\varkappa^H_{A,X,H}=}\left( m_{A}^{\circ }\left( A\circ g\right) \bullet \left( X\circ
%H\right) \right) \zeta _{A,X,H,H}\left( \left( A\bullet X\right) \circ
%\Delta _{H}^{\bullet }\right) :\left( A\bullet X\right) \circ H\rightarrow
%A\bullet \left( X\circ H\right)
\end{eqnarray*}
\end{definition}

Of course, one could consider the left-hand versions of the morphisms defined above. For example, given a
morphism $g:C\rightarrow H$  in $\Comon(\Cc^\bullet)$, one can define
\begin{equation*}\label{formgr}
g^{l} :=\left( m_{H}^{\circ }\left( H\circ g\right) \bullet \left(
\varepsilon _{H}^{\bullet }\circ C\right) \right) \Delta _{H\circ
C}^{\bullet }:H\circ C\rightarrow H\bullet \left( J\circ C\right).
\end{equation*}
%and, given a morphism  $g:H\rightarrow A$ in $\Mon(\Cc^\circ)$, one can define
%\begin{equation*}
%g_{l}:=m^{\circ}_{H\bullet A}((H\bullet g)\Delta_{H}^{\bullet}\circ(u_{H}^{\circ}\bullet A)):H\circ(I\bullet A)\to H\bullet A
%\end{equation*}
In particular,  $(\id_H)^l:H\circ H\to H\bullet(J\circ H)$ is the \emph{left fusion morphism} defined in \cite[\S 4.7]{Street}, whose invertibility defines what Street calls \emph{left Hopf monoids}. Thus, we can adopt a similar terminology and call $(\id_H)^r$ and $(\id_H)_r$ the \emph{fusion morphisms}. As we will see, these morphisms will play a role in relevant cases, see \cref{pro:Hopfmonbraidedcat} and  \cref{pro:Bohm}.

\begin{remark}
\label{rmk:Bohm}
    Note that, by using the co-Galois maps, we can write%
    \begin{eqnarray*}
g^{r} &=&\left( \left( C\circ \varepsilon _{H}^{\bullet }\right) \bullet
m_{H}^{\circ }\left( g\circ H\right) \right) \Delta _{C\circ H}^{\bullet } \\
&=&\left( \left( C\circ \varepsilon _{H}^{\bullet }\right) \bullet
m_{H}^{\circ }\left( g\circ H\right) \right) \zeta _{C,C,H,H}\left( \Delta
_{C}^{\bullet }\circ \Delta _{H}^{\bullet }\right) \nonumber \\
&=&\left( \left( C\circ J\right) \bullet m_{H}^{\circ }\right) \left( \left(
C\circ \varepsilon _{H}^{\bullet }\right) \bullet \left( g\circ H\right)
\right) \zeta _{C,C,H,H}\left( \Delta _{C}^{\bullet }\circ \Delta
_{H}^{\bullet }\right) \nonumber \\
&=&\left( \left( C\circ J\right) \bullet m_{H}^{\circ }\right) \zeta
_{C,H,J,H}\left( \left( C\bullet g\right) \circ \left( \varepsilon
_{H}^{\bullet }\bullet H\right) \right) \left( \Delta _{C}^{\bullet }\circ
\Delta _{H}^{\bullet }\right) \nonumber \\
&=&\left( \left( C\circ J\right) \bullet m_{H}^{\circ }\right) \zeta
_{C,H,J,H}\left( \left( C\bullet g\right) \Delta _{C}^{\bullet }\circ \left(
l_{H}^{\bullet }\right) ^{-1}\right) \nonumber=\varsigma^H _{C,J,H}(C\circ \left(
l_{H}^{\bullet }\right) ^{-1})
\end{eqnarray*}%
so that
\begin{equation}
 g^{r}=\varsigma^H _{C,J,H}(C\circ \left(
l_{H}^{\bullet }\right) ^{-1}),
\label{equivalentgr}
\end{equation}
where we are taking  $\rho^\bullet=\left( C\bullet g\right) \Delta^\bullet _{C}:C\to C\bullet H.$
\begin{invisible}
Similarly, one gets
\[g^{l}=\left( m_{H}^{\circ }\bullet \left( J\circ C\right) \right) \zeta
_{H,J,H,C}\left( \left( r_{H}^{\bullet }\right) ^{-1}\circ \left( \left(
g\bullet C\right) \Delta _{C}^{\bullet }\right) \right)=\varkappa _{H,J,C}^{H}(\left( r_{H}^{\bullet }\right) ^{-1}\circ C) .\]
In the following computation $A$ should be $H$
\begin{eqnarray*}
g^{l} &=&\left( m_{A}^{\circ }\left( A\circ g\right) \bullet \left(
\varepsilon _{A}^{\bullet }\circ C\right) \right) \Delta _{A\circ
C}^{\bullet } \\
&=&\left( m_{A}^{\circ }\left( A\circ g\right) \bullet \left( \varepsilon
_{A}^{\bullet }\circ C\right) \right) \zeta _{A,A,C,C}\left( \Delta
_{A}^{\bullet }\circ \Delta _{C}^{\bullet }\right)  \\
&=&\left( m_{A}^{\circ }\bullet \left( J\circ C\right) \right) \left( \left(
A\circ g\right) \bullet \left( \varepsilon _{A}^{\bullet }\circ C\right)
\right) \zeta _{A,A,C,C}\left( \Delta _{A}^{\bullet }\circ \Delta
_{C}^{\bullet }\right)  \\
&=&\left( m_{A}^{\circ }\bullet \left( J\circ C\right) \right) \zeta
_{A,J,A,C}\left( \left( A\bullet \varepsilon _{A}^{\bullet }\right) \circ
\left( g\bullet C\right) \right) \left( \Delta _{A}^{\bullet }\circ \Delta
_{C}^{\bullet }\right)  \\
&=&\left( m_{A}^{\circ }\bullet \left( J\circ C\right) \right) \zeta
_{A,J,A,C}\left( \left( r_{A}^{\bullet }\right) ^{-1}\circ \left( \left(
g\bullet C\right) \Delta _{C}^{\bullet }\right) \right) .
\end{eqnarray*}%
\end{invisible}
%\rd{[Added: 2025/11/13]}
Dually, by using the Galois maps, one proves that
\begin{equation}
 g_{r} := (A\bullet l_H^\circ )\varkappa _{A,I,H}^{H},\label{equivalentformgr}
\end{equation} where we are taking $\mu_A^\circ=m_A^\circ (A\circ g).$
\begin{invisible}
We have
 \begin{eqnarray*}
g_{r} &:=&m_{A\bullet H}^{\circ }\left( \left( A\bullet u_{H}^{\circ
}\right) \circ \left( g\bullet H\right) \Delta _{H}^{\bullet }\right)  \\
&=&\left( m_{A}^{\circ }\bullet m_{H}^{\circ }\right) \zeta _{A,H,A,H}\left(
\left( A\bullet u_{H}^{\circ }\right) \circ \left( g\bullet H\right) \right)
\left( \left( A\bullet I\right) \circ \Delta _{H}^{\bullet }\right) \nonumber \\
&=&\left( m_{A}^{\circ }\bullet m_{H}^{\circ }\right) \left( \left( A\circ
g\right) \bullet \left( u_{H}^{\circ }\circ H\right) \right) \zeta
_{A,I,H,H}\left( \left( A\bullet I\right) \circ \Delta _{H}^{\bullet
}\right) \nonumber \\
&=&\left( m_{A}^{\circ }\left( A\circ g\right) \bullet l_{H}^{\circ }\right)
\zeta _{A,I,H,H}\left( \left( A\bullet I\right) \circ \Delta _{H}^{\bullet
}\right)\\&=& (A\bullet l^\circ_H)\varkappa^H_{A,I,H}. \nonumber
\end{eqnarray*}%
\end{invisible}
As a consequence, $g_r$ and $g^r$ are invertible if so are the Galois and co-Galois maps.
%This means that $g_r$ is $\beta_A$ as in \cite[(3.2))]{BCZ} in the particular case $\gamma=m_{A}^{\circ }\left( A\circ g\right):A\circ H\to A$.

%Note that requiring the invertibility of the natural transformation $\beta$ was presented in \cite[Remark 3.12]{BCZ} as a possible way to define Hopf monoids in a duoidal category.
\begin{invisible}
Similarly, one gets
\[
\begin{split}
g_{l}&=m^{\circ}_{H\bullet A}((H\bullet g)\Delta_{H}^{\bullet}\circ(u_{H}^{\circ}\bullet A))\\&=(m^{\circ}_{H}\bullet m^{\circ}_{A})\zeta_{H,A,H,A}((H\bullet g)\circ(u_{H}^{\circ}\bullet A))(\Delta_{H}^{\bullet}\circ(I\bullet A)\\&=(m^{\circ}_{H}\bullet m^{\circ}_{A})((H\circ u_{H}^{\circ})\bullet(g\circ A))\zeta_{H,H,I,A}(\Delta_{H}^{\bullet}\circ(I\bullet A))\\&=(r_{H}^{\circ}\bullet m^{\circ}_{A}(g\circ A))\zeta_{H,H,I,A}(\Delta_{H}^{\bullet}\circ(I\bullet A))\\&=(r^\circ_H\bullet A)\varsigma^H_{H, I, A}.
\end{split}
\]
\end{invisible}
\end{remark}

%As we will see in \cref{pro:Hopfmonbraidedcat}, Hopf monoids in braided monoidal categories can be recovered as an instance of the general notion introduced above. For the moment, however, we proceed by establishing results in the broader setting. 
We now derive some identities involving the morphisms introduced above.

\begin{lemma}
\label{lem:grnat}
Let $\left( \mathcal{C},\circ ,I,\bullet ,J\right) $ be a duoidal category. Let $H$ be an object in $\Bimon \left( \mathcal{C},\circ,\bullet\right) $ and let $%
g:C\rightarrow H$ be a morphism in $\Comon \left( \mathcal{C}%
^{\bullet }\right) .$

If $f:H\rightarrow H^{\prime }$ is a morphism in $\Bimon \left(
\mathcal{C},\circ,\bullet\right) $, then%
\begin{equation}\label{form:grnat1}
\left( \left( C\circ J\right) \bullet f\right) g^{r}=\left( fg\right)
^{r}\left( C\circ f\right) .
\end{equation}%
\begin{invisible}
\begin{equation*}\label{form:glnat1}
\left( f\bullet \left( J\circ C\right) \right) g^{l} =\left( fg\right)
^{l}\left( f\circ C\right) ,
\end{equation*}
\end{invisible}
If $f:C^{\prime }\rightarrow C$ is a morphism in $\Comon \left(
\mathcal{C}^{\bullet }\right) ,$ then
\begin{equation}\label{form:grnat2}
g^{r}\left( f\circ H\right)  =\left( \left( f\circ J\right) \bullet
H\right) \left( gf\right) ^{r}.
\end{equation}
\begin{invisible}
\begin{equation*}\label{form:glnat2}
g^{l}\left( H\circ f\right)  =\left( H\bullet \left( J\circ f\right)
\right) \left( gf\right) ^{l},
\end{equation*}
\end{invisible}
Moreover, let $g:H\to A$ be a morphism in $\mathsf{Mon}(\Cc^{\circ})$.

If $f:H'\to H$ is in $\mathsf{Bimon}(\Cc,\circ,\bullet)$, then
\begin{equation}\label{formgrnat1}
   g_{r}((A\bullet I)\circ f)  =(A\bullet f)(gf)_{r}.
\end{equation}
If $f:A\to A'$ is a morphism in $\Mon(\Cc^\circ)$,
then
\begin{equation}\label{formgrnat2}
(f\bullet H)g_{r}= (fg)_{r}((f\bullet I)\circ H).
\end{equation}
\end{lemma}

\begin{proof}
If $f:H\rightarrow H^{\prime }$ is a morphism in $\Bimon \left(
\mathcal{C}\right) ,$ then it is a morphism of left $H$-modules from $(H,m^\circ_H)$ to $(H',m^\circ_{H'}(f\circ H'))$ in $\Cc^\circ$ so that, by \eqref{equivalentgr}, naturality of the co-Galois map and the easily checked identity $\varsigma^H _{C,J,H'}=\varsigma^{H'} _{C,J,H'}$, one gets \eqref{form:grnat1}.
\begin{invisible} We compute 
\[
\begin{split}
    \left( \left( C\circ J\right) \bullet f\right) g^{r}&\overset{\eqref{equivalentgr}}{=}\left( \left( C\circ J\right) \bullet f\right)\varsigma^H _{C,J,H}(C\circ \left(
l_{H}^{\bullet }\right) ^{-1})=\varsigma^H _{C,J,H'}\left( C\circ (J\bullet f\right))(C\circ \left(
l_{H}^{\bullet }\right) ^{-1})\\&\hspace{0.13cm}=\varsigma^{H'} _{C,J,H'}(C\circ \left(
l_{H'}^{\bullet }\right) ^{-1})\left( C\circ f\right)\overset{\eqref{equivalentgr}}{=}(fg)^{r}(C\circ f).
\end{split}
\]
First note that \begin{eqnarray*}
\varsigma _{C,J,H^{\prime }}^{H} &=&\left( \left( C\circ J\right) \bullet
m_{H^{\prime }}^{\circ }\left( f\circ H^{\prime }\right) \right) \zeta
_{C,H,J,H^{\prime }}\left( \rho _{C}^{\bullet }\circ \left( J\bullet
H^{\prime }\right) \right)  \\
&=&\left( \left( C\circ J\right) \bullet m_{H^{\prime }}^{\circ } \right) \zeta _{C,H^{\prime },J,H^{\prime }}\left(
\left( C\bullet f\right) \rho _{C}^{\bullet }\circ \left( J\bullet H^{\prime
}\right) \right) =\varsigma _{C,J,H^{\prime }}^{H^{\prime }}.
\end{eqnarray*}
We keep also the elementary proof of the formula:
\begin{eqnarray*}
\left( \left( C\circ J\right) \bullet f\right) g^{r} &=&\left( \left( C\circ
J\right) \bullet f\right) \left( \left( C\circ \varepsilon _{H}^{\bullet
}\right) \bullet m_{H}^{\circ }\left( g\circ H\right) \right) \Delta^{\bullet}
_{C\circ H} \\
&=&\left( \left( C\circ \varepsilon _{H^{\prime }}^{\bullet }f\right)
\bullet m_{H^{\prime }}^{\circ }\left( fg\circ f\right) \right) \Delta^{\bullet}
_{C\circ H} \\
&=&\left( \left( C\circ \varepsilon _{H^{\prime }}^{\bullet }\right) \bullet
m_{H^{\prime }}^{\circ }\left( fg\circ H^{\prime }\right) \right) \Delta^{\bullet}
_{C\circ H^{\prime }}\left( C\circ f\right)  \\
&=&\left( fg\right) ^{r}\left( C\circ f\right)
\end{eqnarray*}
so we get \eqref{form:grnat1}.
\end{invisible}
If $f:C^{\prime }\rightarrow C$ is a morphism in $\Comon \left(
\mathcal{C}^{\bullet }\right) ,$ then $f$ is a morphism of right $H$-comodules from $(C',(C'\bullet gf)\Delta_{C'}^\bullet)$ to $(C,(C\bullet g)\Delta_{C}^\bullet)$ so that, by \eqref{equivalentgr} and the naturality of the co-Galois map one gets \eqref{form:grnat2}.
\begin{invisible}
We have
\[
\begin{split}
    g^{r}\left( f\circ H\right)&\overset{\eqref{equivalentgr}}{=}\varsigma^H _{C,J,H}(C\circ \left(
l_{H}^{\bullet }\right) ^{-1})\left( f\circ H\right)=\varsigma^H _{C,J,H}\left( f\circ(J\bullet H)\right)(C'\circ \left(
l_{H}^{\bullet }\right) ^{-1})\\&\hspace{0.13cm}=\left(( f\circ J)\bullet H\right)\varsigma^H _{C',J,H}(C'\circ \left(
l_{H}^{\bullet }\right) ^{-1})\overset{\eqref{equivalentgr}}{=}\left( \left( f\circ J\right) \bullet
H\right) \left( gf\right) ^{r}.
\end{split}
\]
We also keep the original elementary proof:
\begin{eqnarray*}
g^{r}\left( f\circ H\right)  &=&\left( \left( C\circ \varepsilon
_{H}^{\bullet }\right) \bullet m_{H}^{\circ }\left( g\circ H\right) \right)
\Delta _{C\circ H}^{\bullet }\left( f\circ H\right)  \\
&=&\left( \left( C\circ \varepsilon _{H}^{\bullet }\right) \bullet
m_{H}^{\circ }\left( g\circ H\right) \right) \left( \left( f\circ H\right)
\bullet \left( f\circ H\right) \right) \Delta _{C^{\prime }\circ H}^{\bullet
} \\
&=&\left( \left( f\circ J\right) \bullet H\right) \left( \left( C^{\prime
}\circ \varepsilon _{H}^{\bullet }\right) \bullet m_{H}^{\circ }\left(
gf\circ H\right) \right) \Delta _{C^{\prime }\circ H}^{\bullet } \\
&=&\left( \left( f\circ J\right) \bullet H\right) \left( gf\right) ^{r},
\end{eqnarray*}
so we get \eqref{form:grnat2}.
\end{invisible}
By dual arguments one can show \eqref{formgrnat1} and \eqref{formgrnat2}, using the naturality of the Galois map.
\begin{invisible}
Given $g:H\to A$ in $\mathsf{Mon}(\Cc^{\circ})$ and $f:H'\to H$ in $\mathsf{Bimon}(\Cc,\circ,\bullet)$, we compute
\[
\begin{split}
    g_{r}((A\bullet I)\circ f)&=m_{A\bullet A}^{\circ }\left( \left( A\bullet u_{H}^{\circ }\right)
\circ \left( g\bullet H\right) \Delta _{H}^{\bullet }\right)((A\bullet I)\circ f)\\&=m_{A\bullet H}^{\circ }\left( \left( A\bullet fu_{H'}^{\circ }\right)
\circ \left( gf\bullet f\right) \Delta _{H'}^{\bullet }\right)\\&=(m^{\circ}_{A}\bullet m^{\circ}_{H})\zeta_{A,H,A,H}((A\bullet f)\circ(A\bullet f))((A\bullet u_{H'}^{\circ})\circ(gf\bullet H')\Delta^{\bullet}_{H'})\\&=(m^{\circ}_{A}\bullet m^{\circ}_{H})((A\circ A)\bullet(f\circ f))\zeta_{A,H',A,H'}((A\bullet u_{H'}^{\circ})\circ(gf\bullet H')\Delta^{\bullet}_{H'})\\&=(A\bullet f)(m^{\circ}_{A}\bullet m^{\circ}_{H'})\zeta_{A,H',A,H'}((A\bullet u_{H'}^{\circ})\circ(gf\bullet H')\Delta^{\bullet}_{H'})\\&=(A\bullet f)m_{A\bullet H'}^{\circ}((A\bullet u_{H'}^{\circ})\circ(gf\bullet H')\Delta^{\bullet}_{H'})\\&=(A\bullet f)(gf)_{r},
\end{split}
\]
so we get \eqref{formgrnat1}.
Let $f:A\to A'$ be a morphism in $\Mon(\Cc^\circ)$. We also compute
\[
\begin{split}
    (f\bullet H)g_{r}&=(f\bullet H)(m_{A}^{\circ}\bullet m_{H}^{\circ})\zeta_{A,H,A,H}\left( \left( A\bullet u_{H}^{\circ }\right)
\circ \left( g\bullet H\right) \Delta _{H}^{\bullet }\right)\\&=(m_{A'}^{\circ}\bullet m_{H}^{\circ})((f\circ f)\bullet( H\circ H))\zeta_{A,H,A,H}\left( \left( A\bullet u_{H}^{\circ }\right)
\circ \left( g\bullet H\right) \Delta _{H}^{\bullet }\right)\\&=(m_{A'}^{\circ}\bullet m_{H}^{\circ})\zeta_{A',H,A',H}((f\bullet u_{H}^{\circ})\circ(fg\bullet H)\Delta_{H}^{\bullet})\\&=m_{A'\bullet H}^{\circ}((A'\bullet u_{H}^{\circ})\circ(fg\bullet H)\Delta_{H}^{\bullet})((f\bullet I)\circ H)\\&=(fg)_{r}((f\bullet I)\circ H),
\end{split}
\]
hence we obtain \eqref{formgrnat2}.
\end{invisible}
\end{proof}

%As a fall-out of \cref{lem:grnat}, we get that $\Hopfmon(\Cc,\circ,\bullet)$ is a replete subcategory of $\Bimon(\Cc,\circ,\bullet)$. Indeed, given an isomorphism $f:H\to H'$ in $\Bimon(\Cc,\circ,\bullet)$ with $H'$ a Hopf monoid, then \eqref{form:grnat1} entails that $g^r$ is invertible for every $g:C\to H$ in $\Comon(\Cc^\bullet)$ as so is $(fg)^r$. Thus $H$ obeys \eqref{C1}. Similarly, by employing \eqref{formgrnat1} one proves that $H$ obeys \eqref{C2} and, consequently, it is a Hopf monoid.

\begin{remark}%[\rd{added:2025/11/13}]
Given a bimonoid $H$ and a morphism $g:C\rightarrow H$  in $\Comon(\Cc^\bullet)$, we compute%
\begin{align*}
&\left( g^{r}\bullet H\right) \left( \left( C\circ H\right) \bullet
m_{H}^{\circ }\left( g\circ H\right) \right) \Delta _{C\circ H}^{\bullet } \\
&=\left( \left( C\circ J\right) \bullet H\bullet m_{H}^{\circ }\left(
g\circ H\right) \right) \left( g^{r}\bullet \left( C\circ H\right) \right)
\Delta _{C\circ H}^{\bullet } \\
&=\left( \left( C\circ J\right) \bullet H\bullet m_{H}^{\circ }\left(
g\circ H\right) \right) \left( \left( C\circ \varepsilon _{H}^{\bullet
}\right) \bullet m_{H}^{\circ }\left( g\circ H\right) \bullet \left( C\circ
H\right) \right) \left( \Delta _{C\circ H}^{\bullet }\bullet \left( C\circ
H\right) \right) \Delta _{C\circ H}^{\bullet } \\
&=\left( \left( C\circ \varepsilon _{H}^{\bullet }\right) \bullet
m_{H}^{\circ }\left( g\circ H\right) \bullet m_{H}^{\circ }\left( g\circ
H\right) \right) \left( \left( C\circ H\right) \bullet \Delta _{C\circ
H}^{\bullet }\right) \Delta _{C\circ H}^{\bullet } \\
&=\left( \left( C\circ J\right) \bullet \Delta _{H}^{\bullet }\right)
\left( \left( C\circ \varepsilon _{H}^{\bullet }\right) \bullet m_{H}^{\circ
}\left( g\circ H\right) \right) \Delta _{C\circ H}^{\bullet } 
=\left( \left( C\circ J\right) \bullet \Delta _{H}^{\bullet }\right) g^{r}
\end{align*}%
so that we get another useful identity, namely
\begin{equation}
\left( g^{r}\bullet H\right) \left( \left( C\circ H\right) \bullet
m_{H}^{\circ }\left( g\circ H\right) \right) \Delta _{C\circ H}^{\bullet
}=\left( \left( C\circ J\right) \bullet \Delta _{H}^{\bullet }\right) g^{r}.
\label{form:Deltagl}
\end{equation}%
\end{remark}

\subsection{Factorization of bimonoids and Hopf monoids}

In order to describe the points of bimonoids, in the setting of the forthcoming \cref{thm:main}, we need an analogue of \cite[Theorem 7.2.3]{Majid-book}, characterizing the factorization of bimonoids in duoidal categories (see the first part of the proof therein).

\begin{theorem}
\label{thm:psi}
Let  $\xymatrixcolsep{1pc}\xymatrix{K\ar[r]^-{i}& A&\ar[l]_-{j} H}$ be  a cospan in $\Bimon (\Cc,\circ,\bullet)$ with $\varphi:= m_A^\circ(i\circ j): K\circ H\to A$ invertible. If we set $\psi:=\varphi^{-1} m_A^\circ(j\circ i):H\circ K \to K\circ H$, then $K\circ H$ inherits a unique bimonoid structure $K\circ_\psi H$ %(the \emph{smash product})
  such that $\varphi$ is a morphism of bimonoids, namely
  \begin{align}
 m^\circ_{K\circ_\psi H}=(m^\circ_K\circ m^\circ_H)(K\circ\psi\circ H),&\quad u^\circ_{K\circ_\psi H}=(u_K^\circ\circ u_H^\circ)(m^\circ_I)^{-1}, \label{eq:alg}\\\Delta^\bullet_{K\circ_\psi H}=\zeta_{K,K,H,H} (\Delta^\bullet_K\circ\Delta^\bullet_H) ,&\quad \varepsilon^\bullet_{K\circ_\psi H}=m^\circ_J(\varepsilon^\bullet_K\circ \varepsilon^\bullet_H).\label{eq:coalg}
  \end{align}
Given $H,K$ in $\Bimon (\Cc,\circ,\bullet)$ and a morphism $\psi:H\circ K\to K\circ H$ in $\Cc$, then $K\circ_\psi H$ is a bimonoid with structures \eqref{eq:alg}, \eqref{eq:coalg} if and only if the following conditions hold:
\begin{align}
 \psi(m_H^\circ\circ K)&=(K\circ m_H^\circ)(\psi\circ H)(H\circ \psi),\quad \psi(u_H^\circ\circ K)(l_K^\circ)^{-1}=(K\circ u_H^\circ)(r_K^\circ)^{-1},\label{eq:psi1}\\
 \psi(H\circ m_K^\circ)&=(m_K^\circ\circ H)(K\circ \psi)(\psi\circ K),\quad
 \psi(H\circ u_K^\circ)(r_H^\circ)^{-1}=(u_K^\circ\circ H)(l_H^\circ)^{-1},\label{eq:psi2}\\
\psi&\in \Comon (\mathcal{C}^\bullet)\label{eq:psi3}.
\end{align}
Thus, $K\circ_\psi H$ is a twisted tensor product of the monoids $H,K$ in $(\Comon(\Cc^\bullet),\circ)$ with twisting morphism $\psi$.
\end{theorem}

\begin{proof}
In view of \cite[Proposition 2.3]{BD-CrossII}, $K\circ H$ becomes a monoid in a unique way such that $\varphi$ is a morphism of monoids. Explicitly, $K\circ H$
inherits the structure of cross product monoid $K\circ_{\psi}^\phi H$ where $\psi$ is as in the statement and
$\phi:=\varphi^{-1} m_A^\circ(j\circ j):H\circ H \to K\circ H.$ Since $j$ in our setting is a morphism of monoids, we get $\phi=\varphi^{-1}j m_H^\circ$. Now, from
\[
\begin{split}
 \varphi(u_{K}^{\circ}\circ H)&=m_{A}^{\circ}(i\circ j)(u_{K}^{\circ}\circ H)=m_{A}^{\circ}(u_{A}^{\circ}\circ j)=m_{A}^{\circ}(u_{A}^{\circ}\circ A)(I\circ j)=l_{A}^{\circ}(I\circ j)=j l_{H}^{\circ}
\end{split}
\]
we get 
$\varphi^{-1}j=(u^\circ_K\circ H)(l_{H}^{\circ})^{-1}$.
Thus
\[\phi=\varphi^{-1}j m_H^\circ{=}(u^\circ_K\circ H)(l_{H}^{\circ})^{-1}m_H^\circ
=(u^\circ_K\circ H)(I\circ m_H^\circ)(l_{H\circ H}^{\circ})^{-1}
=(u^\circ_K\circ m_H^\circ)(l_{H\circ H}^{\circ})^{-1}.\]
By  \cite[Proposition 2.5]{BD-CrossII}, $K\circ_\psi^\phi H=K\circ_\psi H$ is a smash product monoid with multiplication and unit as in the statement. As observed in the proof of \cite[Proposition 2.6]{BD-CrossI}, the necessary and sufficient conditions
for $K\circ_\psi H$ being an algebra with structure as in \eqref{eq:alg} are given by the equations \eqref{eq:psi1} and \eqref{eq:psi2}.

\begin{invisible}
   Assume $\varphi$ is an isomorphism in $\mathcal{C}$. Define $\psi:=\varphi^{-1} m_A^\circ(j\circ i):H\circ K \to K\circ H$. We transport the structure of $A$ on $K\circ H$ via $\varphi$. Define
\[
m^\circ_{K\circ H}:=\varphi^{-1}m_A^\circ (\varphi\circ\varphi):K\circ H\circ K\circ H \to K\circ H.
\]
We have
\[
\begin{split}
\varphi(m^\circ_K\circ m_H^\circ)(K\circ\psi \circ H)&=m_A^\circ (i\circ j)(m^\circ_K\circ m^\circ_H)(K\circ\psi\circ H)
\\&=m_A^\circ(m_A^\circ\circ m_A^\circ )(i\circ i\circ j\circ j)(K\circ\psi\circ H)\\&=m_A^\circ (A\circ m_A^\circ)(i\circ A\circ j)(K\circ\varphi\psi\circ H)\\&=
m_A^\circ (A\circ m_A^\circ)(A\circ m_A^\circ \circ A)(i\circ j\circ i\circ j)\\&=m_A^\circ (m_A^\circ\circ m_A^\circ)(i\circ j\circ i\circ j)=m_A^\circ(\varphi\circ \varphi)=\varphi m^\circ_{K\circ H},
\end{split}
\]
so $m^\circ_{K\circ H}=(m^\circ_K\circ m^\circ_H)(K\circ\psi\circ H)$.

Define $u^\circ_{K\circ H}:=\varphi^{-1}u_A^\circ :I\to K\circ H$. Then, $\varphi(u^\circ_K\circ u^\circ_H)(m^\circ_I)^{-1}=m_A^\circ (i\circ j)(u^\circ_K\circ u^\circ_H)(m^\circ_I)^{-1}=m_A^\circ (u_A^\circ\circ u_A^\circ)(m_I^\circ)^{-1}=u_A^\circ$, so $u^\circ_{K\circ H}=(u_K^\circ\circ u_H^\circ)(m^\circ_I)^{-1}$.
\end{invisible}

Now define $\Delta^\bullet_{K\circ H}:=(\varphi^{-1}\bullet \varphi^{-1})\Delta_A^\bullet \varphi:K\circ H\to (K\circ H)\bullet (K\circ H)$. We have
\[
\begin{split}
	(\varphi\bullet \varphi)\zeta_{K,K,H,H}(\Delta_K^\bullet\circ\Delta^\bullet_H)&=(m_A^\circ\bullet m_A^\circ) ((i\circ j)\bullet (i\circ j))\zeta _{K,K,H,H}(\Delta_K^\bullet\circ\Delta_H^\bullet)\\&=(m_A^\circ\bullet m_A^\circ)\zeta_{A,A,A,A}((i\bullet i)\circ (j\bullet j))(\Delta^\bullet_K\circ\Delta^\bullet_H)\\&
	=(m_A^\circ\bullet m_A^\circ)\zeta_{A,A,A,A} (\Delta^\bullet_A\circ\Delta_A^\bullet)(i\circ j)=\Delta_A^\bullet m_A^\circ (i\circ j)=\Delta_A^\bullet\varphi,
\end{split}
\]
hence $\Delta^\bullet_{K\circ H}=\zeta_{K,K,H,H} (\Delta^\bullet_K\circ\Delta^\bullet_H)$. Define $\varepsilon^\bullet_{K\circ H}:=\varepsilon_A^\bullet\varphi: K\circ H\to J$. We have $\varepsilon_A^\bullet\varphi=\varepsilon_A^\bullet m_A^\circ (i\circ j)=m^\circ_J(\varepsilon_A^\bullet\circ\varepsilon_A^\bullet)(i\circ j)=m_J^\circ(\varepsilon_K^\bullet\circ\varepsilon^\bullet_H)$, so $\varepsilon^\bullet_{K\circ H}=m^\circ_J(\varepsilon^\bullet_K\circ \varepsilon^\bullet_H)$.

Note that the comonoid structure of $K\circ_{\psi} H$ is exactly the one coming from the fact that $K$ and $ H$ are in $\Comon (\mathcal{C}^\bullet)$ and that the latter is a monoidal category through $\circ$, cf. \cite[Proposition 6.35]{Aguiar}.
It remains to rewrite the compatibility conditions between the monoid and comonoid structures in terms of equations involving $\psi$. Since $u_K^\circ$ and $u_H^\circ$ are morphisms in $(\Comon (\mathcal{C}^\bullet),\circ,I)$, it is always true that $u_K^\circ\circ u_H^\circ$ is a morphism in  $\Comon (\mathcal{C}^\bullet)$. Thus we just have to understand when $m^\circ_{K\circ_\psi H}$ is a morphism of comonoids. If $m^\circ_{K\circ_\psi H}$ is a morphism of comonoids, then so is
\begin{align*}
m^\circ_{K\circ_\psi H}
&(u_K^\circ\circ H\circ K\circ u_H^\circ)((l_H^\circ)^{-1}\circ (r_K^\circ)^{-1})\\
&=(m^\circ_K\circ m^\circ_H)(K\circ\psi\circ H)(u_K^\circ\circ H\circ K\circ u_H^\circ)((l_H^\circ)^{-1}\circ (r_K^\circ)^{-1})
\\
&=(m^\circ_K\circ m^\circ_H)(u_K^\circ\circ {K}\circ  {H}\circ u_H^\circ)(I\circ\psi\circ I)((l_H^\circ)^{-1}\circ (r_K^\circ)^{-1})=\psi
\end{align*} as it is a composition of morphisms in the monoidal category $(\Comon (\mathcal{C}^\bullet),\circ, I)$.
Similarly, if $\psi$ is a morphism of comonoids, then so is the composition $m^\circ_{K\circ_\psi H}=(m^\circ_K\circ m^\circ_H)(K\circ\psi\circ H).$
Therefore, $m^\circ_{K\circ_\psi H}$ is a morphism of comonoids if and only if so is $\psi$.
\end{proof}

We now explore the way in which an antipode shapes the context of  \cref{thm:psi}.

\begin{theorem}%[\rd{modified:2026/03/13}]
\label{thm:psitriang}
Let $(\Cc,\circ,I,\bullet,J)$ be a duoidal category with a reversion.
Given $H,K$ in $\Hopfmon(\Cc,\circ,\bullet)$ and a morphism $\psi :H\circ
K\rightarrow K\circ H$ in $\Cc$ obeying the conditions \eqref{eq:psi1}, \eqref{eq:psi2} and \eqref{eq:psi3}, then $K\circ _{\psi }H$ is a Hopf monoid with antipode $\Sigma
_{K\circ _{\psi }H}=\psi ^{\revsn }\phi _{K,H}^{\circ }\left( \Sigma _{K}\circ
\Sigma _{H}\right) .$

Assume further that there is a morphism of comonoids $%
\varepsilon _{K}^{\circ }:K\rightarrow I$ and set
\[
\chi  :=\left( K\circ \varepsilon _{H}^{\bullet }\right) \psi
:H\circ K\rightarrow K\circ J, \qquad
\triangleleft  :=l_{H}^{\circ }\left( \varepsilon _{K}^{\circ }\circ
H\right) \psi :H\circ K\rightarrow H.
\]
Then $\psi =\left( \left( u_{H}^{\circ }\varepsilon _{K}^{\circ }\right)
^{r}\right) ^{-1}\left( \chi \bullet \triangleleft \right) \Delta
_{H\circ K}^{\bullet }$, so it is completely determined by the morphisms $\chi, \triangleleft$.

If $\varepsilon^\circ_K$ is indeed a morphism of bimonoids, the following identities hold true:
\begin{eqnarray}
\chi \left( m_{H}^{\circ }\circ K\right)  &=&\left( K\circ m_{J}^{\circ
}\right) \left( \chi \circ J\right) \left( H\circ \chi \right) ,\qquad \chi
\left( u_{H}^{\circ }\circ K\right)(l^\circ_K)^{-1} =(K\circ u_{J}^{\circ })(r^\circ_K)^{-1},  \label{eq:chi1}
\\
\chi \left( H\circ m_{K}^{\circ }\right)  &=&\left( m_{K}^{\circ }\circ
J\right) \left( K\circ \chi \right) \left( \psi \circ K\right) ,\qquad \chi
\left( H\circ u_{K}^{\circ }\right) (r^\circ_H)^{-1}=(u_{K}^{\circ }\circ \varepsilon
_{H}^{\bullet })(l^\circ_H)^{-1},  \label{eq:chi2} \\
\triangleleft \left( m_{H}^{\circ }\circ K\right)  &=&m_{H}^{\circ }\left(
\triangleleft \circ H\right) \left( H\circ \psi \right) ,\qquad
\triangleleft \left( u_{H}^{\circ }\circ K\right) (l^\circ_K)^{-1}=u_{H}^{\circ }\varepsilon
_{K}^{\circ }  \label{eq:triangl1} \\
\triangleleft \left( H\circ m_{K}^{\circ }\right)  &=&\triangleleft \left(
\triangleleft \circ K\right) ,\qquad \triangleleft \left( H\circ
u_{K}^{\circ }\right) (r^\circ_H)^{-1}=\mathrm{Id}_{H}.  \label{eq:triangl2}
\end{eqnarray}
In particular $(H,\triangleleft)$ is a right $K$-module.
\end{theorem}

\begin{proof}
By \cref{thm:psi}, we already know that $K\circ_{\psi}H$ is in $\Bimon(\Cc,\circ,\bullet)$ with structures \eqref{eq:alg} and \eqref{eq:coalg}. We compute
\begin{eqnarray*}
&&\Sigma _{K\circ _{\psi }H}\star \mathrm{Id}_{K\circ _{\psi }H} \\
&=&\left( J\circ \left( I\bullet m_{K\circ _{\psi }H}^{\circ }\right)
\right) \psi _{K\circ H,K\circ H}\left( \Sigma _{K\circ _{\psi }H}\bullet
\left( K\circ H\right) \right) \Delta _{K\circ _{\psi }H}^{\bullet } \\
&=&
\left( J\circ \left( I\bullet m_{K\circ _{\psi }H}^{\circ }\right) \right)
\psi _{K\circ H,K\circ H}\left( \psi ^{\revsn}\bullet \left( K\circ H\right)
\right) \left( \phi _{K,H}^{\circ }\bullet \left( K\circ H\right) \right)
\\
&&\left( \left( \Sigma _{K}\circ \Sigma _{H}\right) \bullet \left( K\circ
H\right) \right) \zeta _{K,K,H,H}\left( \Delta _{K}^{\bullet }\circ \Delta
_{H}^{\bullet }\right)\\
&=&
\left( J\circ \left( I\bullet m_{K\circ _{\psi }H}^{\circ }\right) \right)
\left( J\circ \left( I\bullet \left( \psi \circ K\circ H\right) \right)
\right) \underbracket[0.140ex]{\psi _{H\circ K,K\circ H}\left( \phi
_{K,H}^{\circ }\bullet \left( K\circ H\right) \right) \zeta _{K^{\revsn
},K,H^{\revsn },H}} \\
&&\left( \left( \Sigma _{K}\bullet K\right) \circ \left( \Sigma _{H}\bullet
H\right) \right) \left( \Delta _{K}^{\bullet }\circ \Delta _{H}^{\bullet
}\right)\\
&\overset{\eqref{form:BohmLem4.2}}{=}&
\left( J\circ \left( I\bullet \underbracket[0.140ex]{m_{K\circ _{\psi }H}^{\circ }}\right) \right)
\left( J\circ \left( I\bullet \left( \psi \circ K\circ H\right) \right)
\right) \left( m_{J}^{\circ }\circ \left( I\bullet \left( H\circ K\circ
K\circ H\right) \right) \right)  \\
&&\left( J\circ \psi _{H,K\circ K\circ H}\right) \left( J\circ \left( l_{H^{%
\revsn}}^{\circ }\bullet \left( K\circ K\circ H\right) \right) \right)
\left( J\circ \zeta _{I,K\circ K,H^{\revsn},H}\right)  \\
&&\left( \psi _{K,K}\left( \Sigma _{K}\bullet K\right) \Delta _{K}^{\bullet
}\circ \left( \Sigma _{H}\bullet H\right) \Delta _{H}^{\bullet }\right)\\
% &=&\left[
% \begin{array}{c}
% \left( J\circ \left( I\bullet \left( m_{K}^{\circ }\circ m_{H}^{\circ
% }\right) \left( K\circ \psi \circ H\right) \right) \right) \left( J\circ
% \left( I\bullet \left( \psi \circ K\circ H\right) \right) \right) \left(
% m_{J}^{\circ }\circ \left( I\bullet \left( H\circ K\circ K\circ H\right)
% \right) \right)  \\
% \left( J\circ \psi _{H,K\circ K\circ H}\right) \left( J\circ \left( l_{H^{%
% \revsn}}^{\circ }\bullet \left( K\circ K\circ H\right) \right) \right)
% \left( J\circ \zeta _{I,K\circ K,H^{\revsn},H}\right)  \\
% \left( \psi _{K,K}\left( \Sigma _{K}\bullet K\right) \Delta _{K}^{\bullet
% }\circ \left( \Sigma _{H}\bullet H\right) \Delta _{K}^{\bullet }\right)
% \end{array}%
% \right]  \\
&=&
\left( J\circ \left( I\bullet \left( K\circ m_{H}^{\circ }\right) \right)
\right) \left( m_{J}^{\circ }\circ \left( I\bullet \left( \underbracket[0.140ex]{\left(
m_{K}^{\circ }\circ H\right) \left( K\circ \psi \right) \left( \psi \circ
K\right) }\circ H\right) \right) \right)  \\
 &&\left( J\circ \psi _{H,K\circ K\circ H}\right)
\left( J\circ \left( l_{H^{\revsn}}^{\circ }\bullet \left( K\circ K\circ
H\right) \right) \right)  \\
&&\left( J\circ \zeta _{I,K\circ K,H^{\revsn},H}\right) \left( \psi
_{K,K}\left( \Sigma _{K}\bullet K\right) \Delta _{K}^{\bullet }\circ \left(
\Sigma _{H}\bullet H\right) \Delta _{H}^{\bullet }\right)\\
&\overset{\eqref{eq:psi2}}{=}&
\left( J\circ \left( I\bullet \left( K\circ m_{H}^{\circ }\right) \right)
\right) \left( m_{J}^{\circ }\circ \left( I\bullet \left( \psi \left( H\circ m_{K}^{\circ
}\right) \circ H\right) \right) \right)  \left( J\circ \psi _{H,K\circ K\circ H}\right) \\
 &&\left( J\circ \left( l_{H^{%
\revsn}}^{\circ }\bullet \left( K\circ K\circ H\right) \right) \right)
\left( J\circ \zeta _{I,K\circ K,H^{\revsn},H}\right)  
\left( \psi _{K,K}\left( \Sigma _{K}\bullet K\right) \Delta _{K}^{\bullet
}\circ \left( \Sigma _{H}\bullet H\right) \Delta _{H}^{\bullet }\right)\\
% &=&\left[
% \begin{array}{c}
% \left( J\circ \left( I\bullet \left( K\circ m_{H}^{\circ }\right) \right)
% \right) \left( J\circ \left( I\bullet \left( \psi \circ H\right) \right)
% \right) \left( m_{J}^{\circ }\circ \left( I\bullet \left( H\circ K\circ
% H\right) \right) \right)  \\
% \left( J\circ \psi _{H,K\circ H}\right) \left( J\circ \left( l_{H^{\revsn%
% }}^{\circ }\bullet \left( K\circ H\right) \right) \right) \left( J\circ
% \zeta _{I,K,H^{\revsn},H}\right)  \\
% \left( \underbracket[0.140ex]{\left( J\circ \left( I\bullet m_{K}^{\circ
% }\right) \right) \psi _{K,K}\left( \Sigma _{K}\bullet K\right) \Delta
% _{K}^{\bullet }}\circ \left( \Sigma _{H}\bullet H\right) \Delta
% _{K}^{\bullet }\right)
% \end{array}%
% \right]  \\
&=&
\left( J\circ \left( I\bullet \left( K\circ m_{H}^{\circ }\right) \right)
\right) \left( m_{J}^{\circ }\circ \left( I\bullet \left( \psi \circ H\right) \right)
\right)   \left( J\circ \psi _{H,K\circ H}\right) \\
&&\left( J\circ \left( l_{H^{\revsn%
}}^{\circ }\bullet \left( K\circ H\right) \right) \right) \left( J\circ
\zeta _{I,K,H^{\revsn},H}\right)  
\left( (\Sigma _{K}\star \mathrm{Id}_{K})\circ \left( \Sigma _{H}\bullet
H\right) \Delta _{H}^{\bullet }\right)\\
&=&
\left( J\circ \left( I\bullet \left( K\circ m_{H}^{\circ }\right) \right)
\right) \left( m_{J}^{\circ }\circ \left( I\bullet \left( \psi \circ H\right) \right)
\right) \left( J\circ \psi _{H,K\circ H}\right) \left( J\circ \left( l_{H^{\revsn
}}^{\circ }\bullet \left( K\circ H\right) \right) \right)  \\
&& \left( J\circ
\zeta _{I,K,H^{\revsn},H}\right)  
\left( (J\circ \left( I\bullet u_{K}^{\circ }\right) \Delta _{I}^{\bullet
}\circ \left( H^{\revsn}\bullet H\right) \right) \left( \left( r_{J}^{\circ
}\right) ^{-1}\varepsilon _{H}^{\bullet }\circ \left( \Sigma _{H}\bullet
H\right) \Delta _{H}^{\bullet }\right)\\
&=&
\left( J\circ \left( I\bullet \left( K\circ m_{H}^{\circ }\right) \right)
\right) \left( m_{J}^{\circ }\circ \left( I\bullet \left( \underbracket[0.140ex]{\psi \left( H\circ u_{K}^{\circ
}\right)} \circ H\right) \right) \right)\left( J\circ \psi _{H,I\circ H}\right)  \\
 && \left( J\circ \left(
l_{H^{\revsn}}^{\circ }\bullet \left( I\circ H\right) \right) \right) \left(
J\circ \zeta _{I,I,H^{\revsn},H}\right)  
\left( (J\circ \Delta _{I}^{\bullet }\circ \left( H^{\revsn}\bullet H\right)
\right) \left( \left( r_{J}^{\circ }\right) ^{-1}\varepsilon _{H}^{\bullet
}\circ \left( \Sigma _{H}\bullet H\right) \Delta _{H}^{\bullet }\right) \\
&\overset{\eqref{eq:psi2}}{=}&
\left( J\circ \left( I\bullet \left( K\circ m_{H}^{\circ }\right) \right)
\right) \left( m_{J}^{\circ }\circ \left( I\bullet \left( \left( u_{K}^{\circ }\circ
H\right) \left( l_{H}^{\circ }\right) ^{-1}r_{H}^{\circ }\circ H\right)
\right) \right) \left( J\circ \psi _{H,I\circ H}\right)  \\
&&\left( J\circ \left(
l_{H^{\revsn}}^{\circ }\bullet \left( I\circ H\right) \right) \right) \left(
J\circ \zeta _{I,I,H^{\revsn},H}\right)  \\
&&\left( (J\circ \Delta _{I}^{\bullet }\circ \left( H^{\revsn}\bullet H\right)
\right) \left( \left( r_{J}^{\circ }\right) ^{-1}\circ \left( H^{\revsn%
}\bullet H\right) \right) \left( \varepsilon _{H}^{\bullet }\circ \left(
\Sigma _{H}\bullet H\right) \Delta _{H}^{\bullet }\right)\\
&=&
\left( J\circ \left( I\bullet \left( u_{K}^{\circ }\circ m_{H}^{\circ
}\right) \right) \right) \left( J\circ \left( I\bullet \left( \left(
l_{H}^{\circ }\right) ^{-1}\circ H\right) \right) \right) 
\left( m_{J}^{\circ }\circ \left( I\bullet \left( H\circ l_{H}^{\circ }\right) \right) \right)\left( J\circ \psi _{H,I\circ H}\right)   \\
&&\left( J\circ \left(
l_{H^{\revsn}}^{\circ }\bullet \left( I\circ H\right) \right) \right) \left(
J\circ \zeta _{I,I,H^{\revsn},H}\right)  \\
&&\left( (J\circ \Delta _{I}^{\bullet }\circ \left( H^{\revsn}\bullet H\right)
\right) \left( J\circ \left( l_{H^{\revsn}\bullet H}^{\circ }\right)
^{-1}\right) \left( \varepsilon _{H}^{\bullet }\circ \left( \Sigma
_{H}\bullet H\right) \Delta _{H}^{\bullet }\right)\\
&=&
\left( J\circ \left( I\bullet \left( u_{K}^{\circ }\circ m_{H}^{\circ
}\right) \right) \right) \left( m_{J}^{\circ }\circ \left( I\bullet \left( l_{H\circ
H}^{\circ }\right) ^{-1}\right) \right)  \left( J\circ \psi
_{H,H}\right) \left( J\circ \underbracket[0.140ex]{\left( l_{H^{\revsn }}^{\circ }\bullet
l_{H}^{\circ }\right) }\right) \\
 &&\left( J\circ \underbracket[0.140ex]{\zeta
_{I,I,H^{\revsn },H}}\right) \left( (J\circ \underbracket[0.140ex]{\Delta
_{I}^{\bullet }\circ \left( H^{\revsn }\bullet H\right) }\right)  
\left( J\circ \underbracket[0.140ex]{\left( l_{H^{\revsn }\bullet H}^{\circ
}\right) ^{-1}}\right) \left( \varepsilon _{H}^{\bullet }\circ \left( \Sigma
_{H}\bullet H\right) \Delta _{H}^{\bullet }\right)\\
&\overset{\eqref{eq:unit1}}{=}&\left( J\circ \left( I\bullet \left(
u_{K}^{\circ }\circ m_{H}^{\circ }\right) \right) \right) \left( m_{J}^{\circ }\circ
\left( I\bullet \left( l_{H\circ H}^{\circ }\right) ^{-1}\right) \right)
 \left( \varepsilon _{H}^{\bullet }\circ \psi _{H,H}\left( \Sigma
_{H}\bullet H\right) \Delta _{H}^{\bullet }\right)  \\
&=&\left( J\circ \left( I\bullet \left( u_{K}^{\circ }\circ H\right) \right)
\right) \left( m_{J}^{\circ }\circ \left( I\bullet \left( l_{H}^{\circ }\right)
^{-1}\right) \right) \left( \varepsilon _{H}^{\bullet }\circ \underbracket[0.140ex]{%
\left( J\circ \left( I\bullet m_{H}^{\circ }\right) \right) \psi
_{H,H}\left( \Sigma _{H}\bullet H\right) \Delta _{H}^{\bullet }}\right)  \\
&=&\left( J\circ \left( I\bullet \left( u_{K}^{\circ }\circ H\right) \right)
\right) \left( m_{J}^{\circ }\circ \left( I\bullet \left( l_{H}^{\circ }\right)
^{-1}\right) \right)  \left( \varepsilon _{H}^{\bullet }\circ (\Sigma _{H}\star \mathrm{Id}%
_{H})\right)  \\
&=&\left( J\circ \left( I\bullet \left( u_{K}^{\circ }\circ H\right) \right)
\right) \left( m_{J}^{\circ }\circ \left( I\bullet \left( l_{H}^{\circ }\right)
^{-1}\right) \right)  \left( \varepsilon _{H}^{\bullet }\circ (J\circ \left( I\bullet
u_{H}^{\circ }\right) )\left( J\circ \Delta _{I}^{\bullet }\right) \left(
r_{J}^{\circ }\right) ^{-1}\varepsilon _{H}^{\bullet }\right)  \\
&=&\left( J\circ \left( I\bullet \left( u_{K}^{\circ }\circ u_{H}^{\circ
}\right) \right) \right) \left( J\circ \left( I\bullet \left( l_{I}^{\circ
}\right) ^{-1}\right) \right) \left( J\circ \Delta _{I}^{\bullet }\right)
\left( m_{J}^{\circ }\circ I\right) \left( J\circ \left( r_{J}^{\circ
}\right) ^{-1}\right) \left( \varepsilon _{H}^{\bullet }\circ \varepsilon
_{H}^{\bullet }\right)  \\
&=&\left( J\circ \left( I\bullet \left( u_{K}^{\circ }\circ u_{H}^{\circ
}\right) \left( m_{I}^{\circ }\right) ^{-1}\right) \Delta _{I}^{\bullet
}\right) \left( m_{J}^{\circ }\circ I\right) \left( r_{J\circ J}^{\circ
}\right) ^{-1}\left( \varepsilon _{H}^{\bullet }\circ \varepsilon
_{H}^{\bullet }\right)  \\
&=&\left( J\circ \left( I\bullet u_{K\circ _{\psi }H}^{\circ }\right) \Delta _{I}^{\bullet
}\right) \left( r_{J}^{\circ }\right) ^{-1}m_{J}^{\circ }\left( \varepsilon
_{K}^{\bullet }\circ \varepsilon _{H}^{\bullet }\right)  \\
&=&\left( J\circ \left( I\bullet u_{K\circ _{\psi }H}^{\circ }\right) \Delta
_{I}^{\bullet }\right) \left( r_{J}^{\circ }\right) ^{-1}\varepsilon
_{K\circ _{\psi }H}^{\bullet }=1_{K\circ _{\psi }H,K\circ _{\psi }H}^{r}.
\end{eqnarray*}
In a similar way, by employing \eqref{form:BohmLem4.2B}, one proves that $\mathrm{Id}_{K\circ _{\psi }H}\star \Sigma
_{K\circ _{\psi }H}=1_{K\circ _{\psi }H,K\circ _{\psi }H}^{l}.$
\begin{invisible}
We compute%
\begin{eqnarray*}
&&\mathrm{Id}_{K\circ _{\psi }H}\star \Sigma _{K\circ _{\psi }H} \\
&=&\left( \left( m_{K\circ _{\psi }H}^{\circ }\bullet I\right) \circ
J\right) \varphi _{K\circ H,K\circ H}\left( \left( K\circ H\right) \bullet
\Sigma _{K\circ _{\psi }H}\right) \Delta _{K\circ _{\psi }H}^{\bullet } \\
&=&\left( \left( m_{K\circ _{\psi }H}^{\circ }\bullet I\right) \circ
J\right) \varphi _{K\circ H,K\circ H}\left( \left( K\circ H\right) \bullet
\psi ^{\revsn }\right) \left( \left( K\circ H\right) \bullet \phi
_{K,H}^{\circ }\right) \left( \left( K\circ H\right) \bullet \left( \Sigma
_{K}\circ \Sigma _{H}\right) \right) \Delta _{K\circ _{\psi }H}^{\bullet } \\
&=&\left( \left( \left( m_{K}^{\circ }\circ m_{H}^{\circ }\right) \left(
K\circ \psi \circ H\right) \bullet I\right) \circ J\right) \left( \left(
\left( K\circ H\circ \psi \right) \bullet I\right) \circ J\right) \varphi
_{K\circ H,H\circ K}\\&&\left( \left( K\circ H\right) \bullet \phi _{K,H}^{\circ
}\right) \left( \left( K\circ H\right) \bullet \left( \Sigma _{K}\circ
\Sigma _{H}\right) \right) \Delta _{K\circ _{\psi }H}^{\bullet } \\
&=&\left( \left( \left( m_{K}^{\circ }\circ H\right) \bullet I\right) \circ
J\right) \left( \left( \left( K\circ \underbracket[0.140ex]{\left( K\circ m_{H}^{\circ
}\right) \left( \psi \circ H\right) \left( H\circ \psi \right) }\right)
\bullet I\right) \circ J\right) \varphi _{K\circ H,H\circ K}\\&&\left( \left(
K\circ H\right) \bullet \phi _{K,H}^{\circ }\right) \left( \left( K\circ
H\right) \bullet \left( \Sigma _{K}\circ \Sigma _{H}\right) \right) \Delta
_{K\circ _{\psi }H}^{\bullet } \\
&=&\left( \left( \left( m_{K}^{\circ }\circ H\right) \bullet I\right) \circ
J\right) \left( \left( \left( K\circ \psi \left( m_{H}^{\circ }\circ
K\right) \right) \bullet I\right) \circ J\right) \varphi _{K\circ H,H\circ
K}\left( \left( K\circ H\right) \bullet \phi _{K,H}^{\circ }\right) \left(
\left( K\circ H\right) \bullet \left( \Sigma _{K}\circ \Sigma _{H}\right)
\right) \Delta _{K\circ _{\psi }H}^{\bullet } \\
&=&\left( \left( \left( m_{K}^{\circ }\circ H\right) \left( K\circ \psi
\right) \bullet I\right) \circ J\right) \left( \left( \left( K\circ
m_{H}^{\circ }\circ K\right) \bullet I\right) \circ J\right) \varphi
_{K\circ H,H\circ K}\\&&\left( \left( K\circ H\right) \bullet \phi _{K,H}^{\circ
}\right) \left( \left( K\circ H\right) \bullet \left( \Sigma _{K}\circ
\Sigma _{H}\right) \right) \zeta _{K,K,H,H}\left( \Delta _{K}^{\bullet
}\circ \Delta _{H}^{\bullet }\right)  \\
&=&\left( \left( \left( m_{K}^{\circ }\circ H\right) \left( K\circ \psi
\right) \bullet I\right) \circ J\right) \left( \left( \left( K\circ
m_{H}^{\circ }\circ K\right) \bullet I\right) \circ J\right) \underbracket[0.140ex]{%
\varphi _{K\circ H,H\circ K}\left( \left( K\circ H\right) \bullet \phi
_{K,H}^{\circ }\right) \zeta _{K,K^{\revsn },H,H^{\revsn }}}\\&&\left( \left(
K\bullet \Sigma _{K}\right) \Delta _{K}^{\bullet }\circ \left( H\bullet
\Sigma _{H}\right) \Delta _{H}^{\bullet }\right)  \\
&\overset{(\ref{form:BohmLem4.2B})}{=}&\left[
\begin{array}{c}
\left( \left( \left( m_{K}^{\circ }\circ H\right) \left( K\circ \psi \right)
\bullet I\right) \circ J\right) \left( \left( \left( K\circ m_{H}^{\circ
}\circ K\right) \bullet I\right) \circ J\right)  \\
\left( \left( \left( K\circ H\circ H\circ K\right) \bullet I\right) \circ
m_{J}^{\circ }\right) \left( \varphi _{K\circ H\circ H,K}\circ J\right)
\left( \left( \left( K\circ H\circ H\right) \bullet r_{K^{\revsn }}^{\circ
}\right) \circ J\right)  \\
\left( \zeta _{K,K^{\revsn },H\circ H,I}\circ J\right) \left( \left( K\bullet
\Sigma _{K}\right) \Delta _{K}^{\bullet }\circ \varphi _{H,H}\left( H\bullet
\Sigma _{H}\right) \Delta _{H}^{\bullet }\right)
\end{array}%
\right]  \\
&=&\left[
\begin{array}{c}
\left( \left( \left( m_{K}^{\circ }\circ H\right) \left( K\circ \psi \right)
\bullet I\right) \circ J\right)  \\
\left( \left( \left( K\circ H\circ K\right) \bullet I\right) \circ
m_{J}^{\circ }\right) \left( \varphi _{K\circ H,K}\circ J\right) \left(
\left( \left( K\circ H\right) \bullet r_{K^{\revsn }}^{\circ }\right) \circ
J\right)  \\
\left( \zeta _{K,K^{\revsn },H,I}\circ J\right) \left( \left( K\bullet \Sigma
_{K}\right) \Delta _{K}^{\bullet }\circ \underbracket[0.140ex]{\left( \left(
m_{H}^{\circ }\bullet I\right) \circ J\right) \varphi _{H,H}\left( H\bullet
\Sigma _{H}\right) \Delta _{H}^{\bullet }}\right)
\end{array}%
\right]  \\
&=&\left[
\begin{array}{c}
\left( \left( \left( m_{K}^{\circ }\circ H\right) \left( K\circ \psi \right)
\bullet I\right) \circ J\right)  \\
\left( \left( \left( K\circ H\circ K\right) \bullet I\right) \circ
m_{J}^{\circ }\right) \left( \varphi _{K\circ H,K}\circ J\right) \left(
\left( \left( K\circ H\right) \bullet r_{K^{\revsn }}^{\circ }\right) \circ
J\right)  \\
\left( \zeta _{K,K^{\revsn },H,I}\circ J\right) \left( \left( K\bullet \Sigma
_{K}\right) \Delta _{K}^{\bullet }\circ \left( \mathrm{Id}_{H}\star \Sigma
_{H}\right) \right)
\end{array}%
\right]  \\
&=&\left[
\begin{array}{c}
\left( \left( \left( m_{K}^{\circ }\circ H\right) \left( K\circ \psi \right)
\bullet I\right) \circ J\right)  \\
\left( \left( \left( K\circ H\circ K\right) \bullet I\right) \circ
m_{J}^{\circ }\right) \left( \varphi _{K\circ H,K}\circ J\right) \left(
\left( \left( K\circ H\right) \bullet r_{K^{\revsn }}^{\circ }\right) \circ
J\right)  \\
\left( \zeta _{K,K^{\revsn },H,I}\circ J\right) \left( \left( K\bullet \Sigma
_{K}\right) \Delta _{K}^{\bullet }\circ \left( \left( u_{H}^{\circ }\bullet
I\right) \Delta _{I}^{\bullet }\circ J\right) \left( l_{J}^{\circ }\right)
^{-1}\varepsilon _{H}^{\bullet }\right)
\end{array}%
\right]  \\
&=&\left[
\begin{array}{c}
\left( \left( \left( m_{K}^{\circ }\circ H\right) \left( K\circ \underbracket[0.140ex]{%
\psi \left( u_{H}^{\circ }\circ K\right) }\right) \bullet I\right) \circ
J\right)  \\
\left( \left( \left( K\circ I\circ K\right) \bullet I\right) \circ
m_{J}^{\circ }\right) \left( \varphi _{K\circ I,K}\circ J\right) \left(
\left( \left( K\circ I\right) \bullet r_{K^{\revsn }}^{\circ }\right) \circ
J\right)  \\
\left( \zeta _{K,K^{\revsn },I,I}\circ J\right) \left( \left( K\bullet K^{\revsn
}\right) \circ \Delta _{I}^{\bullet }\circ J\right) \left( \left( K\bullet
K^{\revsn }\right) \circ \left( l_{J}^{\circ }\right) ^{-1}\right) \left(
\left( K\bullet \Sigma _{K}\right) \Delta _{K}^{\bullet }\circ \varepsilon
_{H}^{\bullet }\right)
\end{array}%
\right]  \\
&\overset{\eqref{eq:psi1}}{=}&\left[
\begin{array}{c}
\left( \left( \left( m_{K}^{\circ }\circ H\right) \left( K\circ \left(
K\circ u_{H}^{\circ }\right) \left( r_{K}^{\circ }\right) ^{-1}l_{K}^{\circ
}\right) \bullet I\right) \circ J\right)  \\
\left( \left( \left( K\circ I\circ K\right) \bullet I\right) \circ
m_{J}^{\circ }\right) \left( \varphi _{K\circ I,K}\circ J\right) \left(
\left( \left( K\circ I\right) \bullet r_{K^{\revsn }}^{\circ }\right) \circ
J\right)  \\
\left( \zeta _{K,K^{\revsn },I,I}\circ J\right) \left( \left( K\bullet K^{\revsn
}\right) \circ \Delta _{I}^{\bullet }\circ J\right) \left( \left(
r_{K\bullet K^{\revsn }}^{\circ }\right) ^{-1}\circ J\right) \left( \left(
K\bullet \Sigma _{K}\right) \Delta _{K}^{\bullet }\circ \varepsilon
_{H}^{\bullet }\right)
\end{array}%
\right]  \\
&=&\left[
\begin{array}{c}
\left( \left( \left( m_{K}^{\circ }\circ H\right) \left( K\circ \left(
K\circ u_{H}^{\circ }\right) \left( r_{K}^{\circ }\right) ^{-1}\right)
\bullet I\right) \circ J\right) \left( \left( \left( K\circ l_{K}^{\circ
}\right) \bullet I\right) \circ J\right)  \\
\left( \left( \left( K\circ I\circ K\right) \bullet I\right) \circ
m_{J}^{\circ }\right) \left( \varphi _{K\circ I,K}\circ J\right) \left(
\left( \left( K\circ I\right) \bullet r_{K^{\revsn }}^{\circ }\right) \circ
J\right)  \\
\left( \zeta _{K,K^{\revsn },I,I}\circ J\right) \left( \left( K\bullet K^{\revsn
}\right) \circ \Delta _{I}^{\bullet }\circ J\right) \left( \left(
r_{K\bullet K^{\revsn }}^{\circ }\right) ^{-1}\circ J\right) \left( \left(
K\bullet \Sigma _{K}\right) \Delta _{K}^{\bullet }\circ \varepsilon
_{H}^{\bullet }\right)
\end{array}%
\right]  \\
&=&\left[
\begin{array}{c}
\left( \left( \left( m_{K}^{\circ }\circ H\right) \left( K\circ K\circ
u_{H}^{\circ }\right) \bullet I\right) \circ J\right) \left( \left( \left(
K\circ \left( r_{K}^{\circ }\right) ^{-1}\right) \bullet I\right) \circ
J\right) \left( \left( \left( r_{K}^{\circ }\circ K\right) \bullet I\right)
\circ J\right)  \\
\left( \left( \left( K\circ I\circ K\right) \bullet I\right) \circ
m_{J}^{\circ }\right) \left( \varphi _{K\circ I,K}\circ J\right) \left(
\left( \left( K\circ I\right) \bullet r_{K^{\revsn }}^{\circ }\right) \circ
J\right)  \\
\left( \zeta _{K,K^{\revsn },I,I}\circ J\right) \left( \left( K\bullet K^{\revsn
}\right) \circ \Delta _{I}^{\bullet }\circ J\right) \left( \left(
r_{K\bullet K^{\revsn }}^{\circ }\right) ^{-1}\circ J\right) \left( \left(
K\bullet \Sigma _{K}\right) \Delta _{K}^{\bullet }\circ \varepsilon
_{H}^{\bullet }\right)
\end{array}%
\right]  \\
&=&\left[
\begin{array}{c}
\left( \left( \left( m_{K}^{\circ }\circ u_{H}^{\circ }\right) \bullet
I\right) \circ J\right) \left( \left( \left( r_{K\circ K}^{\circ }\right)
^{-1}\bullet I\right) \circ J\right) \left( \left( \left( K\circ K\right)
\bullet I\right) \circ m_{J}^{\circ }\right) \left( \varphi _{K,K}\circ
J\right)  \\
\left( \underbracket[0.140ex]{\left( r_{K}^{\circ }\bullet r_{K^{\revsn }}^{\circ
}\right) }\circ J\right) \left( \underbracket[0.140ex]{\zeta _{K,K^{\revsn },I,I}}\circ
J\right) \left( \underbracket[0.140ex]{\left( K\bullet K^{\revsn }\right) \circ \Delta
_{I}^{\bullet }}\circ J\right) \left( \underbracket[0.140ex]{\left( r_{K\bullet K^{\revsn
}}^{\circ }\right) ^{-1}}\circ J\right) \left( \left( K\bullet \Sigma
_{K}\right) \Delta _{K}^{\bullet }\circ \varepsilon _{H}^{\bullet }\right)
\end{array}%
\right]  \\
&\overset{\eqref{eq:unit1}}{=}&\left( \left( \left( m_{K}^{\circ }\circ u_{H}^{\circ }\right) \bullet
I\right) \circ J\right) \left( \left( \left( r_{K\circ K}^{\circ }\right)
^{-1}\bullet I\right) \circ J\right) \left( \left( \left( K\circ K\right)
\bullet I\right) \circ m_{J}^{\circ }\right) \left( \varphi _{K,K}\left(
K\bullet \Sigma _{K}\right) \Delta _{K}^{\bullet }\circ \varepsilon
_{H}^{\bullet }\right)  \\
&=&\left( \left( \left( K\circ u_{H}^{\circ }\right) \bullet I\right) \circ
J\right) \left( \left( \left( r_{K}^{\circ }\right) ^{-1}\bullet I\right)
\circ J\right) \left( \left( K\bullet I\right) \circ m_{J}^{\circ }\right)
\left( \underbracket[0.140ex]{\left( \left( m_{K}^{\circ }\bullet I\right) \circ
J\right) \varphi _{K,K}\left( K\bullet \Sigma _{K}\right) \Delta
_{K}^{\bullet }}\circ \varepsilon _{H}^{\bullet }\right)  \\
&=&\left( \left( \left( K\circ u_{H}^{\circ }\right) \bullet I\right) \circ
J\right) \left( \left( \left( r_{K}^{\circ }\right) ^{-1}\bullet I\right)
\circ J\right) \left( \left( K\bullet I\right) \circ m_{J}^{\circ }\right)
\left( \left( \mathrm{Id}_{K}\star \Sigma _{K}\right) \circ \varepsilon
_{H}^{\bullet }\right)  \\
&=&\left( \left( \left( K\circ u_{H}^{\circ }\right) \bullet I\right) \circ
J\right) \left( \left( \left( r_{K}^{\circ }\right) ^{-1}\bullet I\right)
\circ J\right) \left( \left( K\bullet I\right) \circ m_{J}^{\circ }\right)
\left( \left( \left( u_{K}^{\circ }\bullet I\right) \Delta _{I}^{\bullet
}\circ J\right) \left( l_{J}^{\circ }\right) ^{-1}\varepsilon _{K}^{\bullet
}\circ \varepsilon _{H}^{\bullet }\right)  \\
&=&\left( \left( \left( u_{K}^{\circ }\circ u_{H}^{\circ }\right) \bullet
I\right) \circ J\right) \left( \left( \left( r_{I}^{\circ }\right)
^{-1}\bullet I\right) \circ J\right) \left( \left( I\bullet I\right) \circ
m_{J}^{\circ }\right) \left( \left( \Delta _{I}^{\bullet }\circ J\right)
\left( l_{J}^{\circ }\right) ^{-1}\varepsilon _{K}^{\bullet }\circ
\varepsilon _{H}^{\bullet }\right)  \\
&=&\left( \left( \left( u_{K}^{\circ }\circ u_{H}^{\circ }\right) \bullet
I\right) \circ J\right) \left( \left( \left( r_{I}^{\circ }\right)
^{-1}\bullet I\right) \circ J\right) \left( \left( I\bullet I\right) \circ
m_{J}^{\circ }\right) \left( \Delta _{I}^{\bullet }\circ J\circ J\right)
\left( \left( l_{J}^{\circ }\right) ^{-1}\circ J\right) \left( \varepsilon
_{K}^{\bullet }\circ \varepsilon _{H}^{\bullet }\right)  \\
&=&\left( \left( \left( u_{K}^{\circ }\circ u_{H}^{\circ }\right) \bullet
I\right) \circ J\right) \left( \left( \left( m_{I}^{\circ }\right)
^{-1}\bullet I\right) \circ J\right) \left( \Delta _{I}^{\bullet }\circ
J\right) \left( I\circ m_{J}^{\circ }\right) \left( l_{J\circ J}^{\circ
}\right) ^{-1}\left( \varepsilon _{K}^{\bullet }\circ \varepsilon
_{H}^{\bullet }\right)  \\
&=&\left( \left( \left( u_{K}^{\circ }\circ u_{H}^{\circ }\right) \left(
m_{I}^{\circ }\right) ^{-1}\bullet I\right) \Delta _{I}^{\bullet }\circ
J\right) \left( l_{J}^{\circ }\right) ^{-1}m_{J}^{\circ }\left( \varepsilon
_{K}^{\bullet }\circ \varepsilon _{H}^{\bullet }\right)  \\
&=&\left( \left( u_{K\circ _{\psi }H}^{\circ }\bullet I\right) \Delta
_{I}^{\bullet }\circ J\right) \left( l_{J}^{\circ }\right) ^{-1}\varepsilon
_{K\circ _{\psi }H}^{\bullet }.
\end{eqnarray*}
\end{invisible} 
Therefore, $K\circ_{\psi}H$ is in $\mathsf{Hopf}(\Cc,\circ,\bullet)$.

Assuming now the existence of a morphism of comonoids $\varepsilon^\circ_K:K\to I$, then $u^\circ_H\varepsilon^\circ_K:K\to H$ is a morphism of comonoids into the Hopf monoid $H$ so that $\left( u_{H}^{\circ }\varepsilon _{K}^{\circ }\right)
^{r}$ is invertible and we get
\begin{eqnarray*}
\psi  &=&\left( \left( u_{H}^{\circ }\varepsilon _{K}^{\circ }\right)
^{r}\right) ^{-1}\underbracket[0.140ex]{\left( u_{H}^{\circ }\varepsilon _{K}^{\circ
}\right) ^{r}}\psi  
=\left( \left( u_{H}^{\circ }\varepsilon _{K}^{\circ }\right) ^{r}\right)
^{-1}\left( \left( K\circ \varepsilon _{H}^{\bullet }\right) \bullet
m_{H}^{\circ }\left( u_{H}^{\circ }\varepsilon _{K}^{\circ }\circ H\right)
\right) \Delta _{K\circ H}^{\bullet }\psi  \\
&=&\left( \left( u_{H}^{\circ }\varepsilon
_{K}^{\circ }\right) ^{r}\right) ^{-1}\left( \left( K\circ \varepsilon
_{H}^{\bullet }\right) \bullet l_{H}^{\circ }\left( \varepsilon _{K}^{\circ
}\circ H\right) \right) \underbracket[0.140ex]{\Delta _{K\circ H}^{\bullet }\psi}  \\
&\overset{\eqref{eq:psi3}}{=}&\left( \left( u_{H}^{\circ }\varepsilon _{K}^{\circ }\right) ^{r}\right)
^{-1}\left( \underbracket[0.140ex]{\left( K\circ \varepsilon _{H}^{\bullet }\right) \psi} \bullet
\underbracket[0.140ex]{l_{H}^{\circ }\left( \varepsilon _{K}^{\circ }\circ H\right) \psi} \right)
\Delta _{H\circ K}^{\bullet } =
\left( \left( u_{H}^{\circ }\varepsilon _{K}^{\circ }\right) ^{r}\right)
^{-1}\left( \chi \bullet \triangleleft \right) \Delta _{H\circ
K}^{\bullet }.
\end{eqnarray*}
The last part of the statement follows by applying \eqref{eq:psi1} and \eqref{eq:psi2} together with the definition of $\chi$ and $\triangleleft$.
\begin{invisible}
    We have%
\begin{eqnarray*}
\chi \left( m_{H}^{\circ }\circ K\right)  &=&\left( K\circ \varepsilon
_{H}^{\bullet }\right) \psi \left( m_{H}^{\circ }\circ K\right)  \\
&\overset{\eqref{eq:psi1}}{=}&\left( K\circ \varepsilon _{H}^{\bullet }\right)
\left( K\circ m_{H}^{\circ }\right) \left( \psi \circ H\right) \left( H\circ
\psi \right)  \\
&=&\left( K\circ m_{J}^{\circ }\right) \left( K\circ \varepsilon
_{H}^{\bullet }\circ \varepsilon _{H}^{\bullet }\right) \left( \psi \circ
H\right) \left( H\circ \psi \right)  \\
&=&\left( K\circ m_{J}^{\circ }\right) \left( \chi \circ J\right) \left(
H\circ K\circ \varepsilon _{H}^{\bullet }\right) \left( H\circ \psi \right)
=\left( K\circ m_{J}^{\circ }\right) \left( \chi \circ J\right) \left(
H\circ \chi \right)
\end{eqnarray*}%
\begin{equation*}
\chi \left( u_{H}^{\circ }\circ K\right) \left( l_{K}^{\circ }\right)
^{-1}=\left( K\circ \varepsilon _{H}^{\bullet }\right) \psi \left(
u_{H}^{\circ }\circ K\right) \left( l_{K}^{\circ }\right) ^{-1}\overset{\eqref%
{eq:psi1}}{=}\left( K\circ \varepsilon _{H}^{\bullet }\right) \left( K\circ
u_{H}^{\circ }\right) \left( r_{K}^{\circ }\right) ^{-1}=\left( K\circ
u_{J}^{\circ }\right) \left( r_{K}^{\circ }\right) ^{-1}
\end{equation*}%
\begin{eqnarray*}
\chi \left( H\circ m_{K}^{\circ }\right)  &=&\left( K\circ \varepsilon
_{H}^{\bullet }\right) \psi \left( H\circ m_{K}^{\circ }\right)  \\
&\overset{\eqref{eq:psi2}}{=}&\left( K\circ \varepsilon _{H}^{\bullet }\right)
\left( m_{K}^{\circ }\circ H\right) \left( K\circ \psi \right) \left( \psi
\circ K\right)  \\
&=&\left( m_{K}^{\circ }\circ J\right) \left( K\circ K\circ \varepsilon
_{H}^{\bullet }\right) \left( K\circ \psi \right) \left( \psi \circ K\right)
\\
&=&\left( m_{K}^{\circ }\circ J\right) \left( K\circ \chi \right) \left(
\psi \circ K\right)
\end{eqnarray*}%
\begin{equation*}
\chi \left( H\circ u_{K}^{\circ }\right) \left( r_{H}^{\circ }\right)
^{-1}=\left( K\circ \varepsilon _{H}^{\bullet }\right) \psi \left( H\circ
u_{K}^{\circ }\right) \left( r_{H}^{\circ }\right) ^{-1}\overset{\eqref%
{eq:psi2}}{=}\left( K\circ \varepsilon _{H}^{\bullet }\right) \left(
u_{K}^{\circ }\circ H\right) \left( l_{H}^{\circ }\right) ^{-1}=\left(
u_{K}^{\circ }\circ \varepsilon _{H}^{\bullet }\right) \left( l_{H}^{\circ
}\right) ^{-1}.
\end{equation*}%
Similarly, we have

\begin{eqnarray*}
\triangleleft \left( m_{H}^{\circ }\circ K\right)  &=&l_{H}^{\circ }\left(
\varepsilon _{K}^{\circ }\circ H\right) \psi \left( m_{H}^{\circ }\circ
K\right)  \\
&\overset{\eqref{eq:psi1}}{=}&l_{H}^{\circ }\left( \varepsilon _{K}^{\circ
}\circ H\right) \left( K\circ m_{H}^{\circ }\right) \left( \psi \circ
H\right) \left( H\circ \psi \right)  \\
&=&l_{H}^{\circ }\left( I\circ m_{H}^{\circ }\right) \left( \varepsilon
_{K}^{\circ }\circ H\circ H\right) \left( \psi \circ H\right) \left( H\circ
\psi \right)  \\
&=&m_{H}^{\circ }l_{H\circ H}^{\circ }\left( \varepsilon _{K}^{\circ }\circ
H\circ H\right) \left( \psi \circ H\right) \left( H\circ \psi \right)  \\
&=&m_{H}^{\circ }\left( l_{H}^{\circ }\circ H\right) \left( \varepsilon
_{K}^{\circ }\circ H\circ H\right) \left( \psi \circ H\right) \left( H\circ
\psi \right)  \\
&=&m_{H}^{\circ }\left( \triangleleft \circ H\right) \left( H\circ \psi
\right)
\end{eqnarray*}%
\begin{equation*}
\triangleleft \left( u_{H}^{\circ }\circ K\right) \left( l_{K}^{\circ
}\right) ^{-1}=l_{H}^{\circ }\left( \varepsilon _{K}^{\circ }\circ H\right)
\psi \left( u_{H}^{\circ }\circ K\right) \left( l_{K}^{\circ }\right) ^{-1}%
\overset{\eqref{eq:psi1}}{=}l_{H}^{\circ }\left( \varepsilon _{K}^{\circ
}\circ H\right) \left( K\circ u_{H}^{\circ }\right) \left( r_{K}^{\circ
}\right) ^{-1}=u_{H}^{\circ }\varepsilon _{K}^{\circ }
\end{equation*}%
If $\varepsilon _{K}^{\circ }$ is a morphism of bimonoids, then
\begin{eqnarray*}
\triangleleft \left( H\circ m_{K}^{\circ }\right)  &=&l_{H}^{\circ }\left(
\varepsilon _{K}^{\circ }\circ H\right) \psi \left( H\circ m_{K}^{\circ
}\right)  \\
&\overset{\eqref{eq:psi2}}{=}&l_{H}^{\circ }\left( \varepsilon _{K}^{\circ
}\circ H\right) \left( m_{K}^{\circ }\circ H\right) \left( K\circ \psi
\right) \left( \psi \circ K\right)  \\
&=&l_{H}^{\circ }\left( m_{I}^{\circ }\circ H\right) \left( \varepsilon
_{K}^{\circ }\circ \varepsilon _{K}^{\circ }\circ H\right) \left( K\circ
\psi \right) \left( \psi \circ K\right)  \\
&=&l_{H}^{\circ }\left( I\circ l_{H}^{\circ }\right) \left( I\circ
\varepsilon _{K}^{\circ }\circ H\right) \left( I\circ \psi \right) \left(
\varepsilon _{K}^{\circ }\circ H\circ K\right) \left( \psi \circ K\right)  \\
&=&l_{H}^{\circ }\left( I\circ \triangleleft \right) \left( \varepsilon
_{K}^{\circ }\circ H\circ K\right) \left( \psi \circ K\right)  \\
&=&\triangleleft l_{H\circ K}^{\circ }\left( \varepsilon _{K}^{\circ }\circ
H\circ K\right) \left( \psi \circ K\right)  \\
&=&\triangleleft \left( l_{H}^{\circ }\circ K\right) \left( \varepsilon
_{K}^{\circ }\circ H\circ K\right) \left( \psi \circ K\right) =\triangleleft
\left( \triangleleft \circ K\right)
\end{eqnarray*}%
\begin{equation*}
\triangleleft \left( H\circ u_{K}^{\circ }\right) \left( r_{H}^{\circ
}\right) ^{-1}=l_{H}^{\circ }\left( \varepsilon _{K}^{\circ }\circ H\right)
\psi \left( H\circ u_{K}^{\circ }\right) \left( r_{H}^{\circ }\right) ^{-1}%
\overset{\eqref{eq:psi2}}{=}l_{H}^{\circ }\left( \varepsilon _{K}^{\circ
}\circ H\right) \left( u_{K}^{\circ }\circ H\right) \left( l_{H}^{\circ
}\right) ^{-1}=\mathrm{Id}_{H}.
\end{equation*}%
\end{invisible}
\end{proof}

In the following, we consider duoidal categories which are braided with respect to one of the two tensor products.

\section{Results on tensor-braided duoidal categories} 
\label{sec:tens-br}
We recall, from \cite[Definition 6.5]{Aguiar}, the notion of $\bullet$\,-\textit{braided} duoidal category. This is a duoidal category $(\mathcal{C}, \circ, I, \bullet, J)$ with $(\Cc,\bullet,J)$ braided such that the braiding $\sigma_{X,Y}:X\bullet Y\to Y\bullet X$, for $X,Y$ in $\Cc$, satisfies the following equalities:
\[
\sigma_{A\circ C,B\circ D}\zeta_{A,B,C,D}=\zeta_{B,A,D,C}(\sigma_{A,B}\circ\sigma_{C,D}),\qquad \sigma_{I,I}\Delta_{I}^{\bullet}=\Delta_{I}^{\bullet}.
\]
In such a duoidal category one can consider a cocommutative bimonoid, i.e.\ a bimonoid $(H,m_H^\circ,u_H^\circ,\Delta_H^\bullet,\varepsilon_H^\bullet)$ with $(H,\Delta_H^\bullet,\varepsilon_H^\bullet)$ in $\mathsf{Comon}_{\mathrm{coc}}(\Cc^{\bullet})$ i.e. such that $\sigma_{H,H}\Delta_H^\bullet=\Delta_H^\bullet$, see \cite[Definition 6.31]{Aguiar}.

Then one gets the category $\cBimon(\mathcal{C}, \circ, \bullet)$ of cocommutative bimonoids.%, or simply denoted by $\cBimon(\mathcal{C})$.

Cocommutative bimonoids which are Hopf monoids constitute a full subcategory of $\cBimon(\Cc,\circ,\bullet)$ that we denote  by $\cHopfmon(\Cc,\circ,\bullet)$. %or simply by $\cHopfmon(\Cc)$.

\begin{remark}\label{rmk:Bimonmonoidal}
By \cite[Proposition 6.37]{Aguiar}, if $(\Cc,\circ,\bullet)$ is a $\bullet $-braided duoidal
category then $\left(\Mon\left( \Cc^\circ\right) ,\bullet ,J\right) $ is braided monoidal. Since $\Bimon\left( \mathcal{C},\circ ,%
\mathcal{\bullet }\right) =\Comon\left(\Mon\left( \mathcal{C}^{\circ }\right)
,\bullet \right) $ we get that $\left(\Bimon\left( \Cc,\circ ,%
\bullet \right) ,\bullet ,J\right) $ is monoidal, see e.g. \cite[page 9]{Aguiar}. Explicitly, if $H
$ and $K$ are bimonoids, then so is $H\bullet K$ with the following structures:%
\begin{eqnarray*}
m_{H\bullet K}^{\circ } &:&
\xymatrix@C=1.5cm{
\left( H\bullet K\right) \circ \left( H\bullet
K\right) \ar[r]^-{\zeta _{H,K,H,K}}&\left( H\circ H\right)
\bullet \left( K\circ K\right) \ar[r]^-{m_{H}^{\circ }\bullet m_{K}^{\circ }}%
&H\bullet K} \\
u_{H\bullet K}^{\circ } &:&
\xymatrix@C=1.3cm{I\ar[r]^-{\Delta _{I}^{\bullet }}&I\bullet I\ar[r]^-{u_{H}^{\circ }\bullet u_{K}^{\circ }}&H\bullet K} \\
\Delta _{H\bullet K}^{\bullet } &:&
\xymatrix@C=1.5cm{H\bullet K\ar[r]^-{\Delta _{H}^{\bullet
}\bullet \Delta _{K}^{\bullet }}&H\bullet H\bullet K\bullet
K\ar[r]^-{H\bullet \sigma _{H,K}\bullet K}&\left( H\bullet
K\right) \bullet \left( H\bullet K\right)}  \\
\varepsilon _{H\bullet K}^{\bullet } &:&
\xymatrix@C=1.3cm{H\bullet K\ar[r]^-{\varepsilon
_{H}^{\bullet }\bullet \varepsilon _{K}^{\bullet }}&
J\bullet J\ar[r]^-{m_{J}^{\bullet }}&J}.
\end{eqnarray*}
\end{remark}

We are now able to prove that, in a $\bullet$-braided duoidal category with a reversion, $A\bullet B$ is a Hopf
monoid whenever $A$ and $B$ are Hopf monoids.

\begin{proposition}%[\rd{update:2026/03/10}]
\label{prop:Hopfbul}
Let $\left( \mathcal{C},\circ ,I,\bullet ,J\right) $  be a $\bullet$-braided duoidal category with a reversion.
If $A$ and $B$ are Hopf monoids, then so is $A\bullet B$.
As a consequence, $(\Hopfmon(\Cc,\circ,\bullet),\bullet,J)$ is a monoidal category.
\end{proposition}

\begin{proof}
By \cref{rmk:Bimonmonoidal}, given $A,B$ in $\mathsf{Hopf}(\Cc,\circ,\bullet)$, we already know that $A\bullet B$ is in $\Bimon(\Cc,\circ,\bullet)$. We define
\begin{equation*}
\Sigma _{A\bullet B}:\xymatrix{
A \bullet B
  \ar[r]^{\sigma_{A,B}}
& B \bullet A
  \ar[r]^{\Sigma_B \bullet \Sigma_A}
& B^{\revsn} \bullet A^{\revsn}
  \ar[r]^{\phi_{B,A}^{\bullet}}
& (A \bullet B)^{\revsn}}
\end{equation*}%
and we check that this is an antipode.
\begin{eqnarray*}
\mathrm{Id}_{A\bullet B}\star \Sigma _{A\bullet B} &=&\left( \left(
m_{A\bullet B}^{\circ }\bullet I\right) \circ J\right) \varphi _{A\bullet
B,A\bullet B}\left( A\bullet B\bullet \Sigma _{A\bullet B}\right) \Delta
_{A\bullet B}^{\bullet } \\
&=&
\left( \left( m_{A}^{\circ }\bullet m_{B}^{\circ }\bullet I\right) \circ
J\right) \left( \left( \zeta _{A,B,A,B}\bullet I\right) \circ J\right)
\varphi _{A\bullet B,A\bullet B}\left( A\bullet B\bullet \phi
_{B,A}^{\bullet }\right)  \\
&&\left( A\bullet B\bullet \Sigma _{B}\bullet \Sigma _{A}\right) \left(
A\bullet B\bullet \sigma _{A,B}\right) \left( A\bullet \sigma _{A,B}\bullet
B\right) \left( \Delta _{A}^{\bullet }\bullet \Delta _{B}^{\bullet }\right)  \\
&=&
\left( \left( m_{A}^{\circ }\bullet m_{B}^{\circ }\bullet I\right) \circ
J\right) \underbracket[0.140ex]{\left( \left( \zeta _{A,B,A,B}\bullet I\right) \circ
J\right) \varphi _{A\bullet B,A\bullet B}\left( A\bullet B\bullet \phi
_{B,A}^{\bullet }\right) } \\
&&\left( A\bullet B\bullet \Sigma _{B}\bullet \Sigma _{A}\right) \left(
A\bullet \sigma _{A,B\bullet B}\right) \left( \Delta _{A}^{\bullet }\bullet
\Delta _{B}^{\bullet }\right)\\
&\overset{\eqref{form:vartheta1}}{=}&
\left( \left( m_{A}^{\circ }\bullet m_{B}^{\circ }\bullet I\right) \circ
J\right) \vartheta _{A,B\circ B,A}^{\prime }\left( A\bullet \varphi
_{B,B}\bullet A^{\revsn }\right)  \\
&&\left( A\bullet \left( B\bullet \Sigma _{B}\right) \Delta _{B}^{\bullet
}\bullet \Sigma _{A}\right) \left( A\bullet \sigma _{A,B}\right) \left(
\Delta _{A}^{\bullet }\bullet B\right)\\
&=&
\left( \left( m_{A}^{\circ }\bullet B\bullet I\right) \circ J\right)
\vartheta _{A,B,A}^{\prime }\left( A\bullet \left( \underbracket[0.140ex]{\left(m_{B}^{\circ
}\bullet I\right) \circ J}\right) \bullet A^{\revsn }\right) \left( A\bullet
\underbracket[0.140ex]{\varphi _{B,B}}\bullet A^{\revsn }\right)  \\
&&\left( A\bullet \underbracket[0.140ex]{\left( B\bullet \Sigma _{B}\right) \Delta _{B}^{\bullet
}}\bullet \Sigma _{A}\right) \left( A\bullet \sigma _{A,B}\right) \left(
\Delta _{A}^{\bullet }\bullet B\right)\\
&=&\left( \left( m_{A}^{\circ }\bullet B\bullet I\right) \circ J\right)
\vartheta _{A,B,A}^{\prime }\left( A\bullet \left( \mathrm{Id}_{B}\star
\Sigma _{B}\right) \bullet \Sigma _{A}\right) \left( A\bullet \sigma
_{A,B}\right) \left( \Delta _{A}^{\bullet }\bullet B\right)  \\
&=&
\left( \left( m_{A}^{\circ }\bullet B\bullet I\right) \circ J\right)
\vartheta _{A,B,A}^{\prime }\left( A\bullet \left( \left( u_{B}^{\circ
}\bullet I\right) \circ J\right) \bullet A^{\revsn }\right)  \\
&&\left( A\bullet \left( \Delta _{I}^{\bullet }\circ J\right) \bullet A^{\revsn
}\right) \left( A\bullet \left( l_{J}^{\circ }\right) ^{-1}\varepsilon
_{B}^{\bullet }\bullet \Sigma _{A}\right) \left( A\bullet \sigma
_{A,B}\right) \left( \Delta _{A}^{\bullet }\bullet B\right)\\
&=&\left( \left( m_{A}^{\circ }\bullet u_{B}^{\circ }\bullet I\right) \circ
J\right) \underbracket[0.140ex]{\vartheta _{A,I,A}^{\prime }\left( A\bullet \left(
\Delta _{I}^{\bullet }\circ J\right) \bullet A^{\revsn }\right) }\\
&&\left(A\bullet \left( l_{J}^{\circ }\right) ^{-1}\bullet A^{\revsn }\right) \left(
A\bullet \sigma _{A^{\revsn },J}\right) \left( \left( A\bullet \Sigma
_{A}\right) \Delta _{A}^{\bullet }\bullet \varepsilon _{B}^{\bullet }\right)\\
&\overset{\eqref{form:vartheta1I}}{=}&\left( \left( m_{A}^{\circ }\bullet
u_{B}^{\circ }\bullet I\right) \circ J\right) \left( \left( \left( A\circ
A\right) \bullet \Delta _{I}^{\bullet }\right) \circ J\right) \varphi
_{A,A}\left( A\bullet l_{A^{\revsn }}^{\bullet }\right)\\ 
&&\left( A\bullet \left(
l_{A^{\revsn }}^{\bullet }\right) ^{-1}r_{A^{\revsn }}^{\bullet }\right) \left(
\left( A\bullet \Sigma _{A}\right) \Delta _{A}^{\bullet }\bullet J\right)
\left( A\bullet \varepsilon _{B}^{\bullet }\right)\\
&=&\left( \left( m_{A}^{\circ }\bullet \left( u_{B}^{\circ }\bullet I\right)
\Delta _{I}^{\bullet }\right) \circ J\right) \varphi _{A,A}r_{A\bullet
A^{\revsn }}^{\bullet }\left( \left( A\bullet \Sigma _{A}\right) \Delta
_{A}^{\bullet }\bullet J\right) \left( A\bullet \varepsilon _{B}^{\bullet
}\right)  \\
&=&\left( \left( A\bullet \left( u_{B}^{\circ }\bullet I\right) \Delta
_{I}^{\bullet }\right) \circ J\right) \underbracket[0.140ex]{\left( \left( m_{A}^{\circ
}\bullet I\right) \circ J\right) \varphi _{A,A}\left( A\bullet \Sigma
_{A}\right) \Delta _{A}^{\bullet }}r_{A}^{\bullet }\left( A\bullet
\varepsilon _{B}^{\bullet }\right)  \\
&=&\left( \left( A\bullet \left( u_{B}^{\circ }\bullet I\right) \Delta
_{I}^{\bullet }\right) \circ J\right) \left( \mathrm{Id}_{A}\star \Sigma
_{A}\right) r_{A}^{\bullet }\left( A\bullet \varepsilon _{B}^{\bullet
}\right)  \\
&=&\left( \left( A\bullet \left( u_{B}^{\circ }\bullet I\right) \Delta
_{I}^{\bullet }\right) \circ J\right) \left( \left( u_{A}^{\circ }\bullet
I\right) \Delta _{I}^{\bullet }\circ J\right) \left( l_{J}^{\circ }\right)
^{-1}\varepsilon _{A}^{\bullet }r_{A}^{\bullet }\left( A\bullet \varepsilon
_{B}^{\bullet }\right)  \\
&=&\left( \left( \left( u_{A}^{\circ }\bullet u_{B}^{\circ }\right) \Delta
_{I}^{\bullet }\bullet I\right) \circ J\right) \left( \Delta _{I}^{\bullet
}\circ J\right) \left( l_{J}^{\circ }\right) ^{-1}r_{J}^{\bullet }\left(
\varepsilon _{A}^{\bullet }\bullet \varepsilon _{B}^{\bullet }\right)  \\
&=&\left( \left( u_{A\bullet B}^{\circ }\bullet I\right) \circ J\right)
\left( \Delta _{I}^{\bullet }\circ J\right) \left( l_{J}^{\circ }\right)
^{-1}\varepsilon _{A\bullet B}^{\bullet }=1_{A\bullet B,A\bullet B}^{l}
\end{eqnarray*}%
In a similar way, using \eqref{form:vartheta} and \eqref{form:varthetaI}, one proves that
$\Sigma _{A\bullet B}\star \mathrm{Id}_{A\bullet B}= 1_{A\bullet B,A\bullet B}^{r}$.
\begin{invisible}
\begin{eqnarray*}
\Sigma _{A\bullet B}\star \mathrm{Id}_{A\bullet B} &=&\left( J\circ \left(
I\bullet m_{A\bullet B}^{\circ }\right) \right) \psi _{A\bullet B,A\bullet
B}\left( \Sigma _{A\bullet B}\bullet A\bullet B\right) \Delta _{A\bullet
B}^{\bullet } \\
&=&\left[
\begin{array}{c}
\left( J\circ \left( I\bullet m_{A}^{\circ }\bullet m_{B}^{\circ }\right)
\right) \left( J\circ \left( I\bullet \zeta _{A,B,A,B}\right) \right) \psi
_{A\bullet B,A\bullet B}\left( \phi _{B,A}^{\bullet }\bullet A\bullet
B\right)  \\
\left( \Sigma _{B}\bullet \Sigma _{A}\bullet A\bullet B\right) \left( \sigma
_{A,B}\bullet A\bullet B\right) \left( A\bullet \sigma _{A,B}\bullet
B\right) \left( \Delta _{A}^{\bullet }\bullet \Delta _{B}^{\bullet }\right)
\end{array}%
\right]  \\
&=&\left[
\begin{array}{c}
\left( J\circ \left( I\bullet m_{A}^{\circ }\bullet m_{B}^{\circ }\right)
\right) \underbracket[0.140ex]{\left( J\circ \left( I\bullet \zeta _{A,B,A,B}\right)
\right) \psi _{A\bullet B,A\bullet B}\left( \phi _{B,A}^{\bullet }\bullet
A\bullet B\right) } \\
\left( \Sigma _{B}\bullet \Sigma _{A}\bullet A\bullet B\right) \left( \sigma
_{A\bullet A,B}\bullet B\right) \left( \Delta _{A}^{\bullet }\bullet \Delta
_{B}^{\bullet }\right)
\end{array}%
\right]  \\
&\overset{\eqref{form:vartheta}}{=}&\left[
\begin{array}{c}
\left( J\circ \left( I\bullet m_{A}^{\circ }\bullet m_{B}^{\circ }\right)
\right) \vartheta _{B,A\circ A,B}\left( B^{\revsn }\bullet \psi _{A,A}\bullet
B\right)  \\
\left( \Sigma _{B}\bullet \Sigma _{A}\bullet A\bullet B\right) \left(
B\bullet \Delta _{A}^{\bullet }\bullet B\right) \left( \sigma _{A,B}\bullet
B\right) \left( A\bullet \Delta _{B}^{\bullet }\right)
\end{array}%
\right]  \\
&=&\left[
\begin{array}{c}
\left( J\circ \left( I\bullet A\bullet m_{B}^{\circ }\right) \right)
\vartheta _{B,A,B}\left( B^{\revsn }\bullet \left( J\circ \left( I\bullet
m_{A}^{\circ }\right) \right) \bullet B\right) \left( B^{\revsn }\bullet \psi
_{A,A}\bullet B\right)  \\
\left( \Sigma _{B}\bullet \Sigma _{A}\bullet A\bullet B\right) \left(
B\bullet \Delta _{A}^{\bullet }\bullet B\right) \left( \sigma _{A,B}\bullet
B\right) \left( A\bullet \Delta _{B}^{\bullet }\right)
\end{array}%
\right]  \\
&=&\left( J\circ \left( I\bullet A\bullet m_{B}^{\circ }\right) \right)
\vartheta _{B,A,B}\left( \Sigma _{B}\bullet \left( \Sigma _{A}\star \mathrm{Id%
}_{A}\right) \bullet B\right) \left( \sigma _{A,B}\bullet B\right) \left(
A\bullet \Delta _{B}^{\bullet }\right)  \\
&=&\left[
\begin{array}{c}
\left( J\circ \left( I\bullet A\bullet m_{B}^{\circ }\right) \right)
\vartheta _{B,A,B}\left( B^{\revsn }\bullet \left( J\circ \left( I\bullet
u_{A}^{\circ }\right) \right) \bullet B\right)  \\
\left( B^{\revsn }\bullet \left( J\circ \Delta _{I}^{\bullet }\right) \bullet
B\right) \left( B^{\revsn }\bullet \left( r_{J}^{\circ }\right)
^{-1}\varepsilon _{A}^{\bullet }\bullet B\right) \left( \sigma _{A,B^{\revsn
}}\bullet B\right) \left( A\bullet \left( \Sigma _{B}\bullet B\right) \Delta
_{B}^{\bullet }\right)
\end{array}%
\right]  \\
&=&\left( J\circ \left( I\bullet u_{A}^{\circ }\bullet m_{B}^{\circ }\right)
\right) \underbracket[0.140ex]{\vartheta _{B,I,B}\left( B^{\revsn }\bullet \left( J\circ
\Delta _{I}^{\bullet }\right) \bullet B\right) }\left( B^{\revsn }\bullet
\left( r_{J}^{\circ }\right) ^{-1}\bullet B\right) \left( \sigma _{J,B^{\revsn
}}\bullet B\right) \left( \varepsilon _{A}^{\bullet }\bullet \left( \Sigma
_{B}\bullet B\right) \Delta _{B}^{\bullet }\right)  \\
&\overset{\eqref{form:varthetaI}}{=}&\left( J\circ \left( I\bullet
u_{A}^{\circ }\bullet m_{B}^{\circ }\right) \right) \left( J\circ \left(
\Delta _{I}^{\bullet }\bullet \left( B\circ B\right) \right) \right) \psi
_{B,B}\left( r_{B^{\revsn }}^{\bullet }\bullet B\right) \left( \left(
r_{B^{\revsn }}^{\bullet }\right) ^{-1}l_{B^{\revsn }}^{\bullet }\bullet
B\right) \left( \varepsilon _{A}^{\bullet }\bullet \left( \Sigma _{B}\bullet
B\right) \Delta _{B}^{\bullet }\right)  \\
&=&\left( J\circ \left( I\bullet u_{A}^{\circ }\bullet m_{B}^{\circ }\right)
\right) \left( J\circ \left( \Delta _{I}^{\bullet }\bullet \left( B\circ
B\right) \right) \right) \psi _{B,B}l_{B^{\revsn }\bullet B}^{\bullet }\left(
J\bullet \left( \Sigma _{B}\bullet B\right) \Delta _{B}^{\bullet }\right)
\left( \varepsilon _{A}^{\bullet }\bullet B\right)  \\
&=&\left( J\circ \left( \left( I\bullet u_{A}^{\circ }\right) \Delta
_{I}^{\bullet }\bullet B\right) \right) \left( J\circ \left( I\bullet
m_{B}^{\circ }\right) \right) \psi _{B,B}\left( \Sigma _{B}\bullet B\right)
\Delta _{B}^{\bullet }l_{B}^{\bullet }\left( \varepsilon _{A}^{\bullet
}\bullet B\right)  \\
&=&\left( J\circ \left( \left( I\bullet u_{A}^{\circ }\right) \Delta
_{I}^{\bullet }\bullet B\right) \right) \left( \Sigma _{B}\star \mathrm{Id}%
_{B}\right) l_{B}^{\bullet }\left( \varepsilon _{A}^{\bullet }\bullet
B\right)  \\
&=&\left( J\circ \left( \left( I\bullet u_{A}^{\circ }\right) \Delta
_{I}^{\bullet }\bullet B\right) \right) \left( J\circ \left( I\bullet
u_{B}^{\circ }\right) \Delta _{I}^{\bullet }\right) \left( r_{J}^{\circ
}\right) ^{-1}\varepsilon _{B}^{\bullet }l_{B}^{\bullet }\left( \varepsilon
_{A}^{\bullet }\bullet B\right)  \\
&=&\left( J\circ \left( I\bullet \left( u_{A}^{\circ }\bullet u_{B}^{\circ
}\right) \Delta _{I}^{\bullet }\right) \right) \left( J\circ \Delta
_{I}^{\bullet }\right) \left( r_{J}^{\circ }\right) ^{-1}l_{J}^{\bullet
}\left( \varepsilon _{A}^{\bullet }\bullet \varepsilon _{B}^{\bullet
}\right)  \\
&=&\left( J\circ \left( I\bullet u_{A\bullet B}^{\circ }\right) \right)
\left( J\circ \Delta _{I}^{\bullet }\right) \left( r_{J}^{\circ }\right)
^{-1}\varepsilon _{A\bullet B}^{\bullet }=1_{A\bullet B,A\bullet B}^{r}.
\end{eqnarray*}
\end{invisible}
We already proved in \cref{lem:IHopf} that $J$ is a Hopf monoid. The last part of the statement follows.
\end{proof}

\subsection{On pullbacks}

To construct kernels of projections, we first analyze appropriate pullbacks. 

Recall that the intersection of two subobjects of the same object can be realized as the pullback
of the inclusions, see \cite[Proposition 4.2.3]{BorI94}.

As in \cite[\S1]{P19}, given a monoidal category $(\Cc,\bullet,J)$, by saying that \emph{$\bullet$ preserves limits (of some type) in $\Cc$}, we mean that the functors $X\bullet
\left( -\right):\Cc\to \Cc $ and $\left( -\right) \bullet X:\Cc\to \Cc$ preserve these limits for every object $X$ in $\Cc$. We begin with a result in the context of braided monoidal categories.

\begin{lemma}%[\rd{Update:2026/02/27}]
\label{lem:ddpulb}Let $\left( \mathcal{C},\bullet ,J,\sigma\right) $ be a braided
monoidal category. Assume that $\mathcal{C}$ has and $\bullet$ preserves binary intersections.
Then every cospan
$\xymatrix{A\ar[r]^f&C&\,B\ar@{>->}[l]_m}$ in $\cComon \left( \mathcal{C}\right) ,$ with $m$ a monomorphism
in $\mathcal{C}$, has a pullback%
\begin{equation}
\label{def:Kpullbcspn}
\xymatrix@R=.7cm{K\pulb\ar[r]^k\ar[d]_{g}&A\ar[d]^f\\B\,\ar@{>->}[r]^m&C}
\end{equation}
in $\Comon \left( \mathcal{C}\right) $ and in $\cComon \left( \mathcal{C}\right) $ with $k$ monomorphism in $%
\mathcal{C}$. More precisely, $K$ arises as the following pullback in $\Cc$ and in $\Comon(\Cc)$
\begin{equation}
\label{def:Kpullb}
\xymatrix@R=.7cm{K\pulb\ar[r]^k\ar[d]_{\lambda}&A\,\ar[d]^{(A\bullet f)\Delta_A}\\A\bullet B\,\ar@{>->}[r]^{A\bullet m}&A\bullet C}
\end{equation}%
Moreover, we have the identities  $k=r_A(A\bullet\varepsilon_B)\lambda$, $g:=l_{B}\left( \varepsilon _{A}\bullet
B\right) \lambda $,
$\lambda=(k\bullet g)\Delta_K$, and the following diagram is an equalizer in $\Cc$ and in $\Comon(\Cc)$.
\begin{equation}
\label{def:Kequa}
    \xymatrix@C=2.5cm{
K \ar[r]^{\lambda} & A\bullet B\ar@<0.5ex>[r]^-{(A\bullet f)\Delta_A\bullet B} \ar@<-0.5ex>[r]_-{A\bullet \left(
m\bullet B\right) \Delta _{B}} & A\bullet C\bullet B}
\end{equation} 
\end{lemma}

\begin{proof}
Since $m$ is a monomorphism in $\Cc$ and $A\bullet \left( -\right) $ preserves
intersections, we get that $A\bullet m$ is a monomorphism (cf.\ \cite[Proposition
2.5.6]{BorI94}). Moreover $\rho _{A}:=\left( A\bullet f\right)
\Delta _{A}:A\to A\bullet C$ is a  split monomorphism in $\Cc$ as $r_{A}\left( A\bullet
\varepsilon _{C}\right) \rho _{A}=r_{A}\left( A\bullet \varepsilon _{C}\right) \left( A\bullet f\right)
\Delta _{A}=r_{A}\left( A\bullet \varepsilon
_{A}\right) \Delta _{A}=\mathrm{Id}_{A}$.

By assumption we can consider the pullback in $\mathcal{C}$ given in \eqref{def:Kpullb}
as it is a binary intersection. Observe that $k$ is a monomorphism in $\Cc$ since monomorphisms are always stable under pullbacks and $A\bullet m$ is a
monomorphism; indeed, for every object $X$ in $\Cc$, the functor $X\bullet (-)$ preserves monomorphisms as it preserves the kernel pair of any monomorphism and this is a binary intersection, see \cite[Proposition 2.5.6]{BorI94}.

One easily checks that the squares on the right in the following diagram are serially commutative.
\[
\xymatrix@C=2.5cm{
K\ar[d]^{k} \ar[r]^{\lambda} & A\bullet B\ar[d]^{A\bullet  m}\ar@<0.5ex>[r]^-{\rho _{A}\bullet B} \ar@<-0.5ex>[r]_-{A\bullet \left(
m\bullet B\right) \Delta _{B}} & A\bullet C\bullet B \ar[d]^{A\bullet C\bullet m}\\
A \ar[r]^{\rho _{A}} & A\bullet  C \ar@<0.5ex>[r]^-{\rho _{A}\bullet C} \ar@<-0.5ex>[r]_-{A\bullet \Delta_{C}} & A\bullet C\bullet C
}
\]
Since $A\bullet C\bullet m$ is a monomorphism, the squares on the right are serially commutative, the lower row is an equalizer in $\Cc$ and the left square is a pullback, we deduce that the upper row, i.e. \eqref{def:Kequa}, is an equalizer in $\Cc$.
\begin{invisible}
See Proposition 2.5 in \href{https://ncatlab.org/nlab/show/pasting+law+for+pullbacks}{ncatlab}. Anyway the proof is easy.
\end{invisible} In fact, it is a coreflexive equalizer as the parallel pair has a common retraction given by $r_A(A\bullet\varepsilon_C)\bullet B:A\bullet C\bullet B\to A\bullet B$. Now, the lower row is a split fork whence an absolute equalizer.
\begin{invisible}
Take $e:=\rho_A$, $f:=\rho_A\bullet C$, $g:=A\bullet\Delta_C$,
$s:=r_A(A\bullet \varepsilon_C)$ and $t:=A\bullet r_C(C\bullet\varepsilon_C)$. Then $se=\id$, $tg=\id$ and $tf=es$.
\end{invisible}
As a consequence, it is preserved by the functors $X\bullet (-)$ and $(-)\bullet X$. Since these functors preserve monomorphisms and binary intersections (so also the pullback in the left square), by the same argument used above we deduce that the upper row is an equalizer in $\Cc$ which is preserved by the same functors. Note also that the parallel pairs of the equalizer are morphisms in $\Comon(\Cc)$ as $A$ and $B$ are cocommutative and hence $\Delta_A$ and $\Delta_B$ are morphisms of comonoids. Therefore, by the same proof of \cite[Proposition A.9]{AM-MMCat}, the forgetful functor $\Comon(\Cc)\to \Cc$ creates this equalizer. In other words, $K$ carries a structure of comonoid $(K,\Delta_K,\varepsilon_K)$ that makes \eqref{def:Kequa} an equalizer in $\Comon(\Cc)$.

Again, since $A\bullet C\bullet m$ is a monomorphism, the squares on the right are serially commutative, the lower row is an equalizer in $\Comon(\Cc)$ and the upper row is an equalizer in $\Comon(\Cc)$, we deduce that the left square, i.e.\ \eqref{def:Kpullb}, is a pullback in $\Comon(\Cc)$, see the dual of \cite[5.8 Lemma]{Barr}.

We check that $K$ is cocommutative:
\begin{equation*}
\left( k\bullet k\right) \sigma _{K,K}\Delta _{K}=\sigma
_{A,A}\left( k\bullet k\right) \Delta _{K}=\sigma _{A,A}\Delta
_{A}k=\Delta _{A}k=\left( k\bullet k\right) \Delta
_{K}
\end{equation*}%
and hence $\sigma _{K,K}\Delta _{K}=\Delta _{K}.$

Set $g:=l_{B}\left( \varepsilon _{A}\bullet
B\right) \lambda :K\to B$. Then, we obtain%
\begin{eqnarray*}
mg &=&ml_{B}\left( \varepsilon _{A}\bullet
B\right) \lambda =l_{C}\left( \varepsilon _{A}\bullet
C\right) \left( A\bullet m\right) \lambda =l_{C}\left(
\varepsilon _{A}\bullet C\right) \rho _{A}k \\
&=&l_{C}\left( \varepsilon _{A}\bullet C\right) \left(
A\bullet f\right) \Delta _{A}k=fl_{A}\left(
\varepsilon _{A}\bullet A\right) \Delta _{A}k=fk
\end{eqnarray*}
so that we have the commutative diagram \eqref{def:Kpullbcspn}.
We prove that it is a pullback in $\Comon \left( \mathcal{C}\right) .$

First, note that $g$ is a morphism of comonoids as so are $l_B(\varepsilon_A\bullet B)$ and $\lambda$.

Let $u:D\rightarrow A,v:D\rightarrow B$ be morphisms in $\mathsf{Comon}(\Cc)$ such that $%
mv=fu.$ Then%
\begin{equation*}
\rho _{A}u=\left( A\bullet f\right) \Delta _{A}u=\left( A\bullet f\right) \left( u\bullet u\right) \Delta _{D}=\left( A\bullet m\right) \left( u\bullet v\right) \Delta _{D}.
\end{equation*}%
Since \eqref{def:Kpullb} is a pullback in $\Comon(\Cc)$, there is a unique morphism $w:D\rightarrow K$ of comonoids with
\begin{equation*}
kw=u\qquad \text{and}\qquad \lambda w=\left( u\bullet v\right) \Delta
_{D}.
\end{equation*}%
Then
\[
gw=l_{B}\left( \varepsilon _{A}\bullet
B\right) \lambda w=l_{B}\left( \varepsilon _{A}\bullet
B\right) \left( u\bullet v\right) \Delta _{D}=vl_{D}\left( \varepsilon _{D}\bullet D\right) \Delta _{D}=v.
\]
To conclude, we show that $w$ is the unique morphism in $\mathsf{Comon}(\Cc)$ such that $gw =v$ and $kw =u$. Suppose that there is another morphism $w^{\prime }:D\rightarrow K$ in $\mathsf{Comon}(\Cc)$ such
that $kw^{\prime }=u$ and $gw^{\prime }=v$. Note that
\[\lambda=(r_A(A\bullet \varepsilon_B)\bullet l_B(\varepsilon_A\bullet B))\Delta_{A\bullet B}\lambda
=(r_A(A\bullet \varepsilon_B)\lambda\bullet l_B(\varepsilon_A\bullet B)\lambda)\Delta_K=(k\bullet g)\Delta_K\]
where we used that $r_A(A\bullet \varepsilon_B)\lambda=r_A(A\bullet \varepsilon_C)(A\bullet m)\lambda=r_A(A\bullet \varepsilon_C)\rho_A k=k.$
Therefore
\begin{equation*}
\left( u\bullet v\right) \Delta _{D}=\left( k\bullet g\right) \left( w^{\prime }\bullet w^{\prime }\right) \Delta _{D}=\underbracket[0.140ex]{\left( k\bullet g\right) \Delta _{K}}w^{\prime }=\lambda w^{\prime }.
\end{equation*}%
Since $kw^{\prime }=u$ and $\lambda w^{\prime }=\left( u\bullet v\right)
\Delta _{D}$, by uniqueness in the definition of $w$, we get $%
w^{\prime }=w.$ We have so proved that \eqref{def:Kpullb} is a pullback in $\Comon(\Cc)$. Since all objects in this diagram are in $\cComon(\Cc)$, it is a fortiori a pullback in $\cComon(\Cc)$.
\end{proof}

We now apply \cref{lem:ddpulb} to construct similar pullbacks of bimonoids in a $\bullet$-braided duoidal category.

\begin{corollary}%[\rd{update:27/02/2026}]
\label{coro:ddpulb}
Let $\left( \mathcal{C},\circ ,I,\bullet ,J\right) $ be a $\bullet $-braided duoidal category. Assume that $\mathcal{C}$ has and $\bullet$ preserves binary intersections.
For every cospan
$\xymatrix{A\ar[r]^f&C&\,B\ar@{>->}[l]_m}$ in $\cBimon(\Cc,\circ,\bullet)$  with $m$ monomorphism in $\Cc$,
then the diagram \eqref{def:Kpullbcspn} is a pullback in $\Bimon (\Cc,\circ,\bullet)$ and in $\cBimon(\Cc,\circ,\bullet)$,
the diagram \eqref{def:Kpullb} is a pullback
 in $\Bimon(\Cc,\circ,\bullet)$ and the diagram \eqref{def:Kequa} is an equalizer  in $\Bimon(\Cc,\circ,\bullet)$.

Moreover, if there is a reversion on $\left( \mathcal{C}%
,\circ ,I,\bullet ,J\right) $ which preserves the equalizer \eqref{def:Kequa} in $\Cc$, the morphisms $\left( \lambda\bullet I\right) \circ J$ and $J\circ \left( I\bullet \lambda\right)$ are monomorphisms in $\Cc$ and the above cospan is in $\cHopfmon(\Cc,\circ,\bullet)$, then we can replace $\Bimon$ with $\Hopfmon$ above.
%\rd{[Nel caso dei bimoduli abbiamo tutti i ker di morfismi Hopf. Per recuperare questo, ho aggiunto l'ultima parte dell'enunciato. Nel caso dei bimoduli, la reversion è un'equivalenza, in particolare un aggiunto destro, e quindi preserva gli equalizzatori. Il fatto che $\left( \lambda\bullet I\right) \circ J$ and $J\circ \left( I\bullet \lambda\right)$ siano mono non è chiaro.]}
\end{corollary}

\begin{proof}
By \cref{lem:ddpulb},
\eqref{def:Kpullbcspn} is a pullback in $\Comon \left( \mathcal{C}\right) $ and in $\cComon(\Cc)$,
\eqref{def:Kpullb} is a pullback
 in $\Comon(\Cc)$ and \eqref{def:Kequa} is an equalizer in $\Cc$ and in $\Comon(\Cc)$.

Since the functor $\Bimon \left(\Cc,\circ,\bullet
\right) \cong \Mon\left( \Comon \left( \mathcal{C}^{\bullet
}\right) ,\circ \right) \rightarrow \Comon \left(\Cc
^{\bullet }\right) $ creates limits, we get that \eqref{def:Kpullbcspn} and \eqref{def:Kpullb} are pullbacks in $\Bimon (\Cc,\circ,\bullet)$, and \eqref{def:Kequa} is an equalizer in $\Bimon (\Cc,\circ,\bullet)$.
Since all the objects in \eqref{def:Kpullbcspn} are cocommutative, we get that this is also a pullback in $\cBimon(\Cc,\circ,\bullet)$.

Let us prove the last part of the statement. Consider now the following diagram.%
\[
\xymatrix@C=2.5cm{
K\ar@{.>}[d]^{\Sigma_K} \ar[r]^{\lambda} & A\bullet B\ar[d]^{\Sigma_{A\bullet B}}\ar@<0.5ex>[r]^-{(A\bullet f)\Delta _{A}^{\bullet }\bullet B} \ar@<-0.5ex>[r]_-{A\bullet \left(
m\bullet B\right) \Delta _{B}^{\bullet }} & A\bullet C\bullet B \ar[d]^{\Sigma_{(A\bullet C)\bullet B}}\\
K^\revsn \ar[r]^{\lambda^\revsn} & (A\bullet  B)^\revsn \ar@<0.5ex>[r]^-{((A\bullet f)\Delta _{A}^{\bullet }\bullet B)^\revsn} \ar@<-0.5ex>[r]_-{(A\bullet \left(
m\bullet B\right) \Delta _{B}^{\bullet })^\revsn} & (A\bullet C\bullet B)^\revsn
}
\]
where the solid vertical arrows are the antipodes arising from \cref{prop:Hopfbul}. %together with \cref{lem:IHopf}. 
Note that the upper row is the equalizer \eqref{def:Kequa} while the
lower row is obtained by applying the reversion to \eqref{def:Kequa} and it is an equalizer in $\Cc$ by hypothesis.
We have shown that the upper row is an equalizer in $\Bimon(\Cc,\circ,\bullet)$. Thus  the corresponding parallel pair consists of morphisms of bimonoids whence compatible with the antipode, in view of \cref{coro:convcomp}.
Note also that the associativity constraint $(A\bullet H)\bullet I\to A\bullet (H\bullet I)$ is a morphism of bimonoids and hence it is compatible with the antipode too, i.e. $\Sigma_{(A\bullet H)\bullet I}=\Sigma_{A\bullet (H\bullet I)}.$
This means that the squares on the right in the diagram above are serially
commutative. Therefore, there exists a unique morphism $\Sigma
_{K}:K\rightarrow K^{\revsn }$ in $\Cc$ such that $\lambda ^{\revsn }\Sigma _{K}=\Sigma
_{A\bullet I}\lambda .$
Since also $\lambda:K\rightarrow A\bullet I$ is a morphism of bimonoids and $\left( \lambda\bullet I\right) \circ J$ and $J\circ \left( I\bullet \lambda\right)$ are monomorphisms in $\Cc$, by \cref{lem:retract} we have that $\Sigma _{K}$ is an antipode.
\end{proof}

%\rd{[We need to prove that $K$ is a Hopf monoid in case $f$ is a morphism of Hopf monoids.]}
%\as{[Suppose that $f:A\to H$ is in $\Hopfmon(\Cc)$, so in particular $A$ and $H$ are in $\Hopfmon(\Cc)$. In order to prove that $K$ is in $\Hopfmon(\Cc)$, we have to show that, for any morphism $h:C\to K$ in $\mathsf{Comon}(\Cc^{\bullet})$ we have that $h^{r}$ is an isomorphism in $\Cc$ and for any morphism $h':K\to E$ in $\mathsf{Mon}(\Cc^{\circ})$ then $h'_{r}$ is an isomorphism in $\Cc$. Since $kh:C\to H$ is a morphism in $\mathsf{Comon}(\Cc^{\bullet})$ and $A$ is a Hopf monoid, then $(kh)^{r}$ is an isomorphism in $\Cc$. We recall that $(kh)^{r}(C\circ k)=((C\circ J)\bullet k)h^{r}$ and $k$ is a monomorphism in $\Cc$.]}
%\begin{invisible} ?$\xymatrix{(C\circ J)\bullet K \ar[r]^-{(C\circ J)\bullet k}& (C\circ J)\bullet A\ar[r]^-{((kh)^r)^{-1}}& C\circ A\ar[r]^-{?} &C\circ K}$
%\end{invisible}
%\medskip

\subsection{Points of Hopf monoids}
\noindent The notion of a point is a fundamental ingredient in the concept of protomodularity. Recall that a \emph{point} $\xymatrix{A\ar@<.5ex>[r]^-\pi&\ar@<.5ex>[l]^-\sigma H}$ in a category $\mathcal{C}$ consists of two morphisms $\pi:A\to H$ and $\sigma:H\to A$ in $\Cc$ such that $\pi\sigma=\id_H$, i.e.\ a split epimorphism $\pi$ in $\Cc$ with section $\sigma$. A morphism between points \begin{tikzcd}
	A & H
	\arrow[shift left, from=1-1, to=1-2,"\pi"]
	\arrow[shift left, from=1-2, to=1-1, "\sigma"]
\end{tikzcd} and \begin{tikzcd}
	A' & H'
	\arrow[shift left, from=1-1, to=1-2, "\pi'"]
	\arrow[shift left, from=1-2, to=1-1, "\sigma'"]
\end{tikzcd}
in $\Cc$ is a pair $(f,g)$ of morphisms $f:A\to A'$ and $g:H\to H'$ in $\Cc$ such that the following diagrams are serially commutative:
\[
\xymatrix@R=0.6cm{A\ar[d]_-{f}\ar@<.5ex>[r]^-\pi&\ar@<.5ex>[l]^-\sigma H\ar[d]^-{g}\\
A'\ar@<.5ex>[r]^-{\pi'} & H'\ar@<.5ex>[l]^-{\sigma'}
}
\]
The category of points in $\Cc$ and their morphisms is denoted by $\mathrm{Pt}(\Cc)$. \medskip

We now state our first main result, showing that, under suitable assumptions, a point of bimonoids admits a factorization in accordance with \cref{thm:psi} involving the kernel of the projection. Note that the kernel of a morphism $f:A\to B$ in a category with an initial object $\initial$ can be realized as the pullback along $f$ of the unique morphism $\initial\to B$.

\begin{theorem}[Factorization of points]\label{thm:main}
%[\rd{added: 2025/11/13}]
Let $\left( \mathcal{C},\circ ,I,\bullet ,J\right) $ be a $\bullet $-braided duoidal category. Assume that $\mathcal{C}$ has and $\bullet$ preserves binary intersections. Consider a point $\xymatrix{A\ar@<.5ex>[r]^-\pi&\ar@<.5ex>[l]^-\sigma H}$ in $\Bimon \left(\Cc,\circ,\bullet\right) $
with $u_{H}^{\circ }:I\rightarrow H$ monomorphism in $\Cc$ and $A$ cocommutative. Then, $\pi $
has a kernel $\left( K,k:K\rightarrow A\right) $ that is the pullback
\begin{equation*}
\xymatrix{K\pulb\ar[r]^k\ar[d]_{\varepsilon_K^\circ}&A\ar[d]^\pi\\I\,\ar@{>->}[r]^{u_H^\circ}&H}
\end{equation*}%
 in $\Bimon(\Cc,\circ,\bullet)$, in  $\cBimon(\Cc,\circ,\bullet)$ and  in $\Comon(\Cc^{\bullet }) $.

Assume further that there is a reversion on $\left( \mathcal{C}%
,\circ ,I,\bullet ,J\right) $.\\ 
If $H$ is a Hopf monoid, then $\varphi :=m_{A}^{\circ }\left( k\circ \sigma \right)
:K\circ H\rightarrow A$ is an isomorphism in $\Cc$. \\
If $A$ is a Hopf monoid, then so is $K$, which turns out to be the kernel of $\pi$ in $\cHopfmon(\Cc,\circ,\bullet)$.
\end{theorem}

\begin{proof}
By definition of $\bullet$-braided duoidal category, $I$ is $\bullet$-cocommutative. Moreover, since $A$ is cocommutative, so is $H$ as $\sigma_{H,H}\Delta^\bullet_H=(\pi\bullet\pi) \sigma_{A,A}\Delta_A^{\bullet}\sigma=(\pi\bullet\pi) \Delta_A^{\bullet}\sigma=\Delta^{\bullet}_H.$ Therefore, the cospan $\xymatrix{A\ar[r]^\pi&H&\,I\ar[l]_{u_H^\circ}}$
is in $\cBimon(\Cc,\circ,\bullet)$.
Since $u_{H}^{\circ }$ is a monomorphism in $\Cc$, by \cref{coro:ddpulb} this cospan
has a pullback in $\Bimon (\Cc,\circ,\bullet)$ and $\cBimon (\Cc,\circ,\bullet)$ as stated with $k$ monomorphism in $%
\mathcal{C}$ and $K$ cocommutative. Moreover, this pullback is preserved by
the forgetful functor $\Bimon \left( \mathcal{C},\circ,\bullet\right) \rightarrow
\Comon \left( \mathcal{C}^{\bullet }\right) .$
This means that $\left( K,k\right) $ is the kernel of $\pi $ in $\Bimon\left( \mathcal{C},\circ,\bullet\right) $ as $I$ is the initial object in this
category.

Assume now there is a reversion. If $H$ is a Hopf monoid, let us prove that the following diagram, which is \eqref{def:Kequa} in our setting and hence an equalizer in $\Bimon(\Cc,\circ,\bullet)$ in view of \cref{coro:ddpulb}, turns out to be a split equalizer in $\Cc$.
\begin{equation}
\label{def:spliteq}
\xymatrix@C=2cm{
K \ar[r]^{\lambda} & A\bullet I \ar@<0.5ex>[r]^-{\left( A\bullet \pi
\right) \Delta _{A}^{\bullet }\bullet I} \ar@<-0.5ex>[r]_-{A\bullet \left(
u_{H}^{\circ }\bullet I\right) \Delta _{I}^{\bullet }} & A\bullet H\bullet I
}
\end{equation}
Recall from \cref{lem:ddpulb} that we have the identity
\begin{equation}
\label{form:lambepsr}
    \left( k\bullet \varepsilon _{K}^{\circ }\right) \Delta _{K}^{\bullet
}=\lambda ,
\end{equation} that will be used afterword.
Set $\rho_A^\bullet:= \left( A\bullet \pi
\right) \Delta _{A}^{\bullet }:A\to A\bullet H$.
% One easily checks that the squares on the right in the following diagram are serially commutative.
% \[
% \xymatrix@C=2.5cm{
% K\ar[d]^{k} \ar[r]^{\lambda} & A\bullet I\ar[d]^{A\bullet  u_{H}^\circ}\ar@<0.5ex>[r]^-{\rho _{A}^{\bullet }\bullet I} \ar@<-0.5ex>[r]_-{A\bullet \left(
% u_{H}^{\circ }\bullet I\right) \Delta _{I}^{\bullet }} & A\bullet H\bullet I \ar[d]^{A\bullet H\bullet u_{H}^\circ}\\
% A \ar[r]^{\rho _{A}^{\bullet }} & A\bullet  H \ar@<0.5ex>[r]^-{\rho _{A}^{\bullet }\bullet H} \ar@<-0.5ex>[r]_-{A\bullet \Delta_{H}^{\bullet}} & A\bullet H\bullet H
% }
% \] Note that, for every object $X$ in $\Cc$, the functor $X\bullet (-)$ preserves monomorphisms as it preserves the kernel pair of any monomorphism and this is a binary intersection, see \cite[Proposition 2.5.6]{BorI94}.
% Thus, the vertical arrow $A\bullet H\bullet u_{H}^\circ$ is a monomorphism. Therefore, since the lower row of the diagram is a fork, we deduce that upper row of the diagram is a fork as well. Let as prove it is in fact a split fork.
Consider
\begin{equation*}
\pi ^{r}=\left( \left( A\circ \varepsilon _{H}^{\bullet }\right) \bullet
m_{H}^{\circ }\left( \pi \circ H\right) \right) \Delta _{A\circ H}^{\bullet
}:A\circ H\rightarrow \left( A\circ J\right) \bullet H
\end{equation*}%
that, since $H$ is a Hopf monoid, has inverse
$\left( \pi ^{r}\right) ^{-1}:\left( A\circ J\right) \bullet H\rightarrow
A\circ H$.
Define%
\[\xi:\quad\xymatrix@C=1.3cm{
A\bullet H\ar[r]^-{\left( r_{A}^{\circ }\right) ^{-1}\bullet H}&\left( A\circ I\right) \bullet H\ar[r]^-{\left( A\circ
\varepsilon _{I}^{\bullet }\right) \bullet H}&\left( A\circ
J\right) \bullet H
\ar[r]^-{\left( \pi ^{r}\right) ^{-1}}&
A\circ H\ar[r]^{A\circ \sigma }&A\circ A\ar[r]^{
m_{A}^{\circ }}&A
}\]
Let us prove that $\xi \left( A\bullet u_{H}^{\circ }\right) :A\bullet
I\rightarrow A$
factors through $K.$ We compute%
\begin{eqnarray*}
\rho _{A}^{\bullet }\xi  &=&\left( A\bullet \pi \right) \Delta _{A}^{\bullet
}m_{A}^{\circ }\left( A\circ \sigma \right) \left( \pi ^{r}\right)
^{-1}\left( \left( A\circ \varepsilon _{I}^{\bullet }\right) \bullet
H\right) \left( \left( r_{A}^{\circ }\right) ^{-1}\bullet H\right)  \\
&=&\left( A\bullet \pi \right) \left( m_{A}^{\circ }\bullet m_{A}^{\circ
}\right) \Delta _{A\circ A}^{\bullet }\left( A\circ \sigma \right) \left(
\pi ^{r}\right) ^{-1}\left( \left( A\circ \varepsilon _{I}^{\bullet }\right)
\bullet H\right) \left( \left( r_{A}^{\circ }\right) ^{-1}\bullet H\right)
\\
&=&\left( A\bullet \pi \right) \left( m_{A}^{\circ }\bullet m_{A}^{\circ
}\right) \left( \left( A\circ \sigma \right) \bullet \left( A\circ \sigma
\right) \right) \Delta _{A\circ H}^{\bullet }\left( \pi ^{r}\right)
^{-1}\left( \left( A\circ \varepsilon _{I}^{\bullet }\right) \bullet
H\right) \left( \left( r_{A}^{\circ }\right) ^{-1}\bullet H\right)  \\
&=&\left( m_{A}^{\circ }\left( A\circ \sigma \right) \bullet m_{H}^{\circ
}\left( \pi \circ H\right) \right) \Delta _{A\circ H}^{\bullet }\left( \pi
^{r}\right) ^{-1}\left( \left( A\circ \varepsilon _{I}^{\bullet }\right)
\bullet H\right) \left( \left( r_{A}^{\circ }\right) ^{-1}\bullet H\right)
\\
&=&\left( m_{A}^{\circ }\left( A\circ \sigma \right) \bullet H\right)
\underbracket[0.140ex]{\left( \left( A\circ H\right) \bullet m_{H}^{\circ }\left( \pi
\circ H\right) \right) \Delta _{A\circ H}^{\bullet }\left( \pi ^{r}\right)
^{-1}}\left( \left( A\circ \varepsilon _{I}^{\bullet }\right) \bullet
H\right) \left( \left( r_{A}^{\circ }\right) ^{-1}\bullet H\right)  \\
&\overset{\eqref{form:Deltagl}}{=}&\left( m_{A}^{\circ }\left( A\circ \sigma
\right) \bullet H\right) \left( \left( \pi ^{r}\right) ^{-1}\bullet H\right)
\left( \left( A\circ J\right) \bullet \Delta _{H}^{\bullet }\right) \left(
\left( A\circ \varepsilon _{I}^{\bullet }\right) \bullet H\right) \left(
\left( r_{A}^{\circ }\right) ^{-1}\bullet H\right)  \\
&=&\left( m_{A}^{\circ }\bullet H\right) \left( A\circ \sigma \bullet
H\right) \left( \left( \pi ^{r}\right) ^{-1}\bullet H\right) \left( \left(
A\circ \varepsilon _{I}^{\bullet }\right) \bullet H\bullet H\right) \left(
\left( r_{A}^{\circ }\right) ^{-1}\bullet H\bullet H\right) \left( A\bullet
\Delta _{H}^{\bullet }\right)  \\
&=&\left( \xi \bullet H\right) \left( A\bullet \Delta _{H}^{\bullet }\right)
\end{eqnarray*}%
so that $\rho _{A}^{\bullet }\xi =\left( \xi \bullet H\right) \left( A\bullet \Delta
_{H}^{\bullet }\right)$ i.e.\ $\xi$ is a morphism of right $H$-comodules.
Hence%
\begin{eqnarray*}
\rho _{A}^{\bullet }\xi \left( A\bullet u_{H}^{\circ }\right)  &=&\left( \xi
\bullet H\right) \left( A\bullet \Delta _{H}^{\bullet }\right) \left(
A\bullet u_{H}^{\circ }\right) =\left( \xi \bullet H\right) \left( A\bullet
u_{H}^{\circ }\bullet u_{H}^{\circ }\right) \left( A\bullet \Delta
_{I}^{\bullet }\right)  \\
&=&\left( A\bullet u_{H}^{\circ }\right) \left( \xi \left( A\bullet
u_{H}^{\circ }\right) \bullet I\right) \left( A\bullet \Delta _{I}^{\bullet
}\right) .
\end{eqnarray*}%
By definition of $K$ there is a unique morphism $\theta :A\bullet
I\rightarrow K$ in $\Cc$ such that%
\begin{equation}
\label{def:theta}
k\theta =\xi \left( A\bullet u_{H}^{\circ }\right) \qquad \text{and}\qquad
\lambda \theta =\left( \xi \left( A\bullet u_{H}^{\circ }\right) \bullet
I\right) \left( A\bullet \Delta _{I}^{\bullet }\right).
\end{equation}%
Let us check that $\xi \rho _{A}^{\bullet }=\mathrm{Id}_{A}.$ Since
\begin{eqnarray*}
\pi ^{r}\left( A\circ u_{H}^{\circ }\right) \left( r_{A}^{\circ }\right)
^{-1} &=&\left( \left( A\circ \varepsilon _{H}^{\bullet }\right) \bullet
m_{H}^{\circ }\left( \pi \circ H\right) \right) \Delta _{A\circ H}^{\bullet
}\left( A\circ u_{H}^{\circ }\right) \left( r_{A}^{\circ }\right) ^{-1} \\
&=&\left( \left( A\circ \varepsilon _{H}^{\bullet }\right) \bullet
m_{H}^{\circ }\left( \pi \circ H\right) \right) \left( \left( A\circ
u_{H}^{\circ }\right) \left( r_{A}^{\circ }\right) ^{-1}\bullet \left(
A\circ u_{H}^{\circ }\right) \left( r_{A}^{\circ }\right) ^{-1}\right)
\Delta _{A}^{\bullet } \\
&=&\left( \left( A\circ \varepsilon _{I}^{\bullet }\right) \bullet H\right)
\left( \left( r_{A}^{\circ }\right) ^{-1}\bullet H\right) \left( A\bullet
m_{H}^{\circ }\left( \pi \circ u_{H}^{\circ }\right) \left( r_{A}^{\circ
}\right) ^{-1}\right) \Delta _{A}^{\bullet } \\
&=&\left( \left( A\circ \varepsilon _{I}^{\bullet }\right) \bullet H\right)
\left( \left( r_{A}^{\circ }\right) ^{-1}\bullet H\right) \left( A\bullet
\pi \right) \Delta _{A}^{\bullet }=\left( \left( A\circ \varepsilon
_{I}^{\bullet }\right) \bullet H\right) \left( \left( r_{A}^{\circ }\right)
^{-1}\bullet H\right) \rho _{A}^{\bullet }
\end{eqnarray*}%
we get
\begin{eqnarray*}
\xi \rho _{A}^{\bullet } &=&m_{A}^{\circ }\left( A\circ \sigma \right)
\underbracket[0.140ex]{\left( \pi ^r\right) ^{-1}\left( \left( A\circ \varepsilon
_{I}^{\bullet }\right) \bullet H\right) \left( \left( r_{A}^{\circ }\right)
^{-1}\bullet H\right) \rho _{A}^{\bullet }} \\
&=&m_{A}^{\circ }\left( A\circ \sigma \right) \left( A\circ u_{H}^{\circ
}\right) \left( r_{A}^{\circ }\right) ^{-1}=m_{A}^{\circ }\left( A\circ
u_{A}^{\circ }\right) \left( r_{A}^{\circ }\right) ^{-1}=\mathrm{Id}_{A}.
\end{eqnarray*}%
As a consequence, we get%
\begin{equation*}
k\theta \lambda \overset{\eqref{def:theta}}{=}\xi \left( A\bullet u_{H}^{\circ }\right) \lambda =\xi (k\bullet \pi k)\Delta^\bullet_K=\xi(A\bullet\pi)\Delta^\bullet_A k=\xi \rho
_{A}^{\bullet }k=k
\end{equation*}%
and hence $\theta \lambda =\mathrm{Id}$ as $k$ is a monomorphism. Since%
\begin{eqnarray*}
\left( \xi \bullet I\right) \left( \left( A\bullet \pi \right) \Delta
_{A}^{\bullet }\bullet I\right)  &=&\left( \xi \bullet I\right) \left( \rho
_{A}^{\bullet }\bullet I\right) =\mathrm{Id}_{A}, \\
\left( \xi \bullet I\right) \left( A\bullet \left( u_{H}^{\circ }\bullet
I\right) \Delta _{I}^{\bullet }\right)  &=&\left( \xi \bullet I\right)
\left( A\bullet u_{H}^{\circ }\bullet I\right) \left( A\bullet \Delta
_{I}^{\bullet }\right) =\lambda \theta ,
\end{eqnarray*}%
we conclude that the fork in \eqref{def:spliteq} is a split equalizer in $\Cc$.

Thus, it is an absolute equalizer, i.e.\ it is preserved by any functor, in
particular by $\left( -\right) \circ H.$

Therefore, we get the upper fork in the following diagram is a split equalizer.
\[
\xymatrix@C=2.5cm{
K\circ H\ar[d]^{\varphi} \ar[r]^{\lambda\circ H} & (A\bullet I)\circ H \ar[d]^{\sigma_r}\ar@<0.5ex>[r]^-{(\rho _{A}^{\bullet }\bullet I)\circ H} \ar@<-0.5ex>[r]_-{(A\bullet \left(
u_{H}^{\circ }\bullet I\right) \Delta _{I}^{\bullet })\circ H} & (A\bullet H\bullet I)\circ H \ar[d]^{\left( \left( \sigma
\bullet H\right) \Delta _{H}^{\bullet }\right) _{r}}\\
A \ar[r]^{\rho _{A}^{\bullet }} & A\bullet  H \ar@<0.5ex>[r]^-{\rho _{A}^{\bullet }\bullet H} \ar@<-0.5ex>[r]_-{A\bullet \Delta_{H}^{\bullet}} & A\bullet H\bullet H
}
\]
Also the lower fork is a split equalizer.
\begin{invisible}
 We can prove it by hand.
Indeed set $\epsilon_X:=r_X^\bullet(X\bullet\varepsilon_H^\bullet)$ and define $s:=\epsilon_A$ and $t:=A\bullet\epsilon_H$. Then $s\rho_A^\bullet=\id_A$, $\rho_A^\bullet s=t(\rho_A\bullet H)$ and $t(A\bullet \Delta_H^\bullet)=\id_{A\bullet H}$.
In alternative, apply the dual of \cite[Diagram 4.19]{BorII94} to the adjunction $L\dashv R$ where $L:(\Cc^\bullet)^H\to \Cc$ is  the forgetful functor and $R(X)=(X\bullet H,X\bullet \Delta_H^\bullet)$. The unit is defined by $L\eta_{(X,\rho)}=\rho$ and the counit $\epsilon_X:LR(X)=X\bullet H\to X$ is $r_X^\bullet (X\bullet \varepsilon_H^\bullet)$. \\
\end{invisible}
Let us prove that, in the diagram
above, the square on the left commutes and the squares on the right are serially commutative. Concerning the square on the left, we compute
\begin{eqnarray*}
\rho _{A}^{\bullet }\varphi  &=&\rho _{A}^{\bullet }m_{A}^{\circ }\left(
k\circ \sigma \right)  \\
&=&\left( A\bullet \pi \right) \Delta _{A}^{\bullet }m_{A}^{\circ }\left(
k\circ \sigma \right)  \\
% &=&\left( A\bullet \pi \right) \left( m_{A}^{\circ }\bullet m_{A}^{\circ
% }\right) \Delta _{A\circ A}^{\bullet }\left( k\circ \sigma \right)  \\
&=&\left( A\bullet \pi \right) \left( m_{A}^{\circ }\left( k\circ \sigma
\right) \bullet m_{A}^{\circ }\left( k\circ \sigma \right) \right) \Delta
_{K\circ H}^{\bullet } \\
&=&\left( m_{A}^{\circ }\left( k\circ \sigma \right) \bullet m_{H}^{\circ
}\left( \pi k\circ H\right) \right) \Delta _{K\circ H}^{\bullet } \\
&=&\left( m_{A}^{\circ }\left( k\circ \sigma \right) \bullet m_{H}^{\circ
}\left( u_{H}^{\circ }\varepsilon _{K}^{\circ }\circ H\right) \right) \zeta
_{K,K,H,H}\left( \Delta _{K}^{\bullet }\circ \Delta _{H}^{\bullet }\right)
\\
&=&\left( m_{A}^{\circ }\bullet m_{H}^{\circ }\right) \left( \left( k\circ
\sigma \right) \bullet \left( u_{H}^{\circ }\varepsilon _{K}^{\circ }\circ
H\right) \right) \zeta _{K,K,H,H}\left( \Delta _{K}^{\bullet }\circ \Delta
_{H}^{\bullet }\right)  \\
&=&\left( m_{A}^{\circ }\bullet m_{H}^{\circ }\right) \zeta _{A,H,A,H}\left(
\left( k\bullet u_{H}^{\circ }\varepsilon _{K}^{\circ }\right) \circ \left(
\sigma \bullet H\right) \right) \left( \Delta _{K}^{\bullet }\circ \Delta
_{H}^{\bullet }\right)  \\
&=&\underbracket[0.140ex]{m_{A\bullet H}^{\circ }\left( \left( A\bullet u_{H}^{\circ
}\right) \circ \left( \sigma \bullet H\right) \Delta _{H}^{\bullet }\right) }%
\left( \underbracket[0.140ex]{\left( k\bullet \varepsilon _{K}^{\circ }\right) \Delta
_{K}^{\bullet }}\circ H\right)  \overset{\eqref{form:lambepsr}}=\sigma _{r}\left( \lambda \circ H\right),
\end{eqnarray*}%
so that $\rho _{A}^{\bullet }\varphi =\sigma _{r}\left( \lambda \circ
H\right) .$ 
Concerning the squares on the right, since $\rho _{A}^{\bullet }\sigma =\left( A\bullet \pi
\right) \Delta _{A}^{\bullet }\sigma =\left( \sigma \bullet \pi \sigma
\right) \Delta _{H}^{\bullet }=\left( \sigma \bullet H\right) \Delta
_{H}^{\bullet }$, we get
\begin{eqnarray*}
\left( \left( \sigma \bullet H\right) \Delta _{H}^{\bullet }\right)
_{r}\left( \left( \rho _{A}^{\bullet }\bullet I\right) \circ H\right)
&=&m_{\left( A\bullet H\right) \bullet H}^{\circ }\left( \left( A\bullet
H\bullet u_{H}^{\circ }\right) \circ \left( \underbracket[0.140ex]{\left( \sigma
\bullet H\right) \Delta _{H}^{\bullet }}\bullet H\right) \Delta
_{H}^{\bullet }\right) \left( \left( \rho _{A}^{\bullet }\bullet I\right)
\circ H\right)  \\
&=&m_{\left( A\bullet H\right) \bullet H}^{\circ }\left( \left( \rho
_{A}^{\bullet }\bullet u_{H}^{\circ }\right) \circ \left( \rho _{A}^{\bullet
}\sigma \bullet H\right) \Delta _{H}^{\bullet }\right)  \\
% &=&m_{\left( A\bullet H\right) \bullet H}^{\circ }\left( \left( \rho
% _{A}^{\bullet }\bullet H\right) \circ \left( \rho _{A}^{\bullet }\bullet
% H\right) \right) \left( \left( A\bullet u_{H}^{\circ }\right) \circ \left(
% \sigma \bullet H\right) \Delta _{H}^{\bullet }\right)  \\
&=&\left( \rho _{A}^{\bullet }\bullet H\right) m_{A\bullet H}^{\circ }\left(
\left( A\bullet u_{H}^{\circ }\right) \circ \left( \sigma \bullet H\right)
\Delta _{H}^{\bullet }\right)  =\left( \rho _{A}^{\bullet }\bullet H\right) \sigma _{r}
\end{eqnarray*}%
and
\begin{eqnarray*}
&&\left( \left( \sigma \bullet H\right) \Delta _{H}^{\bullet }\right)
_{r}\left( \left( A\bullet \left( u_{H}^{\circ }\bullet I\right) \Delta
_{I}^{\bullet }\right) \circ H\right)  \\
&=&m_{\left( A\bullet H\right) \bullet H}^{\circ }\left( \left( A\bullet
H\bullet u_{H}^{\circ }\right) \circ \left( \left( \sigma \bullet H\right)
\Delta _{H}^{\bullet }\bullet H\right) \Delta _{H}^{\bullet }\right) \left(
\left( A\bullet \left( u_{H}^{\circ }\bullet I\right) \Delta _{I}^{\bullet
}\right) \circ H\right)  \\
&=&m_{\left( A\bullet H\right) \bullet H}^{\circ }\left( \left( A\bullet
\left( u_{H}^{\circ }\bullet u_{H}^{\circ }\right) \Delta _{I}^{\bullet
}\right) \circ \left( \sigma \bullet H\bullet H\right) \left( \Delta
_{H}^{\bullet }\bullet H\right) \Delta _{H}^{\bullet }\right)  \\
&=&m_{A\bullet \left( H\bullet H\right) }^{\circ }\left( \left( A\bullet
\Delta _{H}^{\bullet }u_{H}^{\circ }\right) \circ \left( \sigma \bullet
\Delta _{H}^{\bullet }\right) \Delta _{H}^{\bullet }\right)  \\
% &=&m_{A\bullet \left( H\bullet H\right) }^{\circ }\left( \left( A\bullet
% \Delta _{H}^{\bullet }\right) \circ \left( A\bullet \Delta _{H}^{\bullet
% }\right) \right) \left( \left( A\bullet u_{H}^{\circ }\right) \circ \left(
% \sigma \bullet H\right) \Delta _{H}^{\bullet }\right)  \\
&=&\left( A\bullet \Delta _{H}^{\bullet }\right) m_{A\bullet H}^{\circ
}\left( \left( A\bullet u_{H}^{\circ }\right) \circ \left( \sigma \bullet
H\right) \Delta _{H}^{\bullet }\right)  =\left( A\bullet \Delta _{H}^{\bullet }\right) \sigma _{r}.
\end{eqnarray*}

Now, since $H$ is a Hopf
monoid and $\sigma :H\rightarrow A$ and $\left( \sigma \bullet H\right)
\Delta _{H}^{\bullet }:H\rightarrow A\bullet H$ are morphisms of $\circ$-monoids, we
get that $\sigma _{r}$ and $\left( \left( \sigma \bullet H\right) \Delta
_{H}^{\bullet }\right) _{r}$ are invertible as so are the Galois maps in view of \cref{thm:Galinv}. As a consequence, also $\varphi $ is
invertible by uniqueness of the equalizer.

From $\rho _{A}^{\bullet }\varphi =\sigma _{r}\left( \lambda \circ H\right) $
we get $\left( \sigma _{r}\right) ^{-1}\rho _{A}^{\bullet }=\left( \lambda
\circ H\right) \varphi ^{-1}$. Composition by $\theta \circ H$ gives that
$\varphi ^{-1}=\left( \theta \circ H\right) \left( \sigma _{r}\right)
^{-1}\rho _{A}^{\bullet }$, providing an explicit formula for the inverse.

If $A$ has an antipode, then also its retract $H$ has one, in view of \cref{lem:retract}.
The fact that, in this setting, also $K$ has an antipode follows by the last part of \cref{coro:ddpulb}
 once observed that \eqref{def:spliteq} is an absolute equalizer in $\Cc$, whence preserved by
 the reversion, and that both
$\left(\lambda\bullet I\right) \circ J$ and $J\circ \left( I\bullet \lambda\right)$ are split-monomorphisms.
\end{proof}

Incorporating Hopf monoid factorization results into the previous theorem, the given factorization of points turns out to depend on suitable morphisms $\chi$ and $\triangleleft$.

\begin{corollary}
\label{coro:main}
In the setting
 of \cref{thm:main} with $H$ Hopf monoid, $\varphi :K\circ _{\psi }H\rightarrow A$
becomes an isomorphism in $\cBimon(\Cc,\circ,\bullet)$, where $\psi :=\varphi ^{-1}m_{A}^{\circ
}\left( \sigma \circ k\right) :H\circ K\rightarrow K\circ H$. If $A$ is a Hopf monoid, then so is $K\circ H$ with antipode $\Sigma
_{K\circ H}:=\psi ^{\revsn }\phi _{K,H}^{\circ }\left( \Sigma _{K}\circ \Sigma
_{H}\right) .$
 If we define
\begin{eqnarray*}
\Upsilon  &:&\xymatrix@C=2cm{
  A\ar[r]^-{\mathrm{Id}_{A}\star \sigma^{\revsn}\,\Sigma_{H}\pi}
  &  (A\bullet I)\circ J
    \ar[r]^{\theta\circ J}
  &  K\circ J}
, \\
\chi  &:&\xymatrix@C=1cm{
  H\circ K \ar[r]^{\sigma \circ k}
  &  A\circ A
    \ar[r]^{m_A^{\circ}}
  &  A    \ar[r]^-{\Upsilon}
  &  K\circ J}
,\\
\triangleleft &:& \xymatrix{H\circ K\ar[r]^{H\circ\varepsilon_K^\circ}&H\circ I\ar[r]^{r^\circ_H}& H},
\end{eqnarray*} we get that $\psi:H\circ K\to K\circ H$ equals the composition
\[\xymatrix@C=1.5cm{H\circ K\ar[r]^-{\Delta
_{H\circ K}^{\bullet }}&(H\circ K)\bullet (H\circ K)\ar[r]^-{\chi \bullet \triangleleft}&(K\circ J)\bullet H
\ar[r]^-{\left( \left( u_{H}^{\circ }\varepsilon _{K}^{\circ }\right)
^{r}\right) ^{-1}}&K\circ H
}\]
so that it is completely determined by the morphisms $%
\chi:H\circ K\to K\circ J $ and $\triangleleft :H\circ K\to H.$  Moreover, we have $\varphi ^{-1}=\left( \left( u_{H}^{\circ }\varepsilon _{K}^{\circ }\right)
^{r}\right) ^{-1}\left( \Upsilon \bullet \pi \right) \Delta _{A}^{\bullet }$.
\end{corollary}

\begin{proof}
By \cref{thm:main}, $\varphi:=m_{A}^{\circ }\left( k\circ \sigma \right) :K\circ H\rightarrow A$ is an isomorphism in $\Cc$. Thus, we can apply \cref{thm:psi} to the cospan
$\xymatrixcolsep{1pc}\xymatrix{K\ar[r]^-{k}& A&\ar[l]_-{\sigma} H}$. As a consequence $\varphi :K\circ _{\psi }H\rightarrow A$
becomes an isomorphism in $\cBimon(\Cc,\circ,\bullet)$, where $\psi :=\varphi ^{-1}m_{A}^{\circ
}\left( \sigma \circ k\right) :H\circ K\rightarrow K\circ H.$

If $A$ is a Hopf monoid,  \cref{thm:main} ensures that so are $H$ and $K$. Therefore, we can apply \cref{thm:psitriang} to get that $K\circ _{\psi }H$ has the desired antipode. The same result tells us that $\psi =\left( \left( u_{H}^{\circ }\varepsilon _{K}^{\circ }\right)
^{r}\right) ^{-1}\left( \chi \bullet \triangleleft \right) \Delta
_{H\circ K}^{\bullet }$ where $\chi=\left( K\circ \varepsilon _{H}^{\bullet }\right) \psi$
 and $\triangleleft  :=l_{H}^{\circ }\left( \varepsilon _{K}^{\circ }\circ
H\right) \psi$. We want to check that these morphisms are as claimed.
Since $\pi \varphi =\pi m_{A}^{\circ }\left( k\circ \sigma \right)
=m_{H}^{\circ }\left( \pi k\circ \pi \sigma \right) =m_{H}^{\circ }\left(
u_{H}^{\circ }\varepsilon _{K}^{\circ }\circ H\right) =l_{H}^{\circ }\left(
\varepsilon _{K}^{\circ }\circ H\right) ,$ we get $l_{H}^{\circ }\left(
\varepsilon _{K}^{\circ }\circ H\right) \varphi ^{-1}=\pi $ and hence
\begin{eqnarray*}
\triangleleft  &=&l_{H}^{\circ }\left( \varepsilon _{K}^{\circ }\circ
H\right) \psi =\underbracket[0.140ex]{l_{H}^{\circ }\left( \varepsilon _{K}^{\circ
}\circ H\right) \varphi ^{-1}}m_{A}^{\circ }\left( \sigma \circ k\right)  \\
&=&\pi m_{A}^{\circ }\left( \sigma \circ k\right) =m_{H}^{\circ }\left( \pi
\sigma \circ \pi k\right) =m_{H}^{\circ }\left( H\circ u_{H}^{\circ
}\varepsilon _{K}^{\circ }\right) =r_{H}^{\circ }\left( H\circ \varepsilon
_{K}^{\circ }\right) .
\end{eqnarray*} Therefore, $\triangleleft$ has the desired form.
Moreover, we compute
\begin{eqnarray*}
&&\left( \left( A\bullet I\right) \circ \varepsilon _{H}^{\bullet }\right)
\left( \sigma _{r}\right) ^{-1}\rho _{A}^{\bullet } \\
&\overset{\eqref{equivalentformgr}}{=}&\left( \left( A\bullet I\right) \circ
\varepsilon _{H}^{\bullet }\right) \overline{\varkappa }_{A,I,H}^{H}\left(
A\bullet \left( l_{H}^{\circ }\right) ^{-1}\right) \left( A\bullet \pi
\right) \Delta _{A}^{\bullet } \\
&=&\left( \left( A\bullet I\right) \circ \varepsilon _{H}^{\bullet }\right)
\left( \left( m_{A}^{\circ }\left( A\circ \sigma \right) \bullet I\right)
\circ H\right) \delta _{A,I,H,H}\left( A\bullet \left( I\circ \left( \Sigma
_{H}\bullet H\right) \Delta _{H}^{\bullet }\right) \right) \left( A\bullet
\left( l_{H}^{\circ }\right) ^{-1}\right) \left( A\bullet \pi \right) \Delta
_{A}^{\bullet } \\
&=&\left( \left( m_{A}^{\circ }\left( A\circ \sigma \right) \bullet I\right)
\circ J\right) \delta _{A,I,H,J}\left( A\bullet \left( I\circ \left( \Sigma
_{H}\bullet \varepsilon _{H}^{\bullet }\right) \Delta _{H}^{\bullet }\right)
\right) \left( A\bullet \left( l_{H}^{\circ }\right) ^{-1}\right) \left(
A\bullet \pi \right) \Delta _{A}^{\bullet } \\
&=&\left( \left( m_{A}^{\circ }\left( A\circ \sigma \right) \bullet I\right)
\circ J\right) \delta _{A,I,H,J}\left( A\bullet \left( I\circ \left( \Sigma
_{H}\bullet J\right) \left( r_{H}^{\bullet }\right) ^{-1}\right) \right)
\left( A\bullet \left( l_{H}^{\circ }\right) ^{-1}\right) \left( A\bullet
\pi \right) \Delta _{A}^{\bullet } \\
&=&\left( \left( m_{A}^{\circ }\bullet I\right) \circ J\right) \delta
_{A,I,A,J}\left( A\bullet \left( I\circ \left( r_{A^{\revsn }}^{\bullet
}\right) ^{-1}\right) \left( l_{A^{\revsn }}^{\circ }\right) ^{-1}\right)
\left( A\bullet \sigma ^{\revsn }\Sigma _{H}\pi \right) \Delta _{A}^{\bullet }
\\
&=&\left( \left( m_{A}^{\circ }\bullet I\right) \circ J\right) \varphi
_{A,A}\left( A\bullet \sigma ^{\revsn }\Sigma _{H}\pi \right) \Delta
_{A}^{\bullet }=\mathrm{Id}_{A}\star \sigma ^{\revsn }\Sigma _{H}\pi
\end{eqnarray*}
so that $\left( \left( A\bullet I\right) \circ \varepsilon _{H}^{\bullet
}\right) \left( \sigma _{r}\right) ^{-1}\rho _{A}^{\bullet }=\mathrm{Id}%
_{A}\star \sigma ^{\revsn }\Sigma _{H}\pi .$ As a consequence, since in \cref{thm:main} we proved
that $\varphi ^{-1}=\left( \theta \circ H\right) \left( \sigma _{r}\right)
^{-1}\rho _{A}^{\bullet },$ we obtain%
\begin{equation*}
\left( K\circ \varepsilon _{H}^{\bullet }\right) \varphi ^{-1}=\left(
K\circ \varepsilon _{H}^{\bullet }\right) \left( \theta \circ H\right)
\left( \sigma _{r}\right) ^{-1}\rho _{A}^{\bullet }
=\left( \theta \circ J\right) {\left( \left( A\bullet I\right)
\circ \varepsilon _{H}^{\bullet }\right) \left( \sigma _{r}\right) ^{-1}\rho
_{A}^{\bullet }}
=\left( \theta \circ J\right) \left( \mathrm{Id}_{A}\star \sigma ^{\revsn
}\Sigma _{H}\pi \right),
\end{equation*}%
so that $\left( K\circ \varepsilon _{H}^{\bullet }\right) \varphi ^{-1}=\Upsilon
.$ Hence $\chi =\left( K\circ \varepsilon _{H}^{\bullet }\right) \psi =\left( K\circ \varepsilon _{H}^{\bullet }\right) \varphi ^{-1}
m_{A}^{\circ }\left( \sigma \circ k\right) =\Upsilon m_{A}^{\circ }\left( \sigma
\circ k\right) $ as desired. Finally,
\begin{eqnarray*}
\varphi ^{-1} &=&\left( \left( u_{H}^{\circ }\varepsilon _{K}^{\circ
}\right) ^{r}\right) ^{-1}\left( u_{H}^{\circ }\varepsilon _{K}^{\circ
}\right) ^{r}\varphi ^{-1} \\
&=&\left( \left( u_{H}^{\circ }\varepsilon _{K}^{\circ }\right) ^{r}\right)
^{-1}\left( \left( K\circ \varepsilon _{H}^{\bullet }\right) \bullet
m_{H}^{\circ }\left( u_{H}^{\circ }\varepsilon _{K}^{\circ }\circ H\right)
\right) \Delta _{K\circ H}^{\bullet }\varphi ^{-1} \\
&=&\left( \left( u_{H}^{\circ }\varepsilon _{K}^{\circ }\right) ^{r}\right)
^{-1}\left( \left( K\circ \varepsilon _{H}^{\bullet }\right) \varphi
^{-1}\bullet m_{H}^{\circ }\left( u_{H}^{\circ }\varepsilon _{K}^{\circ
}\circ H\right) \varphi ^{-1}\right) \Delta _{A}^{\bullet } \\
&=&\left( \left( u_{H}^{\circ }\varepsilon _{K}^{\circ }\right) ^{r}\right)
^{-1}\left( \underbracket[0.140ex]{\left( K\circ \varepsilon _{H}^{\bullet }\right)
\varphi ^{-1}}\bullet \underbracket[0.140ex]{l_{H}^{\circ }\left( \varepsilon _{K}^{\circ }\circ
H\right) \varphi ^{-1}}\right) \Delta _{A}^{\bullet } 
=\left( \left( u_{H}^{\circ }\varepsilon _{K}^{\circ }\right) ^{r}\right)
^{-1}\left( \Upsilon \bullet \pi \right) \Delta _{A}^{\bullet }
\end{eqnarray*}
which proves the last assertion in the statement.
\end{proof}

By means of \cref{thm:main}, we are ready to prove the following analogue of the Split Short Five Lemma in the framework of duoidal categories.

\begin{lemma}[Split Short Five]\label{lem:split-short}
In the most restrictive setting of \cref{thm:main}, where for $A'$ and $H'$ we make the same assumptions as we did for $A$ and $H$, consider the following diagram
\[\xymatrix@R=.7cm@C=1.5cm{K\ar[d]^f\ar[r]^-k
&A\ar[d]^-\phi\ar@<.5ex>[r]^-\pi&\ar@<.5ex>[l]^-\sigma H\ar[d]^g\\ K'\ar[r]^-{k'}
&A'\ar@<.5ex>[r]^-{\pi'}&\ar@<.5ex>[l]^-{\sigma'} H'}\]
where the right square is a morphism of points while $\left( K,k\right) $
and $\left( K^{\prime },k^{\prime }\right) $ are the kernels of $\pi $ and $%
\pi ^{\prime }$ respectively. If $f$ and $g$ are both isomorphisms in $\Bimon(\Cc,\circ,\bullet)$, then so is $\phi .$
\end{lemma}

\begin{proof}
By \cref{thm:main}, we have isomorphisms $\varphi :K\circ H\rightarrow A$ and $\varphi
^{\prime }:K^{\prime }\circ H^{\prime }\rightarrow A^{\prime }$ in $\Cc$. We compute
$\varphi ^{\prime }\left( f\circ g\right) =m_{A^{\prime }}^{\circ }\left(
k^{\prime }\circ \sigma ^{\prime }\right) \left( f\circ g\right)
=m_{A^{\prime }}^{\circ }\left( \phi \circ \phi \right) \left( k\circ \sigma
\right) =\phi m_{A}^{\circ }\left( k\circ \sigma \right) =\phi \varphi$.
Thus, $\phi =\varphi ^{\prime }\left( f\circ g\right) \varphi ^{-1}$ is an isomorphism in $\Cc$.
\end{proof}

The next goal is to prove that every point admits pullbacks along any morphism. To do this, we first need the two following preliminary results.

\begin{lemma}
\label{lem:urinv}
Given $H$ in $\Bimon\left( \Cc,\circ,\bullet\right) $, we have $%
\left( u_{H}^{\circ }\right) ^{r}=\left( \left( l_{J}^{\circ
}\right) ^{-1}\bullet H\right) \left( l_{H}^{\bullet }\right)
^{-1}l_{H}^{\circ }$. Thus, $\left( u_{H}^{\circ }\right) ^{r}$ is invertible.
\end{lemma}

\begin{proof}
By \eqref{equivalentgr}, we have $\left( u_{H}^{\circ }\right) ^{r} =\varsigma_{I,J,H}^H(I\circ (l_H^\bullet)^{-1})$ where $\rho_I^\bullet=(I\bullet u_H^\circ)\Delta_I^\bullet$. One checks that $\varsigma_{I,J,H}^H=\varsigma_{I,J,H}^I$ where in the latter one takes $\mu_H^\circ=l_ H^\circ$.
\begin{invisible}
$\varsigma_{I,J,H}^H= ((I\circ J)\bullet m_ H^\circ)\zeta_{I,H,J,H}((I\bullet u_H^\circ)\Delta_I^\bullet\circ (J\bullet H))
=((I\circ J)\bullet m_ H^\circ(u_H^\circ\circ H))\zeta_{I,I,J,H}(\Delta_I^\bullet\circ (J\bullet H))
=((I\circ J)\bullet l_ H^\circ)\zeta_{I,I,J,H}(\Delta_I^\bullet\circ (J\bullet H))=\varsigma_{I,J,H}^I.$
\end{invisible}
As a consequence, in view of \cref{lem:inv}, we get
\[\left( u_{H}^{\circ }\right) ^{r} =\varsigma_{I,J,H}^I(I\circ (l_H^\bullet)^{-1})=\left( \left( l_{J}^{\circ }\right) ^{-1}\bullet H\right) l_{J\bullet
H}^{\circ }\left( I\circ \left( l_{H}^{\bullet }\right) ^{-1}\right)
=\left( \left( l_{J}^{\circ }\right) ^{-1}\bullet H\right) \left(
l_{H}^{\bullet }\right) ^{-1}l_{H}^{\circ }.\qedhere\]
% We compute
% \begin{eqnarray*}
% \left( u_{H}^{\circ }\right) ^{r} &=&\left( \left( I\circ \varepsilon
% _{H}^{\bullet }\right) \bullet m_{H}^{\circ }\left( u_{H}^{\circ }\circ
% H^{\prime }\right) \right) \Delta _{I\circ H}^{\bullet } \\
% &=&\left( \left( I\circ \varepsilon _{H^{\prime }}^{\bullet }\right) \bullet
% l_{H}^{\circ }\right) \Delta _{I\circ H}^{\bullet } \\
% &=&\left( \left( I\circ J\right) \bullet l_{H}^{\circ }\right) \left( \left(
% I\circ \varepsilon _{H}^{\bullet }\right) \bullet \left( I\circ H\right)
% \right) \zeta _{I,I,H,H}\left( \Delta _{I}^{\bullet }\circ \Delta
% _{H}^{\bullet }\right)  \\
% &=&\left( \left( I\circ J\right) \bullet l_{H}^{\circ }\right) \zeta
% _{I,I,J,H}\left( \left( I\bullet I\right) \circ \left( \varepsilon
% _{H}^{\bullet }\bullet H\right) \right) \left( \Delta _{I}^{\bullet }\circ
% \Delta _{H}^{\bullet }\right)  \\
% &=&\left( \left( I\circ J\right) \bullet l_{H}^{\circ }\right) \zeta
% _{I,I,J,H}\left( \Delta _{I}^{\bullet }\circ \left( J\bullet H\right)
% \right) \left( I\circ \left( l_{H}^{\bullet }\right) ^{-1}\right)  \\
% &=&\left( \left( l_{J}^{\circ }\right) ^{-1}\bullet H\right) l_{J\bullet
% H}^{\circ }\left( I\circ \left( l_{H}^{\bullet }\right) ^{-1}\right)  \\
% &=&\left( \left( l_{J}^{\circ }\right) ^{-1}\bullet H\right) \left(
% l_{H}^{\bullet }\right) ^{-1}l_{H}^{\circ }
% \end{eqnarray*}
\end{proof}

\begin{lemma}
\label{lem:grcomon}
Let $\left( \mathcal{C},\circ ,I,\bullet ,J\right) $ be a $\bullet $-braided duoidal
category. If $H$ is in $\cBimon\left( \mathcal{C,\circ, \bullet }\right) $ and $g:C\rightarrow H$ is in $\cComon\left( \mathcal{C}%
^{\bullet }\right) $, then $g^{r}:C\circ H\rightarrow \left( C\circ J\right)
\bullet H$ is in $\Comon\left( \mathcal{C}^{\bullet }\right) $.
\end{lemma}

\begin{proof}
Since $\mathcal{C}^{\bullet }$ is braided, we know that $\left( C\circ
J\right) \bullet H$ is a comonoid. By the dual of \cite[Proposition 6.37]{Aguiar}, we know that $\left( \cComon\left( \Cc^{\bullet }\right)
,\circ \right) $ is monoidal so that $C\circ H$ is a cocommutative comonoid.
Since $\left( \Comon\left( \mathcal{C}^{\bullet }\right) ,\circ \right) $ is
monoidal, we have that $C\circ \varepsilon _{H}^{\bullet }$ and $g\circ H$
are in $\Comon\left( \mathcal{C}^{\bullet }\right) .$ 

As a consequence, for $D:=C\circ H$, we have
\begin{align*}
\left( D \bullet \sigma _{D,D}\bullet
D\right) \left( \Delta _{D}^{\bullet }\bullet \Delta _{D}^{\bullet }\right) \Delta _{D}^{\bullet }
&=
\left( D \bullet \sigma _{D,D}\Delta _{D}^{\bullet }\bullet
D\right) \left( \Delta _{D}^{\bullet }\bullet D\right) \Delta _{D}^{\bullet }\\
&=
\left( D \bullet \Delta _{D}^{\bullet }\bullet
D\right) \left( \Delta _{D}^{\bullet }\bullet D\right) \Delta _{D}^{\bullet }
=
\left( \Delta _{D}^{\bullet }\bullet \Delta _{D}^{\bullet }\right) \Delta _{D}^{\bullet }.\end{align*}
Therefore, we have
\begin{eqnarray*}
&&\Delta _{\left( C\circ J\right) \bullet H}^{\bullet }g^{r} \\
&=&\left( \left(
C\circ J\right) \bullet \sigma _{C\circ J,H}\bullet H\right) \left( \Delta
_{C\circ J}^{\bullet }\bullet \Delta _{H}^{\bullet }\right) \left( \left(
C\circ \varepsilon _{H}^{\bullet }\right) \bullet m_{H}^{\circ }\left(
g\circ H\right) \right) \Delta _{C\circ H}^{\bullet } \\
&=&\left( \left( C\circ J\right) \bullet \sigma _{C\circ J,H}\bullet
H\right) \left( \left( C\circ \varepsilon _{H}^{\bullet }\right) \bullet
\left( C\circ \varepsilon _{H}^{\bullet }\right) \bullet m_{H}^{\circ
}\left( g\circ H\right) \bullet m_{H}^{\circ }\left( g\circ H\right) \right)
\left( \Delta _{D}^{\bullet }\bullet \Delta _{D}^{\bullet
}\right) \Delta _{D}^{\bullet } \\
&=&\left( \left( C\circ \varepsilon _{H}^{\bullet }\right) \bullet
m_{H}^{\circ }\left( g\circ H\right) \bullet \left( C\circ \varepsilon
_{H}^{\bullet }\right) \bullet m_{H}^{\circ }\left( g\circ H\right) \right)
\underbracket[0.140ex]{\left(D \bullet \sigma _{D,D}\bullet
D\right) \left( \Delta _{D}^{\bullet }\bullet \Delta _{D}^{\bullet }\right) \Delta _{D}^{\bullet }} \\
&=&\left( \left( C\circ \varepsilon _{H}^{\bullet }\right) \bullet
m_{H}^{\circ }\left( g\circ H\right) \bullet \left( C\circ \varepsilon
_{H}^{\bullet }\right) \bullet m_{H}^{\circ }\left( g\circ H\right) \right)
\left( \Delta _{D}^{\bullet }\bullet \Delta _{D}^{\bullet
}\right) \Delta _{D}^{\bullet }
=\left( g^{r}\bullet g^{r}\right) \Delta _{C\circ H}^{\bullet }
\end{eqnarray*}%
and%
\begin{eqnarray*}
\varepsilon _{\left( C\circ J\right) \bullet H}^{\bullet }g^{r}
&=&m_{J}^{\bullet }\left( \varepsilon _{C\circ J}^{\bullet }\bullet
\varepsilon _{H}^{\bullet }\right) \left( \left( C\circ \varepsilon
_{H}^{\bullet }\right) \bullet m_{H}^{\circ }\left( g\circ H\right) \right)
\Delta _{C\circ H}^{\bullet } \\
&=&m_{J}^{\bullet }\left( \varepsilon _{C\circ H}^{\bullet }\bullet
\varepsilon _{C\circ H}^{\bullet }\right) \Delta _{C\circ H}^{\bullet
}=m_{J}^{\bullet }\Delta _{J}^{\bullet }\varepsilon _{C\circ H}^{\bullet
}=\varepsilon _{C\circ H}^{\bullet },
\end{eqnarray*}
hence $g^r$ is in $\Comon(\Cc^{\bullet})$.
\end{proof}

%\begin{invisible}
%\rd{[Il seguente è la versione duoidale di \cref{thm:spliepistability}.]}
%\end{invisible}

\begin{theorem}[Pullback of points along morphisms]
\label{thm:ddlstab}
Let $H$ and $H'$ be Hopf monoids. Then each point $\xymatrix{A\ar@<.5ex>[r]^-\pi&\ar@<.5ex>[l]^-\sigma H}$  as in \cref{thm:main}
has a pullback  along every morphism $f:H'\to H$ in $\cBimon \left(\Cc,\circ,\bullet\right) $.
\end{theorem}

\begin{proof} Since $H$ is a Hopf monoid, by \cref{coro:main} we have the isomorphism $\varphi:=m_A^\circ(k\circ \sigma):K\circ_\psi H\to A$ in $\cBimon(\Cc,\circ,\bullet)$.
    Set $p_{2}:=l_{H}^{\circ }\left( \varepsilon _{K}^{\circ }\circ H\right)
:K\circ H\rightarrow H$ and $i_{2}:=\left( u_{K}^{\circ }\circ H\right)
\left( l_{H}^{\circ }\right) ^{-1}:H\to K\circ H$. Then
\begin{eqnarray*}
\pi \varphi  &=&\pi m_{A}^{\circ }\left( k\circ \sigma \right) =m_{H}^{\circ
}\left( \pi k\circ \pi \sigma \right) =m_{H}^{\circ }\left( u_{H}^{\circ
}\varepsilon _{K}^{\circ }\circ H\right) =l_{H}^{\circ }\left( \varepsilon
_{K}^{\circ }\circ H\right) =p_{2} \\
\varphi i_{2} &=&m_{A}^{\circ }\left( k\circ \sigma \right) \left(
u_{K}^{\circ }\circ H\right) \left( l_{H}^{\circ }\right) ^{-1}=m_{A}^{\circ
}\left( u_{A}^{\circ }\circ A\right) \left( l_{A}^{\circ }\right)
^{-1}\sigma =\sigma .
\end{eqnarray*}%
As a consequence, $p_2=\pi\varphi$ and $i_2=\varphi^{-1}\sigma$ are morphisms of bimonoids and the point $\xymatrix{A\ar@<.5ex>[r]^-\pi&\ar@<.5ex>[l]^-\sigma H}$ is isomorphic to $\xymatrix{K\circ_\psi H\ar@<.5ex>[r]^-{p_2}&\ar@<.5ex>[l]^-{i_2} H}$ in $\Bimon \left(\Cc,\circ,\bullet\right) $ so that we can work with the latter.

Set $p_{2}':=l_{H'}^{\circ }\left( \varepsilon _{K}^{\circ }\circ H'\right)
:K\circ H'\rightarrow H'$ and let us prove that the following commutative diagram is a pullback in $\cComon\left(
\Cc^{\bullet }\right) $.%
\[\xymatrix{K\circ H^{\prime } \ar[d]_{K\circ f}\ar[r]^-{p'_2}& H'\ar[d]^f\\K\circ H\ar[r]^-{p_2}& H}\]
Let $\delta _{1}:D\rightarrow K\circ H$ and $\delta _{2}:D\rightarrow
H^{\prime }$ be morphisms in $\Comon(\Cc,\circ,\bullet)$ such that $p_{2}\delta
_{1}=f\delta _{2}$. Then, since $H'$ is a Hopf monoid and $u_{H^{\prime }}^{\circ }\varepsilon _{K}^{\circ }:K\to H'$ is a morphism of comonoids, we can define%
\[\delta: \xymatrix@C=2cm{D\ar[r]^-{\left( \delta _{1}\bullet \delta _{2}\right) \Delta
_{D}^{\bullet }}&\left( K\circ H\right) \bullet H^{\prime }\ar[r]^-{\left( K\circ \varepsilon _{H}^{\bullet }\right) \bullet H^{\prime }}&\left( K\circ J\right) \bullet H^{\prime }
\ar[r]^-{\left(
\left( u_{H^{\prime }}^{\circ }\varepsilon _{K}^{\circ }\right) ^{r}\right)
^{-1}}&K\circ H^{\prime }
}.\]
By \cref{lem:urinv}, we have that $\left( u_{H^{\prime }}^{\circ }\right)
^{r}=\left( \left( l_{J}^{\circ }\right) ^{-1}\bullet H^{\prime }\right)
\left( l_{H^{\prime }}^{\bullet }\right) ^{-1}l_{H^{\prime }}^{\circ },$ so
that $l_{H^{\prime }}^{\circ }\left( \left( u_{H^{\prime }}^{\circ }\right)
^{r}\right) ^{-1}=l_{H^{\prime }}^{\bullet }\left( l_{J}^{\circ }\bullet
H^{\prime }\right) .$
As a consequence, since $l_{J}^{\circ }\left( \varepsilon _{K}^{\circ }\circ J\right)
=m_{J}^{\circ }\left( u_{J}^{\circ }\varepsilon _{K}^{\circ }\circ J\right)
=m_{J}^{\circ }\left( \varepsilon _{I}^{\bullet }\varepsilon _{K}^{\circ
}\circ \varepsilon _{J}^{\bullet }\right) =m_{J}^{\circ }\left( \varepsilon
_{K}^{\bullet }\circ \varepsilon _{J}^{\bullet }\right) =\varepsilon
_{K\circ J}^{\bullet }$, we get
\begin{eqnarray*}
p_{2}^{\prime }\left( \left( u_{H^{\prime }}^{\circ }\varepsilon _{K}^{\circ
}\right) ^{r}\right) ^{-1} &=&l_{H^{\prime }}^{\circ }%
\underbracket[0.140ex]{\left( \varepsilon _{K}^{\circ }\circ H^{\prime
}\right) \left( \left( u_{H^{\prime }}^{\circ }\varepsilon _{K}^{\circ
}\right) ^{r}\right) ^{-1}}\overset{\eqref{form:grnat2}}{=}\underbracket[0.140ex]{l_{H^{\prime }}^{\circ }\left( \left( u_{H^{\prime
}}^{\circ }\right) ^{r}\right) ^{-1}}\left( \left( \varepsilon _{K}^{\circ
}\circ J\right) \bullet H^{\prime }\right)\\ &=&l_{H^{\prime }}^{\bullet }\left(
\underbracket[0.140ex]{l_{J}^{\circ }\left( \varepsilon _{K}^{\circ }\circ J\right)} \bullet
H^{\prime }\right)  =l_{H^{\prime }}^{\bullet }\left( \varepsilon _{K\circ J}^{\bullet
}\bullet H^{\prime }\right)
\end{eqnarray*}
and hence
\begin{eqnarray*}
p_{2}^{\prime }\delta  &=&\underbracket[0.140ex]{p_{2}^{\prime } \left( \left( u_{H^{\prime
}}^{\circ }\varepsilon _{K}^{\circ }\right) ^{r}\right) ^{-1}}\left( \left(
K\circ \varepsilon _{H}^{\bullet }\right) \bullet H^{\prime }\right) \left(
\delta _{1}\bullet \delta _{2}\right) \Delta _{D}^{\bullet } \\
&=&l_{H^{\prime }}^{\bullet }\left( \varepsilon _{K\circ
J}^{\bullet }\bullet H^{\prime }\right) \left( \left( K\circ
\varepsilon _{H}^{\bullet }\right) \bullet H^{\prime }\right) \left( \delta
_{1}\bullet \delta _{2}\right) \Delta _{D}^{\bullet }  \\
&=&\delta _{2}l_{D}^{\bullet }\left( \varepsilon _{K\circ H}^{\bullet
}\bullet D\right) \left( \delta _{1}\bullet D\right) \Delta _{D}^{\bullet }
=\delta _{2}l_{D}^{\bullet }\left( \varepsilon _{D}^{\bullet }\bullet
D\right) \Delta _{D}^{\bullet }=\delta _{2}.
\end{eqnarray*}%
Moreover%
\begin{eqnarray*}
\left( K\circ f\right) \delta  &=&\underbracket[0.140ex]{\left( K\circ f\right) \left(
\left( u_{H^{\prime }}^{\circ }\varepsilon _{K}^{\circ }\right) ^{r}\right)
^{-1}}\left( \left( K\circ \varepsilon _{H}^{\bullet }\right) \bullet
H^{\prime }\right) \left( \delta _{1}\bullet \delta _{2}\right) \Delta
_{D}^{\bullet } \\
&\overset{\eqref{form:grnat1}}{=}&\left( \left( fu_{H^{\prime }}^{\circ
}\varepsilon _{K}^{\circ }\right) ^{r}\right) ^{-1}\left( \left( K\circ
J\right) \bullet f\right) \left( \left( K\circ \varepsilon _{H}^{\bullet
}\right) \bullet H^{\prime }\right) \left( \delta _{1}\bullet \delta
_{2}\right) \Delta _{D}^{\bullet } \\
&=&\left( \left( u_{H}^{\circ }\varepsilon _{K}^{\circ }\right) ^{r}\right)
^{-1}\left( \left( K\circ \varepsilon _{H}^{\bullet }\right) \bullet
H\right) \left( \delta _{1}\bullet f\delta _{2}\right) \Delta _{D}^{\bullet }
\\
&=&\left( \left( u_{H}^{\circ }\varepsilon _{K}^{\circ }\right) ^{r}\right)
^{-1}\left( \left( K\circ \varepsilon _{H}^{\bullet }\right) \bullet
H\right) \left( \delta _{1}\bullet p_{2}\delta _{1}\right) \Delta
_{D}^{\bullet } \\
&=&\left( \left( u_{H}^{\circ }\varepsilon _{K}^{\circ }\right) ^{r}\right)
^{-1}\left( \left( K\circ \varepsilon _{H}^{\bullet }\right) \bullet
l_{H}^{\circ }\left( \varepsilon _{K}^{\circ }\circ H\right) \right) \Delta
_{K\circ H}^{\bullet }\delta _{1} \\
&=&\left( \left( u_{H}^{\circ }\varepsilon _{K}^{\circ }\right) ^{r}\right)
^{-1}\left( \left( K\circ \varepsilon _{H}^{\bullet }\right) \bullet
m_{H}^{\circ }\left( u_{H}^{\circ }\varepsilon _{K}^{\circ }\circ H\right)
\right) \Delta _{K\circ H}^{\bullet }\delta _{1} \\
&=&\left( \left( u_{H}^{\circ }\varepsilon _{K}^{\circ }\right) ^{r}\right)
^{-1}\left( u_{H}^{\circ }\varepsilon _{K}^{\circ }\right) ^{r}\delta
_{1}=\delta _{1}.
\end{eqnarray*}%
Let us check that $\delta $ is a morphism in $\cComon(\Cc^\bullet)$. Since $\Cc^{\bullet }$
is braided, then $\left( \Comon(\Cc^\bullet) ,\bullet
\right) $ is monoidal and hence $\delta _{1}\bullet \delta _{2}$ is in $\Comon(\Cc^\bullet).$ Since $D$ is
cocommutative, we have that $\Delta _{D}^{\bullet }$ is a morphism in $\Comon(\Cc^\bullet).$ Moreover $\left( K\circ
\varepsilon _{H}^{\bullet }\right) \bullet H^{\prime }$ is in $\Comon(\Cc^\bullet)$ as $K\circ \varepsilon _{H}^{\bullet }$ is therein because $\left(
\Comon\left( \Cc^{\bullet }\right) ,\circ \right) $ is monoidal. Finally, by \cref{lem:grcomon}, we know
that  $\left( u_{H^{\prime }}^{\circ }\varepsilon _{K}^{\circ }\right)
^{r}:K\circ H^{\prime }\rightarrow \left( K\circ J\right) \bullet H^{\prime }
$ is in $\Comon\left( \mathcal{C}^{\bullet }\right) .$ Hence $\delta $ is a
morphism of comonoids as a composition of morphisms of comonoids. Note that $K\circ H'$ is cocommutative as both $K$ and $H'$ are cocommutative and $(\cComon(\Cc^\bullet),\circ)$ is monoidal.
Thus $\delta$ is a morphism in $\cComon(\Cc^\bullet)$.

Let us prove the uniqueness of $\delta .$ Let $\delta ^{\prime
}:D\rightarrow K\circ H^{\prime }$ be a morphism of comonoids such that $%
p_{2}^{\prime }\delta' =\delta _{2}$ and $\left( K\circ f\right) \delta'
=\delta _{1}.$ Then%
\begin{eqnarray*}
\delta ^{\prime } &=&\left( \left( u_{H^{\prime }}^{\circ }\varepsilon
_{K}^{\circ }\right) ^{r}\right) ^{-1}\left( u_{H^{\prime }}^{\circ
}\varepsilon _{K}^{\circ }\right) ^{r}\delta ^{\prime } \\
&=&\left( \left( u_{H^{\prime }}^{\circ }\varepsilon _{K}^{\circ }\right)
^{r}\right) ^{-1}\left( \left( K\circ \varepsilon _{H^{\prime }}^{\bullet
}\right) \bullet m_{H^{\prime }}^{\circ }\left( u_{H^{\prime }}^{\circ
}\varepsilon _{K}^{\circ }\circ H^{\prime }\right) \right) \Delta _{K\circ
H^{\prime }}^{\bullet }\delta ^{\prime } \\
&=&\left( \left( u_{H^{\prime }}^{\circ }\varepsilon _{K}^{\circ }\right)
^{r}\right) ^{-1}\left( \left( K\circ \varepsilon _{H^{\prime }}^{\bullet
}\right) \bullet l_{H^{\prime }}^{\circ }\left( \varepsilon _{K}^{\circ
}\circ H^{\prime }\right) \right) \left( \delta ^{\prime }\bullet \delta
^{\prime }\right) \Delta _{D}^{\bullet } \\
&=&\left( \left( u_{H^{\prime }}^{\circ }\varepsilon _{K}^{\circ }\right)
^{r}\right) ^{-1}\left( \left( K\circ \varepsilon _{H}^{\bullet }\right)
\bullet H^{\prime }\right) \left( \left( K\circ f\right) \delta ^{\prime
}\bullet p_{2}^{\prime }\delta ^{\prime }\right) \Delta _{D}^{\bullet } \\
&=&\left( \left( u_{H^{\prime }}^{\circ }\varepsilon _{K}^{\circ }\right)
^{r}\right) ^{-1}\left( \left( K\circ \varepsilon _{H}^{\bullet }\right)
\bullet H^{\prime }\right) \left( \delta _{1}\bullet \delta _{2}\right)
\Delta _{D}^{\bullet }=\delta .
\end{eqnarray*}%
We have so proved that the diagram above is a pullback in $\cComon\left(
\mathcal{C}^{\bullet }\right) .$ Since the forgetful functor $%
\cBimon\left( \mathcal{C},\circ ,\bullet \right) =\Mon\left(
\cComon\left( \mathcal{C}^{\bullet }\right) ,\circ \right) \rightarrow
\cComon\left( \mathcal{C}^{\bullet }\right) $ creates limits, we can
endow $K\circ H^{\prime }$ with a structure of cocommutative bimonoid such
that the diagram above is also a pullback in $\cBimon\left( \mathcal{C}%
,\circ ,\bullet \right) .$
Now, the morphism $u^{\circ}_{K}$ is in $\mathsf{Comon}(\Cc^{\bullet})$ since $K$ is in $\mathsf{Bimon}(\Cc,\circ,\bullet)$ and then, since $(\mathsf{Comon}(\Cc^{\bullet}),\circ)$ is a monoidal category, we have that $i_2':=\left( u_{K}^{\circ }\circ H^{\prime }\right) \left(
l_{H^{\prime }}^{\circ }\right) ^{-1}:H'\to K\circ H'$ is in $\mathsf{Comon}(\Cc^{\bullet})$. We have
\begin{eqnarray*}
p_{2}^{\prime }i_2'&=&l_{H^{\prime }}^{\circ }\left(
\varepsilon _{K}^{\circ }\circ H^{\prime }\right) \left( u_{K}^{\circ }\circ
H^{\prime }\right) \left( l_{H^{\prime }}^{\circ }\right) ^{-1}=\mathrm{Id}%
_{H^{\prime }} \\
\left( K\circ f\right) i_2'&=&\left( K\circ f\right) \left(
u_{K}^{\circ }\circ H^{\prime }\right) \left( l_{H^{\prime }}^{\circ
}\right) ^{-1}=\left( u_{K}^{\circ }\circ H\right) \left( l_{H}^{\circ
}\right) ^{-1}f=i_{2}f.
\end{eqnarray*}
As a consequence, $i_{2}^{\prime }$ is the unique morphism of comonoids such that $p_{2}^{\prime }i_{2}^{\prime }=\mathrm{Id}%
_{H^{\prime }}$ and $\left( K\circ f\right) i_{2}^{\prime }=i_{2}f.$ Since our pullback is also a pullback of bimonoids and $\id_{H'}$ and $i_2f$ are morphisms of bimonoids, we get that $i_2'$ must be a morphism of bimonoids too. Thus we get that $\xymatrix{K\circ H'\ar@<.5ex>[r]^-{p_2'}&\ar@<.5ex>[l]^-{i_2'} H'}$ is a point in $\Bimon \left(\Cc,\circ,\bullet\right) $.
\begin{invisible}
Moreover, we have
\[
\begin{split}
(K\circ f)m^{\circ}_{K\circ H'}(i'_{2}\circ i'_{2})&=m^{\circ}_{K\circ H}(K\circ f\circ K\circ f)(i'_{2}\circ i'_{2})\\&=m^{\circ}_{K\circ H}(K\circ f\circ K\circ f)((u_{K}^{\circ}\circ H')(l_{H'}^{\circ})^{-1}\circ (u_{K}^{\circ}\circ H')(l_{H'}^{\circ})^{-1})\\&=m^{\circ}_{K\circ H}(u_{K}^{\circ}\circ H\circ u_{K}^{\circ}\circ H)((I\circ f)(l^{\circ}_{H'})^{-1}\circ (I\circ f)(l^{\circ}_{H'})^{-1})\\&=m^{\circ}_{K\circ H}(u_{K}^{\circ}\circ H\circ u_{K}^{\circ}\circ H)((l^{\circ}_{H})^{-1}\circ (l^{\circ}_{H})^{-1})(f\circ f)\\&=m^{\circ}_{K\circ H}(i_{2}\circ i_{2})(f\circ f)=i_{2}m^{\circ}_{H}(f\circ f)=i_{2}fm^{\circ}_{H'}\\&=(K\circ f)i_{2}'m^{\circ}_{H'}
\end{split}
\]
and also
\[
p'_{2}m^{\circ}_{K\circ H'}(i'_{2}\circ i'_{2})=m^{\circ}_{H'}(p'_{2}\circ p'_{2})(i'_{2}\circ i'_{2})=m^{\circ}_{H'}=p'_{2}i'_{2}m_{H'}^{\circ}
\]
hence by the universal property of the pullback we obtain that $m_{K\circ H'}^{\circ}(i'_{2}\circ i'_{2})=i'_{2}m_{H'}^{\circ}$. Moreover, we have
\[
\begin{split}
(K\circ f)u^{\circ}_{K\circ H'}&=u^{\circ}_{K\circ H}=(u_{K}^{\circ}\circ u_{H}^{\circ})(m_{I}^{\circ})^{-1}=(K\circ f)(u_{K}^{\circ}\circ H')(I\circ u_{H'}^{\circ})(m_{I}^{\circ})^{-1}\\&=(K\circ f)(u_{K}^{\circ}\circ H')(l_{H'}^{\circ})^{-1}u_{H'}^{\circ}=(K\circ f)i'_{2}u_{H'}^{\circ}
\end{split}
\]
and also
\[
\begin{split}
p'_{2}u^{\circ}_{K\circ H'}&=u^{\circ}_{H'}=p'_{2}i'_{2}u_{H'}^{\circ}
\end{split}
\]
then by the universal property of the pullback, we get $u^{\circ}_{K\circ H'}=i'_{2}u_{H'}^{\circ}$, hence $i'_{2}$ is in $\mathsf{Mon}(\Cc^{\circ})$.
\end{invisible}
\end{proof}

\begin{remark}%[\rd{revised:2026/03/13}]
\label{rmk:pulbantip}
 In the setting of \cref{thm:ddlstab}, we have that $fu_{H'}^\circ=u_{H}^\circ$. Since by assumption we have that $u_{H}^\circ$ is a monomorphism in $\Cc$, then so is $u_{H'}^\circ$. Therefore, also the point $\xymatrix{K\circ H'\ar@<.5ex>[r]^-{p'_2}&\ar@<.5ex>[l]^-{i'_2} H'}$ verifies the assumptions of \cref{thm:main}. We want to take a deeper look at the structure of $K\circ H'$. To this aim we first identify the kernel of the projection $p_2'$.
 Set $i_{1}:=\left( K\circ u_{H}^{\circ }\right) \left( r_{K}^{\circ }\right)
^{-1}:K\rightarrow K\circ H.$ Then
\begin{equation*}
\varphi i_{1}=m_{A}^{\circ }\left( k\circ \sigma \right) \left( K\circ
u_{H}^{\circ }\right) \left( r_{K}^{\circ }\right) ^{-1}=m_{A}^{\circ
}\left( A\circ u_{A}^{\circ }\right) \left( r_{A}^{\circ }\right) ^{-1}k=k.
\end{equation*}%
In view of the isomorphism between the points $\xymatrix{A\ar@<.5ex>[r]^-\pi&\ar@<.5ex>[l]^-\sigma H}$ and $\xymatrix{K\circ_\psi H\ar@<.5ex>[r]^-{p_2}&\ar@<.5ex>[l]^-{i_2} H}$ in $\Bimon \left(\Cc,\circ,\bullet\right) $, the above equality means that $\left( K,i_{1}\right) $ is the kernel of $p_{2}$ in $\cBimon\left(
\mathcal{C},\circ ,\bullet \right) .$ 

Set $i_{1}^{\prime }:=\left( K\circ
u_{H^{\prime }}^{\circ }\right) \left( r_{K}^{\circ }\right)
^{-1}:K\rightarrow K\circ H^{\prime }$ and note that
\begin{align*}
p_2'i_1'&=l_H^\circ(\varepsilon^\circ_K\circ H')\left( K\circ
u_{H^{\prime }}^{\circ }\right) \left( r_{K}^{\circ }\right)
^{-1}=u_{H^{\prime }}^{\circ }\varepsilon^\circ_K,\\
\left( K\circ f\right) i_{1}^{\prime }&=\left( K\circ f\right) \left( K\circ
u_{H^{\prime }}^{\circ }\right) \left( r_{K}^{\circ }\right) ^{-1}=\left(
K\circ u_{H}^{\circ }\right) \left( r_{K}^{\circ }\right) ^{-1}=i_{1}.
\end{align*}%
Then, as we did in the case of $i_2'$, we deduce that $i_1'$ is the unique morphism of comonoids such that $p_2'i_1'=u_{H^{\prime }}^{\circ }\varepsilon^\circ_K$ and $\left( K\circ f\right) i_{1}^{\prime }=i_{1}$ and that in fact it is a morphism of bimonoids.
Consider the following diagram
    \[\xymatrix{K\ar[r]^-{i_1'}\ar[d]^{\varepsilon_K^\circ}&K\circ H'\ar[r]^{K\circ f}\ar[d]^{p_2'}&K\circ H\ar[d]^{p_2}\\ I\ar[r]^{u_{H'}^\circ}& H'\ar[r]^{f}&H}\]
By pullback pasting law, since the right square is a pullback and the total
rectangle is the pullback expressing that $(K,i_1)$ is the kernel of $p_2$, we conclude that the left square is a pullback.
This means that $\left( K,i_{1}^{\prime }\right) $ is the kernel of $p_{2}^{\prime
}$ in $\cBimon\left( \mathcal{C},\circ ,\bullet \right) .$
By the first part of \cref{coro:main}, we have the isomorphism of bimonoids $\varphi':=m_{K\circ H^{\prime
}}^{\circ }\left( i_{1}^{\prime }\circ i_{2}^{\prime }\right):K\circ_{\psi'} H'\to K\circ
H'$ where $\psi' :=(\varphi') ^{-1}m_{K\circ H'}^{\circ
}\left( i_2' \circ i_1'\right) :H'\circ K\rightarrow K\circ H'$. 
If $A$ is a Hopf monoid, so is $K$ by \cref{thm:main}. Thus, since $H'$ is a Hopf monoid too,
by \cref{thm:psitriang}, we get that $K\circ H'$ has antipode $\Sigma
_{K\circ H'}:=(\psi') ^{\revsn }\phi _{K,H'}^{\circ }\left( \Sigma _{K}\circ \Sigma
_{H'}\right) .$
Moreover
\begin{eqnarray*}
p_{2}^{\prime }\varphi ^{\prime } &=&p_{2}^{\prime }m_{K\circ H^{\prime
}}^{\circ }\left( i_{1}^{\prime }\circ i_{2}^{\prime }\right) =m_{H^{\prime
}}^{\circ }\left( p_{2}^{\prime }\circ p_{2}^{\prime }\right) \left(
i_{1}^{\prime }\circ i_{2}^{\prime }\right) =m_{H^{\prime }}^{\circ }\left(
u^{\circ}_{H^{\prime }}\varepsilon _{K}^{\circ }\circ H^{\prime }\right)
=l_{H^{\prime }}^{\circ }\left( \varepsilon _{K}^{\circ }\circ H^{\prime
}\right) =p_{2}^{\prime } \\
\left( K\circ f\right) \varphi ^{\prime } &=&\left( K\circ f\right)
m_{K\circ H^{\prime }}^{\circ }\left( i_{1}^{\prime }\circ i_{2}^{\prime
}\right) =m_{K\circ H}^{\circ }\left( K\circ f\circ K\circ f\right) \left(
i_{1}^{\prime }\circ i_{2}^{\prime }\right) =m_{K\circ H}^{\circ }\left(
i_{1}\circ i_{2}\right) \left( K\circ f\right) \\&=&m_{K\circ H}^{\circ }\left(
\varphi ^{-1}k\circ \varphi ^{-1}\sigma \right) \left( K\circ f\right)
=\varphi ^{-1}m_{A}^{\circ }\left( k\circ \sigma \right) \left( K\circ
f\right) =\varphi ^{-1}\varphi \left( K\circ f\right) =K\circ f
\end{eqnarray*}%
so that $\varphi'=\id$ and hence the Hopf monoid structure of $K\circ H^{\prime }$ is precisely the one of $K\circ _{\psi ^{\prime }}H^{\prime }$.
\end{remark}

\subsection{Protomodularity}

The results obtained so far on duoidal categories with a reversion will enable us to establish protomodularity for a broad class of cocommutative Hopf monoids. In particular, we will apply these findings in the context of quasitriangular Hopf algebras.

We recall some preliminary notions and results about protomodular categories, which have been introduced in \cite[Definition 4]{Bourn}, see also \cite[Definition 3.1.3]{BB04}.

Given a point \begin{tikzcd}
	A& H
	\arrow[shift left, from=1-1, to=1-2,"\pi"]
	\arrow[shift left, from=1-2, to=1-1,"\sigma"]
\end{tikzcd}
in a category $\Cc$, the codomain assignment gives rise to a fibration $F:\mathrm{Pt}(\Cc)\to\Cc$, the so-called \textit{fibration of points}. A category $\Cc$ is \textit{protomodular} when:
\begin{itemize}
    \item[1)] $\Cc$ has pullbacks of split epimorphisms along any morphism;
    \item[2)] all the inverse image functors of the fibration $F:\mathrm{Pt}(\Cc)\to\Cc$ of points reflect
isomorphisms.
\end{itemize}

It is known that, if $\Cc$ is pointed (i.e.\ it has
a zero object) and finitely complete, the protomodularity of $\Cc$ is equivalent to the fact that the Split Short Five Lemma holds in $\Cc$, see e.g. \cite[Proposition 3.1.2]{BB04}. When $\Cc$ only admits an initial object $\initial$, for any object $A$ in $\Cc$, one can consider the change of base along the initial morphism $\initial\to A$. This defines a kernel functor $\mathcal{K}_{A}:\mathrm{Pt}_{A}(\Cc)\to\Cc$. In the presence of an initial object, the protomodularity condition can be equivalently described  by requiring that just kernel functors are conservative, together with the fact that $\Cc$ has pullbacks of split epimorphisms along any morphism in $\Cc$. \medskip\newline
We now define our candidate protomodular category.

\begin{definition}
 Given a $\bullet$-braided duoidal category $\left( \mathcal{C},\circ ,I,\bullet ,J\right) $ with a reversion, we define the full subcategory $\cHopfmon^{\mathrm{m}}(\Cc,\circ, \bullet )$ of $\cHopfmon(\Cc,\circ,\bullet)$ whose objects are cocommutative Hopf monoids $H$ whose unit morphism $u_{H}^{\circ}:I\to H$ (that is the initial morphism) is a  monomorphism in $\Cc$.
\end{definition}
%\rd{[Possibile notazione alternativa: la si potrebbe anche indicare con $\cHopfmon^\Mm(\Cc,\bullet,\circ)$ intendendo che $u_H^\circ\in\Mm$ dove $\Mm$ è la classe dei morphismi in $\cHopfmon(\Cc,\circ,\bullet)$ che sono mono in $\Cc$.]}

We now show that, under suitable assumptions on $\Cc$, the above category is protomodular.

% \begin{lemma}\label{lem:split-Hopf}
%     Let $\left( \mathcal{C},\circ ,I,\bullet ,J\right) $ be a $\bullet$-braided duoidal category with a reversion. Assume that $\mathcal{C}$ has and $\bullet$ preserves binary intersections.
% Then, the Split Short Five Lemma holds in $\cHopfmon^{\mathrm{m}}(\Cc,\circ, \bullet )$.
% \end{lemma}
% \begin{lemma}\label{lem:pulb-hopf}
%     Let $\left( \mathcal{C},\circ ,I,\bullet ,J\right) $ be a $\bullet $-braided duoidal category with a reversion. Then, each point in $\cHopfmon^{\mathrm{m}}(\Cc,\circ, \bullet )$ has pullback along any morphism.
% \end{lemma}

\begin{theorem}
\label{thm:protomon}
  Let $\left( \mathcal{C},\circ ,I,\bullet ,J\right) $ be a $\bullet $-braided duoidal category with a reversion.  Assume that $\mathcal{C}$ has and $\bullet$ preserves binary intersections. Then, the category $\cHopfmon^{\mathrm{m}}(\Cc,\circ, \bullet )$ is protomodular.
\end{theorem}

\begin{proof}
First, we prove that the Split Short Five Lemma holds in $\cHopfmon^{\mathrm{m}}(\Cc,\circ, \bullet )$.
Consider the following commutative diagram in $\mathsf{Hopf}^\mathrm{m}_{\mathrm{coc}}(\Cc, \circ,\bullet)$
    \[\xymatrix@R=.6cm@C=1.5cm{K\ar[d]^f\ar[r]^-k
&A\ar[d]^\phi\ar@<.5ex>[r]^-\pi&\ar@<.5ex>[l]^-\sigma H\ar[d]^g\\ K'\ar[r]^-{k'}
&A'\ar@<.5ex>[r]^-{\pi'}&\ar@<.5ex>[l]^-{\sigma'} H'}\]
where $K$ and $K'$ are the kernels of $\pi$ and $\pi'$ in $\cHopfmon(\Cc,\circ,\bullet)$ (whence in $\cHopfmon^{\mathrm{m}}(\Cc,\circ,\bullet)$), respectively, as in \cref{thm:main}. If $f$ and $g$ are isomorphisms in $\cHopfmon^{\mathrm{m}}(\Cc,\circ, \bullet )$, hence in $\mathsf{Bimon}(\Cc,\circ,\bullet)$, then $\phi$ is an isomorphism in $\mathsf{Bimon}(\Cc,\circ,\bullet)$ by \cref{lem:split-short}, hence an isomorphism in  $\cHopfmon^{\mathrm{m}}(\Cc,\circ, \bullet )$.

Now, we prove that each point in $\cHopfmon^{\mathrm{m}}(\Cc,\circ, \bullet )$ has pullback along any morphism. According to \cref{thm:ddlstab}, a point  $\xymatrix{A\ar@<.5ex>[r]^-\pi&\ar@<.5ex>[l]^-\sigma H}$ in $\cHopfmon^{\mathrm{m}}(\Cc,\circ, \bullet )$, admits a pullback in $\cBimon(\Cc,\circ,\bullet)$ along any morphism $f:H'\to H$ in $\cHopfmon^{\mathrm{m}}(\Cc,\circ, \bullet )$. In view of \cref{rmk:pulbantip}, the pullback bimonoid $K\circ H'$ has an antipode and hence what we have obtained is in fact a pullback in $\cHopfmon^{\mathrm{m}}(\Cc,\circ, \bullet )$.

Therefore, the category $\cHopfmon^{\mathrm{m}}(\Cc,\circ, \bullet )$ is protomodular.
\end{proof}

\section{Examples and Applications}\label{sec:exapp}

\subsection{Duoidal categories from braided categories} 
In light of \cite[\S 5]{Bohm-Lack}, it is reasonable to expect that a braided monoidal category can be viewed as a duoidal category equipped with a reversion. This is the main focus of the present subsection.

\begin{example}\label{ex:braidedasdouoidal}
By \cite[Proposition 6.10]{Aguiar}, a braided monoidal category $(\Cc,\circ,I,\sigma)$ gives rise to a duoidal category $(\Cc,\circ,I,\circ,I)$ defining
\[
\zeta_{A,B,C,D}=%(a^{\circ}_{A,C,B\circ D})^{-1}(A\circ a^{\circ}_{C,B,D})(
A\circ\sigma_{B,C}\circ D%)(A\circ (a^{\circ}_{B,C,D})^{-1})a^{\circ}_{A,B,C\circ D}
,\quad \Delta_{I}=(l_{I})^{-1},\quad m_{I}=l_{I},\quad \varepsilon_{I}=\mathrm{Id}.
\]
In this case, the structure morphisms \eqref{def:zetaduoidal} and \eqref{def:deltamepsilonduoidal} are isomorphisms in $\Cc$, so the duoidal category is called \textit{strong}, see \cite[Definition 6.3]{Aguiar}.

The category $(\Cc,\circ^\rev,\sigma^\rev)$ is braided where $\sigma^\rev_{X,Y}:=\sigma_{Y,X}$ is the transpose braiding, see e.g.\ \cite[\S 1.1.2.]{Aguiar}. By applying \cite[Example 1.59]{BCPvO} to this braided category, we get that $(\id_\Cc,\phi):(\Cc,\circ^\rev)\to (\Cc,\circ)$ is a strong monoidal functor where $\phi_0=\id$ and $\phi_{X,Y}:=\sigma_{X,Y}$. By using the braid relations, it is straightforward to check that $(\id_\Cc,\phi^\circ,\phi^\bullet)$ is a double lax monoidal functor, where $\phi^\circ=\phi^\bullet:=\phi$.
\begin{invisible}
    We check the interchange law%
\begin{eqnarray*}
&&\phi _{C\circ A,D\circ B}^{\bullet }\left( \phi _{A,C}^{\circ }\bullet \phi
_{B,D}^{\circ }\right) \zeta _{A,B,C,D} =\underbracket[0.140ex]{\sigma _{C\circ
A,D\circ B}}\left( \sigma _{A,C}\circ \sigma _{B,D}\right) \left( A\circ
\sigma _{B,C}\circ D\right)  \\
&=&\left( \sigma _{C,D\circ B}\circ A\right) \underbracket[0.140ex]{\left( C\circ
\sigma _{A,D\circ B}\right) \left( \sigma _{A,C}\circ \sigma _{B,D}\right) }%
\left( A\circ \sigma _{B,C}\circ D\right)  \\
&=&\left( \sigma _{C,D\circ B}\circ A\right) \left( C\circ \sigma
_{B,D}\circ A\right) \left( C\circ \underbracket[0.140ex]{\sigma _{A,B\circ D}}\right)
\left( \sigma _{A,C}\circ B\circ D\right) \left( A\circ \sigma _{B,C}\circ
D\right)  \\
&=&\left( \sigma _{C,D\circ B}\circ A\right) \left( C\circ \sigma
_{B,D}\circ A\right) \left( C\circ B\circ \sigma _{A,D}\right) \left(
\underbracket[0.140ex]{C\circ \sigma _{A,B}}\circ D\right) \left( \underbracket[0.140ex]{\sigma
_{A,C}\circ B}\circ D\right) \left( \underbracket[0.140ex]{A\circ \sigma _{B,C}}\circ
D\right)  \\
&=&\left( \sigma _{C,D\circ B}\circ A\right) \underbracket[0.140ex]{\left( C\circ
\sigma _{B,D}\circ A\right) \left( C\circ B\circ \sigma _{A,D}\right) \left(
\sigma _{B,C}\circ A\circ D\right) \left( B\circ \sigma _{A,C}\circ D\right)
}\left( \sigma _{A,B}\circ C\circ D\right)  \\
&=&\left( \underbracket[0.140ex]{\sigma _{C,D\circ B}}\circ A\right) \sigma _{B\circ
A,C\circ D}\left( \sigma _{A,B}\circ C\circ D\right)  \\
&=&\left( D\circ \sigma _{C,B}\circ A\right) \underbracket[0.140ex]{\left( \sigma
_{C,D}\circ B\circ A\right) \sigma _{B\circ A,C\circ D}}\left( \sigma
_{A,B}\circ C\circ D\right)  \\
&=&\left( D\circ \sigma _{C,B}\circ A\right) \sigma _{B\circ A,D\circ
C}\left( \sigma _{A,B}\circ \sigma _{C,D}\right)  \\
&=&\zeta _{D,C,B,A}\phi _{B\bullet A,D\bullet C}^{\circ }\left( \phi
_{A,B}^{\bullet }\circ \phi _{C,D}^{\bullet }\right)
\end{eqnarray*}%
and unitality%
\begin{eqnarray*}
\phi _{I,I}^{\bullet }\left( \phi _{0}^{\circ }\bullet \phi _{0}^{\circ
}\right) \Delta _{I}^{\bullet } &=&\sigma _{I,I}\left( \mathrm{Id}\circ
\mathrm{Id}\right) \Delta _{I}=r_{I}^{-1}l_{I}l_{I}^{-1}=r_{I}^{-1}=\Delta
_{I}=\Delta _{I}^{\bullet }\phi _{0}^{\circ }, \\
m_{J}^{\circ }\phi _{J,J}^{\circ }\left( \phi _{0}^{\bullet }\circ \phi
_{0}^{\bullet }\right)  &=&m_{I}\sigma _{I,I}\left( \mathrm{Id}\circ \mathrm{%
Id}\right) =r_{I}r_{I}^{-1}l_{I}=l_{I}=\phi _{0}^{\bullet }m_{J}^{\circ }. \\
\phi _{0}^{\bullet }\varepsilon _{I}^{\bullet } &=&\mathrm{Id}=\varepsilon
_{I}^{\bullet }\phi _{0}^{\circ }.
\end{eqnarray*}
\end{invisible}
In order to show that this leads to a reversion on the duoidal category $(\Cc,\circ,\circ)$, define
\begin{eqnarray*}
\gamma _{X,Y,Z,T} &=&\zeta _{X,Y,Z,T}=X\circ \sigma _{Y,Z}\circ T:X\circ
Y\circ Z\circ T\rightarrow X\circ Z\circ Y\circ T, \\
\delta _{X,Y,Z,T} &=&\zeta _{X,Y,Z,T}:X\circ Y\circ Z\circ T\rightarrow
X\circ Z\circ Y\circ T, \\
\widetilde{\gamma }_{X,Y,Z,T} &=&\sigma_{X\circ Y,Z}\circ T:X\circ Y\circ Z\circ
T\rightarrow Z\circ X\circ Y\circ T, \\
\widetilde{\delta }_{X,Y,Z,T} &=&X\circ \sigma _{Y,Z\circ T}:X\circ Y\circ
Z\circ T\rightarrow X\circ Z\circ T\circ Y.
\end{eqnarray*}%
It is straightforward to verify that all diagrams in the definition of reversion commute, as a direct consequence of the properties of the braiding. \begin{invisible}
One could expect that this commutativity could follow by coherence theorem for braided categories, see \cite[Theorem 2.5]{Joyal-Street}. Unfortunately this theorem does not works as one may think as for instance $\sigma_{Y,X}\sigma_{X,Y}$ and $\id_{X\circ Y}$ have the same domain and codomain and are build using the braiding but do not coincide unless we have a symmetry. Thus we collect here all the equalities that are clearly commutative.
\begin{itemize}
\item \eqref{form:psi1} is $\left( X\circ \sigma _{Y,Z\circ U}\circ T\circ
V\right) \left( X\circ Y\circ Z\circ \sigma _{T,U}\circ V\right) =\left(
X\circ Z\circ \sigma _{Y\circ T,U}\circ V\right) \left( X\circ \sigma
_{Y,Z}\circ T\circ U\circ V\right) $.

\item \eqref{form:psi2} is $\left( X\circ \sigma _{Y,T}\circ U\circ Z\circ
V\right) \left( X\circ Y\circ \sigma _{Z,T\circ U}\circ V\right) =\left(
X\circ T\circ Y\circ \sigma _{Z,U}\circ V\right) \left( X\circ \sigma
_{Y\circ Z,T}\circ U\circ V\right) .$

\item \eqref{form:varphi1} is $\left( X\circ Z\circ \sigma _{Y\circ T,U}\circ
V\right) \left( X\circ \sigma _{Y,Z}\circ T\circ U\circ V\right) =\left(
X\circ \sigma _{Y,Z\circ U}\circ T\circ V\right) \left( X\circ Y\circ Z\circ
\sigma _{T,U}\circ V\right) $ i.e. \eqref{form:psi1}.

\item \eqref{form:varphi2} is $\left( X\circ T\circ Y\circ \sigma _{Z,U}\circ
V\right) \left( X\circ \sigma _{Y\circ Z,T}\circ U\circ V\right) =\left(
X\circ \sigma _{Y,T}\circ U\circ Z\circ V\right) \left( X\circ Y\circ \sigma
_{Z,T\circ U}\circ V\right) $ i.e. \eqref{form:psi2}.

\item \eqref{form:gammadelta} is $\left( X\circ \sigma _{Y,T}\circ U\circ
Z\circ V\right) \left( X\circ Y\circ \sigma _{Z,T\circ U}\circ V\right)
=\left( X\circ T\circ Y\circ \sigma _{Z,U}\circ V\right) \left( X\circ
\sigma _{Y\circ Z,T}\circ U\circ V\right) $ i.e. \eqref{form:psi2}.

\item \eqref{form:phi0gamma} is $l_{I\circ I}^{\circ }\left( I\circ I\circ
l_{I}^{\circ }\right) \left( I\circ \sigma _{I,I}\circ I\right) \left(
\Delta _{I}\circ I\circ I\right) =l_{I}^{\circ }\circ I$ which is true by
coherence theorem for monoidal categories as these morphisms are unit
constraints.

\item \eqref{form:phi0delta} is also true by coherence theorem for monoidal
categories as the morphisms involved are unit constraints.

\item \eqref{form:phipsi2} is 
$\left( X\circ Z\circ \sigma _{Y,I\circ I}\circ
U\right) \left( X\circ \sigma _{Y,Z}\circ I\circ I\circ U\right) \left(
X\circ Y\circ r_{Z}^{-1}\circ I\circ U\right) \\=\left( X\circ Z\circ I\circ
\sigma _{Y,I}\circ U\right) \left( X\circ Z\circ l_{Y}^{-1}\circ
l_{U}^{-1}\right) \left( X\circ l_{Z\circ Y\circ U}\right) \left( X\circ
\sigma _{Z\circ Y,I}\circ U\right) \left( X\circ \sigma _{Y,Z}\circ I\circ
U\right) $ 
which is equivalent to $\left( Z\circ \sigma _{Y,I\circ I}\circ
U\right) \left( \sigma _{Y,Z}\circ I\circ I\circ U\right) \left( Y\circ
r_{Z}^{-1}\circ I\circ U\right)\\ =\left( Z\circ I\circ \sigma _{Y,I}\circ
U\right) \left( Z\circ l_{Y}^{-1}\circ l_{U}^{-1}\right) l_{Z\circ Y\circ
U}\left( \sigma _{Z\circ Y,I}\circ U\right) \left( \sigma _{Y,Z}\circ I\circ
U\right) $ 

i.e. $\left( Z\circ \sigma _{Y,I\circ I}\circ U\right) \left(
\sigma _{Y,Z}\circ I\circ I\circ U\right) \left( r_{Y\circ Z}^{-1}\circ
I\circ U\right) \\=\left( Z\circ I\circ \sigma _{Y,I}\circ U\right) \left(
Z\circ l_{Y}^{-1}\circ l_{U}^{-1}\right) l_{Z\circ Y\circ U}\left( \sigma
_{Z\circ Y,I}\circ U\right) \left( \sigma _{Y,Z}\circ I\circ U\right) $ 

i.e.
$\left( Z\circ \sigma _{Y,I\circ I}\circ U\right) \left( r_{Z\circ
Y}^{-1}\circ I\circ U\right) \left( \sigma _{Y,Z}\circ I\circ U\right)
\\=\left( Z\circ I\circ \sigma _{Y,I}\circ U\right) \left( Z\circ
l_{Y}^{-1}\circ l_{U}^{-1}\right) l_{Z\circ Y\circ U}\left( \sigma _{Z\circ
Y,I}\circ U\right) \left( \sigma _{Y,Z}\circ I\circ U\right) $ 

i.e. $\left(
Z\circ \sigma _{Y,I\circ I}\circ U\right) \left( r_{Z\circ Y}^{-1}\circ
I\circ U\right) =\left( Z\circ I\circ \sigma _{Y,I}\circ U\right) \left(
Z\circ l_{Y}^{-1}\circ l_{U}^{-1}\right) l_{Z\circ Y\circ U}\left( \sigma
_{Z\circ Y,I}\circ U\right) $ which is true by coherence theorem for
monoidal categories as these morphisms are unit constraints.

\item \eqref{form:phivarphi2} is $\left( X\circ \sigma _{I\circ I,Z}\circ Y\circ U\right) \left( X\circ
I\circ I\circ \sigma _{Y,Z}\circ U\right) \left( X\circ I\circ
l_{Y}^{-1}\circ Z\circ U\right)\\ =\left( X\circ \sigma _{I,Z}\circ I\circ
Y\circ U\right) \left( r_{X}^{-1}\circ r_{Z}^{-1}\circ Y\circ U\right)
\left( r_{X\circ Z\circ Y}\circ U\right) \left( X\circ \sigma _{I,Z\circ
Y}\circ U\right) \left( X\circ I\circ \sigma _{Y,Z}\circ U\right) $

which is equivalent to $\left( X\circ \sigma _{I\circ I,Z}\circ Y\right)
\left( X\circ I\circ I\circ \sigma _{Y,Z}\right) \left( X\circ I\circ
l_{Y}^{-1}\circ Z\right) \\=\left( X\circ \sigma _{I,Z}\circ I\circ Y\right)
\left( r_{X}^{-1}\circ r_{Z}^{-1}\circ Y\right) \left( r_{X\circ Z\circ
Y}\right) \left( X\circ \sigma _{I,Z\circ Y}\right) \left( X\circ I\circ
\sigma _{Y,Z}\right) $ 

i.e. $\left( X\circ \sigma _{I\circ I,Z}\circ
Y\right) \left( X\circ I\circ I\circ \sigma _{Y,Z}\right) \left( X\circ
I\circ l_{Y\circ Z}^{-1}\right) \\=\left( X\circ \sigma _{I,Z}\circ I\circ
Y\right) \left( r_{X}^{-1}\circ r_{Z}^{-1}\circ Y\right) \left( r_{X\circ
Z\circ Y}\right) \left( X\circ \sigma _{I,Z\circ Y}\right) \left( X\circ
I\circ \sigma _{Y,Z}\right) $

i.e. $\left( X\circ \sigma _{I\circ I,Z}\circ Y\right) \left( X\circ I\circ
l_{Z\circ Y}^{-1}\right) \left( X\circ I\circ \sigma _{Y,Z}\right) \\=\left(
X\circ \sigma _{I,Z}\circ I\circ Y\right) \left( r_{X}^{-1}\circ
r_{Z}^{-1}\circ Y\right) \left( r_{X\circ Z\circ Y}\right) \left( X\circ
\sigma _{I,Z\circ Y}\right) \left( X\circ I\circ \sigma _{Y,Z}\right) $ 

i.e.
$\left( X\circ \sigma _{I\circ I,Z}\circ Y\right) \left( X\circ I\circ
l_{Z\circ Y}^{-1}\right) =\left( X\circ \sigma _{I,Z}\circ I\circ Y\right)
\left( r_{X}^{-1}\circ r_{Z}^{-1}\circ Y\right) \left( r_{X\circ Z\circ
Y}\right) \left( X\circ \sigma _{I,Z\circ Y}\right) $ which is true by
coherence theorem for monoidal categories as these morphisms are unit
constraints.

\item \eqref{form:varphipsitilde} is $\widetilde{\varphi }_{X\circ Y,Z}\left(
\psi _{X,Y}\bullet Z^{\revsn }\right) =\widetilde{\psi }_{X,Y\circ Z}\left(
X^{\revsn }\bullet \varphi _{Y,Z}\right) $ i.e.
\begin{eqnarray*}
&&r_{X\circ Y\circ Z}\left( l_{X\circ Y\circ Z}\circ I\right) \left( I\circ
\sigma _{I,X\circ Y\circ Z}\right) \left( J\circ \sigma _{X,I}\circ Y\circ
Z\right) \left( l_{X}^{-1}\circ I\circ Y\circ Z\right) \left(
r_{X}^{-1}\circ  Y\circ Z\right)  \\
&=&l_{X\circ Y\circ Z}\left( I\circ r_{X\circ Y\circ Z}\right) \left(
c_{X\circ Y\circ Z,I}\circ I\right) \left( X\circ Y\circ \sigma _{I,Z}\circ
J\right) \left( X\circ Y\circ \left( I\circ r_{Z}^{-1}\right) \right) \left(
X\circ Y\circ l_{Y}^{-1}\right) .
\end{eqnarray*}%
which is true by coherence theorem for monoidal categories as these
morphisms are unit constraints.

\item \eqref{form:gammatildeIIJ} is also true by coherence theorem for
monoidal categories as the morphisms involved are unit constraints.

\item \eqref{form:deltatildeJII} is also true by coherence theorem for
monoidal categories as the morphisms involved are unit constraints.

\item \eqref{form: psi1gamma} is $\left( X\circ T\circ Z\circ \sigma
_{Y,U}\circ V\right) \left( X\circ \sigma _{Z\circ Y,T}\circ U\circ V\right)
\left( X\circ \sigma _{Y,Z}\circ T\circ U\circ V\right) \\=\left( X\circ
\sigma _{Y,T\circ Z\circ U}\circ V\right) \left( X\circ Y\circ \sigma
_{Z,T}\circ U\circ V\right) $.

\item \eqref{form:varphi1delta} is $\left( X\circ \sigma _{Y,U}\circ T\circ
Z\circ V\right) \left( X\circ Y\circ \sigma _{Z,U\circ T}\circ V\right)
\left( X\circ Y\circ Z\circ \sigma _{T,U}\circ V\right) \\=\left( X\circ
\sigma _{Y\circ T\circ Z,U}\circ V\right) \left( X\circ Y\circ \sigma
_{Z,T}\circ U\circ V\right) .$
\item \eqref{form:phipsi2NEW} becomes $(X\circ Z\circ \sigma_{T,U}\circ Y\circ V)(X\circ \sigma_{Y,Z\circ T\circ U}\circ V)\\=(X\circ Z\circ \sigma_{T\circ Y,U}\circ V)(X\circ Z\circ \sigma_{Y,T}\circ U\circ V)(X\circ \sigma_{Y,Z}\circ T\circ U\circ V)$.
\item \eqref{form:phivarphi2NEW} becomes $\left( X\circ U\circ \sigma _{Y,Z}\circ T\circ V\right) \left( X\circ
\sigma _{Y\circ Z\circ T,U}\circ V\right) \\=\left( X\circ \sigma _{Y,U\circ
Z}\circ T\circ V\right) \left( X\circ Y\circ \sigma _{Z,U}\circ T\circ
V\right) \left( X\circ Z\circ \sigma _{T,U}\circ V\right) $.
\item \eqref{form:phi0gamma-bul} and \eqref{form:phi0delta-bul} are true by coherence theorem for
monoidal categories as they only  involve unit constraints.
\end{itemize}
\end{invisible}
\end{example}

Note that, given the duoidal category $(\Cc,\circ,I,\circ,I)$ of \cref{ex:braidedasdouoidal}, in view of \eqref{equivalnecesBimonC}, we have that $\mathsf{Bimon}(\Cc,\circ,\circ)$ is the category $\mathsf{Bimon}(\Cc)$ of bimonoids in the braided monoidal category $\Cc$. 

In this setting, inspired by \cite[Theorem 7.17]{Bohm-Lack}, we prove the following result showing that $\mathsf{Hopf}(\Cc,\circ,\circ)$ is the category $\mathsf{Hopf}(\Cc)$ of Hopf monoids in the braided monoidal category $\Cc$.

\begin{proposition}%[\rd{upgrade:2026/03/02}]
\label{pro:Hopfmonbraidedcat}
    Let $(\Cc,\circ,I,\sigma)$ be a braided monoidal category.
For any object $H$ in $\mathsf{Bimon}(\Cc,\circ,\circ)$, the following conditions are equivalent:
\begin{enumerate}
    \item $H$ is a Hopf monoid in the braided monoidal category $(\Cc,\circ,I,\sigma)$;
     \item $H$ is a Hopf monoid in the duoidal category with a reversion $(\Cc,\circ,I,\circ,I)$ of \cref{ex:braidedasdouoidal};
     \item the Galois maps
 \[
 \varkappa _{P,X,Q}^{H} :=(\mu _{P}\circ X\circ Q)(P\circ\sigma_{X,H}\circ Q)(P\circ X \circ \rho _{Q}):P\circ X\circ Q\to P\circ X\circ Q
 \]
    are invertible, for any object $X$ in $\Cc$, right $H$-module $\left( P,\mu _{P}\right) $ and left $H$-comodule $\left( Q,\rho _{Q}\right)$;
      \item the co-Galois maps
\[
\varsigma _{P,X,Q}^{H}:=( P\circ X\circ \mu _{Q})(P\circ\sigma_{H,X}\circ Q)(\rho
_{P}\circ X\circ Q ):P\circ X\circ Q\to P\circ X\circ Q
\]
    are invertible, for any object $X$, right $H$-comodule $\left( P,\rho _{P}\right) $ and left $H$-module $\left( Q,\mu _{Q}\right) $;
    \item $g_r$ is an isomorphism in $\Cc$ for any morphism $g:H\to A$ in $\Mon(\Cc)$;
    \item $g^r$ is an isomorphism in $\Cc$ for any morphism $g:C\to H$ in $\Comon(\Cc)$;
     \item the fusion morphism $(\id_H)_r$ is invertible;
    \item the fusion morphism $(\id_H)^r$ is  invertible.
\end{enumerate}
As a consequence, the category $\Hopfmon(\Cc,\circ,\circ)$ is the category $\mathsf{Hopf}(\Cc)$.
\end{proposition}

\begin{proof}
    In view of the aforementioned identification of $\mathsf{Bimon}(\Cc,\circ,\circ)$ with $\mathsf{Bimon}(\Cc)$, the category $\Hopfmon(\Cc,\circ,\circ)$ is so the full subcategory of $\Bimon(\Cc)$ whose objects $H$ have an antipode $S_H:H\to H$ in the sense of \cref{def:Hopfmon}.

    $(1)\Rightarrow$ $(2)$ If $S_{H}:H\rightarrow H$ is an antipode in the braided category $(\Cc,\circ,I,\sigma)$, one easily checks that
\[
\begin{split}
    \mathrm{Id}%
_{H}\star S_{H}&=(m_{H}\circ I\circ I)\varphi_{H,H}(H\circ S_{H})\Delta_{H}\\&=(m_{H}\circ I\circ I)(H\circ\sigma_{I,H}\circ I)(H\circ I\circ r_{H}^{-1})(H\circ l_{H}^{-1})(H\circ S_{H})\Delta_{H}\\&=
\left( r_{H}^{-1}\circ I\right) r_{H}^{-1}m_{H}\left( H\circ
S_{H}\right) \Delta _{H}=\left( r_{H}^{-1}\circ I\right)
r_{H}^{-1}u_{H}\varepsilon _{H}=1_{H,H}^{l}
\end{split}
\]
and, similarly, $S_{H}\star
\mathrm{Id}_{H}=1_{H,H}^{r}$. Since $H^\revsn=H$, we conclude that $S_H$ is an antipode in the duoidal sense.

\begin{invisible}
\begin{eqnarray*}
\mathrm{Id}_{H}\star S_{H} &=&\left( m_{H}\circ I\circ I\right) \varphi
_{H,H}\left( H\circ S_{H}\right) \Delta _{H} \\
&=&\left( m_{H}\circ I\circ I\right) \left( H\circ \sigma _{I,H}\circ
I\right) \underbracket[0.140ex]{\left( H\circ I\circ r_{H}^{-1}\right) \left( H\circ
l_{H}^{-1}\right) }\left( H\circ S_{H}\right) \Delta _{H} \\
&=&\left( m_{H}\circ I\circ I\right) \left( H\circ r_{H}^{-1}\circ I\right)
\left( H\circ l_{H}\circ I\right) \left( H\circ l_{H\circ I}^{-1}\right)
\left( H\circ r_{H}^{-1}\right) \left( H\circ S_{H}\right) \Delta _{H} \\
&=&\left( m_{H}\circ I\circ I\right) \left( r_{H\circ H}^{-1}\circ I\right)
\left( H\circ l_{H\circ I}\right) \left( H\circ l_{H\circ I}^{-1}\right)
r_{H\circ H}^{-1}\left( H\circ S_{H}\right) \Delta _{H} \\
&=&\left( m_{H}\circ I\circ I\right) \left( r_{H\circ H}^{-1}\circ I\right)
r_{H\circ H}^{-1}\left( H\circ S_{H}\right) \Delta _{H} \\
&=&\left( r_{H}^{-1}\circ I\right) r_{H}^{-1}m_{H}\left( H\circ S_{H}\right)
\Delta _{H} \\
&=&\left( r_{H}^{-1}\circ I\right) r_{H}^{-1}u_{H}\varepsilon _{H} \\
&=&\left( u_{H}\circ I\circ I\right) \left( r_{I}^{-1}\circ I\right)
r_{I}^{-1}\varepsilon _{H} \\
&=&\left( u_{H}\circ I\circ I\right) \left( \Delta _{I}\circ I\right)
l_{I}^{-1}\varepsilon _{H}=1_{H,H}^{l}
\end{eqnarray*}%
Similarly
\begin{eqnarray*}
S_{H}\star \mathrm{Id}_{H} &=&\left( I\circ I\circ m_{H}\right) \psi
_{H,H}\left( S_{H}\circ H\right) \Delta _{H} \\
&=&\left( I\circ I\circ m_{H}\right) \left( I\circ \sigma _{H,I}\circ
H\right) \underbracket[0.140ex]{\left( l_{H}^{-1}\circ I\circ H\right) \left(
r_{H}^{-1}\circ H\right) }\left( S_{H}\circ H\right) \Delta _{H} \\
&=&\left( I\circ I\circ m_{H}\right) \left( I\circ l_{H}^{-1}\circ H\right)
\left( I\circ r_{H}\circ H\right) \left( r_{I\circ H}^{-1}\circ H\right)
\left( l_{H}^{-1}\circ H\right) \left( S_{H}\circ H\right) \Delta _{H} \\
&=&\left( I\circ I\circ m_{H}\right) \left( I\circ l_{H\circ H}^{-1}\right)
\left( r_{I\circ H}\circ H\right) \left( r_{I\circ H}^{-1}\circ H\right)
\left( l_{H}^{-1}\circ H\right) \left( S_{H}\circ H\right) \Delta _{H} \\
&=&\left( I\circ I\circ m_{H}\right) \left( I\circ l_{H\circ H}^{-1}\right)
l_{H\circ H}^{-1}\left( S_{H}\circ H\right) \Delta _{H} \\
&=&\left( I\circ l_{H}^{-1}\right) l_{H}^{-1}m_{H}\left( S_{H}\circ H\right)
\Delta _{H} \\
&=&\left( I\circ l_{H}^{-1}\right) l_{H}^{-1}u_{H}\varepsilon _{H} \\
&=&\left( I\circ I\circ u_{H}\right) \left( I\circ l_{I}^{-1}\right)
l_{I}^{-1}\varepsilon _{H} \\
&=&\left( I\circ I\circ u_{H}\right) \left( I\circ \Delta _{I}\right)
r_{I}^{-1}\varepsilon _{H}=1_{H,H}^{r}
\end{eqnarray*}
\end{invisible}

    $(2)\Rightarrow$ $(4)$ This is \cref{thm:Galinv}.

     $(4)\Rightarrow$ $(6)$ This follows from \eqref{equivalentgr}.

     $(6)\Rightarrow$ $(8)$ It is obvious.

    $(8)\Rightarrow$ $(1)$ It is easy to check that $(r_H\circ H)(\id_H)^r=(H\circ m_H)(\Delta_H\circ H)$ so that the latter is invertible. We deduce that $H$ has an antipode by the symmetric version of \cite[Theorem 1.8]{Ve}, where we note that the standing assumptions are superfluous for the part of the proof we are interested in.
\begin{invisible}
Given $g:C\rightarrow H$ as above, we  set
\begin{equation*}
\partial \left( g\right):=\left( C\circ m_{H}\right) \left( C\circ g\circ
H\right) \left( \Delta _{C}\circ H\right):C\circ H\to  C\circ H.
\end{equation*}
First, given $g:C\to H$ one easily proves that $(C\circ \partial(g))(\Delta_C\circ H)=(\Delta_C\circ H)\partial(g)$ and $(C\circ m_H)(\partial(g)\circ H)=\partial(g)(C\circ m_H)$. Thus, if $\partial(g)$ is invertible, we get that  $(\Delta_C\circ H)\partial(g)^{-1}=(C\circ \partial(g)^{-1})(\Delta_C\circ H)$ and $\partial(g)^{-1}(C\circ m_H)=(C\circ m_H)(\partial(g)^{-1}\circ H)$. Therefore, if we set $\bar g:=(\varepsilon_C\circ H)\partial(g)^{-1}(C\circ u_H)$, we get (by omitting unit constraints)
\begin{align*}
\partial(g)^{-1}&=(C\circ  \varepsilon_C\circ H)(\Delta_C\circ H)\partial(g)^{-1}(C\circ m_H)(C\circ  u_H\circ H)
\\ &=(C\circ  \varepsilon_C\circ H)(C\circ \partial(g)^{-1})(\Delta_C\circ H)(C\circ m_H)(C\circ  u_H\circ H)
\\ &=(C\circ  \varepsilon_C\circ H)(C\circ \partial(g)^{-1})(C\circ C\circ m_H)(C\circ C\circ   u_H\circ H)(\Delta_C\circ I\circ H)
\\ &=(C\circ  \varepsilon_C\circ H)(C\circ C\circ m_H)(C\circ \partial(g)^{-1}\circ H)(C\circ C\circ   u_H\circ H)(\Delta_C\circ I\circ H)
\\ &=(C\circ I\circ m_H)(C\circ  \varepsilon_C\circ H\circ H)(C\circ \partial(g)^{-1}\circ H)(C\circ C\circ   u_H\circ H)(\Delta_C\circ I\circ H)
\\ &=(C\circ m_H)(C\circ  \bar g\circ H)(\Delta_C\circ H)=\partial(\bar g).
\end{align*}
Thus $\partial(g)^{-1}=\partial(\bar g)$. As a consequence $\partial(g*\bar g)=\partial(g)\partial(\bar g)=\id$. Hence
$g*\bar g=(\varepsilon_C\circ H)\partial(g*\bar g)(C\circ  u_H)=(\varepsilon_C\circ H)(C\circ  u_H)=u_H\varepsilon_C$. Similarly $\bar g*g=u_H\varepsilon_C$. We have so proved that $\partial(g)$ invertible implies $g$ convolution invertible.
In particular, for $g=\id_H$, we get that $\partial(\id_H)$ invertible implies $\id_H$ convolution invertible i.e. there is an antipode $S$. \\
NOTE: Given a morphism $g:C\rightarrow H$ in $\Comon(\Cc),$ the morphism $\partial \left( g\right):C\circ H\to  C\circ H$, employed above,
    is the categorical counterpart of the \emph{Hopf-Galois map} for the right comodule $(C,(C\circ g)\Delta_C)$ defined in \cite[Definition 3.0.1]{Chemla}.
\end{invisible}
Similarly one proves  $(2)\Rightarrow(3)\Rightarrow(5)\Rightarrow(7)\Rightarrow(1)$.
\end{proof}

Let $(\Cc,\circ,I,\sigma)$ be a braided monoidal category. In \cite[Proposition 6.13]{Aguiar} it is shown that the duoidal category  $(\Cc,\circ,I,\circ,I)$ recalled in \cref{ex:braidedasdouoidal} is $\circ$-braided precisely when $\sigma$ is a symmetry.

\begin{corollary}%[\rd{added:2026/03/13}]
Let $(\Cc,\circ,I,\sigma)$ be a symmetric monoidal category. Assume that $\mathcal{C}$ has and $\circ$ preserves binary intersections. Then, the category $\cHopfmon(\Cc)$ of cocommutative Hopf monoids therein is protomodular.
\end{corollary}

\begin{proof}
Since we know that our duoidal category has a reversion, we get that $\cHopfmon^{\mathrm{m}}(\Cc,\circ, \circ )$ is protomodular, by \cref{thm:protomon}. Now, \cref{pro:Hopfmonbraidedcat} entails that $\Hopfmon(\Cc,\circ,\circ)$ is the category $\mathsf{Hopf}(\Cc)$. Moreover, any Hopf monoid $H$ here verifies $\varepsilon_Hu_H=\id$ so that the assumption that $u_H$ is a monomorphism is redundant. As a consequence $\cHopfmon^{\mathrm{m}}(\Cc,\circ, \circ )$ is the category $\cHopfmon(\Cc)$.  \end{proof}

\begin{remark}%[\as{[Added on 25/03/2026]}]
\label{rmk:Hopfcc}
    By \cite[Lemma 1.2.10]{Johnstone-v1} and the subsequent comment, binary intersections allow one to compute coreflexive equalizers (that is equalizers of pairs with a common retraction) and  these equalizers together with finite products suffice to construct all finite limits. In \cite[Proposition 3.11]{SZ} and \cite[Remark 3.4]{SZ}, it is proven that $\cHopfmon(\Cc)$ is protomodular once $(\Cc,\circ,I,\sigma)$ is a symmetric monoidal category with coreflexive equalizers which are preserved by $\circ$. This simply follows from the fact that $\cHopfmon(\Cc)=\mathsf{Grp}(\cComon(\Cc))$, where $\cComon(\Cc)$ is a finitely complete category as $(\cComon(\Cc),\circ,I)$ is cartesian monoidal (so it has finite products) and it has coreflexive equalizers. We point out that, if $(\Cc,\circ,I,\sigma)$ is an abelian symmetric monoidal category which satisfies some faithful (co)flatness conditions, the category $\cHopfmon(\Cc)$ is even semi-abelian by \cite[Theorem 8.3]{SZ}.
\end{remark}

% \rd{[ALERT: possible intersections with Sciandra-Zhenbang's paper.]}
% \as{[Nel caso $(\Cc,\ot,\mathbf{1},\sigma)$ categoria monoidale simmetrica (con equalizzatori) la protomodularità diventa immediata perchè $\cHopf(\Cc)=\Grp(\cComon(\Cc))$ e $\cComon(\Cc)$ è finitamente completa (si usa proprio risultato che dice "$\mathsf{Grp}(\Cc)$ con $\Cc$ finitamente completa è protomodulare"). Abbiamo osservato questo nel nostro lavoro, ma poi ci siamo dedicati completamente alla regolarità, generalizzando il Teorema di Newman. Potremmo poi pensare di provare a generalizzare quest'ultimo al contesto duoidale, per la regolarità.]}
% \rd{[Potremmo dire che nel caso simmetrico la protomodularità segue facilmente come suddetto e accennare al fatto che nel vostro articolo state studiando la semiabelianità.]}
% \as{[Ok, mi occupo di aggiungere qualche riga a riguardo.]}

\subsection{Internal groups revised} In this section, we extend the notion of internal group to an arbitrary monoidal category with finite products which is endowed with a reversion and we explore in this setting the fallouts of the results we have obtained so far.

\begin{example}
\label{exa:duoidfinprod}
Every monoidal category $(\Cc,\circ,I)$ with finite products becomes a $\times$-braided duoidal category $(\Cc,\circ,I,\times,\mathsf{1})$ where $\times$ is the binary product and $\mathsf{1}$ the terminal object, see \cite[Example 6.19 and Remark 6.21]{Aguiar}. %Then it makes sense to define the category of internal groups  in $(\Cc,\circ,I)$ to be $\cHopfmon(\Cc,\circ,\times)$. %We denote this category by $\Grp(\Cc,\circ,I)$.
\end{example}

The previous example suggests the following definition. 
% \rd{[Non so se vogliamo tenerlo. Nel caso, occorre coinvolgere una reversion. Potremmo anche limitarci a un funtore monoidale rispetto a $\circ$ e cartesiano rispetto a $\times$.]} \as{[A me sembra interessante usare questa nozione di gruppo interno generalizzato per recuperare $B/\mathsf{TrHopf}$, personalmente lo terrei (magari dicendo che tale nozione può essere poi approfondita in un altro lavoro ?)]} \rd{[In altro lavoro coinvolgendo Zhenbang?]}
% \as{[Mi sembra possa essere una buona idea]}

\begin{definition}
    Let $(\Cc,\circ, I)$ be a monoidal category with binary product $\times$ and terminal object $\terminal$. When there exists a reversion on the duoidal category $(\Cc,\circ,I,\times,\terminal)$, we define the category $\mathsf{Grp}(\Cc,\circ,I)$ of \textit{internal groups} in $(\Cc,\circ, I)$ as the category $\Hopfmon(\Cc,\circ,\times)$.
\end{definition}

\begin{remark}%[\rd{revised:2026/02/11}]
\label{rmk:internalgroupduoidal}
Since $(\Cc,\times,\terminal)$ is a cartesian monoidal category, we have
\[\mathsf{Comon}_{\mathrm{coc}}(\Cc,\times,\terminal)=\mathsf{Comon}(\Cc,\times,\terminal)\cong \Cc,\] see e.g. \cite[Proposition 2.22]{Ulrich-Myriam}. It follows that $\cBimon(\Cc,\circ,\times)=\mathsf{Bimon}(\Cc,\circ,\times)$  whence
\[
\cBimon(\Cc,\circ,\times)\cong\mathsf{Mon}(\mathsf{Comon}(\Cc,\times,\mathbf{1}),\circ, I)\cong\mathsf{Mon}(\Cc,\circ,I).
\]
Moreover, $\cHopfmon(\Cc,\circ,\times)=\Hopfmon(\Cc,\circ,\times)=\mathsf{Grp}(\Cc,\circ,I).$
\end{remark}

\begin{example}%[\rd{revised:2026/02/11}]
\label{exa:usualGrp}
    In case $(\Cc,\circ,I)=(\Cc,\times,\terminal)$, which is a symmetric monoidal category, see e.g.\ \cite[Example 1.4]{Aguiar}, we are in the setting of \cref{pro:Hopfmonbraidedcat}, where the duoidal category $(\Cc,\times,\terminal,\times,\terminal)$ is considered. Thus, $\cHopfmon(\Cc,\times,\times)=\Hopfmon(\Cc,\times,\times)=\mathsf{Hopf}(\Cc)$. Since $\Cc$ is cartesian, $\mathsf{Hopf}(\Cc)$ coincides with the category $\mathsf{Grp}(\Cc)$ of internal groups in $\Cc$.
\end{example}

Denote by $\mathsf{Grp}^{\mathrm{m}}(\Cc,\circ,I)$ the full subcategory  of $\mathsf{Grp}(\Cc,\circ,I)$ of internal groups $G$ whose unit morphism $u_{G}^{\circ}:I\to G$ is a  monomorphism in $\Cc$. By the foregoing we have $\mathsf{Grp}^{\mathrm{m}}(\Cc,\circ,I)=\cHopfmon^{\mathrm{m}}(\Cc,\circ,\times)$.

\begin{corollary}%[\rd{added:2026/03/18}]
\label{coro:Grproto}
Let $(\Cc,\circ,I)$ be a finitely complete monoidal category. If there is a reversion on the $\times$-braided duoidal category $(\Cc,\circ,I,\times,\terminal)$, then the category $\mathsf{Grp}^{\mathrm{m}}(\Cc,\circ,I)$ is protomodular.
\end{corollary}

\begin{proof} Since limits commute with limits, we have that $\times$ preserves pullbacks whence binary intersections. Thus, \cref{thm:protomon} entails that $\cHopfmon^{\mathrm{m}}(\Cc,\circ, \times)=\mathsf{Grp}^{\mathrm{m}}(\Cc,\circ,I)$ is protomodular.
\end{proof}

% \begin{invisible}
% \as{[Immagino sia noto. Non ho una referenza ora ma dovrei averlo provato qui sotto]}

% Let $\Cc$ be a category with finite products and binary intersections (i.e.\ pullbacks di mono). Then, by \cite[Proposition 2.8.2]{BorI94}, $\Cc$ is finitely complete if and only if it has equalizers. Let $f,g:A\to B$ be morphisms in $\Cc$ and consider the morphism $\langle1_{A},f\rangle:A\to A\times B$ and $\langle1_{A},g\rangle:A\to A\times B$ in $\Cc$ which are monomorphisms since $\pi_{A}\langle1_{A},f\rangle=1_{A}$ and $\pi_{A}\langle1_{A},g\rangle=1_{A}$. Then, consider the pullback $(P,p_{1},p_{2})$ of the pair $(\langle1_{A},f\rangle,\langle1_{A},g\rangle)$ in $\Cc$. Since $\langle p_{1},fp_{1}\rangle=\langle1_{A},f\rangle p_{1}=\langle1_{A},g\rangle p_{2}=\langle p_{2},gp_{2}\rangle$, we obtain $p_{1}=p_{2}$ and $fp_{1}=gp_{2}=gp_{1}$. We prove that $p_{1}:P\to A$ is the equalizer of $(f,g)$ showing that the universal property is satisfied. Let $h:C\to A$ in $\Cc$ such that $fh=gh$. Then, since $\langle1_{A},f\rangle h=\langle h,fh\rangle=\langle h,gh\rangle=\langle1_{A},g\rangle h$, by the universal property of the pullback $(P,p_{1},p_{2})$, there exists a unique morphism $\psi:C\to P$ in $\Cc$ such that $p_{1}\psi=h$.
% \end{invisible}
We note that the assumptions in \cref{coro:Grproto} are optimal. Indeed, we need finite products as in \cref{exa:duoidfinprod} and binary intersections as in  \cref{thm:protomon}.
Moreover, as recalled in \cref{rmk:Hopfcc},
%by \cite[Lemma 1.2.10]{Johnstone-v1} and the subsequent comment, 
binary intersections allow one to compute coreflexive equalizers and  these equalizers together with finite products suffice to construct all finite limits.

It would be interesting to find conditions ensuring the existence of a reversion as in \cref{coro:Grproto}.\medskip

On the one hand, by applying the previous corollary to \cref{exa:usualGrp}, one recovers the protomodularity of usual internal groups. On the other hand, to substantiate our new notion of internal group, we need an example where $\circ$ and $\times$ differ. We will be able to present such an example in \cref{rmk:nntrvntgrps} after discussing the duoidal category of bimodules.

\subsection{The duoidal category of bimodules}
Here we study a particular duoidal category consisting of bimodules over a bialgebra.

\begin{example}
\label{es:bimod} As shown in \cite[Proposition 2.7]{BT25}, see also \cite[Remark 1.6]{Sar21}, the category $\mathcal{C}={}_B\mm_B$ of $B$-bimodules, where $B$ is a $\Bbbk$-bialgebra, is duoidal. In fact,  the monoidal structures are given by $\mathcal{C}^\circ=(\mathcal{C},  \circ =\otimes_B, I=B )$  and $\mathcal{C}^\bullet= (\mathcal{C}, \bullet=\otimes_\Bbbk, J=\Bbbk )$, where the base field $\Bbbk$ is regarded as a $B$-bimodule via the counit $\varepsilon_{B}$ of $B$. The interchange law $\zeta$ is the natural transformation whose component $\zeta_{A,E,C,D}:(A\bullet E)\circ (C\bullet D)\to (A\circ C)\bullet (E\circ D),$ for all $A,E,C,D\in\mathcal{C}$, is defined by
\[
\zeta_{A,E,C,D}\big((a\otimes e)\otimes_B (c\otimes d)\big)=(a\otimes_Bc)\otimes (e\otimes_B d).
\]
Moreover, the morphisms $\Delta_{I}:I\to I\bullet I$, $m_{J}:J\circ J\to J$ and $\varepsilon_{I}:I\to J$ are, respectively, given by $\Delta_{B}$, $m_{\Bbbk}(k\ot_{B}k')=kk'$ and $\varepsilon_{B}$.
\end{example}

\begin{remark}
\label{rmk:bimonbimod}
In view of  \cite[Corollary 2.10]{BT25}, in the setting of \cref{es:bimod},   giving a bimonoid $(H,m_H^\circ,u_H^\circ,\Delta_H^\bullet,\varepsilon_H^\bullet)$ in $\Cc$
is equivalent to giving a bialgebra $H$ together with a bialgebra map $u_{H}^{\circ}:B\rightarrow H$ which defines the $B$-bimodule
structure of $H$ so that $u_{H}^{\circ}\left(
b\right) =b\cdot 1_{H}.$
Through this identification, a morphism of bimonoids turns out to be a bialgebra map $f:H\to H'$ such that $fu_H^\circ=u_{H'}^\circ$.
In other words, the category  $\Bimon(\Cc,\circ,\bullet)$ of bimonoids is isomorphic to the coslice category of bialgebras under $B$ i.e. $B/\Bialg$.
\end{remark}

Next, we construct a reversion on the duoidal category of \cref{es:bimod} for $B$ a $\Bbbk$-Hopf algebra.

\subsubsection{Reversion and antipode}\label{subsec:reversionbim}
In this subsection, $B$ will denote a $\Bbbk$-Hopf algebra.

For every $B$-bimodule $X$, one can define a $B$%
-bimodule $X^{\revsn }$ whose underlying vector space is $X$, in the following
way. If we denote by $x^{\revsn }$ an element $x\in X$ regarded in $X^{\revsn },$
then the $B$-bimodule structure of $X^{\revsn }$ is defined by $bx^{\revsn }b^{\prime }:=\left(
S_{B}\left( b^{\prime }\right) xS_{B}\left( b\right) \right) ^{\revsn }.$
Given a $B$-bimodule morphism $f:X\rightarrow Y,$ we can define $f^{\revsn
}:X^{\revsn }\rightarrow Y^{\revsn }$ by setting $f^{\revsn }\left( x^{\revsn
}\right) :=f\left( x\right) ^{\revsn }.$ This yields a functor $\left(
-\right) ^{\revsn }:\mathcal{C}\rightarrow \mathcal{C}$. Then
it is easy to check that $\left( \left( -\right) ^{\revsn },\phi ^{\circ
},\phi ^{\bullet }\right) \ $ is double lax monoidal, where
\begin{align*}
\phi^\circ _{X,Y} &:X^{\revsn }\circ Y^{\revsn }\rightarrow \left( Y\circ X\right)
^{\revsn },\; x^{\revsn }\circ y^{\revsn }\mapsto \left( y\circ x\right) ^{\revsn
}, \qquad
\phi^\circ _{0} :B\rightarrow B^{\revsn },\; b\mapsto S_{B}\left( b\right)
^{\revsn },\\
\phi _{X,Y}^{\bullet } &:X^{\revsn }\bullet Y^{\revsn }\rightarrow \left(
Y\bullet X\right) ^{\revsn },\; x^{\revsn }\bullet y^{\revsn }\mapsto
\left( y\bullet x\right) ^{\revsn },\qquad
\phi _{0}^{\bullet } :J\rightarrow J^{\revsn },\; k\mapsto k^{\revsn }.
\end{align*}
Note that $\phi^\circ_0$ is  invertible if and only if so is $S_B$, so $((-)^\revsn,\phi^\circ)$ is not strong monoidal, in general.

\begin{invisible}
We check that $\phi^\circ _{X,Y}$ is well-defined:
\begin{eqnarray*}
\phi^\circ _{X,Y}\left( x^{\revsn }b\circ y^{\revsn }\right) &=&\phi^\circ _{X,Y}\left(
\left( S_{B}\left( b\right) x\right) ^{\revsn }\circ y^{\revsn }\right) =\left(
y\circ S_{B}\left( b\right) x\right) ^{\revsn } \\
&=&\left( yS_{B}\left( b\right) \circ x\right) ^{\revsn }=\phi^\circ _{X,Y}\left(
x^{\revsn }\circ \left( yS_{B}\left( b\right) \right) ^{\revsn }\right) =\phi^\circ
_{X,Y}\left( x^{\revsn }\circ by^{\revsn }\right) .
\end{eqnarray*}%
We check that $\phi^\circ _{X,Y}$ is a $B$-bimodule morphism:%
\begin{eqnarray*}
\phi^\circ _{X,Y}\left( b\left( x^{\revsn }\circ y^{\revsn }\right) b^{\prime }\right)
&=&\phi^\circ _{X,Y}\left( bx^{\revsn }\circ y^{\revsn }b^{\prime }\right) =\phi^\circ
_{X,Y}\left( \left( xS_{B}\left( b\right) \right) ^{\revsn }\circ \left(
S_{B}\left( b^{\prime }\right) y\right) ^{\revsn }\right) \\
&=&\left( S_{B}\left( b^{\prime }\right) y\circ xS_{B}\left( b\right)
\right) ^{\revsn }=\left( S_{B}\left( b^{\prime }\right) \left( y\circ
x\right) S_{B}\left( b\right) \right) ^{\revsn } \\
&=&b\left( y\circ x\right) ^{\revsn }b^{\prime }=b\phi^\circ _{X,Y}\left( x^{\revsn
}\circ y^{\revsn }\right) b^{\prime }.
\end{eqnarray*}%
We check that $\phi^\circ _{X,Y}$ is natural. Given $f:X\rightarrow X^{\prime }$
and $g:Y\rightarrow Y^{\prime },$ we have%
\begin{eqnarray*}
\left( g\circ f\right) ^{\revsn }\phi^\circ _{X,Y}\left( x^{\revsn }\circ y^{\revsn
}\right) &=&\left( g\circ f\right) ^{\revsn }\left( \left( y\circ x\right)
^{\revsn }\right) =\left( \left( g\circ f\right) \left( y\circ x\right)
\right) ^{\revsn } \\
&=&\left( g\left( y\right) \circ f\left( x\right) \right) ^{\revsn }=\phi^\circ
_{X^{\prime },Y^{\prime }}\left( f\left( x\right) ^{\revsn }\circ g\left(
y\right) ^{\revsn }\right) =\phi^\circ _{X^{\prime },Y^{\prime }}\left( f^{\revsn
}\left( x^{\revsn }\right) \circ g^{\revsn }\left( y^{\revsn }\right) \right) \\
&=&\phi^\circ _{X^{\prime },Y^{\prime }}\left( f^{\revsn }\circ g^{\revsn }\right)
\left( x^{\revsn }\circ y^{\revsn }\right) .
\end{eqnarray*}%
We check that $\phi^\circ _{0}\ $is a $B$-bimodule morphism:%
\begin{eqnarray*}
\phi^\circ _{0}\left( bxb^{\prime }\right) &=&S_{B}\left( bxb^{\prime }\right)
^{\revsn }=\left( S_{B}\left( b^{\prime }\right) S_{B}\left( x\right)
S_{B}\left( b\right) \right) ^{\revsn } \\
&=&bS_{B}\left( x\right) ^{\revsn }b^{\prime }=b\phi^\circ _{0}\left( x\right)
b^{\prime }.
\end{eqnarray*}%
We check that $\phi^\circ $ is associative:%
\begin{equation*}
\begin{array}{ccc}
X^{\revsn }\circ Y^{\revsn }\circ Z^{\revsn } & \overset{X^{\revsn }\circ \phi^\circ _{Y,Z}%
}{\longrightarrow } & X^{\revsn }\circ \left( Z\circ Y\right) ^{\revsn } \\
\phi^\circ _{X,Y}\circ Z^{\revsn }\downarrow &  & \downarrow \phi^\circ _{X,Z\circ Y} \\
\left( Y\circ X\right) ^{\revsn }\circ Z^{\revsn } & \overset{\phi^\circ _{Y\circ X,Z}}%
{\longrightarrow } & \left( Z\circ Y\circ X\right) ^{\revsn }%
\end{array}%
\qquad
\begin{array}{ccc}
x^{\revsn }\circ y^{\revsn }\circ z^{\revsn } & \overset{X^{\revsn }\circ \phi^\circ _{Y,Z}%
}{\longrightarrow } & x^{\revsn }\circ \left( z\circ y\right) ^{\revsn } \\
\phi^\circ _{X,Y}\circ Z^{\revsn }\downarrow &  & \downarrow \phi^\circ _{X,Z\circ Y} \\
\left( y\circ x\right) ^{\revsn }\circ z^{\revsn } & \overset{\phi^\circ _{Y\circ X,Z}}%
{\longrightarrow } & \left( z\circ y\circ x\right) ^{\revsn }%
\end{array}%
\end{equation*}%
$\phi^\circ $ is left unital:%
\begin{equation*}
\begin{array}{ccc}
B\circ X^{\revsn } & \overset{l_{X^{\revsn }}^{\circ }}{\longrightarrow } &
X^{\revsn } \\
\phi^\circ _{0}\circ X^{\revsn }\downarrow &  & \uparrow \left( r_{X}^{\circ
}\right) ^{\revsn } \\
B^{\revsn }\circ X^{\revsn } & \overset{\phi^\circ _{B,X}}{\longrightarrow } & \left(
X\circ B\right) ^{\revsn }%
\end{array}%
\qquad
\begin{array}{ccc}
b\circ x^{\revsn } & \overset{l_{X^{\revsn }}^{\circ }}{\longrightarrow } &
bx^{\revsn }=\left( xS_{B}\left( b\right) \right) ^{\revsn } \\
\phi^\circ _{0}\circ X^{\revsn }\downarrow &  & \uparrow \left( r_{X}^{\circ
}\right) ^{\revsn } \\
S_{B}\left( b\right) ^{\revsn }\circ x^{\revsn } & \overset{\phi^\circ _{B,X}}{%
\longrightarrow } & \left( x\circ S_{B}\left( b\right) \right) ^{\revsn }%
\end{array}%
\end{equation*}%
$\phi^\circ $ is right unital:%
\begin{equation*}
\begin{array}{ccc}
X^{\revsn }\circ B & \overset{r_{X^{\revsn }}^{\circ }}{\longrightarrow } &
X^{\revsn } \\
X^{\revsn }\circ \phi^\circ _{0}\downarrow &  & \uparrow \left( l_{X}^{\circ
}\right) ^{\revsn } \\
X^{\revsn }\circ B^{\revsn } & \overset{\phi^\circ _{X,B}}{\longrightarrow } & \left(
B\circ X\right) ^{\revsn }%
\end{array}%
\qquad
\begin{array}{ccc}
x^{\revsn }\circ b & \overset{r_{X^{\revsn }}^{\circ }}{\longrightarrow } &
x^{\revsn }b=\left( S_{B}\left( b\right) x\right) ^{\revsn } \\
X^{\revsn }\circ \phi^\circ _{0}\downarrow &  & \uparrow \left( l_{X}^{\circ
}\right) ^{\revsn } \\
x^{\revsn }\circ S_{B}\left( b\right) ^{\revsn } & \overset{\phi^\circ _{X,B}}{%
\longrightarrow } & \left( S_{B}\left( b\right) \circ x\right) ^{\revsn }%
\end{array}%
\end{equation*}
We now consider $\left( \left( -\right) ^{\revsn },\phi ^{\bullet }\right).$\\
$\phi _{X,Y}^{\bullet }$ is a morphism of $B$-bimodules%
\begin{eqnarray*}
\phi _{X,Y}^{\bullet }\left( b\left( x^{\revsn }\bullet y^{\revsn }\right)
b^{\prime }\right) &=&\phi _{X,Y}^{\bullet }\left( b_{1}x^{\revsn
}b_{1}^{\prime }\bullet b_{2}y^{\revsn }b_{2}^{\prime }\right) \\
&=&\phi _{X,Y}^{\bullet }\left( \left( S_{B}\left( b_{1}^{\prime }\right)
xS_{B}\left( b_{1}\right) \right) ^{\revsn }\bullet \left( S_{B}\left(
b_{2}^{\prime }\right) yS_{B}\left( b_{2}\right) \right) ^{\revsn }\right) \\
&=&\left( S_{B}\left( b_{2}^{\prime }\right) yS_{B}\left( b_{2}\right)
\bullet S_{B}\left( b_{1}^{\prime }\right) xS_{B}\left( b_{1}\right) \right)
^{\revsn } \\
&=&\left( S_{B}\left( b^{\prime }\right) \left( y\bullet x\right)
S_{B}\left( b\right) \right) ^{\revsn } \\
&=&b\left( y\bullet x\right) ^{\revsn }b^{\prime }=b\phi _{X,Y}^{\bullet
}\left( x^{\revsn }\bullet y^{\revsn }\right) b^{\prime }
\end{eqnarray*}

We check that $\phi _{X,Y}^{\bullet }$ is natural. Given $f:X\rightarrow
X^{\prime }$ and $g:Y\rightarrow Y^{\prime },$ we have%
\begin{eqnarray*}
\left( g\bullet f\right) ^{\revsn }\phi _{X,Y}^{\bullet }\left( x^{\revsn
}\bullet y^{\revsn }\right) &=&\left( g\bullet f\right) ^{\revsn }\left( \left(
y\bullet x\right) ^{\revsn }\right) =\left( \left( g\bullet f\right) \left(
y\bullet x\right) \right) ^{\revsn } \\
&=&\left( g\left( y\right) \bullet f\left( x\right) \right) ^{\revsn }=\phi
_{X^{\prime },Y^{\prime }}^{\bullet }\left( f\left( x\right) ^{\revsn }\bullet
g\left( y\right) ^{\revsn }\right) =\phi _{X^{\prime },Y^{\prime }}^{\bullet}\left(
f^{\revsn }\left( x^{\revsn }\right) \bullet g^{\revsn }\left( y^{\revsn }\right)
\right) \\
&=&\phi _{X^{\prime },Y^{\prime }}^{\bullet}\left( f^{\revsn }\bullet g^{\revsn }\right)
\left( x^{\revsn }\bullet y^{\revsn }\right) .
\end{eqnarray*}%
Since $bk^{\revsn }b^{\prime }=\left( S_{B}\left( b^{\prime }\right)
kS_{B}\left( b\right) \right) ^{\revsn }=\left( \varepsilon _{B}^{\bullet
}\left( b^{\prime }\right) k\varepsilon _{B}^{\bullet }\left( b\right)
\right) ^{\revsn }=\varepsilon _{B}^{\bullet }\left( b\right) k^{\revsn
}\varepsilon _{B}^{\bullet }\left( b^{\prime }\right) $ it is clear that $%
J^{\revsn }=J$ as $B$-bimodules and hence $\phi _{0}^{\bullet }=\mathrm{Id}%
_{\Bbbk }$ is a morphism of $B$-bimodules.

We check that $\phi^\bullet $ is associative:%
\begin{equation*}
\begin{array}{ccc}
X^{\revsn }\bullet Y^{\revsn }\bullet Z^{\revsn } & \overset{X^{\revsn }\bullet \phi
_{Y,Z}^{\bullet }}{\longrightarrow } & X^{\revsn }\bullet \left( Z\bullet
Y\right) ^{\revsn } \\
\phi _{X,Y}^{\bullet }\bullet Z^{\revsn }\downarrow  &  & \downarrow \phi
_{X,Z\bullet Y}^{\bullet } \\
\left( Y\bullet X\right) ^{\revsn }\bullet Z^{\revsn } & \overset{\phi
_{Y\bullet X,Z}}{\longrightarrow } & \left( Z\bullet Y\bullet X\right)
^{\revsn }%
\end{array}%
\qquad
\begin{array}{ccc}
x^{\revsn }\bullet y^{\revsn }\bullet z^{\revsn } & \overset{X^{\revsn }\bullet \phi
_{Y,Z}^{\bullet }}{\longrightarrow } & x^{\revsn }\bullet \left( z\bullet
y\right) ^{\revsn } \\
\phi _{X,Y}^{\bullet }\bullet Z^{\revsn }\downarrow  &  & \downarrow \phi
_{X,Z\bullet Y}^{\bullet } \\
\left( y\bullet x\right) ^{\revsn }\bullet z^{\revsn } & \overset{\phi
_{Y\bullet X,Z}^{\bullet }}{\longrightarrow } & \left( z\bullet y\bullet
x\right) ^{\revsn }%
\end{array}%
\end{equation*}%
$\phi ^{\bullet }$ is left unital:%
\begin{equation*}
\begin{array}{ccc}
\Bbbk \bullet X^{\revsn } & \overset{l_{X^{\revsn }}^{\bullet }}{\longrightarrow
} & X^{\revsn } \\
\phi _{0}^{\bullet }\bullet X^{\revsn }\downarrow  &  & \uparrow \left(
r_{X}^{\bullet }\right) ^{\revsn } \\
\Bbbk ^{\revsn }\bullet X^{\revsn } & \overset{\phi _{\Bbbk,X}^{\bullet }}{%
\longrightarrow } & \left( X\bullet \Bbbk \right) ^{\revsn }%
\end{array}%
\qquad
\begin{array}{ccc}
k\bullet x^{\revsn } & \overset{l_{X^{\revsn }}^{\bullet }}{\longrightarrow } &
kx^{\revsn }=\left( xk\right) ^{\revsn } \\
\phi _{0}^{\bullet }\bullet X^{\revsn }\downarrow  &  & \uparrow \left(
r_{X}^{\bullet }\right) ^{\revsn } \\
k^{\revsn }\bullet x^{\revsn } & \overset{\phi _{\Bbbk,X}^{\bullet }}{%
\longrightarrow } & \left( x\bullet k\right) ^{\revsn }%
\end{array}%
\end{equation*}%
$\phi ^{\bullet }$ is right unital:%
\begin{equation*}
\begin{array}{ccc}
X^{\revsn }\bullet \Bbbk  & \overset{r_{X^{\revsn }}^{\bullet }}{\longrightarrow
} & X^{\revsn } \\
X^{\revsn }\bullet \phi _{0}^{\bullet }\downarrow  &  & \uparrow \left(
l_{X}^{\bullet }\right) ^{\revsn } \\
X^{\revsn }\bullet \Bbbk ^{\revsn } & \overset{\phi _{X,\Bbbk}^{\bullet }}{%
\longrightarrow } & \left( \Bbbk \bullet X\right) ^{\revsn }%
\end{array}%
\qquad
\begin{array}{ccc}
x^{\revsn }\bullet k & \overset{r_{X^{\revsn }}^{\bullet }}{\longrightarrow } &
x^{\revsn }k=\left( kx\right) ^{\revsn } \\
X^{\revsn }\bullet \phi _{0}^{\bullet }\downarrow  &  & \uparrow \left(
l_{X}^{\bullet }\right) ^{\revsn } \\
x^{\revsn }\bullet k^{\revsn } & \overset{\phi _{X,\Bbbk}^{\bullet }}{%
\longrightarrow } & \left( k\bullet x\right) ^{\revsn }%
\end{array}%
\end{equation*}%
Thus $\left( \left( -\right) ^{\revsn },\phi ^{\bullet }\right) $ is lax
monoidal. Let us see the interchange law:%
\begin{equation*}
\begin{array}{ccc}
\left( A^{\revsn }\bullet B^{\revsn }\right) \circ \left( C^{\revsn }\bullet
D^{\revsn }\right)  & \overset{\zeta _{A^{\revsn },B^{\revsn },C^{\revsn },D^{\revsn }}%
}{\longrightarrow } & \left( A^{\revsn }\circ C^{\revsn }\right) \bullet \left(
B^{\revsn }\circ D^{\revsn }\right)  \\
\phi _{A,B}^{\bullet }\circ \phi _{C,D}^{\bullet }\downarrow  &  &
\downarrow \phi _{A,C}^{\circ }\bullet \phi _{B,D}^{\circ } \\
\left( B\bullet A\right) ^{\revsn }\circ \left( D\bullet C\right) ^{\revsn } &
& \left( C\circ A\right) ^{\revsn }\bullet \left( D\circ B\right) ^{\revsn } \\
\phi _{B\bullet A,D\bullet C}^{\circ }\downarrow &  & \downarrow\phi _{C\circ A,D\circ
B}^{\bullet } \\
\left( \left( D\bullet C\right) \circ \left( B\bullet A\right) \right)
^{\revsn } & \overset{\left( \zeta _{D,C,B,A}\right) ^{\revsn }}{\longrightarrow
} & \left( \left( D\circ B\right) \bullet \left( C\circ A\right) \right)
^{\revsn }%
\end{array}%
\quad
\begin{array}{ccc}
\left( a^{\revsn }\bullet b^{\revsn }\right) \circ \left( c^{\revsn }\bullet
d^{\revsn }\right)  & \overset{\zeta _{A^{\revsn },B^{\revsn },C^{\revsn },D^{\revsn }}%
}{\longrightarrow } & \left( a^{\revsn }\circ c^{\revsn }\right) \bullet \left(
b^{\revsn }\circ d^{\revsn }\right)  \\
\phi _{A,B}^{\bullet }\circ \phi _{C,D}^{\bullet }\downarrow  &  &
\downarrow \phi _{A,C}^{\circ }\bullet \phi _{B,D}^{\circ } \\
\left( b\bullet a\right) ^{\revsn }\circ \left( d\bullet c\right) ^{\revsn } &
& \left( c\circ a\right) ^{\revsn }\bullet \left( d\circ b\right) ^{\revsn } \\
\phi _{B\bullet A,D\bullet C}^{\circ }\downarrow &  & \downarrow\phi _{C\circ A,D\circ
B}^{\bullet } \\
\left( \left( d\bullet c\right) \circ \left( b\bullet a\right) \right)
^{\revsn } & \overset{\left( \zeta _{D,C,B,A}\right) ^{\revsn }}{\longrightarrow
} & \left( \left( d\circ b\right) \bullet \left( c\circ a\right) \right)
^{\revsn }%
\end{array}%
\end{equation*}%
We now check the unitality conditions:%
\begin{equation*}
\begin{array}{ccc}
I & \overset{\phi _{0}^{\circ }}{\longrightarrow }I^{\revsn }\overset{\left(
\Delta _{I}^{\bullet }\right) ^{\revsn }}{\longrightarrow } & \left( I\bullet
I\right) ^{\revsn } \\
\Delta _{I}^{\bullet }\downarrow  &  & \uparrow \phi _{I,I}^{\bullet } \\
I\bullet I & \overset{\phi _{0}^{\circ }\bullet \phi _{0}^{\circ }}{%
\longrightarrow } & I^{\revsn }\bullet I^{\revsn }%
\end{array}%
\qquad
\begin{array}{ccc}
b & \overset{\phi _{0}^{\circ }}{\longrightarrow }\left( S_{B}\left(
b\right) \right) ^{\revsn }\overset{\left( \Delta _{I}^{\bullet }\right)
^{\revsn }}{\longrightarrow } & \left( S_{B}\left( b\right) _{1}\bullet
S_{B}\left( b\right) _{2}\right) ^{\revsn }=\left( S_{B}\left( b_{2}\right)
\bullet S_{B}\left( b_{1}\right) \right) ^{\revsn } \\
\Delta _{I}^{\bullet }\downarrow  &  & \uparrow \phi _{I,I}^{\bullet } \\
b_{1}\bullet b_{2} & \overset{\phi _{0}^{\circ }\bullet \phi _{0}^{\circ }}{%
\longrightarrow } & \left( S_{B}\left( b_{1}\right) \right) ^{\revsn }\bullet
\left( S_{B}\left( b_{2}\right) \right) ^{\revsn }%
\end{array}%
\end{equation*}%
\begin{equation*}
\begin{array}{ccc}
J & \overset{\phi _{0}^{\bullet }}{\longrightarrow }J^{\revsn }\overset{\left(
m_{J}^{\circ }\right) ^{\revsn }}{\longleftarrow } & \left( J\circ J\right)
^{\revsn } \\
m_{J}^{\circ }\uparrow  &  & \uparrow \phi _{I,I}^{\circ } \\
J\circ J & \overset{\phi _{0}^{\bullet }\circ \phi _{0}^{\bullet }}{%
\longrightarrow } & J^{\revsn }\circ J^{\revsn }%
\end{array}%
\qquad
\begin{array}{ccc}
kk^{\prime } & \overset{\phi _{0}^{\bullet }}{\longrightarrow }\left(
kk^{\prime }\right) ^{\revsn }=\left( k^{\prime }k\right) ^{\revsn }\overset{%
\left( m_{J}^{\circ }\right) ^{\revsn }}{\longleftarrow } & \left( k^{\prime
}\circ k\right) ^{\revsn } \\
m_{J}^{\circ }\uparrow  &  & \uparrow \phi _{I,I}^{\circ } \\
k\circ k^{\prime } & \overset{\phi _{0}^{\bullet }\circ \phi _{0}^{\bullet }}%
{\longrightarrow } & k^{\revsn }\circ \left( k^{\prime }\right) ^{\revsn }%
\end{array}%
\end{equation*}%
\begin{equation*}
\begin{array}{ccc}
I^{\revsn } & \overset{\left( \varepsilon _{I}^{\bullet }\right) ^{\revsn }}{%
\longrightarrow } & J^{\revsn } \\
\phi _{0}^{\circ }\uparrow  &  & \uparrow \phi _{0}^{\bullet } \\
I & \overset{\varepsilon _{I}^{\bullet }}{\longrightarrow } & J^{\revsn }%
\end{array}%
\qquad
\begin{array}{ccc}
\left( S_{B}\left( b\right) \right) ^{\revsn } & \overset{\left( \varepsilon
_{I}^{\bullet }\right) ^{\revsn }}{\longrightarrow } & \left( \varepsilon
_{B}^{\bullet }S_{B}\left( b\right) \right) ^{\revsn }=\left( \varepsilon
_{B}^{\bullet }\left( b\right) \right) ^{\revsn } \\
\phi _{0}^{\circ }\uparrow  &  & \uparrow \phi _{0}^{\bullet } \\
b & \overset{\varepsilon _{I}^{\bullet }}{\longrightarrow } & \varepsilon
_{B}^{\bullet }\left( b\right)
\end{array}%
\end{equation*}
\end{invisible}
One can define the following natural transformations:%
\begin{eqnarray*}
\gamma _{X,Y,Z,T} &:&\left( \left( X\bullet Y^{\revsn }\right) \circ Z\right)
\bullet T\rightarrow X\circ \left( Z\bullet \left( Y\circ T\right) \right)
,\quad \left( \left( x\bullet y^{\revsn }\right) \circ z\right) \bullet
t\mapsto x\circ \left( z\bullet \left( y\circ t\right) \right) , \\
\delta _{X,Y,Z,T} &:&X\bullet \left( Y\circ \left( Z^{\revsn }\bullet T\right)
\right) \rightarrow \left( \left( X\circ Z\right) \bullet Y\right) \circ
T,\quad x\bullet \left( y\circ \left( z^{\revsn }\bullet t\right) \right)
\mapsto \left( \left( x\circ z\right) \bullet y\right) \circ t,\\
\widetilde{\gamma }_{X,Y,Z,T} &:&X^{\revsn }\bullet \left( \left( Y\bullet
Z\right) \circ T\right) \rightarrow Z\circ \left( \left( X\circ Y\right)
^{\revsn }\bullet T\right) ,\quad x^{\revsn }\bullet \left( \left( y\bullet
z\right) \circ t\right) \mapsto z\circ \left( \left( x\circ y\right) ^{\revsn
}\bullet t\right),  \\
\widetilde{\delta }_{X,Y,Z,T} &:&\left( X\circ \left( Y\bullet Z\right)
\right) \bullet T^{\revsn }\rightarrow \left( X\bullet \left( Z\circ T\right)
^{\revsn }\right) \circ Y,\quad \left( x\circ \left( y\bullet z\right)
\right) \bullet t^{\revsn }\mapsto \left( x\bullet \left( z\circ t\right)
^{\revsn }\right) \circ y.
\end{eqnarray*}
\begin{invisible}
We see that $\gamma $ is well-defined
\begin{eqnarray*}
\gamma _{X,Y,Z,T}\left( \left( \left( x\bullet y^{\revsn }\right) b\circ
z\right) \bullet t\right) &=&\gamma _{X,Y,Z,T}\left( \left( \left(
xb_{1}\bullet \left( S_{B}\left( b_{2}\right) y\right) ^{\revsn }\right) \circ
z\right) \bullet t\right) \\
&=&xb_{1}\circ \left( z\bullet \left( S_{B}\left( b_{2}\right) y\circ
t\right) \right) \\
&=&x\circ b_{1}\left( z\bullet S_{B}\left( b_{2}\right) \left( y\circ
t\right) \right) \\
&=&x\circ \left( b_{1}z\bullet b_{2}S_{B}\left( b_{3}\right) \left( y\circ
t\right) \right) \\
&=&x\circ \left( bz\bullet \left( y\circ t\right) \right) =\gamma
_{X,Y,Z,T}\left( \left( \left( x\bullet y^{\revsn }\right) \circ bz\right)
\bullet t\right) .
\end{eqnarray*}%
It is also a morphism of $B$-bimodules:%
\begin{eqnarray*}
\gamma _{X,Y,Z,T}\left( b\left( \left( \left( x\bullet y^{\revsn }\right)
\circ z\right) \bullet t\right) b^{\prime }\right) &=&\gamma
_{X,Y,Z,T}\left( \left( \left( b_{1}x\bullet \left( yS_{B}\left(
b_{2}\right) \right) ^{\revsn }\right) \circ zb_{1}^{\prime }\right) \bullet
b_{3}tb_{2}^{\prime }\right) \\
&=&b_{1}x\circ \left( zb_{1}^{\prime }\bullet \left( yS_{B}\left(
b_{2}\right) \circ b_{3}tb_{2}^{\prime }\right) \right) \\
&=&b_{1}x\circ \left( zb_{1}^{\prime }\bullet \left( yS_{B}\left(
b_{2}\right) b_{3}\circ tb_{2}^{\prime }\right) \right) \\
&=&bx\circ \left( zb_{1}^{\prime }\bullet \left( y\circ t\right)
b_{2}^{\prime }\right) =b\gamma _{X,Y,Z,T}\left( \left( \left( x\bullet
y^{\revsn }\right) \circ z\right) \bullet t\right) b^{\prime }.
\end{eqnarray*}%
We see that $\delta $ is well-defined%
\begin{eqnarray*}
\delta _{X,Y,Z,T}\left( x\bullet \left( y\circ b\left( z^{\revsn }\bullet
t\right) \right) \right) &=&\delta _{X,Y,Z,T}\left( x\bullet \left( y\circ
\left( \left( zS_{B}\left( b_{1}\right) \right) ^{\revsn }\bullet
b_{2}t\right) \right) \right) \\
&=&\left( \left( x\circ zS_{B}\left( b_{1}\right) \right) \bullet y\right)
\circ b_{2}t \\
&=&\left( \left( x\circ zS_{B}\left( b_{1}\right) \right) \bullet y\right)
b_{2}\circ t \\
&=&\left( \left( x\circ zS_{B}\left( b_{1}\right) b_{2}\right) \bullet
yb_{3}\right) \circ t \\
&=&\left( \left( x\circ z\right) \bullet yb\right) \circ t=\delta
_{X,Y,Z,T}\left( x\bullet \left( yb\circ \left( z^{\revsn }\bullet t\right)
\right) \right) .
\end{eqnarray*}%
It is a morphism of $B$-bimodules:%
\begin{eqnarray*}
\delta _{X,Y,Z,T}\left( b\left( x\bullet \left( y\circ \left( z^{\revsn
}\bullet t\right) \right) \right) b^{\prime }\right) &=&\delta
_{X,Y,Z,T}\left( b_{1}xb_{1}^{\prime }\bullet \left( b_{2}y\circ \left(
\left( S_{B}\left( b_{2}^{\prime }\right) z\right) ^{\revsn }\bullet
tb_{3}^{\prime }\right) \right) \right) \\
&=&\left( \left( b_{1}xb_{1}^{\prime }\circ S_{B}\left( b_{2}^{\prime
}\right) z\right) \bullet b_{2}y\right) \circ tb_{3}^{\prime } \\
&=&\left( \left( b_{1}x\circ b_{1}^{\prime }S_{B}\left( b_{2}^{\prime
}\right) z\right) \bullet b_{2}y\right) \circ tb_{3}^{\prime } \\
&=&\left( \left( b_{1}x\circ z\right) \bullet b_{2}y\right) \circ tb^{\prime
} \\
&=&b\left( \left( x\circ z\right) \bullet y\right) \circ tb^{\prime
}=b\delta _{X,Y,Z,T}\left( x\bullet \left( y\circ \left( z^{\revsn }\bullet
t\right) \right) \right) b^{\prime }.
\end{eqnarray*}%
We see that $\widetilde{\gamma }_{X,Y,Z,T}$ is well-defined%
\begin{eqnarray*}
\widetilde{\gamma }_{X,Y,Z,T}\left( x^{\revsn }\bullet \left( \left( y\bullet
z\right) b\circ t\right) \right)  &=&\widetilde{\gamma }_{X,Y,Z,T}\left(
x^{\revsn }\bullet \left( \left( yb_{1}\bullet zb_{2}\right) \circ t\right)
\right)  \\
&=&zb_{2}\circ \left( \left( x\circ yb_{1}\right) ^{\revsn }\bullet t\right)
\\
&=&z\circ b_{2}\left( \left( x\circ yb_{1}\right) ^{\revsn }\bullet t\right)
\\
&=&z\circ \left( b_{2}\left( x\circ yb_{1}\right) ^{\revsn }\bullet
b_{3}t\right)  \\
&=&z\circ \left( \left( x\circ yb_{1}S_{B}\left( b_{2}\right) \right) ^{\revsn
}\bullet b_{3}t\right)  \\
&=&z\circ \left( \left( x\circ y\right) ^{\revsn }\bullet bt\right)  \\
&=&\widetilde{\gamma }_{X,Y,Z,T}\left( x^{\revsn }\bullet \left( \left(
y\bullet z\right) \circ bt\right) \right) .
\end{eqnarray*}%
We see that $\widetilde{\gamma }_{X,Y,Z,T}$ is a morphism of $B$-bimodules%
\begin{eqnarray*}
\widetilde{\gamma }_{X,Y,Z,T}\left( b\left( x^{\revsn }\bullet \left( \left(
y\bullet z\right) \circ t\right) \right) b^{\prime }\right)  &=&\widetilde{%
\gamma }_{X,Y,Z,T}\left( b_{1}x^{\revsn }b_{1}^{\prime }\bullet \left( \left(
b_{2}y\bullet b_{3}z\right) \circ tb_{2}^{\prime }\right) \right)  \\
&=&\widetilde{\gamma }_{X,Y,Z,T}\left( \left( S_{B}\left( b_{1}^{\prime
}\right) xS_{B}\left( b_{1}\right) \right) ^{\revsn }\bullet \left( \left(
b_{2}y\bullet b_{3}z\right) \circ tb_{2}^{\prime }\right) \right)  \\
&=&b_{3}z\circ \left( \left( S_{B}\left( b_{1}^{\prime }\right) xS_{B}\left(
b_{1}\right) \circ b_{2}y\right) ^{\revsn }\bullet tb_{2}^{\prime }\right)  \\
&=&bz\circ \left( \left( S_{B}\left( b_{1}^{\prime }\right) x\circ
y\right) ^{\revsn }\bullet tb_{2}^{\prime }\right)  \\
&=&bz\circ \left( \left( x\circ y\right) ^{\revsn }b_{1}^{\prime }\bullet
tb_{2}^{\prime }\right)  \\
&=&b\left( z\circ \left( \left( x\circ y\right) ^{\revsn }\bullet t\right)
\right) b^{\prime } \\
&=&b\widetilde{\gamma }_{X,Y,Z,T}\left( x^{\revsn }\bullet \left( \left(
y\bullet z\right) \circ t\right) \right) b^{\prime }.
\end{eqnarray*}%
We see that $\widetilde{\delta }_{X,Y,Z,T}$ is well-defined%
\begin{eqnarray*}
\widetilde{\delta }_{X,Y,Z,T}\left( \left( x\circ b\left( y\bullet z\right)
\right) \bullet t^{\revsn }\right)  &=&\widetilde{\delta }_{X,Y,Z,T}\left(
\left( x\circ \left( b_{1}y\bullet b_{2}z\right) \right) \bullet t^{\revsn
}\right)  \\
&=&\left( x\bullet \left( b_{2}z\circ t\right) ^{\revsn }\right) \circ b_{1}y
\\
&=&\left( x\bullet \left( b_{2}z\circ t\right) ^{\revsn }\right) b_{1}\circ y
\\
&=&\left( xb_{1}\bullet \left( b_{3}z\circ t\right) ^{\revsn }b_{2}\right)
\circ y \\
&=&\left( xb_{1}\bullet \left( S_{B}\left( b_{2}\right) b_{3}z\circ t\right)
^{\revsn }\right) \circ y \\
&=&\left( xb\bullet \left( z\circ t\right) ^{\revsn }\right) \circ y \\
&=&\widetilde{\delta }_{X,Y,Z,T}\left( \left( xb\circ \left( y\bullet
z\right) \right) \bullet t^{\revsn }\right) .
\end{eqnarray*}%
We see that $\widetilde{\delta }_{X,Y,Z,T}$ is a morphism of $B$-bimodules%
\begin{eqnarray*}
\widetilde{\delta }_{X,Y,Z,T}\left( b\left( \left( x\circ \left( y\bullet
z\right) \right) \bullet t^{\revsn }\right) b^{\prime }\right)  &=&\widetilde{%
\delta }_{X,Y,Z,T}\left( \left( b_{1}x\circ \left( yb_{1}^{\prime }\bullet
zb_{2}^{\prime }\right) \right) \bullet b_{2}t^{\revsn }b_{3}^{\prime }\right)
\\
&=&\widetilde{\delta }_{X,Y,Z,T}\left( \left( b_{1}x\circ \left(
yb_{1}^{\prime }\bullet zb_{2}^{\prime }\right) \right) \bullet \left(
S_{B}\left( b_{3}^{\prime }\right) tS_{B}\left( b_{2}\right) \right) ^{\revsn
}\right)  \\
&=&\left( b_{1}x\bullet \left( zb_{2}^{\prime }\circ S_{B}\left(
b_{3}^{\prime }\right) tS_{B}\left( b_{2}\right) \right) ^{\revsn }\right)
\circ yb_{1}^{\prime } \\
&=&\left( b_{1}x\bullet \left( z\circ b_{2}^{\prime }S_{B}\left(
b_{3}^{\prime }\right) tS_{B}\left( b_{2}\right) \right) ^{\revsn }\right)
\circ yb_{1}^{\prime } \\
&=&\left( b_{1}x\bullet \left( z\circ tS_{B}\left( b_{2}\right) \right)
^{\revsn }\right) \circ yb^{\prime } \\
&=&\left( b_{1}x\bullet b_{2}\left( z\circ t\right) ^{\revsn }\right) \circ
yb^{\prime } \\
&=&b\left( \left( x\bullet \left( z\circ t\right) ^{\revsn }\right) \circ
y\right) b^{\prime } \\
&=&b\widetilde{\delta }_{X,Y,Z,T}\left( \left( x\circ \left( y\bullet
z\right) \right) \bullet t^{\revsn }\right) b^{\prime }.
\end{eqnarray*}
They are clearly natural.
\end{invisible}
Consequently, we have%
\begin{eqnarray*}
\psi _{X,Y} &:&X^{\revsn }\bullet Y\rightarrow J\circ \left( I\bullet \left(
X\circ Y\right) \right) ,\quad x^{\revsn }\bullet y\mapsto 1_{\Bbbk }\circ
\left( 1_{B}\bullet \left( x\circ y\right) \right) , \\
\varphi _{X,Y} &:&X\bullet Y^{\revsn }\rightarrow \left( \left( X\circ
Y\right) \bullet I\right) \circ J,\quad x\bullet y^{\revsn }\mapsto ((x\circ
y)\bullet 1_{B})\circ 1_{\Bbbk }, \\
\widetilde{\psi }_{X,Y} &:&X^{\revsn }\bullet \left( \left( Y\bullet I\right)
\circ J\right) \rightarrow \left( X\circ Y\right) ^{\revsn },\quad x^{\revsn
}\bullet \left( \left( y\bullet b\right) \circ k\right) \mapsto \left(
x\circ yS_{B}\left( b\right) k\right) ^{\revsn }, \\
\widetilde{\varphi }_{X,Y} &:&\left( J\circ \left( I\bullet X\right) \right)
\bullet Y^{\revsn }\rightarrow \left( X\circ Y\right) ^{\revsn },\quad \left(
k\circ \left( b\bullet x\right) \right) \bullet y^{\revsn }\mapsto \left(
kS_{B}\left( b\right) x\circ y\right) ^{\revsn }.
\end{eqnarray*}
It is now straightforward to verify that the diagrams required to have a reversion are commutative.
\begin{invisible}
 They clearly obey \eqref{form:psi1}, \eqref{form:psi2}, \eqref{form:varphi1}, \eqref{form:varphi2} and \eqref{form:gammadelta}.\\
 On elements, \eqref{form:phi0gamma} is  easily verified as
 \begin{equation*}
\left[
\begin{array}{c}
a\bullet b\overset{\left( l_{I}^{\circ }\right) ^{-1}\bullet I}{\longmapsto }%
\left( 1_{B}\circ a\right) \bullet b\overset{\left( \left( I\bullet \phi
_{0}^{\circ }\right) \Delta _{I}^{\bullet }\circ I\right) \bullet I}{%
\longmapsto }\left( \left( 1_{B}\bullet S_{B}\left( 1_{B}\right) ^{\revsn
}\right) \circ a\right) \bullet b \\
=\left( \left( 1_{B}\bullet \left( 1_{B}\right) ^{\revsn }\right) \circ
a\right) \bullet b\overset{\gamma }{\longmapsto }1_{B}\circ \left( a\bullet
\left( 1_{B}\circ b\right) \right) \overset{l_{I\bullet I}^{\circ }\left(
I\circ \left( I\bullet l_{I}^{\circ }\right) \right) }{\longmapsto }a\bullet
b%
\end{array}%
\right] =\mathrm{Id}
\end{equation*}
On elements, \eqref{form:phi0delta}  is easily verified as
\begin{equation*}
\left[
\begin{array}{c}
a\bullet b\overset{I\bullet \left( r_{I}^{\circ }\right) ^{-1}}{\longmapsto }%
a\bullet \left( b\circ 1_{B}\right) \overset{I\bullet \left( I\circ \left(
\phi _{0}^{\circ }\bullet I\right) \Delta _{I}^{\bullet }\right) }{%
\longmapsto }a\bullet \left( b\circ \left( S_{B}\left( 1_{B}\right) ^{\revsn
}\bullet 1_{B}\right) \right)  \\
=a\bullet \left( b\circ \left( \left( 1_{B}\right) ^{\revsn }\bullet
1_{B}\right) \right) \overset{\delta }{\longmapsto }\left( \left( a\circ
1_{B}\right) \bullet b\right) \circ 1_{B}\overset{r_{I\bullet I}^{\circ
}\left( \left( r_{I}^{\circ }\bullet I\right) \circ I\right) }{\longmapsto }%
a\bullet b%
\end{array}%
\right] =\mathrm{Id.}
\end{equation*}
The diagrams \eqref{form:phipsi2}, \eqref{form:phivarphi2}, \eqref{form:varphipsitilde} are clearly commutative.\\
On elements \eqref{form:gammatildeIIJ} is verified as
\begin{equation*}
\left[
\begin{array}{c}
x^{\revsn }\overset{\left[ r_{X^{\revsn }}^{\bullet }\left( X^{\revsn }\bullet
l_{J}^{\circ }\right) \right] ^{-1}}{\longmapsto }x^{\revsn }\bullet \left(
1_{B}\circ 1_{\Bbbk }\right) \overset{X^{\revsn }\bullet \left( \Delta
_{I}^{\bullet }\circ J\right) }{\longmapsto }x^{\revsn }\bullet \left( \left(
1_{B}\bullet 1_{B}\right) \circ 1_{\Bbbk }\right)  \\
\overset{\widetilde{\gamma} }{\longmapsto }1_{B}\circ \left( \left( x\circ 1_{B}\right)
^{\revsn }\bullet 1_{\Bbbk }\right) \overset{\left( r_{X}^{\circ }\right)
^{\revsn }l_{\left( X\circ I\right) ^{\revsn }}^{\circ }\left( I\circ r_{\left(
X\circ I\right) ^{\revsn }}^{\bullet }\right) }{\longmapsto }x^{\revsn }%
\end{array}%
\right] =\mathrm{Id}
\end{equation*}
On elements \eqref{form:deltatildeJII} is verified as
\begin{equation*}
\left[
\begin{array}{c}
x^{\revsn }\overset{\left[ l_{X^{\revsn }}^{\bullet }\left( r_{J}^{\circ
}\bullet X^{\revsn }\right) \right] ^{-1}}{\longmapsto }\left( 1_{\Bbbk }\circ
1_{B}\right) \bullet x^{\revsn }\overset{\left( J\circ \Delta _{I}^{\bullet
}\right) \bullet X^{\revsn }}{\longmapsto }\left( 1_{\Bbbk }\circ \left(
1_{B}\bullet 1_{B}\right) \right) \bullet x^{\revsn } \\
\overset{\widetilde{\delta }}{\longmapsto }\left( 1_{\Bbbk }\bullet \left(
1_{B}\circ x\right) ^{\revsn }\right) \circ 1_{B}\overset{\left( l_{X}^{\circ
}\right) ^{\revsn }r_{\left( I\circ X\right) ^{\revsn }}^{\circ }\left(
l_{\left( I\circ X\right) ^{\revsn }}^{\bullet }\circ I\right) }{\longmapsto }%
x^{\revsn }%
\end{array}%
\right] =\mathrm{Id}
\end{equation*}
The diagrams \eqref{form: psi1gamma}, \eqref{form:varphi1delta}, \eqref{form:phipsi2NEW} and \eqref{form:phivarphi2NEW} are clearly commutative.\\
On elements \eqref{form:phi0gamma-bul} is verified as:
\begin{equation*}
\left(
\begin{array}{c}
k\overset{\left( r_{J}^{\bullet }\left( r_{J}^{\circ }\bullet J\right)
\right) ^{-1}}{\longmapsto }\left( k\circ 1_{B}\right) \bullet 1_{\Bbbk }%
\overset{\left( \left( J\bullet \phi _{0}^{\bullet }\right) \Delta
_{J}^{\bullet }\circ I\right) \bullet J}{\longmapsto }\left( \left( k\bullet
\left( 1_{\Bbbk }\right) ^{\revsn }\right) \circ 1_{B}\right) \bullet 1_{\Bbbk
} \\
\overset{\gamma _{J,J,I,J}}{\longmapsto }k\circ \left( 1_{B}\bullet \left(
1_{\Bbbk }\circ 1_{\Bbbk }\right) \right) \overset{r_{J}^{\circ }\left(
I\circ r_{I}^{\bullet }\right) \left( J\circ \left( I\bullet m_{J}^{\circ
}\right) \right) }{\longmapsto }k%
\end{array}%
\right) =\mathrm{Id.}
\end{equation*}
On elements \eqref{form:phi0delta-bul} is verified as:
\begin{equation*}
\left(
\begin{array}{c}
k\overset{\left( l_{J}^{\bullet }\left( J\bullet l_{J}^{\circ }\right)
\right) ^{-1}}{\longmapsto }1_{\Bbbk }\bullet \left( 1_{B}\circ k\right)
\overset{J\bullet \left( I\circ \left( \phi _{0}^{\bullet }\bullet J\right)
\Delta _{J}^{\bullet }\right) }{\longmapsto }1_{\Bbbk }\bullet \left(
1_{B}\circ \left( \left( 1_{\Bbbk }\right) ^{\revsn }\bullet k\right) \right)
\\
\overset{\delta _{J,I,J,J}}{\longmapsto }\left( \left( 1_{\Bbbk }\circ
1_{\Bbbk }\right) \bullet 1_{B}\right) \circ k\overset{l_{J}^{\circ }\left(
l_{I}^{\bullet }\circ I\right) \left( \left( m_{J}^{\circ }\bullet I\right)
\circ J\right) }{\longmapsto }k%
\end{array}%
\right) =\mathrm{Id.}
\end{equation*}
\end{invisible}

\subsubsection{On Hopf monoids}

Inspired by \cite[Proposition 4.9]{BCZ} and \cite[Theorem 7.17]{Bohm-Lack}, we prove the following result.

\begin{proposition}%[\rd{upgrade:2026/03/03}]
\label{pro:Bohm}
In the setting of
\cref{es:bimod} assume that $B$ is a $\Bbbk$-Hopf algebra.
 Then, for an object $H$ in $\mathsf{Bimon}(\Cc,\circ,\bullet)$,
 the following conditions are equivalent:
\begin{enumerate}
    \item $H$ is a $\Bbbk$-Hopf algebra;
     \item $H$ is a Hopf monoid in the duoidal category $(\Cc,\circ,I,\bullet,J)$;
     \item the Galois maps
 $\varkappa _{P,X,Q}^{H} :(P\bullet X)\circ Q\to P\bullet (X\circ Q)$
    are invertible, for any object $X$ in $\Cc$,  right $H$-module $\left( P,\mu^\circ _{P}\right) $ and left $H$-comodule $\left( Q,\rho^\bullet _{Q}\right)$;
      \item the co-Galois maps
$\varsigma _{P,X,Q}^{H}:P\circ (X\bullet Q)\to (P\circ X)\bullet Q$
    are invertible, for any object $X$, right $H$-comodule $\left( P,\rho^\bullet _{P}\right) $ and left $H$-module $\left( Q,\mu^\circ _{Q}\right) $;
    \item $g_r$ is an isomorphism in $\Cc$ for any morphism $g:H\to A$ in $\Mon(\Cc^\circ)$;
    \item $g^r$ is an isomorphism in $\Cc$ for any morphism $g:C\to H$ in $\Comon(\Cc^\bullet)$;
     \item the fusion morphism $(\id_H)_r$ is invertible;
    \item $g^r$ is an isomorphism in $\Cc$ for $g:C\to H,x\bullet b\mapsto xb$, where $C:=H\bullet I$ is a $B$-bimodule via $b^{\prime }\left( x\bullet b\right) b^{\prime \prime }=b^{\prime
}x\bullet bb^{\prime \prime }$ and a comonoid via $\Delta_C^\bullet(x\bullet b)=(x_1\bullet b_1)\bullet (x_2\bullet b_2)$ and $\varepsilon_C^\bullet(x\bullet b)=\varepsilon_H^\bullet(x)\varepsilon_B^\bullet(b)$.
\end{enumerate}
As a consequence,  $\Hopfmon(\Cc,\circ,\bullet)$ is isomorphic to the coslice category $B/\Hopfmon_\Bbbk$.
\end{proposition}

\begin{proof} In view of \cref{rmk:bimonbimod}, we identify a bimonoid with a bialgebra $H$ together with a bialgebra map $u_H^\circ:B\to H$ which defines its $B$-bimodule structure.

$(1)\Rightarrow(2)$ Given an antipode $S_H:H\to H$, define $\Sigma _{H}:H\rightarrow H^{\revsn },\;
h\mapsto S_{H}\left( h\right) ^{\revsn }$. Note that this is a morphism of $B$-bimodules as, being the $B$-bimodule structure induced by $u_H^\circ$, we get
\[
\begin{split}
S_H(bhb')
&=S_H(u_H^\circ(b)hu_H^\circ(b'))
=S_Hu_H^\circ(b')S_H(h)S_Hu_H^\circ(b)
=u_H^\circ S_B(b')S_H(h)u_H^\circ S_B(b)
\\&=S_B(b')S_H(h)S_B(b)
\end{split}
\]
and hence
$\Sigma_H(bhb')=S_H(bhb')^\revsn
=(S_B(b')S_H(h)S_B(b))^\revsn
=bS_H(h)^\revsn b'=b\Sigma_H(h)b'$.
We have%
\begin{equation*}
\begin{split}
\left( \Sigma _{H}\star \mathrm{Id}_{H}\right) \left( h\right) &=\left(
J\circ \left( I\bullet m_{H}^{\circ }\right) \right) \psi _{H,H}\left(
S_{H}\left( h_{1}\right) ^{\revsn }\bullet h_{2}\right) =\left( J\circ \left( I\bullet m_{H}^{\circ }\right) \right) \left(
1_{\Bbbk }\circ \left( 1_{B}\bullet \left( S_{H}\left( h_{1}\right) \circ
h_{2}\right) \right) \right) \\
&=1_{\Bbbk }\circ \left( 1_{B}\bullet \left( S_{H}\left( h_{1}\right)
h_{2}\right) \right) =1_{\Bbbk }\circ \left( 1_{B}\bullet 1_{H}\right)
\varepsilon _{H}^{\bullet }\left( h\right)=(J\circ (I\bullet u_H^\circ)\Delta_I^\bullet)(r^\circ_J)^{-1}\varepsilon^\bullet_H(h)\end{split}
\end{equation*}%
and
\begin{equation*}
\begin{split}
(\mathrm{Id}_{H}\star \Sigma _{H})(h)& =((m_{H}^{\circ }\bullet I)\circ
J)\varphi _{H,H}(h_{1}\bullet S_{H}(h_{2})^{\revsn })  =((m_{H}^{\circ }\bullet I)\circ J)(((h_{1}\circ S_{H}(h_{2}))\bullet
1_{B})\circ 1_{\Bbbk }) \\
& =((h_{1}S_{H}(h_{2})\bullet 1_{B})\circ 1_{\Bbbk })=\varepsilon
_{H}^{\bullet }(h)(1_{H}\bullet 1_{B})\circ 1_{\Bbbk }  =((u_{H}^{\circ }\bullet I)\Delta _{I}^{\bullet }\circ J)(l_{J}^{\circ })^{-1}\varepsilon _{H}^{\bullet }(h).
\end{split}%
\end{equation*}
Therefore, $\Sigma_H$ is an antipode in the duoidal sense.

$(2)\Rightarrow(4)$ This is \cref{thm:Galinv}.

$(4)\Rightarrow(6)$ This follows by \cref{equivalentgr}.

$(6)\Rightarrow(8)$ It is obvious
as $g:C\to H$ as in the statement is a morphism of $\bullet$-comonoids.
\begin{invisible}
Since $H$ is an object in $\mathsf{Comon}(_{B}\mm_{B},\ot,\Bbbk)$, it is an object in $\mathsf{Comon}(\mathsf{Vec},\ot,\Bbbk)$ such that $\Delta^{\bullet}_{H}$ and $\varepsilon^{\bullet}_{H}$ are $B$-bimodule morphisms. Then, since also $I$ is in $\mathsf{Comon}(\mathsf{Vec},\ot,\Bbbk)$, the object $C=H\bullet I$ is in $\mathsf{Comon}(\mathsf{Vec},\ot,\Bbbk)$ with the tensor product coalgebra structure. In order to have that $C$ is in $\mathsf{Comon}(_{B}\mm_{B},\ot,\Bbbk)$ it remains to prove that $\Delta_{C}^{\bullet}$ and $\varepsilon_{C}^{\bullet}$ are $B$-bilinear morphisms with respect to the $B$-bimodule structure of $C$. The morphisms $\Delta_{C}^{\bullet}$ and $\varepsilon_{C}^{\bullet}$ are left $B$-linear since $\Delta^{\bullet}_{H}$ and $\varepsilon^{\bullet}_{H}$ are left $B$-linear, respectively. Moreover, $\Delta_{C}^{\bullet}$ and $\varepsilon_{C}^{\bullet}$ are right $B$-linear since $B$ is $\Bbbk$-bialgebra.\\
We have
$\Delta^{\bullet}_{H}g(x\bullet b)=\Delta_{H}^{\bullet}(xb)=(x_{1}b_{1})\bullet(x_{2}b_{2})=g(x_{1}\bullet b_{1})\bullet g(x_{2}\bullet b_{2})=(g\bullet g)\Delta_{C}^{\bullet}(x\bullet b)$ using that $\Delta_{H}^{\bullet}$ is right $B$-linear and $\varepsilon_{H}^{\bullet}g(x\bullet b)=\varepsilon_{H}^{\bullet}(xb)=\varepsilon_{H}^{\bullet}(x)\varepsilon_{B}^{\bullet}(b)=\varepsilon_{C}^{\bullet}(x\bullet b)$ using that $\varepsilon_{H}^{\bullet}$ is right $B$-linear. The forgetful functor $R:{}_B\mm_B\to{}_B\mm$ has a left adjoint $L:=(-)\bullet I$. The counit of the adjunction is $\epsilon_Y:Y\bullet I\to Y,y\bullet b\mapsto yb$. Then $g:=\epsilon_H.$
\end{invisible}

$(8)\Rightarrow(1)$
By assumption, we have an isomorphism in $\Cc$ defined as:%
\begin{equation*}
g^{r}:(H\bullet I) \circ H%
\longrightarrow \left((H\bullet I) \circ J\right) \bullet H,\ (x\bullet b)\circ y\mapsto ((x_1\bullet b_1)\circ 1_{\Bbbk})\bullet x_2b_2y.
\end{equation*}%
Clearly, we have the $\Bbbk $-linear isomorphisms
\begin{align*}
f&:(H\bullet I)
\circ H\overset{\cong }{\longrightarrow }H\bullet H,\ \left( x\bullet
b\right) \circ y\mapsto x\bullet by,\\
f'&:\left[ (H\bullet I) \circ J\overset{\cong }{\longrightarrow }H\bullet J\overset{\cong }{%
\longrightarrow }H\right] ,\left( x\bullet b\right) \circ k\mapsto x\varepsilon_B^\bullet(b)k.
\end{align*}
Their composition by $g^r$ yields the $\Bbbk$-linear isomorphism
\begin{align*}
H\bullet H&\overset{f^{-1}}{\longrightarrow }(H\bullet I)\circ H%
\overset{g^r}{\longrightarrow }\left( (H\bullet I) \circ J\right) \bullet H\overset{f'\bullet H}{\longrightarrow }%
H\bullet H,\\
x\bullet y&\overset{f^{-1}}{\mapsto }\left( x\bullet 1_{B}\right) \circ y%
\overset{g^r}{\mapsto }\left( \left( x_{1}\bullet 1_{B}\right) \circ 1_{\Bbbk}\right) \bullet
 x_{2}y \overset{f'\bullet H}{\mapsto }x_{1}\bullet 
x_{2}y.
\end{align*}
Since this is the canonical map, we deduce that $B$ has antipode, again by the symmetric version of \cite[Theorem 1.8]{Ve}.

Similarly, one proves that $(2)\Rightarrow(3)\Rightarrow(5)$. Moreover, $(5)\Rightarrow (7)$ is obvious. It remains to prove $(7)\Rightarrow (1)$. If $(\id_H)_r$ is invertible, consider again the $\Bbbk$-linear isomorphism $f:(H\bullet I)\circ
H\rightarrow H\bullet H$ as before. Note also that the $\Bbbk$-linear map $\mathrm{can}_{I}:H\bullet I\rightarrow
H\bullet I,\ x\bullet b\mapsto xb_{1}\bullet b_{2}$ is invertible with inverse
$x\bullet b\mapsto xS\left( b_{1}\right) \bullet b_{2}$ as $B$ has antipode by assumption and it is a morphism of right $B$-modules so that we can consider $\mathrm{can}_I\circ H$. Then, we have the $\Bbbk$-linear isomorphism
\begin{align*}
H\bullet H&\overset{f^{-1}}{\longrightarrow }(H\bullet I)\circ H\overset{\mathrm{can}%
_{I}\circ H}{\longrightarrow }\left( H\bullet I\right) \circ H\overset{%
\left( \mathrm{Id}_{H}\right) _{r}}{\longrightarrow }H\bullet H,\\
x\bullet y&\overset{f^{-1}}{\longmapsto }\left( x\bullet 1_{B}\right) \circ y%
\overset{\mathrm{can}_{I}\circ H}{\longmapsto }\left( x\bullet 1_{B}\right)
\circ y\overset{\left( \mathrm{Id}_{H}\right) _{r}}{\longmapsto }%
xy_{1}\bullet y_{2}.
\end{align*}
Since the latter is the canonical map, we deduce that $H$ has antipode.
\end{proof}

\subsection{Quasitriangular Hopf algebras}

In order to apply our results on tensor-braided duoidal categories we now turn to the case when $B$ is a quasitriangular Hopf algebra, a notion introduced by Drinfel'd in the seminal paper \cite{Dr87}.
Let us recall the definition and some basic properties; for furthere details, we refer the reader to \cite{Kassel} and \cite{Majid-book}.

Given a bialgebra $H$ and an element $T\in H\otimes H$, we adopt the short notation $T=T^{i}\otimes T_{i
}$, summation understood. If $T$ is invertible we write $T^{-1}=\overline{T}=\overline{T}^{i}\otimes\overline{T}_{i}$. Moreover, we employ the leg notation $T_{12}=T\otimes1_{H}$, $T_{23}=1_{H}\otimes T$, $T_{13}=(\mathrm{Id}_{H}\otimes\tau_{H,H})(T_{12})$, $T_{21}=T^{\mathrm{op}}\otimes1_{H}$, etc.

A bialgebra $H$ is \textbf{quasitriangular} if it admits an invertible element $\Rr\in H\otimes H$, the \textit{universal $\Rr$-matrix} or \textit{quasitriangular structure}, such that $H$ is quasi-cocommutative, i.e.,
$\Delta^\mathrm{op}(\cdot)=\Rr\Delta(\cdot)\Rr^{-1}$,
and the hexagon equations are satisfied i.e. $(\mathrm{Id}_H\otimes\Delta)(\Rr)=\Rr_{13}\Rr_{12}$ and $(\Delta\otimes\mathrm{Id}_H)(\Rr)=\Rr_{13}\Rr_{23}$.
% \begin{equation*}
% (\mathrm{Id}_H\otimes\Delta)(\Rr)=\Rr_{13}\Rr_{12},\qquad
%     (\Delta\otimes\mathrm{Id}_H)(\Rr)=\Rr_{13}\Rr_{23}.
% \end{equation*}
A morphism of quasitriangular bialgebras $f:(A,\Rr)\to(H,\Ss)$ is morphism of bialgebras $f:A\to H$ with $(f\ot f)(\Rr)=\Ss$.
Quasitriangular bialgebras and their morphisms form the category $\QTrBialg$.

We denote by  $\QTrHopf$ its full subcategory consisting of quasitriangular Hopf algebras.\medskip

Clearly, any cocommutative bialgebra is quasitriangular with $\mathcal{R}=1\otimes 1$.

By e.g.\ \cite[Theorem VIII.2.4]{Kassel}, a quasitriangular bialgebra $(H,\Rr)$ satisfies $(\varepsilon\otimes\mathrm{Id}_{H})(\Rr)=1_{H}=(\mathrm{Id}_{H}\otimes\varepsilon)(\Rr)$ and, if $H$ has an antipode $S$, also
$(S\otimes\mathrm{Id})(\Rr)=\Rr^{-1}$ and $(\mathrm{Id}\otimes S)(\Rr^{-1})=\Rr.$

% $%
% H_{\left( \mathcal{S}\right) }=\mathcal{S}_{\left( \ell \right) }\mathcal{S}%
% _{\left( r\right) }$, see \cite[ Definition 12.2.14]{Radford-book}, where%
% \begin{equation*}
% \mathcal{S}_{\left( \ell \right) } =\left\{ \left( \mathrm{Id}\otimes
% f\right) \left( \mathcal{S}\right) \mid f\in H^{\star }\right\},  \qquad
% \mathcal{S}_{\left( r\right) } =\left\{ \left( f\otimes \mathrm{Id}\right)
% \left( \mathcal{S}\right) \mid f\in H^{\star }\right\} .
% \end{equation*}%
% Note that $\mathcal{S}_{\left( \ell \right) }$, $\mathcal{S}_{\left(
% r\right) }$ and $H_{\left( \mathcal{S}\right) }$ are finite-dimensional Hopf
% subalgebras of $H$ (\cite[Proposition 12.2.13]{Radford-book}).

%Writing  $\mathcal{S}=\sum_{i=1}^{t}a_{i}\otimes b_{i}$ with $\left\{a_{1},\ldots ,a_{t}\right\} $ and $\left\{ b_{1},\ldots ,b_{t}\right\} $ linearly independent, let $a_{i}^{\star },b_{i}^{\star }\in H^{\star }$ be such that $a_{i}^{\star }\left( a_{j}\right) =\delta _{ij}$ and $b_{i}^{\star }\left( b_{j}\right) =\delta _{ij}.$ Then $\left\{ a_{1},\ldots ,a_{t}\right\} $ is a basis of $\mathcal{S}_{\left( \ell \right) }$ and $\left\{ b_{1},\ldots ,b_{t}\right\} $ is a basis of $\mathcal{S}_{\left(r\right) }$.
\medskip

Our next goal is to show that the duoidal category of $B$-bimodules becomes tensor-braided when $B$ is a quasitriangular Hopf algebra, and to identify the cocommutative bimonoids in it.

\begin{remark}\label{rmk:duoidaleqbim}%\as{[Added on 20/04/2026]}
    In the setting of \cref{es:bimod}, if $(B,\Ss)$ is a quasitriangular Hopf algebra, then $\left( \left( -\right) ^{\revsn },\phi^\circ \right)$ and $\left( \left( -\right) ^{\revsn },\phi^\bullet \right)$ turn out to be strong monoidal functors. Indeed, thanks to the invertibility of the antipode of a quasitriangular Hopf algebra, see \cite[Proposition 1]{Radford}, the inverses $(\phi^{\circ}_{X,Y})^{-1}:(Y\circ X)^{\revsn }\to X^{\revsn }\circ Y^{\revsn }$, $(y\circ x)^{\revsn }\mapsto x^{\revsn }\circ y^{\revsn }$, for $X,Y$ in $\Cc$, and $(\phi^{\circ}_{0})^{-1}:B^{\revsn}\to B$, $b^{\revsn}\mapsto S^{-1}_{B}(b)$, are well-defined. 
    \begin{invisible}
   Since
\[
\begin{split}
(\phi^{\circ}_{X,Y})^{-1}((yb\circ x)^{\revsn})&=x^{\revsn}\circ(yb)^{\revsn}=x^{\revsn}\circ(yS_{B}S^{-1}_{B}(b))^{\revsn}=x^{\revsn}\circ S^{-1}_{B}(b)y^{\revsn}=x^{\revsn}S^{-1}_{B}(b)\circ y^{\revsn}=(bx)^{\revsn}\circ y^{\revsn}\\&=(\phi^{\circ}_{X,Y})^{-1}((y\circ bx)^{\revsn}),
\end{split}
\]
we have that $(\phi^{\circ}_{X,Y})^{-1}$ is well defined, and it is the $\Bbbk$-linear inverse of $\phi^{\circ}_{X,Y}$, hence the inverse of $\phi^{\circ}_{X,Y}$ in $\Cc$. Clearly, $(\phi^{\circ}_{0})^{-1}$ is the $\Bbbk$-linear inverse of $\phi^{\circ}_{0}$, hence its inverse in $\Cc$, so $((-)^{\revsn},\phi^{\circ})$ is strong monoidal.\textbf{}    
    \end{invisible} 
    Moreover, one can define $(\phi^{\bullet}_{X,Y})^{-1}:(Y\bullet X)^{\revsn }\to X^{\revsn }\bullet Y^{\revsn }$, $(y\bullet x)^{\revsn }\mapsto x^{\revsn }\bullet y^{\revsn }$, and $(\phi^{\bullet}_{0})^{-1}:J^{\revsn}\to J$, $k^{\revsn}\mapsto k$, that make $((-)^{\revsn},\phi^{\bullet})$ a strong monoidal functor. 
    
    In fact, the double lax monoidal functor $\left( \left( -\right) ^{\revsn },\phi^\circ,\phi^{\bullet} \right)$ becomes a duoidal isomorphism of categories. Indeed, one can define a functor $(-)^{\revsnp}:\Cc\to\Cc$ which sends an object $X\in\Cc$ to the $B$-bimodule $X^{\revsnp}$ whose underline vector space is $X$ and $B$-bimodule structure is $bx^{\revsnp}b':=(S^{-1}_{B}(b')xS^{-1}_{B}(b))^{\revsnp}$, and a morphism $f:X\rightarrow Y$ in $\Cc$ to the morphism $f^{\revsnp
}:X^{\revsnp}\rightarrow Y^{\revsnp}$ given by $f^{\revsnp}\left( x^{\revsnp}\right) :=f\left( x\right) ^{\revsnp}$. This provides an inverse to the functor $(-)^{\revsn}:\Cc\to\Cc$.
\end{remark}

\begin{proposition}
\label{pro:bimoncoc}
In the setting of \cref{es:bimod}, if $(B,\Ss)$ is a quasitriangular Hopf algebra, then the duoidal category $(\mathcal{C}, \circ, I, \bullet, J)$ is $\bullet$\,-braided. Moreover, the category  $\cBimon(\Cc,\circ,\bullet)$ is isomorphic to the coslice category of quasitriangular bialgebras under $(B,\Ss)$, i.e.\ $(B,\Ss)/\QTrBialg$.  
\end{proposition}

\begin{proof}
Since $\Ss$ is a quasitriangular structure for $B$, then $\overline{\Ss}^\op$ is a quasitriangular structure for $B$ and $\Ss^\op$ is a quasitriangular structure for $B^\op$, see \cite[Exercise 2.1.3]{Majid-book}. Therefore, by e.g. \cite[Theorem 2.2]{Chen-quasi}, we have that $(\id\ot\tau\ot\id)(\overline{\Ss}^\op\ot\Ss^\op)$ is a quasitriangular structure for $B\ot B^\op$ so that the monoidal category $(_{B\ot B^\op}\mm,\ot,\Bbbk)$ of left $B\ot B^\op$-modules is braided. By identifying left $B\ot B^\op$-modules with $B$-bimodules, we get that the monoidal category $(\Cc,\bullet,J)=({}_B\mathfrak{M}_B,\ot_\Bbbk,\Bbbk)$ is braided, where, for all $X,Y\in\Cc$, the braiding  $\sigma_{X,Y}:X\ot Y\to Y\ot X$ is given by $\sigma_{X,Y}(x\ot y)=\Ss^{-1}(y \ot x)\Ss$.
\begin{invisible}
 Set $\Rr:=(\id\ot\tau\ot\id)(\overline{\Ss}^\op\ot\Ss^\op)
=\overline{\Ss}_i\ot\Ss_i\ot\overline{\Ss}^i\ot\Ss^i.$ 
Then $\Rr^\op(y\ot x)
=(\overline{\Ss}^i\ot\Ss^i\ot \overline{\Ss}_i\ot\Ss_i)(y\ot x)
=(\overline{\Ss}^i\ot\Ss^i)y\ot ( \overline{\Ss}_i\ot\Ss_i)x
=\overline{\Ss}^iy\Ss^i\ot  \overline{\Ss}_ix\Ss_i
=\Ss^{-1}(y\ot x)\Ss.$
\end{invisible}

It is indeed $\bullet$\,-braided, as follows from the following computation:
\begin{align*}
&\zeta _{X,W,Z,Y}\left( \sigma _{W,X}\circ \sigma _{Y,Z}\right) (\left(
a\bullet b\right) \circ \left( c\bullet d\right) ) 
=\zeta _{X,W,Z,Y}\left(\left(
\overline{\mathcal{S}}^{i}b\mathcal{S}^{j}\bullet \overline{\mathcal{S}}_{i}a%
\mathcal{S}_{j}\right) \circ \left( \overline{\mathcal{S}}^{s}d\mathcal{S}%
^{t}\bullet \overline{\mathcal{S}}_{s}c\mathcal{S}_{t}\right) \right) \\
&=\left( \overline{\mathcal{S}}^{i}b\mathcal{S}^{j}\circ \overline{\mathcal{%
S}}^{s}d\mathcal{S}^{t}\right) \bullet \left( \overline{\mathcal{S}}_{i}a%
\mathcal{S}_{j}\circ \overline{\mathcal{S}}_{s}c\mathcal{S}_{t}\right)  
=\left( \overline{\mathcal{S}}^{i}b\circ \mathcal{S}^{j}\overline{\mathcal{%
S}}^{s}d\mathcal{S}^{t}\right) \bullet \left( \overline{\mathcal{S}}%
_{i}a\circ \mathcal{S}_{j}\overline{\mathcal{S}}_{s}c\mathcal{S}_{t}\right)
\\
&=\left( \overline{\mathcal{S}}^{i}b\circ d\mathcal{S}^{t}\right) \bullet
\left( \overline{\mathcal{S}}_{i}a\circ c\mathcal{S}_{t}\right)  
=(\overline{\mathcal{S}}^{i}\left( b\circ d\right) \mathcal{S}^{t})\bullet(
\overline{\mathcal{S}}_{i}\left( a\circ c\right) \mathcal{S}_{t}) \\
&=\sigma _{W\circ Y,X\circ Z}(\left( a\circ c\right) \bullet \left( b\circ
d\right))  
=\sigma _{W\circ Y,X\circ Z}\zeta _{W,X,Y,Z}(\left( a\bullet b\right)
\circ \left( c\bullet d\right) )
\end{align*}%
and $\sigma _{I,I}\Delta _{I}^{\bullet }\left( b\right) =\sigma _{I,I}\left(
b_{1}\otimes b_{2}\right) =\mathcal{S}^{-1}\left( b_{2}\otimes b_{1}\right)
\mathcal{S}=\Delta^\bullet_I \left( b\right) $, by quasi-cocommutativity of $I=B$.

Then, as in \cite[Definition 6.31]{Aguiar}, we can consider a cocommutative bimonoid in this $\bullet$\,-braided category i.e.\ a bimonoid $(H,m_H^\circ,u_H^\circ,\Delta_H^\bullet,\varepsilon_H^\bullet)$ with $(H,\Delta_H^\bullet,\varepsilon_H^\bullet)$ cocommutative.
By \cref{rmk:bimonbimod}, the fact it is a bimonoid is equivalent to say that $H$ is a bialgebra together with a bialgebra map $u_H^\circ:B\to H$ where the $B$-bimodule structure of $H$ is induced by $u_H^\circ$. The cocommutativity means $\sigma_{H,H}\Delta_H^\bullet=\Delta_H^\bullet$ i.e.\ $\Ss^{-1}\Delta_H^{\bullet\mathrm{cop}}(\cdot)\Ss=\Delta_H^\bullet(\cdot)$. Since the $B$-bimodule structure of $H$ is induced by $u_H^\circ$, this means $\Rr^{-1}\Delta_H^{\bullet\mathrm{cop}}(\cdot)\Rr=\Delta_H^\bullet(\cdot)$ where we set $\Rr\coloneqq (u_H^\circ\bullet u_H^\circ)(\Ss)\in H\otimes H.$ Thus $\Rr$ satisfies the quasi-cocommutativity and, since $u_H^\circ$ is a bialgebra map, it also obeys the hexagon equations. Therefore, $(H,\Rr)$ is a quasitriangular bialgebra and $u_H^\circ$ is a morphism of quasitriangular bialgebras. In summary, the category $\cBimon(\Cc,\circ,\bullet)$ identifies with the coslice category $(B,\Ss)/\QTrBialg$. \qedhere
\end{proof}

\begin{remark}%[\rd{added:27/5/2026}] 
In principle we could also consider the quasitriangular structure $(\id\ot\tau\ot\id)(\Ss\ot\Ss^\op)$ on $B\ot B^\op$, from which we get another braiding for $_B\mm_B$, namely  $\sigma'_{X,Y}(x\ot y)=\Ss^{\op}(y \ot x)\Ss$. However, this choice does not seem to lead to a $\bullet$-braided structure in general. 
\begin{invisible}
 \begin{align*}
&\zeta _{X,W,Z,Y}\left( \sigma' _{W,X}\circ \sigma' _{Y,Z}\right) (\left(
a\bullet b\right) \circ \left( c\bullet d\right) ) 
=\zeta _{X,W,Z,Y}\left(\left(
\Ss_{i}b\mathcal{S}^{j}\bullet \Ss^{i}a%
\mathcal{S}_{j}\right) \circ \left( \Ss_{s}d\mathcal{S}%
^{t}\bullet \Ss^{s}c\mathcal{S}_{t}\right) \right) \\
&=\left( \Ss_{i}b\mathcal{S}^{j}\circ \Ss_{s}d\mathcal{S}^{t}\right) \bullet \left( \Ss^{i}a%
\mathcal{S}_{j}\circ \Ss^{s}c\mathcal{S}_{t}\right)  
=\left( \Ss_{i}b\circ \mathcal{S}^{j}\Ss_{s}d\mathcal{S}^{t}\right) \bullet \left( \Ss^{i}a\circ \mathcal{S}_{j}\Ss^{s}c\mathcal{S}_{t}\right)
\\
&=?=\left( \Ss_{i}b\circ d\mathcal{S}^{t}\right) \bullet
\left( \Ss^{i}a\circ c\mathcal{S}_{t}\right)  
=(\Ss_{i}\left( b\circ d\right) \mathcal{S}^{t})\bullet(
\Ss^{i}\left( a\circ c\right) \mathcal{S}_{t}) \\
&=\sigma' _{W\circ Y,X\circ Z}(\left( a\circ c\right) \bullet \left( b\circ
d\right))  
=\sigma' _{W\circ Y,X\circ Z}\zeta _{W,X,Y,Z}(\left( a\bullet b\right)
\circ \left( c\bullet d\right) ).
\end{align*} 
Hence the equality holds if $\Ss\Ss^\op=\id$ i.e. for a triangular bialgebra.
\end{invisible}
\end{remark}

The following result is an immediate consequence of \cref{pro:Bohm} and \cref{pro:bimoncoc}.

\begin{corollary}
\label{cor:quasiascocommutativeHopf}
In the setting of \cref{es:bimod}, if $(B,\Ss)$ is a quasitriangular Hopf algebra, then the category  $\cHopfmon(\Cc,\circ,\bullet)$ of cocommutative Hopf monoids therein is isomorphic to the coslice category of quasitriangular Hopf algebras under $(B,\Ss)$, i.e.\ $(B,\Ss)/\QTrHopf$.
\end{corollary}

We now turn to protomodularity. We denote by $((B,\Ss)/\QTrHopf)^{\mathrm{i}}$ the full subcategory of $(B,\Ss)/\QTrHopf$ whose objects are the pairs  $((H,\Rr),u_H^\circ)$ with  $u_H^\circ:(B,\Ss)\to (H,\Rr)$ injective.

\begin{corollary}%[\rd{added:2026/03/18}]
\label{cor:cosliceQTrHprot}
In the setting of \cref{es:bimod}, if $(B,\Ss)$ is a quasitriangular Hopf algebra, then the category  $\cHopfmon^\mathrm{m}(\Cc,\circ,\bullet)$ is isomorphic to $((B,\Ss)/\QTrHopf)^{\mathrm{i}}$ and is protomodular.
\end{corollary}

\begin{proof}
By the foregoing, the duoidal category $(\Cc,\circ,I,\bullet,J)$  is $\bullet$-braided (see \cref{pro:bimoncoc}) and has a reversion. Clearly $\Cc={}_B\mathfrak{M}_B$ has binary intersections and $\bullet=\otimes_\Bbbk$ preserves them as $X\bullet(-)$ and $(-)\bullet X$ are left exact additive endofunctors of $\Cc$, which is abelian, and hence they preserve finite limits, for any $X$.
\begin{invisible}
See e.g. Proposition 8.6 in Stenstroems book. We also keep here the direct proof.
Note that, although the tensor product does not preserve limits in general,
still it preserves inverse images (which are pullbacks) and so also intersections. To see this,
consider a linear map  $f:X\rightarrow Y$ and  a subspace $Y^{\prime }$ of $%
Y.$ Then $\left( f\otimes \mathrm{Id}_{Z}\right) ^{-1}\left( Y^{\prime
}\otimes Z\right) =$ $f^{-1}\left( Y^{\prime }\right) \otimes Z.$ Indeed,
call $p:Y\rightarrow Y/Y^{\prime }$ the canonical projection. Then $w\in
\left( f\otimes \mathrm{Id}_{Z}\right) ^{-1}\left( Y^{\prime }\otimes
Z\right) $ means $\left( f\otimes \mathrm{Id}_{Z}\right) \left( w\right) \in
Y^{\prime }\otimes Z$ and hence $\left( pf\otimes \mathrm{Id}_{Z}\right)
\left( w\right) =0$ so that $w\in \ker \left( pf\otimes \mathrm{Id}%
_{Z}\right) =\ker \left( pf\right) \otimes
Z=f^{-1}\left( Y^{\prime }\right) \otimes Z.$ Thus $\left( f\otimes \mathrm{%
Id}_{Z}\right) ^{-1}\left( Y^{\prime }\otimes Z\right) \subseteq
f^{-1}\left( Y^{\prime }\right) \otimes Z.$ The other inclusion is always
true.
\end{invisible}
By \cref{thm:protomon}, we get that the category $\cHopfmon^{\mathrm{m}}(\Cc,\circ, \bullet )$ is protomodular. By \cref{cor:quasiascocommutativeHopf}, the category $\cHopfmon(\Cc,\circ,\bullet)$ is isomorphic to the coslice category $(B,\Ss)/\QTrHopf$. Through this isomorphism, $\cHopfmon^{\mathrm{m}}(\Cc,\circ, \bullet )$ identifies with the full subcategory of $(B,\Ss)/\QTrHopf$ whose objects are the pairs  $((H,\Rr),u_H^\circ)$ with  $u_H^\circ:(B,\Ss)\to (H,\Rr)$ a monomorphism in $\Cc$, which means that $u_H^\circ$ is injective. Thus, $\cHopfmon^{\mathrm{m}}(\Cc,\circ, \bullet )$ is isomorphic to $((B,\Ss)/\QTrHopf)^{\mathrm{i}}$.
\end{proof}

\subsection{Triangular Hopf algebras}

Recall that a quasitriangular bialgebra $(B,\Ss)$ is called triangular if it further obeys the condition $\Ss^{-1}=\Ss^{\mathrm{op}}$.

\begin{remark}
\label{rmk:nntrvntgrps}
%[\as{Added on 26/11/2025}]
We recall that the $\bullet$-braiding of $(\Cc,\circ,\bullet)$ as in \cref{cor:quasiascocommutativeHopf} is defined, for all $X,Y\in\Cc$, by $\sigma_{X,Y}(x\bullet y):=\Ss^{-1}(y \bullet x)\Ss$, whose inverse is given by $\sigma_{X,Y}^{-1}(y\bullet x):=\Ss^{\mathrm{op}}(x\bullet y)(\Ss^{-1})^{\mathrm{op}}$. If $(B,\Ss)$ is a triangular Hopf algebra, then $(\Cc,\bullet,J,\sigma)$ becomes  symmetric monoidal and so the category $(\cComon(\Cc^\bullet),\bullet,J)$ becomes cartesian monoidal, see e.g.\ \cite[Corollary 2.24]{Ulrich-Myriam}. Therefore, one can consider the $\bullet$-braided duoidal category $(\cComon(\Cc^\bullet),\circ, I,\bullet,J)$, as an instance of \cref{exa:duoidfinprod}, and, in case this has a reversion, the category $\Grp(\cComon(\Cc^\bullet),\circ,I)=\Hopfmon(\cComon(\Cc^\bullet),\circ,\bullet)$. It is worth noting that here the tensor product $\circ$ is not the binary product $\bullet$ in $\cComon(\Cc^\bullet)$, unlike for ordinary internal groups.
% \rd{[Per parlare di monoidi di Hopf occorre che $(\cComon(\Cc^{\bullet}),\circ,\bullet)$ abbia una reversion. Presumibilmente è la stessa della categoria originale che si solleva ai comonoidi, assieme alle 4 trasformazioni naturali di definizione; ma mi rifiuto di verificarlo :)]}
\end{remark}

We show that the $\bullet$-braided duoidal category $(\cComon(\Cc^\bullet),\circ, I,\bullet,J)$ inherits the reversion $\left( \left( -\right) ^{\revsn },\phi^\circ ,\phi^{\bullet},\gamma ,\delta,\widetilde{\gamma},\widetilde{\delta}
\right) $ given in \cref{subsec:reversionbim}. Since $((-)^{\revsn},\phi^{\bullet}):\left( \mathcal{C},\bullet^\rev \right) \rightarrow \left( \mathcal{C},\bullet \right)$ is a strong monoidal functor (Remark \ref{rmk:duoidaleqbim}), we obtain a functor $(-)^{\revsn}:\Comon(\Cc^{\bullet})\to\Comon(\Cc^{\bullet})$. More precisely, given $(C,\Delta,\varepsilon)$ in $\cComon(\Cc^{\bullet})$, one defines $\Delta^{\bullet}_{C^{\revsn}}:=(\phi^{\bullet}_{C,C})^{-1}(\Delta_{C}^{\bullet})^{\revsn}$ and $\varepsilon^{\bullet}_{C^{\revsn}}:=(\phi^{\bullet}_{0})^{-1}(\varepsilon^{\bullet}_{C})^{\revsn}$ that make $(C^{\revsn},\Delta^{\bullet}_{C^{\revsn}},\varepsilon^{\bullet}_{C^{\revsn}})$ an object in $\Comon(\Cc^{\bullet})$. Moreover, we have
\[
\begin{split}
(\sigma_{Y,X})^{\revsn}\phi^{\bullet}_{X,Y}(x^{\revsn}\bullet y^{\revsn})&=\sigma_{Y,X}(y\bullet x)^{\revsn}=(\Ss^{-1}(x\bullet y)\Ss)^{\revsn}=\phi^{\bullet}_{Y,X}((\overline{\Ss}_{i}y\Ss_{k})^{\revsn}\bullet(\overline{\Ss}^{i}x\Ss^{k})^{\revsn})\\&=\phi^{\bullet}_{Y,X}((S_{B}(\overline{\Ss}_{i})yS_{B}(\Ss_{k}))^{\revsn}\bullet(S_{B}(\overline{\Ss}^{i})xS_{B}(\Ss^{k}))^{\revsn})\\&=\phi^{\bullet}_{Y,X}(\Ss_{k}y^{\revsn}\overline{\Ss}_{i}\bullet\Ss^{k}x^{\revsn}\overline{\Ss}^{i})=\phi^{\bullet}_{Y,X}(\Ss^{\mathrm{op}}(y^{\revsn}\bullet x^{\revsn})(\Ss^{-1})^{\mathrm{op}})\\&\overset{(\dag)}{=}\phi^{\bullet}_{Y,X}(\Ss^{-1}(y^{\revsn}\bullet x^{\revsn})\Ss)=\phi^{\bullet}_{Y,X}\sigma_{X^{\revsn},Y^{\revsn}}(x^{\revsn}\bullet y^{\revsn}),
\end{split}
\]
where $(\dag)$ follows since $\Ss^{-1}=\Ss^{\mathrm{op}}$. Therefore, the monoidal functor $((-)^{\revsn},\phi^{\bullet}):\left( \mathcal{C},\bullet^\rev \right) \rightarrow \left( \mathcal{C},\bullet \right)$ is braided, and we obtain a strong monoidal functor $((-)^{\revsn},\phi^{\bullet}):(\cComon(\Cc^{\bullet}),\bullet^\rev)\to(\cComon(\Cc^{\bullet}),\bullet)$. Moreover, the morphism $\phi^{\circ}_{0}$ is of $\bullet$-comonoids
\begin{invisible}
since
\[
\Delta^{\bullet}_{B^{\revsn}}\phi^{\circ}_{0}(b)=(\phi^{\bullet}_{B,B})^{-1}(\Delta_{B}^{\bullet})^{\revsn}(S_{B}(b)^{\revsn})=(\phi^{\bullet}_{B,B})^{-1}(\Delta_{B}^{\bullet}S_{B}(b)^{\revsn})=S_{B}(b_{1})^{\revsn}\bullet S_{B}(b_{2})^{\revsn}=\phi^{\circ}_{0}(b_{1})\bullet\phi^{\circ}_{0}(b_{2})
\]
and $\varepsilon_{B^{\revsn}}^{\bullet}\phi^{\circ}_{0}(b)=(\phi^{\bullet}_{0})^{-1}(\varepsilon^{\bullet}_{B})^{\revsn}(S_{B}(b)^{\revsn})=\varepsilon_{B}S_{B}(b)=\varepsilon_{B}(b)$
\end{invisible}
as well as 
$\phi^{\circ}_{X,Y}$. 
\begin{invisible}
Indeed, we have
\[
\begin{split}
\Delta^{\bullet}_{(Y\circ X)^{\revsn}}\phi^{\circ}_{X,Y}(x^{\revsn}\circ y^{\revsn})&=(\phi^{\bullet}_{Y\circ X,Y\circ X})^{-1}(\Delta_{Y\circ X}^{\bullet})^{\revsn}((y\circ x)^{\revsn})=(y_{2}\circ x_{2})^{\revsn}\bullet(y_{1}\circ x_{1})^{\revsn}\\&=(\phi^{\circ}_{X,Y}\bullet\phi^{\circ}_{X,Y})((x_{2}^{\revsn}\circ y_{2}^{\revsn})\bullet(x_{1}^{\revsn}\circ y_{1}^{\revsn}))=(\phi^{\circ}_{X,Y}\bullet\phi^{\circ}_{X,Y})\Delta^{\bullet}_{X^{\revsn}\circ Y^{\revsn}}(x^{\revsn}\circ y^{\revsn})
\end{split}
\]
and $\varepsilon^{\bullet}_{(Y\circ X)^{\revsn}}\phi^{\circ}_{X,Y}(x^{\revsn}\circ y^{\revsn})=(\phi^{\bullet}_{0})^{-1}(\varepsilon^{\bullet}_{Y\circ X})^{\revsn}((y\circ x)^{\revsn})=\varepsilon^{\bullet}_{Y\circ X}(y\circ x)=\varepsilon^{\bullet}_{X^{\revsn}\circ Y^{\revsn}}(x^{\revsn}\circ y^{\revsn})$.
\end{invisible}
Thus, we obtain a strong monoidal functor $((-)^{\revsn},\phi^{\circ}):(\cComon(\Cc^{\bullet}),\circ^\rev)\to(\cComon(\Cc^{\bullet}),\circ)$, hence a duoidal equivalence $\left( \left( -\right)^{\revsn },\phi^\circ,\phi^{\bullet} \right):(\cComon(\Cc^{\bullet}),\circ^\rev,\bullet^{\rev})\to(\cComon(\Cc^{\bullet}),\circ,\bullet)$ that lifts the one provided in \cref{rmk:duoidaleqbim}. In order to prove that the reversion of $(\Cc,\circ,\bullet)$ is lifted to $(\cComon(\Cc^{\bullet}),\circ,\bullet)$, it remains to prove that the natural morphisms $\gamma_{X,Y,Z,T},\delta_{X,Y,Z,T},\widetilde{\gamma}_{X,Y,Z,T},\widetilde{\delta}_{X,Y,Z,T}$ given in \cref{subsec:reversionbim} are morphisms of $\bullet$-comonoids, for all $X,Y,Z,T\in\cComon(\Cc^{\bullet})$. This is straightforward.

\medskip We obtain the following result:

\begin{corollary}\label{coro:trHopfasgroups}
%[\as{Added on 26/11/2025}]
    In the setting of \cref{es:bimod}, if $(B,\Ss)$ is a triangular Hopf algebra, then $(B,\Ss)/\TrHopf$ is isomorphic to the category $\Grp(\cComon(\Cc^\bullet),\circ,I)$.
\end{corollary}

\begin{proof}
By \cref{rmk:internalgroupduoidal} and the dual of \cite[Proposition 6.37]{Aguiar}, we have that
 \[
\cBimon(\cComon(\Cc^{\bullet}),\circ,\bullet)\cong\mathsf{Mon}(\cComon(\Cc^{\bullet}),\circ)\cong\cBimon(\Cc,\circ,\bullet)
\]
and then \[\Grp(\cComon(\Cc^\bullet),\circ,I)=\cHopfmon(\cComon(\Cc^{\bullet}),\circ,\bullet)\cong\cHopfmon(\Cc,\circ,\bullet)\cong (B,\Ss)/\QTrHopf,\] where the last isomorphism comes from \cref{cor:quasiascocommutativeHopf}. It remains to observe that $(B,\Ss)/\QTrHopf$ turns out to be $(B,\Ss)/\TrHopf$ for $(B,\Ss)$ triangular. Indeed, given a quasitriangular Hopf algebra $(H,\Rr)$ with a morphism of quasitriangular Hopf algebras $u_{H}^{\circ}:(B,\Ss)\to (H,\Rr)$, then the codomain becomes triangular as so is the domain and $\Rr=(u_{H}^{\circ}\bullet u_{H}^{\circ})(\Ss)$.
\end{proof}

\subsection{Minimal quasitriangular Hopf algebras}
Given a quasitriangular Hopf algebra $(H, \mathcal{S})$, one can define the \textit{minimal quasitriangular Hopf algebra} $H_{\Ss}$, which is a finite-dimensional Hopf
subalgebra of $H$, see \cite[Proposition 12.2.13]{Radford-book}. Minimality is made explicit in \cite[Exercise 12.2.9]{Radford-book}, which states that, if $K$ is a Hopf subalgebra of $H$ with $\Ss\in K\otimes K$, then $H_{\Ss}\subseteq K$.

\begin{remark}%[\rd{added:2026/03/18}]
\label{rmk:coslicemin}
An injective morphism $f:(B,\Ss)\to (A,\Rr)$  of quasitriangular Hopf algebras induces an isomorphism $(B_{\Ss},\Ss)\to (A_{\Rr},\Rr)$, cf.\ \cite[Lemma 1]{Radford-min}. Thus, if we assume $(B,\Ss)$ is minimal, then $f$ induces an isomorphism $(B,\Ss)\to (A_{\Rr},\Rr)$; as a consequence,  $((B,\Ss)/\QTrHopf)^{\mathrm{i}}$ is isomorphic to the category whose objects are pairs $((A,\Rr),\alpha)$ where $(A,\Rr)\in \QTrHopf$ and $\alpha:(B,\Ss)\to (A_{\Rr},\Rr)$ is an isomorphism in $\QTrHopf$. A morphism $\phi:((A,\Rr),\alpha)\to ((A',\Rr'),\alpha')$ is a morphism $\phi:(A,\Rr)\to (A',\Rr')$ such that $\phi i\alpha =i'\alpha'$, where $i:A_{\Rr}\to H$ and $i':A'_{\Rr'}\to H'$ are the canonical injections. The last equality rewrites as $\phi i=i'(\alpha' \alpha^{-1})$ which just means that $\phi$ induces the isomorphism $\alpha'\alpha^{-1}:A_{\Rr}\to A'_{\Rr'}$.
\end{remark}

Given a quasitriangular Hopf subalgebra $(H,\Rr)$ of $(A,\Rr)$, by a projection $\pi:A\to H$ we mean that $\pi$ is a Hopf algebra map with $\pi\sigma=\id$, where $\sigma:H\to A$ is the canonical injection.

\begin{remark}
\label{rmk:pointquasi}
Note that a $\pi$ as above is necessarily a morphism of quasitriangular Hopf algebras as $\Rr\in H\otimes_\Bbbk H$ and $\pi_{\mid H}=\id_H$ imply  $(\pi\otimes_\Bbbk\pi)(\Rr)=\Rr.$ Thus, we get a point $\xymatrix{A\ar@<.5ex>[r]^-\pi&\ar@<.5ex>[l]^-\sigma H}$ in $\QTrHopf$. We know that $A$ can be described as the Radford-Majid bosonization $R\# H=R\otimes_\Bbbk H$, where $R:=A^{\mathrm{co}H}=\{r\in A\mid r_1\ot_\Bbbk\pi(r_2)=r\ot_\Bbbk1_H\}$ is the space of $H$-coinvariant elements in $A$. In the following result, we will see that our results allow for another description of $A$.
\end{remark}

\begin{theorem}%[\rd{added:2026/03/24}]
\label{thm:crospt}
 Let $(A,\Rr)$ be a quasitriangular Hopf algebra with a projection $\pi:A\to H$ onto a quasitriangular Hopf subalgebra $(H,\Rr)$. Set $B:=H_{\Rr}$ and $R:=A^{\mathrm{co}H}$. Then $K:=RB$ is a quasitriangular Hopf subalgebra of $(A,\Rr)$ and $\varphi:K\otimes_B H\to A,w\otimes_B h\mapsto wh$, is an isomorphism of quasitriangular Hopf algebras, where the domain is structured as follows
 \begin{gather*}\label{structures}
 (w\otimes_B h)(w'\otimes_B h')=w (h_1\triangleright w_1')\otimes_B (h_2\triangleleft w_2') h',\qquad 1_{K\otimes_B H}=1_K\otimes_B1_H,\\
 \Delta(w\otimes_B h)=(w_1\otimes_B h_1)\otimes_\Bbbk (w_2\otimes_B h_2),\qquad \varepsilon(w\otimes_B h)=\varepsilon_K(w)\varepsilon_H(h),\\
 S_{K\otimes_B H}(w\otimes_B h)=(S_H(h_2)\triangleright S_K(w_2))\otimes_B (S_H(h_1)\triangleleft S_K(w_1)),\\ \Rr_{K\otimes_B H}=(1_K\otimes_B \Rr^i)\otimes_\Bbbk(1_K\otimes_B\Rr_i),
 \end{gather*}
 for every $w,w'\in K,h,h'\in H,$ where $h\triangleright w:=h_1w_1S_H\pi(w_2)S_H\pi(h_2)\in R$ and
 $h\triangleleft w:=h\pi(w)\in H$.
\end{theorem}

\begin{proof}
By assumption $\pi:A\to H$ is a Hopf algebra map with $\pi\sigma=\id$, where $\sigma:H\to A$ is the canonical injection. As observed in \cref{rmk:pointquasi}, we have a point $\xymatrix{A\ar@<.5ex>[r]^-\pi&\ar@<.5ex>[l]^-\sigma H}$ in $\QTrHopf$. Set $B:=H_{\Rr}$, let $u_H^\circ:(B,\Rr)\to(H,\Rr)$ be the canonical inclusion and $u_A^\circ:=\sigma u_H^\circ:(B,\Rr)\to(A,\Rr)$. Then $\pi u_A^\circ=u_H^\circ$ so that our point is in fact a point in $(B,\Rr)/\QTrHopf$. Since $u_H^\circ$ and $u_A^\circ$ are injections, it is even a point in $((B,\Rr)/\QTrHopf)^{\mathrm{i}}$. Since we proved in \cref{cor:cosliceQTrHprot} that the latter category is protomodular and isomorphic to $\cHopfmon^\mathrm{m}(\Cc,\circ,\bullet)$, where $\Cc={}_B\mathfrak{M}_B,\circ=\ot_B,I=B,\bullet=\ot_\Bbbk, J=\Bbbk$, we can apply \cref{thm:main} to this specific point. Explicitly, let $(K,k:K\to A)$ be the kernel of $\pi$ in $\cHopfmon(\Cc,\circ,\bullet)$ defined by $K:=\{a\in A\mid a_1\bullet\pi(a_2)\in A\bullet I\}$ and where $k$ is the canonical inclusion.
Set $\tau(a):=a_1\sigma S_H\pi(a_2)\in R$ for every $a\in A$. Then, for any $w\in K$, we have $w=\tau(w_1)\sigma\pi(w_2)\in RI$ so that $K\subseteq RI$ and the other inclusion is also true whence $K=RI.$
Moreover, $\varphi=m_A^\circ(k\circ\sigma):K\circ H\to A,w\circ h\mapsto wh$ is an isomorphism. We can use \cref{coro:main} to get a more precise description. Indeed, $K\circ H$ becomes the quasitriangular Hopf algebra $K\circ_\psi H$. In order to write it, consider
$(\pi^r)^{-1}:(A\circ J)\bullet H\to A\circ H$ given by
\[(\pi^r)^{-1}((a\circ 1_\Bbbk)\bullet h)=(A\circ l_H^\circ)\overline{\varsigma}_{A,J,H}^H((a\circ 1_\Bbbk)\bullet h)
=a_1\circ S_H\pi(a_2)h.\] Then we have
$\Upsilon:A\to K\circ J$ as
\begin{align*}
\Upsilon(a)&=(\theta\circ J)(\id_A*\sigma^\revsn\Sigma_H\pi)(a)
=(\theta\circ J)((m_A^\circ\bullet I)\circ J)\varphi_{A,A}(\id_A\bullet \sigma^\revsn\Sigma_H\pi)\Delta^\bullet_A(a)\\
&=(\theta\circ J)((m_A^\circ\bullet I)\circ J)\varphi_{A,A}(a_1\bullet \sigma^\revsn\Sigma_H\pi(a_2))
=(\theta\circ J)((m_A^\circ\bullet I)\circ J)\varphi_{A,A}(a_1\bullet (\sigma S_H\pi(a_2))^\revsn)\\
&=(\theta\circ J)((m_A^\circ\bullet I)\circ J)(((a_1\circ \sigma S_H\pi(a_2))\bullet 1_B)\circ 1_{\Bbbk})=\theta (a_1\sigma S_H\pi(a_2)\bullet 1_B)\circ 1_\Bbbk\\
&=\theta (\tau(a)\bullet 1_B)\circ 1_\Bbbk
=\xi (\tau(a)\bullet 1_B)\circ 1_\Bbbk
=m_A^\circ(A\circ\sigma)(\pi^r)^{-1}((\tau(a)\circ 1_\Bbbk)\bullet 1_B)\circ 1_\Bbbk\\
&=m_A^\circ(A\circ\sigma)(\pi^r)^{-1}((\tau(a)\circ 1_\Bbbk)\bullet 1_B)\circ 1_\Bbbk
=m_A^\circ(A\circ\sigma)(\tau(a)_1\circ S_H\pi(\tau(a)_2))\circ 1_\Bbbk\\
&=m_A^\circ(A\circ\sigma)(\tau(a)\circ S_H(1_H))\circ 1_\Bbbk=\tau(a)\circ 1_\Bbbk.
\end{align*}
Therefore, we get $\chi:H\circ K\to K\circ J$, defined as $\chi(h\circ w)=\Upsilon(hw)=\tau(hw)\circ 1_\Bbbk=(h\triangleright w)\circ 1_\Bbbk$, where we set $h\triangleright w:=\tau(hw)\in R$. We also have $h\triangleleft w:=h\varepsilon_K^\circ(w)=h\pi(w),$ where $\varepsilon_K^\circ(w):=\pi(w)$. Since $((u^\circ_H\varepsilon^\circ_K)^r)^{-1}:(K\circ J)\bullet H\to K\circ H$ is given by 
\[
\begin{split}
((u^\circ_H\varepsilon^\circ_K)^r)^{-1}((w\circ 1_\Bbbk)\bullet h)&=(K\circ l_H^\bullet)\overline{\varsigma}_{K,J,H}^H((w\circ 1_\Bbbk)\bullet h)=(K\circ l_H^\bullet)(w_1\circ(1_\Bbbk\bullet S_Hu^\circ_H\varepsilon_K^\circ(w_2)h))\\&=w_1\circ S_H\pi(w_2)h
\end{split}
\] for every $w\in K,h\in H$, we get
\begin{align*}
\psi(h\circ w)&=\left( \left( u_{H}^{\circ }\varepsilon _{K}^{\circ }\right)
^{r}\right) ^{-1}\left( \chi \bullet \triangleleft \right) \Delta
_{H\circ K}^{\bullet }(h\circ w)
=\left( \left( u_{H}^{\circ }\varepsilon _{K}^{\circ }\right)
^{r}\right) ^{-1}( \chi (h_1\circ w_1)\bullet (h_2\triangleleft w_2))\\
&=\left( \left( u_{H}^{\circ }\varepsilon _{K}^{\circ }\right)
^{r}\right) ^{-1}( ((h_1\triangleright w_1)\circ 1_\Bbbk)\bullet (h_2\triangleleft w_2))
=(h_1\triangleright w_1)_1\circ S_H\pi((h_1\triangleright w_1)_2) (h_2\triangleleft w_2)\\
&=(h_1\triangleright w_1)\circ (h_2\triangleleft w_2)
\end{align*}
where in the last step we used that $h_1\triangleright w_1\in R.$
\begin{invisible}
Moreover,
\begin{align*}
 \varphi ^{-1}(a)&=\left( \left( u_{H}^{\circ }\varepsilon _{K}^{\circ }\right)
^{r}\right) ^{-1}\left( \Upsilon \bullet \pi \right) \Delta _{A}^{\bullet }(a)
=\left( \left( u_{H}^{\circ }\varepsilon _{K}^{\circ }\right)
^{r}\right) ^{-1}\left( \Upsilon(a_1) \bullet \pi(a_2) \right)\\
&=\left( \left( u_{H}^{\circ }\varepsilon _{K}^{\circ }\right)
^{r}\right) ^{-1}\left( (\tau(a_1)\circ 1_\Bbbk) \bullet \pi(a_2) \right)
=\tau(a_1)_1\circ S_H\pi(\tau(a_1)_2)\pi(a_2)
=\tau(a_1)\circ \pi(a_2),
\end{align*} where in the last step we used that $\tau(a_1)\in R$.
\end{invisible}
In view of \cref{thm:psitriang}, we have
\begin{align*}
\Sigma_{K\circ H}(w\circ h)&=
\psi^\revsn\phi^\circ_{K,H}(\Sigma_K\circ \Sigma_H)(w\circ h)
=\psi^\revsn\phi^\circ_{K,H}(S_K(w)^\revsn\circ S_H(h)^\revsn)\\
&=\psi^\revsn((S_H(h)\circ S_K(w))^\revsn)=\big(\psi (S_H(h)\circ S_K(w))\big)^\revsn
\end{align*}
so that the antipode is $S_{K\circ H}(w\circ h)=\psi (S_H(h)\circ S_K(w))=(S_H(h_2)\triangleright S_K(w_2))\circ (S_H(h_1)\triangleleft S_K(w_1)).$
The quasitriangular structure of $K\circ H$ is given by  $(u^\circ_{K\circ H}\bullet u^\circ_{K\circ H})(\Rr)=(1_K\circ \Rr^i)\bullet(1_K\circ\Rr_i).$
\end{proof}

Let $K\bowtie_B H$ be the quasitriangular Hopf algebra structure on $K\otimes_B H$ obtained in \cref{thm:crospt}. Some of the rules that $\triangleright:H\circ K\to K$ and $\triangleleft:H\circ K\to H$ obey can be deduced from \cref{thm:psitriang}. Indeed, \eqref{eq:chi1}, \eqref{eq:chi2}, \eqref{eq:triangl1} and
\eqref{eq:triangl2} amount respectively to, for $h,h'\in H,w,w'\in K$,
\begin{gather*}
(hh')\triangleright w
=h\triangleright (h'\triangleright w),\qquad (1_H\triangleright w)\circ 1_\Bbbk=w\circ 1_\Bbbk,\\
h\triangleright(ww')=(h_1\triangleright w_1)((h_2\triangleleft w_2)\triangleright w'),\qquad h\triangleright 1_K=\varepsilon_H(h)1_K,\\
(hh')\triangleleft w=(h\triangleleft(h'_1\triangleright w_1))(h'_2\triangleleft w_2),\qquad 1_H\triangleleft w=1_H\varepsilon_K(w),\\
h\triangleleft(ww')=(h\triangleleft w)\triangleleft w',\qquad h\triangleleft 1_K= h.
\end{gather*}
For instance, the first equality in \eqref{eq:chi2} says that $h\triangleright(ww')\circ 1_\Bbbk=(h_1\triangleright w_1)((h_2\triangleleft w_2)\triangleright w')\circ1_\Bbbk$ but both these elements are of the form $r\circ 1_\Bbbk$ for some $r\in R$ and $\varphi((u^\circ_H\varepsilon^\circ_K)^r)^{-1}((r\circ 1_\Bbbk)\bullet1_H)=r_1S_H\pi(r_2)1_H=r$ so that we get the desired equality.
\begin{remark}
When $A$ is cocommutative and $\Rr=1\ot1,$ we have $B=\Bbbk$, and $K\bowtie_B H$ reduces to the double cross product $K\bowtie H$ of \cite[Theorem 7.2.2]{Majid-book}. Therefore, our construction can be viewed as an extension of the double cross product of cocommutative Hopf algebras to the quasitriangular setting and the datum $(K,H,\triangleright,\triangleleft)$ as a generalized matched pair.
Although we noted that $A$ can also be described as the Radford-Majid bosonisation $R\#H=R\otimes_\Bbbk H$, the advantage of our approach is that both factors $K,H$ are Hopf subalgebras, rather than having only $H$ as such and $R$ a Hopf monoid in the category of Yetter–Drinfeld modules $^H_H\mathcal{YD}$; the drawback, on the other hand, is that the tensor product is taken over $B$ rather than $\Bbbk$.
\end{remark}
\begin{example}
%[\as{[Added on 27/03/2026]}]
Let $\Bbbk$ be such that $\mathrm{char}(\Bbbk)\not=2$. We consider the $\Bbbk$-algebra $A$ generated by $g,h,x$ modulo the relations 
\[
g^{2}=h^{2}=1,\quad x^{2}=0,\quad gh=hg,\quad gx=-xg,\quad hx=-xh.
\]
This becomes a Hopf algebra with comultiplication and counit determined by
\[
\Delta(g)=g\otimes g,\ \Delta(h)=h\otimes h,\ \Delta(x)=1\otimes x+x\otimes g,\ \varepsilon(g)=\varepsilon(h)=1_{\Bbbk},\ \varepsilon(x)=0,
\]
and antipode given by $S(g)=g$, $S(h)=h$ and $S(x)=-xg$. This Hopf algebra is sometimes denoted by $A_{C_{2}\times C_{2}}$ since its coradical is given by $H:=\Bbbk\langle g,h\rangle$ and $\langle g,h\rangle\cong C_{2}\times C_{2}$, where $C_{2}$ is the cyclic group of order two. As stated in \cite[p.\ 559]{CDR}, $A$ is
isomorphic to $H_{4}\ot\Bbbk C_{2}$, where $H_{4}$ is the Sweedler Hopf algebra.  By \cite[Theorem 1.4]{Wakui}, $A$ has two 1-parameter families of universal $\Rr$-matrices, among which there is the following one:
\[
\Rr_{\lambda}:=\frac{1}{2}(1\ot1+g\ot1+1\ot g-g\ot g)+\frac{\lambda}{2}(x\ot gx+x\ot x+gx\ot gx-gx\ot x),
\]
for $\lambda\in\Bbbk$. We observe that the latter is an exhaustive 1-parameter family of triangular structures for the Sweedler Hopf algebra $H_{4}$ with group-like element $g$, see e.g.\ \cite[Exercise 2.1.7]{Majid-book} (where the opposite coproduct for $H_{4}$ is considered).

Now regard $(H,\Rr_{0})$ as a quasitriangular Hopf subalgebra of $(A,\Rr_{0})$. There is a surjective Hopf algebra map $\pi:A\to H$, $x^{m}g^{a}h^{b}\mapsto\delta_{m,0}g^{a}h^{b}$. %, which provides a surjective morphism of quasitriangular Hopf algebras $\pi:(A,\Rr_{0})\to(H,\Rr_{0})$. 
It is immediate to see that $R=\Bbbk\langle xg\rangle$ while $B:=H_{\Rr_{0}}=\Bbbk \langle g\rangle$. Hence, $K:=RB=H_{4}$. By
\cref{thm:crospt}, $H_{4}$ is a quasitriangular Hopf subalgebra of $(A,\Rr_{0})$ and $\varphi:H_{4}\ot_{\Bbbk\langle g\rangle}\Bbbk\langle g,h\rangle\to A$, $\omega\ot_{\Bbbk\langle g\rangle}h\mapsto \omega h$ is an isomorphism of quasitriangular Hopf algebras, where $H_{4}\ot_{\Bbbk\langle g\rangle}\Bbbk\langle g,h\rangle$ is structured as in \cref{structures}. Let us describe the morphisms $\triangleright$ and $\triangleleft$. For $a\in H,\omega\in H_{4}$, we have $a\triangleright\omega=a_{1}\omega_{1}S\pi(\omega_{2})S\pi(a_{2})$ and $a\triangleleft\omega=a\pi(\omega)$ which, on the basis, becomes  $g^ah^b\triangleright x^mg^c=(-1)^{m(a+b)}x^mg^m$ and $g^ah^b\triangleleft x^mg^c=\delta_{m,0}g^{a+c}h^b$ i.e. \[
\begin{tabular}{c|cccc}
    $\triangleright$& $1$ & $g$ & $x$ & $xg$\\ \hline
    $1$ & $1$ & $1$ & $xg$ & $xg$ \\ 
    $g$ & $1$ & $1$ & $-xg$ & $-xg$ \\ 
    $h$ & $1$ & $1$ & $-xg$ & $-xg$ \\ 
    $gh$ & $1$ & $1$ & $xg$ & $xg$ \\ 
    \end{tabular} 
    \qquad
    \begin{tabular}{c|cccc}
    $\triangleleft$& $1$ & $g$ & $x$ & $xg$\\ \hline
    $1$ & $1$ & $g$ & $0$ & $0$ \\ 
    $g$ & $g$ & $1$ & $0$ & $0$ \\ 
    $h$ & $h$ & $hg$ & $0$ & $0$ \\ 
    $gh$ & $gh$ & $h$ & $0$ & $0$ \\ 
    \end{tabular}
\]
\end{example}

%\subsection{Essentially fibre protomodularity}
%\rd{[\href{https://ncatlab.org/nlab/show/essential+fiber}{nlab}]}

% By analogy with  \cite[Definition 7.3.10]{BorII94}, we introduce the following terminology.

% \begin{definition}
% We say that a functor $U:\mathcal{A}\to \mathcal{B}$ is \emph{essentially fibre protomodular} if, for every object $B\in \mathcal{B}$, the essential fibre $\mathcal{A}_B$ of $U$ at $B$ is a protomodular category.
% \end{definition}

%\rd{[Modified: 2026/03/26]} 
A morphism of quasitriangular Hopf algebras needs not to preserve the minimal quasitriangular Hopf subalgebras. For this reason we let $\mathcal{E}$ be the class of morphisms $f:(A,\Rr)\to (H,\Ss)$ with $f(A_\Rr)\subseteq H_{\Ss}$.
Denote by $\QTrHopf^\mathcal{E}
$ the subcategory of $\QTrHopf$ whose morphisms are the morphisms in the class $\mathcal{E}$.
We also denote by $\MQTrHopf^\mathcal{E}$ the full subcategory of $\QTrHopf^\mathcal{E}$ consisting of minimal quasitriangular Hopf algebras.

The inclusion $F:\MQTrHopf^\mathcal{E}\to \QTrHopf^\mathcal{E}$ has a right adjoint
$U:\QTrHopf^\mathcal{E}\to \MQTrHopf^\mathcal{E}$ that assigns to a quasitriangular Hopf algebra $(H,\Ss)$ its minimal quasitriangular Hopf algebra $(H_{\Ss},\Ss)$ and to a morphism $f:(A,\Rr)\to (H,\Ss)$ in $\mathcal{E}$ the induced morphism $U(f):(A_{\Rr},\Rr)\to (H_\Ss,\Ss)$ between the minimal quasitriangular Hopf algebras, which exists by definition of $\mathcal{E}$. \\
Since $UF=\id$, we can define the unit $\eta:\id\to UF$ to be the identity while the counit $\epsilon:FU\to \id$ is given on a component $(H,\Ss)$ by the canonical inclusion $(H_\Ss,\Ss)\to (H,\Ss)$.
\begin{invisible}
If $f:(A,\Rr)\to (H,\Ss)$ is a morphism in $\QTrHopf^\mathcal{E}$, then by definition it induces $U(f):(A_\Rr,\Rr)\to (H_\Ss,\Ss)$ such that $f\epsilon_{(A,\Rr)}=\epsilon_{(H,\Ss)}FU(f)$ so that $\epsilon$ is natural. To conclude it is an adjunction, we have to check the triangle identities, which are equivalent to $U\epsilon=\id_U$ and $\epsilon F=\id_F$. Now $U\epsilon_{(A,\Rr)}=\id_{(A_\Rr,\Rr)}=\id_{U(A,\Rr)}$. Moreover $\epsilon _{F(B,\Ss)}=\epsilon _{(B,\Ss)}=\id_{(B,\Ss)}=\id _{F(B,\Ss)}$.
\end{invisible}
Note that $\MQTrHopf^\mathcal{E}$ is  a replete subcategory of $\QTrHopf^\mathcal{E}$.
\begin{invisible}
Indeed, if $(B,\Ss)\in\MQTrHopf^\mathcal{E}$,  $(A,\Rr)\in \QTrHopf^\mathcal{E}$ and $f:(B,\Ss)\to (A,\Rr)$ is an isomorphism in $\QTrHopf^\mathcal{E}$, then so is the canonical inclusion $\epsilon_{(A,\Rr)}=f\epsilon_{(B,\Ss)}FU(f^{-1}):FU(A)\to A$, as composition of isomorphisms, and so $(A,\Rr)\in \MQTrHopf^\mathcal{E}$.
\end{invisible}
Therefore, $\MQTrHopf^\mathcal{E}$ is a coreflective subcategory of $\QTrHopf^\mathcal{E}$ and $U$ is the \emph{coreflector}, see e.g. \cite[page 91]{Mac98} or the dual of \cite[Definition 3.5.2]{BorI94}.\medskip

Given a functor $U:\mathcal{A}\to \mathcal{B}$ and an object $B\in \mathcal{B}$, the \emph{essential fibre} of $U$ at $B$ is the subcategory $\mathcal{A}_B$ of $\mathcal{A}$ whose objects are pairs $(A,\alpha)$ where $A\in \mathcal{A}$ and $\alpha:B\to U(A)$ is an isomorphism in $\Bb$, and morphisms $f:(A,\alpha)\to (A',\alpha')$ are morphisms $f:A\to A'$ in $ \mathcal{A}$ such that $U(f)\alpha=\alpha'$.

\begin{theorem}\label{thm:protomodularityforfibre}
Any essential fibre of
$U:\QTrHopf^\mathcal{E}\to \MQTrHopf^\mathcal{E}$ is protomodular.
\end{theorem}

\begin{proof}
Given a minimal quasitriangular Hopf algebra $(B,\Ss)$, we can then consider the essential fibre $\QTrHopf^\mathcal{E}_{(B,\Ss)}$ of $U$ at $(B,\Ss)$ whose objects are pairs $((A,\Rr),\alpha)$ where $(A,\Rr)\in \QTrHopf^\mathcal{E}$ and $\alpha :(B,\Ss)\to (A_{\Rr},\Rr)$ is an isomorphism in $\MQTrHopf^\mathcal{E}$, while a morphism $f:((A,\Rr),\alpha)\to ((A',\Rr'),\alpha')$ is a morphism $f:(A,\Rr)\to (A',\Rr')$
such that $U(f)\alpha=\alpha'.$ Thus, in view of \cref{rmk:coslicemin}, the category $\QTrHopf^\mathcal{E}_{(B,\Ss)}$ is isomorphic to $ ((B,\Ss)/\QTrHopf)^i$ which, in turn, is protomodular by \cref{cor:cosliceQTrHprot}.
\end{proof}

\begin{remark}
Note that $\MQTrHopf^\mathcal{E}$ has $(\Bbbk,1\otimes_\Bbbk1)$ as a terminal object. The corresponding essential fibre $\QTrHopf^\mathcal{E}_{(\Bbbk,1\otimes_\Bbbk1)}$ of $U$ is isomorphic to the category $ (\Bbbk/\QTrHopf)^i=\Bbbk/\QTrHopf$ which is the category
$\cHopf$ of cocommutative Hopf algebras. Indeed, a quasitriangular Hopf algebra $(A,\Rr)$ has a morphism $f:(\Bbbk,1\otimes_\Bbbk1)\to (A,\Rr)$ of quasitriangular Hopf algebras if and only if $f:\Bbbk\to A$ is the unit morphism, which is injective, and $\Rr=(f\otimes_\Bbbk f)(1\otimes_\Bbbk1)=1\otimes_\Bbbk1$ i.e. if and only if it is a cocommutative Hopf algebra. Moreover, a morphism in $\Bbbk/\QTrHopf$ is just a morphism of cocommutative Hopf algebras. Thus we recover $\cHopf$ as the essential fibre of the coreflector $U$ at the terminal object. In particular, this confirms that $\cHopf$ is protomodular (see \cite{GrStVe}). It would be interesting to determine whether the protomodularity of this essential fibre implies that of all the others, possibly as a consequence of some specific properties of the coreflector that, at present, remain elusive.
\end{remark}

\noindent\emph{Acknowledgments.}
The authors thank Gabriella B\"ohm for insightful discussions concerning the notion of reversion on duoidal categories.
This paper was written while the authors were members of the “National Group for Algebraic and Geometric Structures and Their Applications” (GNSAGA-INdAM). The authors were also partially supported by the project funded by the European Union -NextGenerationEU under NRRP, Mission 4 Component 2 CUP D53D23005960006 - Call PRIN 2022 No. 104 of February 2, 2022 of Italian Ministry of University and Research; Project 2022S97PMY \emph{Structures for Quivers, Algebras and Representations} (SQUARE).
A. Sciandra was supported by a postdoctoral fellowship at the Université libre de Bruxelles within the framework of the PDR project ``Reconstruction of modules and algebraic objects from closed and monoidal structures on their representation categories'' funded by the FNRS under the grant number T.0318.25F (PI Joost Vercruysse).

\bibliography{protqtri}

\begin{thebibliography}{10}

\bibitem{Aguiar}
{\sc Aguiar, M., and Mahajan, S.}
\newblock {\em Monoidal functors, species and {H}opf algebras}, vol.~29 of {\em
  CRM Monograph Series}.
\newblock American Mathematical Society, Providence, RI, 2010.
\newblock With forewords by Kenneth Brown and Stephen Chase and Andr\'e{}
  Joyal.

\bibitem{AM-MMCat}
{\sc Ardizzoni, A., and Menini, C.}
\newblock Milnor-{M}oore categories and monadic decomposition.
\newblock {\em J. Algebra 448\/} (2016), 488--563.

\bibitem{Barr}
{\sc Barr, M.}
\newblock Functorial semantics and {HSP} type theorems.
\newblock {\em Algebra Universalis 31}, 2 (1994), 223--251.

\bibitem{Batanin-Markl}
{\sc Batanin, M., and Markl, M.}
\newblock Centers and homotopy centers in enriched monoidal categories.
\newblock {\em Adv. Math. 230}, 4-6 (2012), 1811--1858.

\bibitem{BD-CrossI}
{\sc Bespalov, Y., and Drabant, B.}
\newblock Cross product bialgebras. {I}.
\newblock {\em J. Algebra 219}, 2 (1999), 466--505.

\bibitem{BD-CrossII}
{\sc Bespalov, Y., and Drabant, B.}
\newblock Cross product bialgebras. {II}.
\newblock {\em J. Algebra 240}, 2 (2001), 445--504.

\bibitem{Bohm-Canada}
{\sc B\"ohm, G.}
\newblock Hopf algebroids: an approach via duoidal endohom categories.
\newblock
  \href{https://www.mun.ca/aac/media/production/memorial/academic/faculty-of-science/mathematics-and-statistics/atlantic-algebra-centre/media-library/workshops/nextwork/hart25/Bohm_Canada2025-final.pdf}{HART25.pdf},
  2025.
\newblock Hopf algebras and related topics, Canada.

\bibitem{BCZ}
{\sc B\"ohm, G., Chen, Y., and Zhang, L.}
\newblock On {H}opf monoids in duoidal categories.
\newblock {\em J. Algebra 394\/} (2013), 139--172.

\bibitem{Bohm-Lack}
{\sc B\"ohm, G., and Lack, S.}
\newblock Hopf comonads on naturally {F}robenius map-monoidales.
\newblock {\em J. Pure Appl. Algebra 220}, 6 (2016), 2177--2213.

\bibitem{BorI94}
{\sc Borceux, F.}
\newblock {\em Handbook of categorical algebra. 1}, vol.~50 of {\em
  Encyclopedia of Mathematics and its Applications}.
\newblock Cambridge University Press, Cambridge, 1994.
\newblock Basic category theory.

\bibitem{BB04}
{\sc Borceux, F., and Bourn, D.}
\newblock {\em Mal'cev, protomodular, homological and semi-abelian categories},
  vol.~566 of {\em Mathematics and its Applications}.
\newblock Kluwer Academic Publishers, Dordrecht, 2004.

\bibitem{Bourn}
{\sc Bourn, D.}
\newblock Normalization equivalence, kernel equivalence and affine categories.
\newblock In {\em Category theory ({C}omo, 1990)}, vol.~1488 of {\em Lecture
  Notes in Math.} Springer, Berlin, 1991, pp.~43--62.

\bibitem{BCPvO}
{\sc Bulacu, D., Caenepeel, S., Panaite, F., and Van~Oystaeyen, F.}
\newblock {\em Quasi-{H}opf algebras}, vol.~171 of {\em Encyclopedia of
  Mathematics and its Applications}.
\newblock Cambridge University Press, Cambridge, 2019.
\newblock A categorical approach.

\bibitem{BT25}
{\sc Bulacu, D., Popescu, D., and Torrecillas, B.}
\newblock Double wreath quasi-hopf algebras.
\newblock {\em J. Algebra 662\/} (2025), 1--71.

\bibitem{CDR}
{\sc Caenepeel, S., D\u{a}sc\u{a}lescu, S., and Raianu, c.~S.}
\newblock Classifying pointed {H}opf algebras of dimension 16.
\newblock {\em Comm. Algebra 28}, 2 (2000), 541--568.

\bibitem{Chen-quasi}
{\sc Chen, H.-X.}
\newblock Quasitriangular structures of bicrossed coproducts.
\newblock {\em J. Algebra 204}, 2 (1998), 504--531.

\bibitem{Dr87}
{\sc Drinfel'd, V.~G.}
\newblock Quantum groups.
\newblock In {\em Proceedings of the {I}nternational {C}ongress of
  {M}athematicians, {V}ol. 1, 2 ({B}erkeley, {C}alif., 1986)\/} (1987), Amer.
  Math. Soc., Providence, RI, pp.~798--820.

\bibitem{GrStVe}
{\sc Gran, M., Sterck, F., and Vercruysse, J.}
\newblock A semi-abelian extension of a theorem by {T}akeuchi.
\newblock {\em J. Pure Appl. Algebra 223}, 10 (2019), 4171--4190.

\bibitem{Johnstone-v1}
{\sc Johnstone, P.~T.}
\newblock {\em Sketches of an elephant: a topos theory compendium. {V}ol. 1},
  vol.~43 of {\em Oxford Logic Guides}.
\newblock The Clarendon Press, Oxford University Press, New York, 2002.

\bibitem{Kassel}
{\sc Kassel, C.}
\newblock Quantum groups.
\newblock In {\em Algebra and operator theory ({T}ashkent, 1997)}. Kluwer Acad.
  Publ., Dordrecht, 1998, pp.~213--236.

\bibitem{Ulrich-Myriam}
{\sc Kr\"ahmer, U., and Mahaman, M.}
\newblock Clones from comonoids.
\newblock {\em Rev. Un. Mat. Argentina 68}, 2 (2025), 369--394.

\bibitem{Mac98}
{\sc Mac~Lane, S.}
\newblock {\em Categories for the working mathematician}, second~ed., vol.~5 of
  {\em Graduate Texts in Mathematics}.
\newblock Springer-Verlag, New York, 1998.

\bibitem{Majid-book}
{\sc Majid, S.}
\newblock {\em Foundations of quantum group theory}.
\newblock Cambridge University Press, Cambridge, 1995.

\bibitem{P19}
{\sc Porst, H.-E.}
\newblock Colimits of monoids.
\newblock {\em Theory Appl. Categ. 34\/} (2019), 456--467.

\bibitem{Radford}
{\sc Radford, D.~E.}
\newblock On the antipode of a quasitriangular {H}opf algebra.
\newblock {\em J. Algebra 151}, 1 (1992), 1--11.

\bibitem{Radford-min}
{\sc Radford, D.~E.}
\newblock Minimal quasitriangular {H}opf algebras.
\newblock {\em J. Algebra 157}, 2 (1993), 285--315.

\bibitem{Radford-book}
{\sc Radford, D.~E.}
\newblock {\em Hopf algebras}, vol.~49 of {\em Series on Knots and Everything}.
\newblock World Scientific Publishing Co. Pte. Ltd., Hackensack, NJ, 2012.

\bibitem{Sar21}
{\sc Saracco, P.}
\newblock Antipodes, preantipodes and {F}robenius functors.
\newblock {\em J. Algebra Appl. 20}, 7 (2021), Paper No. 2150124, 32.

\bibitem{SZ}
{\sc Sciandra, A., and Zuo, Z.}
\newblock On the semi-abelianness of cocommutative {H}opf monoids.
\newblock \href{https://arxiv.org/pdf/2603.22200}{2603.22200}, 2026.

\bibitem{Street}
{\sc Street, R.}
\newblock Monoidal categories in, and linking, geometry and algebra.
\newblock {\em Bull. Belg. Math. Soc. Simon Stevin 19}, 5 (2012), 769--821.

\bibitem{Ve}
{\sc Vercruysse, J.}
\newblock {Hopf Algebras—Variant Notions and Reconstruction Theorems}.
\newblock In {\em {Quantum Physics and Linguistics: A Compositional,
  Diagrammatic Discourse}}. Oxford University Press, 02 2013.

\bibitem{Wakui}
{\sc Wakui, M.}
\newblock Various structures associated to the representation categories of
  eight-dimensional nonsemisimple {H}opf algebras.
\newblock {\em Algebr. Represent. Theory 7}, 5 (2004), 491--515.

\end{thebibliography}
\bibliographystyle{acm}

\end{document}